\documentclass[3p]{elsarticle}

\usepackage{amsthm}
\usepackage{amssymb}
\usepackage{amsmath}
\usepackage{amsfonts}
\usepackage{epsfig}
\usepackage[mathscr]{eucal}
\usepackage[ruled,vlined]{algorithm2e}
\usepackage{algorithmic}
\usepackage{setspace,color}
\usepackage{tikz}
\usepackage{hyperref}
\usepackage{cleveref}
\usepackage{numcompress}
\usepackage[caption=false]{subfig}
\usepackage{multirow}

\usetikzlibrary{snakes}
\usetikzlibrary{arrows,shapes}

\numberwithin{equation}{section}

\crefname{section}{section}{sections}
\crefname{subsection}{subsection}{subsections}

\Crefname{figure}{Figure}{Figures}

\crefformat{equation}{\textup{#2(#1)#3}}
\crefrangeformat{equation}{\textup{#3(#1)#4--#5(#2)#6}}
\crefmultiformat{equation}{\textup{#2(#1)#3}}{ and \textup{#2(#1)#3}}
{, \textup{#2(#1)#3}}{, and \textup{#2(#1)#3}}
\crefrangemultiformat{equation}{\textup{#3(#1)#4--#5(#2)#6}}%
{ and \textup{#3(#1)#4--#5(#2)#6}}{, \textup{#3(#1)#4--#5(#2)#6}}{, and \textup{#3(#1)#4--#5(#2)#6}}

\Crefformat{equation}{#2Equation~\textup{(#1)}#3}
\Crefrangeformat{equation}{Equations~\textup{#3(#1)#4--#5(#2)#6}}
\Crefmultiformat{equation}{Equations~\textup{#2(#1)#3}}{ and \textup{#2(#1)#3}}
{, \textup{#2(#1)#3}}{, and \textup{#2(#1)#3}}
\Crefrangemultiformat{equation}{Equations~\textup{#3(#1)#4--#5(#2)#6}}%
{ and \textup{#3(#1)#4--#5(#2)#6}}{, \textup{#3(#1)#4--#5(#2)#6}}{, and \textup{#3(#1)#4--#5(#2)#6}}

\crefdefaultlabelformat{#2\textup{#1}#3}

\DeclareMathAlphabet{\mathpzc}{OT1}{pzc}{m}{it}

\def\R{{\mathbb R}}
\def\C{{\mathbb C}}
\def\N{{\mathbb N}}

\def\Z{{\mathbb Z}}

\def\KK{{\mathbb K}}
\def\MM{{\mathbb M}}

\def\OO{{\mathbb O}}

\def\TGV{\mathrm{TGV}}

\def\rank{\mathrm{rank}}

\def\rd{\mathrm{d}}
\def\Om{\Omega}

\def\f{\frac}
\def\p{\partial}

\def\na{\nabla}
\def\la{\langle}
\def\ra{\rangle}

\def\a{{\boldsymbol a}}

\def\bb{{\boldsymbol b}}

\def\bc{{\boldsymbol c}}

\def\bsd{{\boldsymbol d}}

\def\bsf{{\boldsymbol f}}

\def\bi{{\mathbf i}}

\def\bk{{\boldsymbol k}}
\def\bsl{{\boldsymbol l}}

\def\bm{\boldsymbol{m}}
\def\bn{\boldsymbol{n}}
\def\bp{\boldsymbol{p}}
\def\bq{\boldsymbol{q}}

\def\bu{\boldsymbol{u}}
\def\bsv{\boldsymbol{v}}

\def\x{\boldsymbol{x}}

\def\bsI{{\boldsymbol I}}

\def\bX{{\boldsymbol X}}
\def\bY{{\boldsymbol Y}}

\def\bZ{{\boldsymbol Z}}
\def\mA{{\mathcal A}}

\def\mD{{\mathcal D}}
\def\mE{{\mathcal E}}

\def\mH{{\mathcal H}}
\def\mI{{\mathcal I}}

\def\mK{{\mathcal K}}

\def\mR{{\mathcal R}}
\def\mS{{\mathcal S}}
\def\mT{{\mathcal T}}

\def\mV{{\mathcal V}}
\def\mW{{\mathcal W}}

\def\bmA{{\boldsymbol \mA}}

\def\bmD{{\boldsymbol \mD}}
\def\bmE{{\boldsymbol \mE}}

\def\bmH{{\boldsymbol \mH}}
\def\bmI{{\boldsymbol \mI}}

\def\bmK{{\boldsymbol \mK}}

\def\bmR{{\boldsymbol \mR}}
\def\bmS{{\boldsymbol \mS}}
\def\bmT{{\boldsymbol \mT}}

\def\bmW{{\boldsymbol \mW}}

\def\msF{{\mathscr F}}

\def\msH{{\mathscr H}}

\def\msR{{\mathscr R}}
\def\msS{{\mathscr S}}
\def\msT{{\mathscr T}}

\def\msV{{\mathscr V}}

\def\bmsF{{\boldsymbol \msF}}

\def\aal{\boldsymbol{\alpha}}

\def\dde{\boldsymbol{\delta}}
\def\gga{\boldsymbol{\gamma}}

\def\ssi{\boldsymbol{\sigma}}

\def\nnu{\boldsymbol{\nu}}

\def\bvphi{\boldsymbol{\varphi}}

\def\bi{\begin{itemize}} \def\ei{\end{itemize}}
\def\be{\begin{eqnarray*}}
\def\ee{\end{eqnarray*}}

\def\0{{\mathbf 0}}

\newcommand{\beq}{\begin{equation}}
\newcommand{\eeq}{\end{equation}}

\def\xxi{{\boldsymbol{\xi}}}

\def\zze{{\boldsymbol\zeta}}

\def\Ph{\boldsymbol{\Phi}}
\def\Sig{\boldsymbol{\Sigma}}

\def\wt{\widetilde}
\def\wh{\widehat}

\def\Na{\boldsymbol \nabla}

\newcommand{\eps}{\varepsilon}
\def\la{\langle}
\def\ra{\rangle}

\def\XXint#1#2#3{{\setbox0=\hbox{$#1{#2#3}{\int}$ }
\vcenter{\hbox{$#2#3$ }}\kern-.55\wd0}}

\newtheorem{thm}{Theorem}[section]
\newtheorem{proposition}[thm]{Proposition}

\newtheorem{theorem}[thm]{Theorem}

\newdefinition{definition}[thm]{Definition}
\newdefinition{rmk}[thm]{Remark}
\newdefinition{notation}[thm]{Notation}
\newdefinition{example}[thm]{Example}

\newproof{pf}{Proof}

\newcommand{\argmin}{\operatornamewithlimits{argmin}}

\title{Image Restoration: Structured Low Rank Matrix Framework for Piecewise Smooth Functions and Beyond}

\author[HKUST]{Jian-Feng Cai\fnref{JFC}}
\ead{jfcai@ust.hk}
\author[TJU]{Jae Kyu Choi\fnref{JKC}\corref{cor}}
\ead{jaycjk@tongji.edu.cn}
\author[HKUST]{Jingyang Li\fnref{JFC}}
\ead{jlieb@connect.ust.hk}
\author[FUDAN]{Ke Wei\fnref{KW}\corref{cor}}
\ead{kewei@fudan.edu.cn}

\cortext[cor]{Corresponding author}

\address[HKUST]{Department of Mathematics, The Hong Kong University of Science and Technology, Hong Kong, China}
\address[TJU]{School of Mathematical Sciences, Tongji University, Shanghai, 200092 China}
\address[FUDAN]{School of Data Science, Fudan University, Shanghai, 200433 China}

\fntext[JFC]{J. F. Cai and J. Li are supported by the Hong Kong Research Grant Council (HKRGC) GRF 16306317 and 16309219.}
\fntext[JKC]{J. K. Choi is supported in part by the National Natural Science Foundation of China Youth Program 11901436.}
\fntext[KW]{K. Wei is supported by the National Natural Science Foundation of China Youth Program grant 11801088 and the Shanghai Sailing Program 18YF1401600.}

\begin{document}

\begin{abstract}
Recently, mapping a signal/image into a low rank Hankel/Toeplitz matrix has become an emerging alternative to the traditional sparse regularization, due to its ability to alleviate the basis mismatch between the true support in the continuous domain and the discrete grid. In this paper, we introduce a novel structured low rank matrix framework to restore piecewise smooth functions. Inspired by the total generalized variation to use sparse higher order derivatives, we derive that the Fourier samples of higher order derivatives satisfy an annihilation relation, resulting in a low rank multi-fold Hankel matrix. We further observe that the SVD of a low rank Hankel matrix corresponds to a tight wavelet frame system which can represent the image with sparse coefficients. Based on this observation, we also propose a wavelet frame analysis approach based continuous domain regularization model for the piecewise smooth image restoration. Finally, numerical results on image restoration tasks are presented as a proof-of-concept study to demonstrate that the proposed approach is compared favorably against several popular discrete regularization approaches and structured low rank matrix approaches.
\end{abstract}


\maketitle

\pagestyle{myheadings}
\thispagestyle{plain}
\markboth{Jian-Feng Cai, Jae Kyu Choi, Jingyang Li, and Ke Wei}{Image Restoration: SLRM Framework for Piecewise Smooth Functions and Beyond}

\section{Introduction}

Image restoration, including image denoising, deblurring, inpainting, computed tomography, etc., is one of the central problems in imaging science. It aims at recovering an image of high-quality from a given measurement which is degraded during the process of imaging, acquisition, and communication. The image restoration problem is typically modeled as the following linear inverse problem:
\begin{align}\label{Linear_IP}
\bsf=\bmA\bu+\zze,
\end{align}
where $\bsf$ is the degraded measurement or the observed image, $\zze$ is a certain additive noise, and $\bmA$ is some linear operator which takes different forms for different image restoration problems.

Since the operator $\bmA$ is in general ill-conditioned or non-invertible, it is in general to use a regularization on the images to be restored. Most widely used examples include the variational approaches including total variation (TV) \cite{L.I.Rudin1992} and its nonlocal variants \cite{X.Zhang2010}, the inf-convolution \cite{A.Chambolle1997}, the total generalized variation (TGV) \cite{K.Bredies2014,K.Bredies2010}, and the combined first and second order total variation \cite{M.Bergounioux2010,K.Papafitsoros2014}. Apart from the variational approaches, the applied harmonic analysis approaches including curvelets \cite{E.Candes2006}, Gabor frames \cite{H.Ji2017}, shearlets \cite{G.Kutyniok2011}, complex tight framelets \cite{B.Han2014}, and wavelet frames \cite{J.F.Cai2009,J.F.Cai2008,J.F.Cai2012,R.H.Chan2003} are also widely used in the literature. The common concept of these methods is to use a sparse regularization on discrete images under a discrete linear transformation to regularize smooth image components while preserving image singularities such as edges, ridges, and corners. However, in many applications, the true singularities lie in the continuous domain, and the discretization error will lead to the \emph{basis mismatch} \cite{Y.Chi2010,G.Ongie2016,J.Ying2017} between the true singularities and the discrete image grid. Such a basis mismatch would destroy the sparse structure of the image, and thus can degrade the restoration quality \cite{J.Ying2017}.

Recently, continuous domain regularization is emerging as a powerful alternative to the discrete domain sparse regularization \cite{B.N.Bhaskar2013,E.J.Candes2014,Y.Chen2014,G.Ongie2018}. By such an ``off-the-grid'' approach, we can exploit the sparsity prior in continuous domain, which enables us to alleviate the basis mismatch due to the discretization \cite{G.Ongie2018}. To the best of our knowledge, this off-the-grid regularization stems from the Prony's method \cite{Prony1795} which corresponds a superposition of a few sinusoids to a structured low rank matrix for the Dirac stream retrieval. Hence, we can adopt the so-called \emph{structured low rank matrix (SLRM) approach} \cite{Y.Chen2014,G.Ongie2017,J.C.Ye2017} for the restoration of spectrally sparse signal whose Fourier transform is a Dirac stream \cite{J.F.Cai2016,J.F.Cai2018a,J.F.Cai2019}. Apart from the spectrally sparse signal restoration, the SLRM can also be used to restore Fourier samples of a one dimensional piecewise constant signal \cite{T.Blu2008,M.Vetterli2002} as in this case the Fourier samples of a derivative becomes the superposition of a few sinusoids. However, even though the SLRM can be easily applied to the case of isolated singularities \cite{E.J.Candes2013,W.Xu2014}, the extension to the multi dimensional (piecewise smooth) image restoration is not straightforward \cite{G.Ongie2016}. Since the image singularities such as the edges and ridges in general form a continuous curve in a two dimensional domain, it is in general challenging to construct a structured low rank matrix from the Fourier samples.

In this paper, we introduce a new structured low rank matrix framework for the piecewise smooth image restoration. Our framework is inspired by the $k$th order total generalized variation ($\TGV^k$) \cite{K.Bredies2014,K.Bredies2010}, which is known to be effective in restoring the piecewise polynomial image with sharp edges \cite{W.Guo2014}. Specifically, following the SLRM framework for the piecewise constant image \cite{G.Ongie2016}, we assume that the image singularities (including both jumps and hidden jumps) are located in the zero level set of a band-limited periodic function (called the \emph{annihilating polynomial}). Then we can derive that the gradient can be decomposed into another vector field and the residual, and the Fourier samples of the residual and the symmetric gradient of this vector field can be annihilated by the convolution with the Fourier coefficients of the annihilating polynomial (called the \emph{annihilating filter}). From these annihilation relations, we deduce that the multi-fold Hankel matrices generated from the Fourier samples of a piecewise smooth image are low rank, which in turn enables a balance between the derivatives of order $1,\ldots,k$ via a combined rank minimization of multi-fold Hankel matrices.

As a by-product of the proposed structured low rank matrix framework, we further introduce a wavelet frame based sparse regularization model for the piecewise image restoration via the continuous domain regularization. Briefly speaking, if we can associate a signal/image with a low rank Hankel matrix, its right singular vectors form tight frame filter banks under which the canonical coefficients have a group sparsity according to the index of filters \cite{J.F.Cai2020}. Then motivated by \cite{D.Guo2018}, we assume that the right singular vectors are estimated from a pre-restoration process, and we propose a \emph{sparse regularization} via the wavelet frame analysis approach (e.g. \cite{J.F.Cai2009/10}) as an image restoration model via the continuous domain regularization. Notice that there are several wavelet frame based approaches for the piecewise smooth image restoration \cite{J.F.Cai2016a,J.K.Choi2020,H.Ji2016} in the literature. However, while these existing approaches focus on the sparse approximation of a discrete image, our approach comes from the relaxation of the structured low rank matrices generated by the Fourier samples for the continuous domain regularization.

\subsection{SLRM for piecewise smooth image restoration}\label{SLMAPS}

Our SLRM framework for the piecewise smooth function is mostly related to the recent extension of the SLRM framework to the two dimensional functions in \cite{G.Ongie2018,G.Ongie2015a,G.Ongie2015,G.Ongie2016,H.Pan2014}. Briefly speaking, if the singularity curves of a piecewise constant/holomorphic function, i.e. the supports of the first order (real/complex) derivatives of a target image, lie in the zero level set of an annihilating polynomial, the Fourier transform of the derivatives can be annihilated by the convolution with the annihilating filter. This annihilation relation in turn corresponds the Fourier samples of derivatives to the structured low rank matrices. Based on this framework, the SLRM framework for the piecewise constant image restoration is proposed and studied in \cite{G.Ongie2015a,G.Ongie2015,G.Ongie2016,G.Ongie2017}, together with a restoration guarantee \cite{G.Ongie2018}.

This annihilation relation of the gradient can be easily extended to the higher order derivatives whenever the jump discontinuities of the corresponding derivatives are located in the zero level set of a trigonometric polynomial \cite{G.Ongie2015a}. For instance, we can derive an annihilation relation for the piecewise linear function by considering the Fourier transform of the second order derivatives. Based on this idea, the authors in \cite{Y.Hu2019} proposed a so-called generalized structured low rank (GSLR) approach for the piecewise smooth image restoration. Inspired by the $2$-fold inf-convolution \cite{A.Chambolle1997}, the GSLR approach restores a superposition of a piecewise constant layer and a piecewise linear layer whose first and second order derivatives correspond to low rank Hankel matrices in the frequency domain respectively. By balancing the first order and second order derivatives, the GSLR approach has demonstrated significant improvements in the piecewise smooth image restoration tasks over the existing approaches.

Note that it is not difficult to extend the previous GSLR framework to generic piecewise smooth functions. More precisely, by considering the annihilation relation of higher order derivatives, we can extend the GSLR to the so-called $k$-fold inf-convolution ($k\geq3$) for functions with higher regularity. However, in a piecewise smooth function, there could exist singularities on which the derivatives have jump discontinuities. In this case, the annihilating polynomial describing the image singularities cannot annihilate the $k$th order derivatives \cite[Proposition 4]{G.Ongie2015a}. Since we then need to consider the power of the annihilating polynomial, or equivalently, increase the size of the filter to guarantee the annihilation relation for higher order derivatives, the GSLR framework can degrade the low rank (multi-fold) Hankel matrix structure of the higher order derivatives.

Unlike the GSLR framework in \cite{Y.Hu2019} based on the inf-convolution, the proposed SLRM framework is based on the total generalized variation framework \cite{K.Bredies2010}. Briefly speaking, we decompose the gradient into the residual whose Fourier samples correspond to the structured low rank matrix and a vector field whose Fourier samples of the symmetric gradient correspond to the structured low rank matrix. In other words, since our proposed SLRM framework uses the Fourier samples of the successive first order derivatives, it does not require the powers of the annihilating polynomial describing the image singularities, which enables us to obtain better low rank Hankel matrix structures corresponding to the piecewise smooth functions. In addition, since the TGV takes the inf-convolution as a special case, it can also be verified that the proposed SLRM framework is more general than the GSLR framework for the piecewise smooth functions.

\subsection{Organization and notation of paper}

The rest of this paper is organized as follows. In \cref{OurFramework}, we present the proposed structured low rank matrix framework for the piecewise smooth functions. We first describe the one dimensional case to see the idea clearly, then we extend to the idea to the two dimensional framework. In \cref{ImageRestoration}, we present an image restoration model based on the wavelet frame as an application of the proposed structured low rank matrix framework, followed by an alternating minimization algorithm, and some numerical results are presented to demonstrate the performance of our framework in the piecewise smooth image restoration in \cref{Experiments}. Finally, \cref{Conclusion} concludes this paper with a few future directions. All technical proofs will be postponed to appendices.

Throughout this paper, all two dimensional images and two dimensional $k$ tensors defined on the discrete grid will be denoted by the bold faced lower case letters. Note that a two dimensional discrete image and a two dimensional discrete $k$ tensor can also be identified with a vector and a $2^k$-tuple of vectors (as well as a sequence or a $2^k$-tuple of sequences supported on the grid) respectively whenever convenient. All matrices will be denoted by the bold faced upper case letters, and the $m$th row and the $n$th column of a matrix $\bZ$ will be denoted by $\bZ^{(m,:)}$ and $\bZ^{(:,n)}$, respectively. Denote by
\begin{align}\label{ImageGrid}
\OO=\left\{-\lfloor N/2\rfloor,\ldots,\lfloor(N-1)/2\rfloor\right\}^2
\end{align}
with $N\in\N$, the set of $N\times N$ grid. The space of complex valued functions on $\OO$ and the space of complex $k$-tensor valued functions on $\OO$ are denoted by $\msV\simeq\C^{|\OO|}$ and $\msV_k\simeq\C^{|\OO|\times{2^k}}$, respectively. Notice that $\msV=\msV_0$. Given two rectangular grids $\KK$ and $\MM$, we define
\begin{align*}
\KK:\MM=\left\{\bk\in\KK:\bk+\MM\subseteq\KK\right\}=\left\{\bk\in\KK:\bk+\bm\in\KK~\text{for all}~\bm\in\MM\right\}.
\end{align*}
Operators on both images and tensors are denoted as the bold faced caligraphic letters. For instance, letting $\bsv\in\msV$ and $\KK$ be a rectangular $K_1\times K_2$ grid, the corresponding Hankel matrix $\bmH\bsv$ is an $M_1\times M_2$ matrix ($M_1=|\OO:\KK|$ and $M_2=|\KK|$) generated by concatenating $K_1\times K_2$ patches of $\bsv$ into row vectors. Note that, in the sense of multi-indices, we have
\begin{align*}
\left(\bmH\bsv\right)(\bk,\bm)=\bsv(\bk+\bm),~~~~~\bk\in\OO:\KK,~~\text{and}~~\bm\in\KK.
\end{align*}
With a slight abuse of notation, we also use
\begin{align}\label{kFoldHankel}
\bmH\bq=\left[\begin{array}{ccc}
\left(\bmH\bq_1\right)^T&\cdots&\left(\bmH\bq_{2^k}\right)^T
\end{array}\right]^T\in\C^{kM_1\times M_2}
\end{align}
to denote the $k$-fold Hankel matrix constructed from $\bq=\left(\bq_1,\ldots,\bq_{2^k}\right)\in\msV_k$.

\section{Structured low rank matrix framework for piecewise smooth functions}\label{OurFramework}

In this section, we introduce our structured low rank matrix framework for piecewise smooth functions. For simplicity, we consider the piecewise linear function throughout this paper. Note, however, it is not difficult to extend the proposed framework into the general piecewise smooth functions.

\subsection{SLRM framework for one dimensional signals}\label{1DFramework}

We first establish the structured low rank matrix framework of the following one dimensional piecewise linear model
\begin{align}\label{signalModel}
u(x)=\sum_{j=1}^{K-1}\left(\alpha_jx+\beta_j\right)1_{[x_j,x_{j+1})}(x),
\end{align}
from its Fourier sample
\begin{align*}
\wh{u}(k)=\msF(u)(k)=\int_{-\infty}^{\infty}u(x)e^{-2\pi ikx}\rd x,~~~~~k\in\left\{-\left\lfloor N/2\right\rfloor,\ldots,\left\lfloor(N-1)/2\right\rfloor\right\},~~~N\in\N.
\end{align*}
In \cref{signalModel}, $\alpha_j,\beta_j\in\C$, $1_{[x_j,x_{j+1})}$ denotes the characteristic function on the interval $[x_j,x_{j+1})$: $1_{[x_j,x_{j+1})}(x)=1$ if $x\in[x_j,x_{j+1})$, and $0$ otherwise, and $-1/2<x_1<x_2<\cdots<x_K<1/2$ are the location of singularities, i.e. by the singularities we mean they include both the jumps and the hidden jumps (jumps of derivatives).

In the sense of distribution, the derivative $u'$ satisfies
\begin{align}\label{uDeri}
u'(x)=\sum_{j=1}^{K-1}\alpha_j1_{[x_j,x_{j+1})}(x)+\sum_{j=1}^K\left[\msT_j(u)(x_j)\delta(x-x_j)\right]
\end{align}
where $\msT_j(u)(x_j)=\left(\alpha_j-\alpha_{j-1}\right)x_j+\left(\beta_j-\beta_{j-1}\right)$ with $\alpha_0=\alpha_K=\beta_0=\beta_K=0$. Letting
\begin{align}\label{p1D}
p(x)=\sum_{j=1}^{K-1}\alpha_j1_{[x_j,x_{j+1})}(x),
\end{align}
we have
\begin{align*}
u'(x)-p(x)=\sum_{j=1}^K\msT_j(u)(x_j)\delta(x-x_j).
\end{align*}
Hence, the Fourier transform of $u'-p$ is a linear combination of complex sinusoids:
\begin{align}
\msF(u'-p)(\xi)=\sum_{j=1}^K\msT_j(u)(x_j)e^{-2\pi ix_j\xi},~~~~~\xi\in\R.
\end{align}
In addition, since the derivative of $p$ is also a Dirac stream:
\begin{align}
p'(x)=\sum_{j=1}^K\left(\alpha_j-\alpha_{j-1}\right)\delta(x-x_j),
\end{align}
its Fourier transform is expressed as
\begin{align}
\msF(p')(\xi)=\sum_{j=1}^K\left(\alpha_j-\alpha_{j-1}\right)e^{-2\pi ix_j\xi},~~~~~\xi\in\R.
\end{align}
i.e. a linear combination of complex sinusoids.

We introduce the following trigonometric polynomial
\begin{align}\label{1DAnnihilatingPolynomial}
\varphi(x)=\prod_{j=1}^K\left(e^{-2\pi ix}-e^{-2\pi ix_j}\right):=\sum_{k=0}^K\a(k)e^{-2\pi ikx}.
\end{align}
Since $\varphi(x_j)=0$ for $j=1,\ldots,K$, it follows that $\varphi\left(u'-p\right)=0$ and $\varphi p'=0$ in the sense of distribution. In the frequency domain, since we have
\begin{align*}
\wh{\varphi}(\xi)=\msF(\varphi)(\xi)=\sum_{k=0}^K\a(k)\delta(\xi+k),
\end{align*}
it follows that
\begin{align}
\left(\msF(u'-p)\ast\wh{\varphi}\right)(\xi)&=\sum_{k=0}^K\msF(u'-p)(\xi+k)\a(k)=0,\label{1DAnni1}\\
\left(\msF(p')\ast\wh{\varphi}\right)(\xi)&=\sum_{k=0}^K\msF(p')(\xi+k)\a(k)=0,\label{1DAnni2}
\end{align}
for $\xi\in\R$. Therefore, the Fourier transforms of both $u'-p$ and $p'$ are annihilated by the convolution under the Fourier coefficients of $\varphi$.

In many practical cases, we consider the contiguous Fourier samples (i.e. the samples on a regular grid $\left[N\right]:=\left\{-\lfloor N/2\rfloor,\ldots,\lfloor(N-1)/2\rfloor\right\}$), so \cref{1DAnni1,1DAnni2} become the following systems of linear equations
\begin{align}
\sum_{k=0}^K\msF(u'-p)(l+k)\a(k)&=0,\label{1DAnniDis1}\\
\sum_{k=0}^K\msF(p')(l+k)\a(k)&=0.\label{1DAnniDis2}
\end{align}
In the matrix-vector multiplication form, \cref{1DAnniDis1,1DAnniDis2} lead to
\begin{align*}
\bmH\left(\msF\left(u'-p\right)\big|_{[N]}\right)\a=\0~~~\text{and}~~~\bmH\left(\msF\left(p'\right)\big|_{[N]}\right)\a=\0,
\end{align*}
which shows that two $(N-K+1)\times K$ Hankel matrices $\bmH\left(\msF\left(u'-p\right)\big|_{[N]}\right)$ and $\bmH\left(\msF\left(p'\right)\big|_{[N]}\right)$ have nontrivial nullspaces. In particular, for $M>K$, we can see that
\begin{align*}
e^{-2\pi imx}\varphi(x)=\sum_{k=m}^{K+m}\a(k-m)e^{-2\pi ikx},~~~~~m=0,\ldots,M-K
\end{align*}
is also an annihilating polynomial, which shows that, if $K\leq(N-M+1)\wedge M$, we have
\begin{align*}
\rank\left(\bmH\left(\msF\left(u'-p\right)\big|_{[N]}\right)\right)\leq K~~~\text{and}~~~\rank\left(\bmH\left(\msF\left(p'\right)\big|_{[N]}\right)\right)\leq K.
\end{align*}
Therefore, both $\bmH\left(\msF\left(u'-p\right)\big|_{[N]}\right)\in\C^{(N-M+1)\times M}$ and $\bmH\left(\msF\left(p'\right)\big|_{[N]}\right)\in\C^{(N-M+1)\times M}$ are \emph{rank-deficient,} which shows that it is possible to convert the piecewise regularity in the continuous domain into the low rank Hankel matrices corresponding to the discrete Fourier samples.

\subsection{SLRM framework for two dimensional images}\label{2DFramework}

Now we aim to establish the structured low rank matrix framework for the following two dimensional piecewise linear function model
\begin{align}\label{uModel}
u(\x)=\sum_{j=1}^Ju_j(\x)1_{\Om_j}(\x):=\sum_{j=1}^J\left(\aal_j^T\x+\beta_j\right)1_{\Om_j}(\x),~~~~~~~\x\in\R^2,
\end{align}
from its Fourier samples
\begin{align}\label{Forward}
\wh{u}(\bk)=\msF(u)(\bk)=\int_{\R^2}u(\x)e^{-2\pi i\bk\cdot\x}\rd\x,~~~~~~\bk\in\OO,
\end{align}
where the sample grid $\OO$ is defined as \cref{ImageGrid}. In \cref{uModel}, $\aal_j=\left(\alpha_{j1},\alpha_{j2}\right)\in\C^2$, $\beta_j\in\C$ and $1_{\Om_j}$ denotes the characteristic function on a domain $\Om_j$; $1_{\Om_j}(\x)=1$ if $\x\in\Om_j$, and $0$ otherwise. Without loss of generality, we assume that $\Om_j$ lies in $[-1/2,1/2)^2$ for simplicity. (However, it is not difficult to generalize the setting into an arbitrary rectangular region $[-L_1/2,L_1/2)\times[-L_2/2,L_2/2)$.) We further assume that \cref{uModel} is expressed with the smallest number of characteristic functions such that $\Om_j$'s are pairwise disjoint. Then the singularities of $u$, including both jump discontinuities of $u$ and hidden jump discontinuities of $u$ (jumps of derivatives), agree with $\Gamma=\bigcup_{j=1}^J\p\Om_j$, which will be called the \emph{singularity set} of $u$ throughout this paper.

Generally, it is difficult to directly establish the SLRM framework without any further information on the singularity set $\Gamma$. Inspired by the two dimensional FRI framework for the piecewise constant function \cite{G.Ongie2015a,G.Ongie2015,G.Ongie2016,G.Ongie2017}, we assume that there exists a finite rectangular and symmetric\footnote{Throughout this paper, we only consider a finite and symmetric grid for the trigonometric polynomials.} index set $\KK$ such that
\begin{align}\label{MajorAssumption}
\Gamma\subseteq\left\{\x\in\R^2:\varphi(\x)=0\right\}~~~~~\text{with}~~~~~\varphi(\x)=\sum_{\bk\in\KK}\a(\bk)e^{-2\pi i\bk\cdot\x}.
\end{align}
Throughout this paper, we call any function $\varphi(\x)$ in the form of \cref{MajorAssumption} the \emph{trigonometric polynomial}, and the zero level set $\left\{\x\in\R^2:\varphi(\x)=0\right\}$ the \emph{trigonometric curve}. For a trigonometric polynomial $\varphi$ in \cref{MajorAssumption}, the degree of $\varphi$ is defined as an ordered pair of degrees in each coordinate, and is denoted by $\deg(\varphi)$. In particular, $\varphi$ in \cref{MajorAssumption} with the smallest degree is called the \emph{minimal polynomial}. The algebraic properties of trigonometric polynomials and curves, including the existence of the minimal polynomial, are studied in \cite{G.Ongie2016}, to which interested readers can refer for details.

Under this setting, we present \cref{Th1} to establish the following \emph{(linear) annihilation relation} for the piecewise linear function. The proof can be found in \ref{ProofTh1}.

\begin{theorem}\label{Th1} Let $u(\x)$ be defined as in \cref{uModel} where the singularity set $\Gamma$ satisfies \cref{MajorAssumption}. Then for some vector field $p=(p_1,p_2)$, we have
\begin{align}
\sum_{\bk\in\KK}\msF\left(\na u-p\right)(\xxi+\bk)\a(\bk)&=0,\label{FirstAnnihil}\\
\sum_{\bk\in\KK}\msF\left(\na_sp\right)(\xxi+\bk)\a(\bk)&=0,\label{SecondAnnihil}
\end{align}
where $\xxi\in\R^2$. Here, $\na_s$ is a symmetric gradient defined for $p=(p_1,p_2)$ as
\begin{align}\label{SymGrad}
\na_sp=\f{1}{2}\left(\na p+\na p^T\right)=\left[\begin{array}{cc}
\p_1p_1&\displaystyle{\f{1}{2}\left(\p_2p_1+\p_1p_2\right)}\\
\displaystyle{\f{1}{2}\left(\p_2p_1+\p_1p_2\right)}&\p_2p_2
\end{array}\right]
\end{align}
and the Fourier transform $\msF$ is performed to each component of $\na u-p$ and $\na_sp$.
\end{theorem}

Based on \cref{Th1}, we call $\varphi$ satisfying \cref{FirstAnnihil,SecondAnnihil} an \emph{annihilating polynomial} (for $u$ under a vector field $p$), and the Fourier coefficients $\a:=\left\{\a(\bk):\bk\in\KK\right\}$ an \emph{annihilating filter}. Note that the gradient of $u$ is decomposed into $\na u-p$ and $p$ such that the Fourier transforms of both $\na u-p$ and $\na_sp$ are annihilated by the same annihilating filter $\a$. In addition, since the proof of \cref{Th1} tells us that $p$ is the piecewise constant vector field in $\na u$, \cref{SecondAnnihil} can be viewed as an extension of annihilation relation for piecewise constant functions to piecewise constant vector fields.

\begin{rmk}\label{RK1} According to the proof, \cref{Th1} can also be written as follows: for a function $u$ defined as in \cref{uModel}, $\na u$ can be decomposed into
\begin{align}\label{GraduDecompose}
\na u=p+\rd\nnu
\end{align}
with a piecewise constant vector field $p$ and a Radon vector measure $\nnu$ supported on $\Gamma$ such that
\begin{align}
\sum_{\bk\in\KK}\msF\left(\nnu\right)(\xxi+\bk)\a(\bk)&=0,\label{FirstAnnihilRe}\\
\sum_{\bk\in\KK}\msF\left(\na_sp\right)(\xxi+\bk)\a(\bk)&=0.\label{SecondAnnihilRe}
\end{align}
We mention that the decomposition \cref{GraduDecompose} is well-defined in the sense that $\msF(\na_s\nnu)$ and $\msF(\na_sp)$ do not share the minimal annihilating polynomial. To see this, let $\varphi$ in \cref{MajorAssumption} be the minimal polynomial for $\Gamma$. For $F\in C_0^{\infty}(\R^2,\mathrm{Sym}^2(\R^2))$, we have
\begin{align*}
\la\varphi\na_s\nnu,F\ra=\int_{\Gamma}\left[\na_s^T\left(\varphi F\right)\right]\cdot\rd\nnu=-\int_{\Gamma}\left[F\left(\na\varphi\right)\right]\cdot\rd\nnu-\int_{\Gamma}\left(\varphi\na_s^TF\right)\cdot\rd\nnu.
\end{align*}
Since $\varphi$ is the minimal polynomial for a trigonometric curve $\Gamma$, $\na\varphi\neq0$ a.e. on $\Gamma$ \cite[Proposition A.4]{G.Ongie2016}. Hence, $\varphi\na_s\nnu\neq0$, or equivalently, $\msF(\na_s\nnu)$ does not satisfy the annihilation relation with $\a$. In fact, $\na_s\nnu$ is annihilated by $\varphi^2$, and since $\varphi$ is the minimal polynomial for $\Gamma$, $\varphi^2$ is the minimal polynomial for annihilating $\na_s\nnu$.
\end{rmk}

We claim that the linear annihilation relations \cref{FirstAnnihil} and \cref{SecondAnnihil} can be seen as balancing the annihilations of first and second order derivatives. To see this, we firstly consider an extreme case; a piecewise constant function case. In this case, it would be favorable to choose $p=0$, so that \cref{FirstAnnihil} reduces to
\begin{align}\label{Extreme1}
\sum_{\bk\in\KK}\msF\left(\na u\right)(\xxi+\bk)\a(\bk)=0,
\end{align}
which is the annihilation relation for $\na u$ in \cite{G.Ongie2015a,G.Ongie2015,G.Ongie2016,G.Ongie2017}. For the second extreme case, we note that
\begin{align*}
\na^2u=\na_s\left(\na u\right)=\na_s\left(\na u-p\right)+\na_sp.
\end{align*}
In addition, when $\na u$ is a piecewise constant vector field (i.e. $u$ is continuous on $\Gamma$), we can simply choose $p=\na u$. Then \cref{SecondAnnihil} becomes
\begin{align}\label{Extreme2}
\sum_{\bk\in\KK}\msF\left(\na_sp\right)(\xxi+\bk)\a(\bk)=\sum_{\bk\in\KK}\msF\left(\na^2u\right)(\xxi+\bk)\a(\bk)=0,
\end{align}
the annihilation relation for the second derivatives of $u$. Hence, the vector field $p$ in turn balances \cref{Extreme1,Extreme2} to annihilate the jumps of both $u$ and $\na u$ under the annihilating polynomial $\varphi$. Hence, \cref{FirstAnnihil,SecondAnnihil} are closely related to the following TGV:
\begin{align}\label{TGV2}
\mathrm{TGV}(u)=\inf_{u,p}\gamma_1\left\|\na u-p\right\|_1+\gamma_2\left\|\na_s p\right\|_1,
\end{align}
which in turn measures the jump discontinuities of $u$ and $p$ in our setting. In particular, if $u$ is continuous on $\p\Om_j$ for some $j$, this $\p\Om_j$ will not be reflected in \cref{FirstAnnihil}.

We also note that if we restrict $p$ in the range of $\na$, i.e. $p=\na u_2$, then by letting $u_1=u-u_2$, we can rewrite \cref{FirstAnnihil,SecondAnnihil} as
\begin{align}
\sum_{\bk\in\KK}\msF\left(\na u_1\right)(\xxi+\bk)\a(\bk)&=0\label{PCAnnihil}\\
\sum_{\bk\in\KK}\msF\left(\na^2u_2\right)(\xxi+\bk)\a(\bk)&=0\label{PLAnnihil},
\end{align}
the GSLR framework in \cite{Y.Hu2019}, which promotes a similar balance through a decomposition $u=u_1+u_2$ with a piecewise constant $u_1$ and a piecewise linear $u_2$. We can easily see that \cref{PCAnnihil,PLAnnihil} are closely related to the following inf-convolution
\begin{align}\label{InfConv2}
\left(J_1\square J_2\right)(u):=\inf_{u=u_1+u_2}\gamma_1\|\na u_1\|_1+\gamma_2\|\na^2u_2\|_1.
\end{align}
Noting that when $p=\na u_2$, then by letting $u_1=u-u_2$, the TGV \cref{TGV2} becomes the above inf-convolution, which means that the TGV takes the inf-convolution as a special case. Likewise, we can see that \cref{FirstAnnihil,SecondAnnihil} in \cref{Th1} take \cref{PCAnnihil,PLAnnihil} in \cite{Y.Hu2019} as a special case.

\begin{rmk}\label{RK2} We further mention that \cref{Th1} is different from the GSLR framework in \cite{Y.Hu2019}. To see this, we assume that the singularity set $\Gamma$ satisfies \cref{MajorAssumption} with the minimal polynomial $\varphi$, i.e. $\KK$ is the smallest support of the Fourier coefficients $\a$. Since $\na\varphi\neq0$ a.e. on $\Gamma$ \cite[Proposition A.4]{G.Ongie2016}, $\varphi$ may not be able to annihilate the second order derivatives of a piecewise linear function in general \cite{G.Ongie2015a}. Hence, \cref{PCAnnihil,PLAnnihil} prefer to decompose $u=u_1+u_2$ where $u_1$ is a piecewise constant function and $u_2$ is a piecewise linear spline (a continuous piecewise linear function). Meanwhile, the proposed framework is established by annihilating the jumps for each derivative successively, which enables to cover a broader range of piecewise linear functions. Hence, our framework is more generalized than the GSLR in \cite{Y.Hu2019}.
\end{rmk}

Next, we will present two examples to illustrate this.

\begin{example}\label{Ex1} We consider one dimensional examples. Let $u(x)=\left(x+5/4\right)1_{[-1/4,0)}(x)+\left(5/4-x\right)1_{[0,1/4)}(x)$. Then we can see that
\begin{align*}
\varphi(x)=\left(e^{-2\pi ix}-e^{\pi i/2}\right)\left(e^{-2\pi ix}-1\right)\left(e^{-2\pi ix}-e^{-\pi i/2}\right):=\sum_{k=0}^3\a(k)e^{-2\pi ikx}.
\end{align*}
is a minimal polynomial for the singularity set $\left\{-1/4,0,1/4\right\}$. Since
\begin{align*}
u'(x)=1_{[-1/4,0)}(x)-1_{[0,1/4)}(x)+\delta(x+1/4)-\delta(x-1/4)
\end{align*}
we can choose $p=1_{[-1/4,0)}-1_{[0,1/4)}$, i.e.
\begin{align*}
p'(x)=\delta(x+1/4)-2\delta(x)+\delta(x-1/4).
\end{align*}
Hence, it follows that
\begin{align*}
\sum_{k=0}^3\msF(u'-p)(\xi+k)\a(k)=0~~~~~\text{and}~~~~~\sum_{k=0}^3\msF(p')(\xi+k)\a(k)=0.
\end{align*}
By letting $u_1(x)=1_{[-1/4,1/4)}(x)$ and $u_2(x)=\left(x+1/4\right)1_{[-1/4,0)}(x)+\left(1/4-x\right)1_{[0,1/4)}(x)$, we also have
\begin{align*}
u_1'(x)&=\delta(x+1/4)-\delta(x-1/4),\\
u_2''(x)&=\delta(x+1/4)-2\delta(x)+\delta(x-1/4),
\end{align*}
which leads to
\begin{align*}
\sum_{k=0}^3\msF(u_1')(\xi+k)\a(k)=0~~~~~\text{and}~~~~~\sum_{k=0}^3\msF(u_2'')(\xi+k)\a(k)=0.
\end{align*}
In other words, both frameworks can establish annihilation relations if a given function can be decomposed into a piecewise constant function and a piecewise linear spline.
\end{example}

\begin{example}\label{Ex2} Let $u(x)=x1_{[-1/4,1/4)}(x)$ with
\begin{align*}
\varphi(x)=\left(e^{-2\pi ix}-e^{\pi i/2}\right)\left(e^{-2\pi ix}-e^{-\pi i/2}\right):=\sum_{k=0}^2\a(k)e^{-2\pi ikx}.
\end{align*}
being a minimal polynomial for the singularity set $\left\{-1/4,1/4\right\}$. In this case, we have
\begin{align*}
u'(x)=1_{[-1/4,1/4)}(x)-\delta(x+1/4)-\delta(x-1/4).
\end{align*}
If we choose $p=1_{[-1/4,1/4)}$, then we have
\begin{align*}
p'(x)=\delta(x+1/4)-\delta(x-1/4),
\end{align*}
which leads to
\begin{align*}
\sum_{k=0}^2\msF(u'-p)(\xi+k)\a(k)=0~~~~~\text{and}~~~~~\sum_{k=0}^2\msF(p')(\xi+k)\a(k)=0.
\end{align*}
However, we cannot find a piecewise constant $u_1$ and a piecewise linear spline $u_2$ such that $u=u_1+u_2$,
\begin{align*}
\sum_{k=0}^2\msF(u_1')(\xi+k)\a(k)=0~~~~~\text{and}~~~~~\sum_{k=0}^2\msF(u_2'')(\xi+k)\a(k)=0.
\end{align*}
Notice that we have
\begin{align*}
u''(x)=\delta(x+1/4)-\delta(x-1/4)-\delta'(x+1/4)-\delta'(x-1/4),
\end{align*}
and $\varphi'(x)\neq0$ for $x=\pm1/4$, so $u''\varphi\neq0$. In fact, we have $u''\varphi^2=0$, so that
\begin{align*}
\sum_{k=0}^4\msF(u'')(\xi+k)\bb(k)=0,
\end{align*}
where
\begin{align*}
\varphi^2(x)=\left(e^{-2\pi ix}-e^{\pi i/2}\right)^2\left(e^{-2\pi ix}-e^{-\pi i/2}\right)^2:=\sum_{k=0}^4\bb(k)e^{-2\pi ikx}.
\end{align*}
Hence, these two examples again illustrate that the proposed framework is more generalized than the GSLR framework in \cite{Y.Hu2019}.
\end{example}

In our setting, the Fourier transform of $u$ is sampled on the grid $\OO$ in \cref{ImageGrid} with $N\in\N$ large enough to guarantee a high image resolution. Hence, \cref{FirstAnnihil,SecondAnnihil} become the following finite systems of linear equations
\begin{align}
\sum_{\bk\in\KK}\msF\left(\na u-p\right)(\bm+\bk)\a(\bk)&=0,\label{FirstSystem}\\
\sum_{\bk\in\KK}\msF\left(\na_sp\right)(\bm+\bk)\a(\bk)&=0,\label{SecondSystem}
\end{align}
where $\bm\in\OO:\KK$. In the matrix-vector multiplication form, we have
\begin{align}
\bmH\left(\left(\msF(\na u-p)\right)\big|_{\OO}\right)\a&=\0,\label{FirstRankDeficient}\\
\bmH\left(\left(\msF(\na_sp)\right)\big|_{\OO}\right)\a&=\0.\label{SecondRankDeficient}
\end{align}
Hence, both $\bmH\left(\left(\msF(\na u-p)\right)\big|_{\OO}\right)$ and $\bmH\left(\left(\msF(\na_sp)\right)\big|_{\OO}\right)$ have nontrivial nullspaces. In addition, when the filter support $\KK'$ defining two multi-fold Hankel matrices is sufficiently large, both of them have nontrivial nullspaces as well. To see this, let $\varphi$ be the minimal polynomial for the singularity set with coefficients $\a$ supported on $\KK$. Then for any trigonometric polynomial $\eta$ with coefficients $\bb$ supported on $\KK':\KK$ such that $\a\ast\bb$ is supported on $\KK'$, it follows that
\begin{align*}
\left(\eta\varphi\right)(\x):=\sum_{\bk\in\KK'}\left(\a\ast\bb\right)(\bk)e^{-2\pi i\bk\cdot\x}
\end{align*}
is an annihilating polynomial and $\a\ast\bb$ is an annihilating filter, due to the associativity of convolution. In particular, by letting
\begin{align*}
\bb(\bk)=\dde(\bk-\bm)=\left\{\begin{array}{cl}
1&\text{if}~\bk=\bm\vspace{0.4em}\\
0&\text{otherwise}
\end{array}\right.~~~~~\bk\in\KK,~~\text{and}~~\bm\in\KK':\KK,
\end{align*}
we can see that
\begin{align*}
e^{-2\pi i\bm\cdot\x}\varphi(\x)=\sum_{\bk\in\bm+\KK}\a(\bk-\bm)e^{-2\pi i\bk\cdot\x},~~~~~\bm\in\KK':\KK,
\end{align*}
is also an annihilating polynomial, or equivalently, the translation $\a(\cdot-\bm)$ also satisfies \cref{FirstRankDeficient,SecondRankDeficient}. Based on this observation, we present proposition \ref{Prop1} to demonstrate the \emph{low rank} properties of multi-fold Hankel matrices, which will establish the correspondence between the Hankel matrices and the complexity of the singularity set of $u$.

\begin{proposition}\label{Prop1} Let $u(\x)$ be defined as in \cref{uModel} with the singularity set $\Gamma=\bigcup_{j=1}^J\p\Om_j$ satisfying \cref{MajorAssumption}, where $\varphi$ is the minimal polynomial with coefficients on $\KK$. For an assumed filter support $\KK'$ strictly containing $\KK$, we have
\begin{align}
\rank\left(\bmH\left(\left(\msF(\na u-p)\right)\big|_{\OO}\right)\right)&\leq|\KK'|-|\KK':\KK|,\label{FirstUpperBound}\\
\rank\left(\bmH\left(\left(\msF(\na_sp)\right)\big|_{\OO}\right)\right)&\leq|\KK'|-|\KK':\KK|.\label{SecondUpperBound}
\end{align}
Hence, both $\bmH\left(\left(\msF(\na u-p)\right)\big|_{\OO}\right)$ and $\bmH\left(\left(\msF(\na_sp)\right)\big|_{\OO}\right)$ are rank deficient.
\end{proposition}

\begin{figure}[ht]
\centering
\subfloat[GSLR framework]{\label{GSLRFramework}\includegraphics[width=1\textwidth]{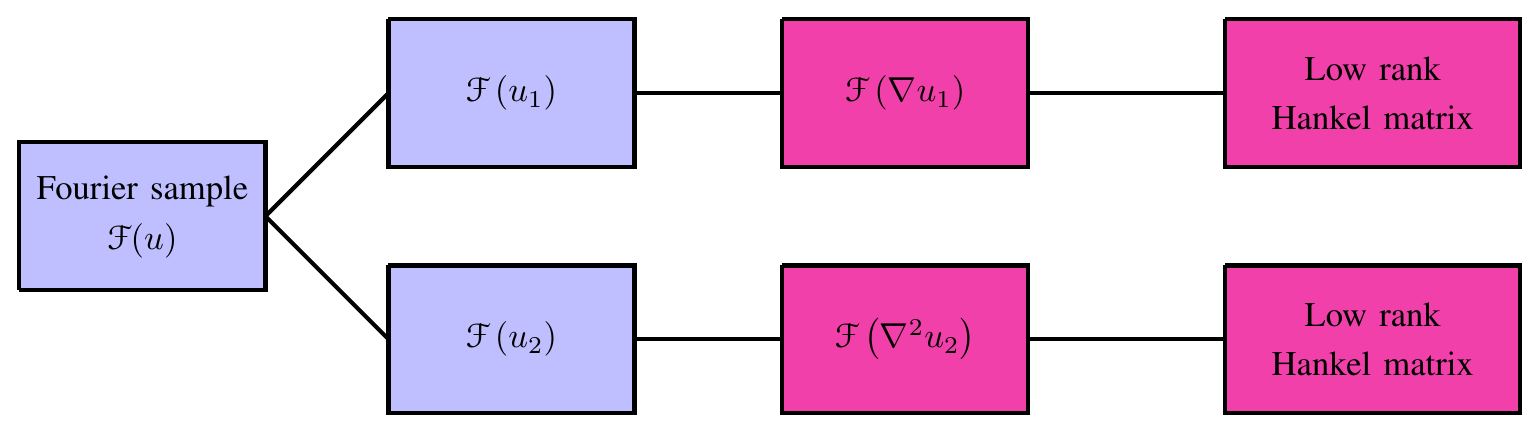}}\\
\subfloat[Proposed SLRM framework]{\label{ProposedSLRM}\includegraphics[width=1\textwidth]{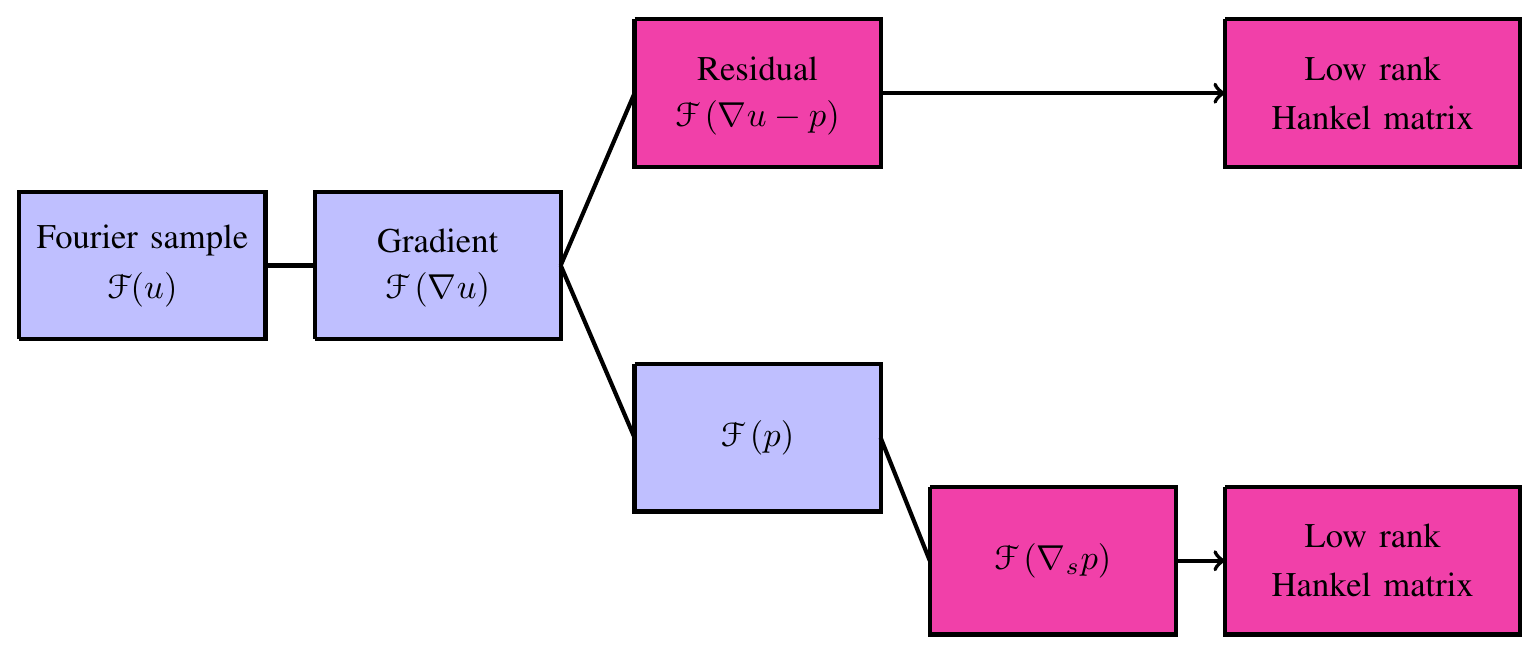}}
\caption{Schematic diagrams for comparison between the GSLR framework in \cite{Y.Hu2019} and the proposed SLRM framework.}\label{SLRMComparison}
\end{figure}

In summary, for $u(\x)$ defined as in \cref{uModel}, the two multi-fold Hankel matrices
\begin{align*}
\bmH\left(\left(\msF(\na u-p)\right)\big|_{\OO}\right)~~~\text{and}~~~\bmH\left(\left(\msF(\na_sp)\right)\big|_{\OO}\right)
\end{align*}
are of low-rank, which enables to convert the piecewise regularity in the continuous domain into the low rank multi-fold Hankel matrices corresponding to the discrete Fourier samples. Notice that the GSLR framework in \cite{Y.Hu2019} promotes the decomposition $u=u_1+u_2$ where the Fourier samples of $\na u_1$ and $\na^2u_2$ correspond to the low rank Hankel matrices. In contrast, the proposed  SLRM framework is established by decomposing the gradient of $u$ into $\na u-p$ and $p$ such that the Fourier samples of $\na u-p$ and $\na_sp$ correspond to the low rank Hankel matrices. See \cref{SLRMComparison} for the schematic illustrations.

\section{Application to image restoration}\label{ImageRestoration}

\subsection{Continuous domain regularization for piecewise smooth image restoration}\label{ImageRestorationModel}

Let $\bsf\in\msV$ be a degraded measurement modeled as
\begin{align}\label{Linear_IP2}
\bsf=\bmA\bsv+\zze,
\end{align}
where $\bsv=\msF(u)\big|_{\OO}$ with $u$ defined as in \cref{uModel}, and $\zze$ is some measurement error\footnote{Here, with a slight abuse of notation, we assume the linear operator $\bmA$ acts on the Fourier samples in what follows.}. According to \cref{OurFramework}, the (multi-fold) Hankel matrices corresponding to
\begin{align*}
\msF(\na u-p)(\xxi)&=\left(2\pi i\xi_1\wh{u}(\xxi)-\wh{p}_1(\xxi),2\pi i\xi_2\wh{u}-\wh{p}_2(\xxi)(\xxi)\right),\\
\msF(\na_sp)(\xxi)&=\left[\begin{array}{cc}
2\pi i\xi_1\wh{p}_1(\xxi)&\pi i\left(\xi_2\wh{p}_1(\xxi)+\xi_1\wh{p}_2(\xxi)\right)\vspace{0.45em}\\
\pi i\left(\xi_2\wh{p}_1(\xxi)+\xi_1\wh{p}_2(\xxi)\right)&2\pi i\xi_2\wh{p}_2(\xxi)
\end{array}\right],
\end{align*}
are low rank. Hence, we can consider
\begin{align}\label{RankMinimization}
\begin{split}
\min_{\bsv,\bq}\rank\left(\bmH\left(\bmD\bsv-\bq\right)\right)+\gamma\rank\left(\bmH\left(\bmE\bq\right)\right)~~~~~\text{subject to}~~\bsf=\bmA\bsv+\zze
\end{split}
\end{align}
to restore $\bsv\in\msV$, as a continuous domain regularization for the piecewise smooth image restoration. Here, $\bq=\msF(p)\big|_{\OO}\in\msV_1$, and $\bmD:\msV\to\msV_1$ and $\bmE:\msV_1\to\msV_2$ are defined as
\begin{align}
\left(\bmD\bsv\right)(\bk)&=\left(2\pi ik_1\bsv(\bk),2\pi ik_2\bsv(\bk)\right),\label{Ddef}\\
\left(\bmE\bq\right)(\bk)&=\left[\begin{array}{cc}
2\pi ik_1\bq_1(\bk)&\pi i\left(k_2\bq_1(\bk)+k_1\bq_2(\bk)\right)\vspace{0.45em}\\
\pi i\left(k_2\bq_1(\bk)+k_1\bq_2(\bk)\right)&2\pi ik_2\bq_2(\bk)
\end{array}\right],\label{Edef}
\end{align}
for $\bk=(k_1,k_2)\in\OO$, respectively.

When $\bq=\0$, \cref{RankMinimization} reduces to
\begin{align}\label{TVRankMinimization}
\min_{\bsv}~\rank\left(\bmH\left(\bmD_1\bsv\right)\right)~~~~~\text{subject to}~~\bsf=\bmA\bsv+\zze
\end{align}
where $\bmD_1=\bmD$. When $\bq=\bmD\bsv$, \cref{RankMinimization} becomes
\begin{align}\label{2ndTVRankMinimization}
\min_{\bsv}~\rank\left(\bmH\left(\bmD_2\bsv\right)\right)~~~~~\text{subject to}~~\bsf=\bmA\bsv+\zze.
\end{align}
where $\bmD_2:\msV\to\msV_2$ is defined as
\begin{align}\label{D2Def}
\left(\bmD_2\bsv\right)(\bk)=\left(\bmE\left(\bmD\bsv\right)\right)(\bk)=\left[\begin{array}{cc}
-4\pi^2k_1^2\bsv(\bk)&-4\pi^2k_1k_2\bsv(\bk)\vspace{0.45em}\\
-4\pi^2k_1k_2\bsv(\bk)&-4\pi^2k_2^2\bsv(\bk)
\end{array}\right]
\end{align}
for $\bk\in\OO$. Finally, when $\bq=\bmD\bsv_2$ for some $\bsv_2\in\msV$, by letting $\bsv_1=\bsv-\bsv_2$, we obtain
\begin{align}\label{GSLRRankMinimization}
\min_{\bsv_1,\bsv_2}~\rank\left(\bmH\left(\bmD_1\bsv_1\right)\right)+\gamma\rank\left(\bmH\left(\bmD_2\bsv_2\right)\right)~~~~~\text{subject to}~~\bsf=\bmA(\bsv_1+\bsv_2)+\zze,
\end{align}
which is a rank minimization model based on the GSLR framework.

From the two extreme cases \cref{TVRankMinimization,2ndTVRankMinimization}, we can see that \cref{RankMinimization} aims to restore piecewise smooth functions by decomposing $\bmD\bsv$ into $\bmD\bsv-\bq$ and $\bq$, thereby balancing the low rank multi-fold Hankel matrices of the Fourier samples of the first order derivatives and the second order derivatives. In addition, we can also see that \cref{RankMinimization} takes the GSLR rank minimization model \cref{GSLRRankMinimization}, which directly decomposes $\bsv=\bsv_1+\bsv_2$, as a special case. Since our SLRM framework is more generalized than the GSLR framework (see remark \ref{RK2} and examples \ref{Ex1} and \ref{Ex2}), it can be expected that \cref{RankMinimization} is able to restore a wider range of piecewise smooth functions than \cref{GSLRRankMinimization}.

\subsection{From low rank model to tight frame approach}\label{ProposedApproach}

Though \cref{RankMinimization} is an NP-hard problem, there are numerous tractable approaches available, including the convex nuclear norm relaxation (e.g. \cite{J.F.Cai2016,M.Fazel2013}), the iterative reweighted least squares (IRLS) for the Schatten $p$-norm minimization \cite{M.Fornasier2011,Y.Hu2019,K.Mohan2012,G.Ongie2017}, etc. Inspired by the SVD of a low rank Hankel matrix, we propose another relaxation of \cref{RankMinimization}, similar to \cite[Theorem 3.2]{J.F.Cai2020} where the continuous domain regularization for the piecewise constant image restoration is studied. For this purpose, we present the main idea of the image restoration model in \cref{Th2}. The proof is postponed to \ref{ProofTh2}.

We begin with introducing some notation. For a set of $K_1\times K_2$ filters $\a_1,\ldots,\a_{M_2}$ supported on $\KK$, we define $\bmW$ and $\bmW^*$ (the adjoint of $\bmW$) as
\begin{align}
\bmW&=\left[\bmS_{\a_1(-\cdot)}^T,\bmS_{\a_2(-\cdot)}^T,\ldots,\bmS_{\a_{M_2}(-\cdot)}^T\right]^T,\label{OurAnalysis}\\
\bmW^*&=\left[\bmS_{\overline{\a}_1},\bmS_{\overline{\a}_2},\ldots,\bmS_{\overline{\a}_{M_2}}\right],\label{OurSynthesis}
\end{align}
where $\bmS_{\a}$ is a discrete convolution under the periodic boundary condition:
\begin{align*}
\left(\bmS_{\a}\bsv\right)(\bk)=\left(\a\ast\bsv\right)(\bk)=\sum_{\bm\in\Z^2}\a(\bk-\bm)\bsv(\bm).
\end{align*}
In other words, both $\bmW$ and $\bmW^*$ are concatenations of discrete convolutions.

\begin{theorem}\label{Th2} Let $u(\x)$ be defined as in \cref{uModel}, and let $\bsv=\msF(u)\big|_{\OO}\in\msV$ and $\bq=\msF(p)\big|_{\OO}\in\msV_1$ be the Fourier samples. Let $\bmH\left(\bmD\bsv-\bq\right)\in\C^{2M_1\times M_2}$ and $\bmH\left(\bmE\bq\right)\in\C^{4M_1\times M_2}$ be multi-fold Hankel matrices with $\bmD$ and $\bmE$ defined as in \cref{Ddef,Edef}. Assume that $\bmH\left(\bmD\bsv-\bq\right)$ and $\bmH\left(\bmE\bq\right)$ satisfy
\begin{align}\label{Assumption}
\rank\left(\bmH\left(\bmD\bsv-\bq\right)\right)=r_1\ll2M_1\wedge M_2~~\text{and}~~\rank\left(\bmH\left(\bmE\bq\right)\right)=r_2\ll4M_1\wedge M_2.
\end{align}
Considering full SVDs $\bmH\left(\bmD\bsv-\bq\right)=\bX_1\Sig_1\bY_1^*$ and $\bmH\left(\bmE\bq\right)=\bX_2\Sig_2\bY_2^*$, we define $\a_{1l}=M_2^{-1/2}\bY_1^{(:,l)}$ and $\a_{2l}=M_2^{-1/2}\bY_2^{(:,l)}$ by reformulating each column vector into a $K_1\times K_2$ filter supported on $\KK$. Then $\bmW_1$ and $\bmW_2$ defined as \cref{OurAnalysis} by using filters $\left\{\a_{11},\ldots,\a_{1M_2}\right\}$ and $\left\{\a_{21},\ldots,\a_{2M_2}\right\}$ satisfies
\begin{align}
\bmW_1^*\bmW_1\left(\bmD\bsv-\bq\right)&=\sum_{l=1}^{M_2}\bmS_{\overline{\a}_{1l}}\left(\bmS_{\a_{1l}(-\cdot)}\left(\bmD\bsv-\bq\right)\right)=\bmD\bsv-\bq,\label{LowRankHankelTightFrame1}\\
\bmW_2^*\bmW_2\left(\bmE\bq\right)&=\sum_{l=1}^{M_2}\bmS_{\overline{\a}_{2l}}\left(\bmS_{\a_{2l}(-\cdot)}\left(\bmE\bq\right)\right)=\bmE\bq,\label{LowRankHankelTightFrame2}
\end{align}
and for $\bk\in\OO:\KK$, we have
\begin{align}
\left(\bmS_{\a_{1l}}\left(\bmD\bsv-\bq\right)\right)(\bk)&=\0,~~~~~l=r_1+1,\ldots,M_2\label{TightFrameSparse1}\\
\left(\bmS_{\a_{2l}}\left(\bmE\bq\right)\right)(\bk)&=\0,~~~~~l=r_2+1,\ldots,M_2\label{TightFrameSparse2},
\end{align}
where the discrete convolution is performed on each component of $\bmD\bsv-\bq$ and $\bmE\bq$, respectively. Consequently, if $\bmH\left(\bmD\bsv-\bq\right)$ (and $\bmH\left(\bmE\bq\right)$, respectively) is of low rank, then its right singular vectors construct a tight frame under which $\bmD\bsv-\bq$ (and $\bmE\bq$, respectively) is sparsely represented.
\end{theorem}

In words, \cref{Th2} tells us that if the SVDs of multi-fold Hankel matrices are known as oracles, we can explicitly construct tight frames under which $\bmD\bsv-\bq$ and $\bmE\bq$ are sparsely represented, and the sparsity of canonical coefficients can be grouped according to the filters. Hence, motivated by the idea in \cite{D.Guo2018}, assume that $\wt{\bsv}\in\mV$ and $\wt{\bq}\in\mV_1$ are the a-priori estimations of $\bsv$ and $\bq$ with the SVDs
\begin{align*}
\bmH\left(\bmD\wt{\bsv}-\wt{\bq}\right)=\wt{\bX}_1\wt{\Sig}_1\wt{\bY}_1^*~~\text{and}~~\bmH\left(\bmE\wt{\bq}\right)=\wt{\bX}_2\wt{\Sig}_2\wt{\bY}_2^*.
\end{align*}
Then we define the tight frame transforms $\bmW_1$ and $\bmW_2$ in \cref{OurAnalysis} via
\begin{align*}
\a_{1l}=M_2^{-1/2}\wt{\bY}_1^{(:,l)}~~~\text{and}~~~\a_{2l}=M_2^{-1/2}\wt{\bY}_2^{(:,l)},~~~~l=1,\ldots,M_2.
\end{align*}
Since \cref{TightFrameSparse1,TightFrameSparse2} are then approximately true under these $\bmW_1$ and $\bmW_2$, we remove the group sparsity pattern in the canonical coefficients for the better sparse approximation instead. This leads us to solve
\begin{align}\label{ProposedModel}
\min_{\bsv,\bq}\f{1}{2}\left\|\bmA\bsv-\bsf\right\|_2^2+\left\|\gga_1\cdot\bmW_1\left(\bmD\bsv-\bq\right)\right\|_1+\left\|\gga_2\cdot\bmW_2\left(\bmE\bq\right)\right\|_1,
\end{align}
where the $\ell_1$ norms take the form of
\begin{align*}
\left\|\gga_1\cdot\bmW_1\left(\bmD\bsv-\bq\right)\right\|_1&=\sum_{l=1}^{M_2}\gamma_{1l}\left\|\bmS_{\a_{1l}(-\cdot)}\left(\bmD\bsv-\bq\right)\right\|_1\\
\left\|\gga_2\cdot\bmW_2\left(\bmE\bq\right)\right\|_1&=\sum_{l=1}^{M_2}\gamma_{2l}\left\|\bmS_{\a_{2l}(-\cdot)}\left(\bmE\bq\right)\right\|_1
\end{align*}
to reflect the different weights according to the index of filters. Finally, we reflect the singular values of the pre-restored Hankel matrices in the regularization parameters as
\begin{align*}
\gamma_{kl}=\f{\nu_k}{\wt{\Sig}_k^{(l,l)}+\eps}~~~~~k=1,2~~\text{and}~~l=1,\ldots,M_2
\end{align*}
with some $\nu_k>0$ and a small $\eps>0$ to avoid the division by zero. Hence, we relax the sparsity of tight frame coefficients over the entire range of $\bmW_k$'s (not necessarily in groups as in \cref{TightFrameSparse1,TightFrameSparse2}), we expect to achieve more flexibility, thereby leading to the improvements in restoration performance.

Since the wavelet frame based relaxation model \cref{ProposedModel} is inspired by the \cref{RankMinimization} via the SVDs of $\bmH\left(\bmD\bsv-\bq\right)$ and $\bmH\left(\bmE\bq\right)$, we can say that our relaxation model \cref{ProposedModel} restores piecewise smooth functions by decomposing $\bmD\bsv$ into $\bmD\bsv-\bq$ and $\bq$, thereby balancing the low rank Hankel matrices of the Fourier samples of the first order derivatives and that of the second order derivatives. In addition, we further note that, following \cite{J.F.Cai2020}, it is also possible to consider the following data driven tight frame model
\begin{align}\label{PSDDTFModel}
\begin{split}
&~~~\min_{\bsv,\bq,\left\{\bc_j,\bmW_j\right\}_{j=1}^2}\f{1}{2}\left\|\bmA\bsv-\bsf\right\|_2^2+\f{\mu_1}{2}\left\|\bmW_1\left(\bmD\bsv-\bq\right)-\bc_1\right\|_2^2+\left\|\gga\cdot\bc_1\right\|_0\\
&\hspace{15.00em}+\f{\mu_2}{2}\left\|\bmW_2\left(\bmE\bq\right)-\bc_2\right\|_2^2+\left\|\gga_2\cdot\bc_2\right\|_0\\
&~~~~~\text{subject to}~~\bmW_1^*\bmW_1=\bmW_2^*\bmW_2=\bmI,
\end{split}
\end{align}
with the $\ell_0$ norm $\|\gga_k\cdot\bc_k\|_0$ ($k=1,2$) encoding the number of nonzero entries in $\bc_k$'s, to learn tight frames and restore the Fourier sample $\bsv$ simultaneously. Even though it is not clear at this point whether the adaptive tight frame system will lead to better restoration results or not, throughout this paper, we only consider the model \cref{ProposedModel} rather than the data driven tight frame model \cref{PSDDTFModel}, and the reasons are as follows. First of all, given that $\bmW_k$'s are properly estimated, it may not be necessary to further learn them with additional computational costs. Second, since \cref{ProposedModel} is convex whereas \cref{PSDDTFModel} is nonconvex, we can easily expect better behavior and theoretical support for the numerical algorithms. Most importantly, from the viewpoint of wavelet frame based image restoration, the model \cref{ProposedModel} is an \emph{analysis approach} \cite{J.F.Cai2009/10} while the data driven tight frame model \cref{PSDDTFModel} can be classified into the \emph{balanced approach} \cite{J.F.Cai2008,R.H.Chan2003}. (See \ref{PreliminariesTightFrame} for the brief descriptions on the wavelet frame based models.) It is well known that the analysis approach reflects the structure of a target image better than other approaches (e.g. \cite{J.F.Cai2012}). Since we derive the wavelet frame based approach from the structured low rank matrix framework, it can be expected that the convex relaxation model \cref{ProposedModel} will reflect our SLRM frameworks for the piecewise smooth functions better than the data driven tight frame model \cref{PSDDTFModel}.

\subsection{Alternating minimization algorithm}\label{AlternatingMinimizationAlgorithm}

Among numerous algorithms which can solve the convex model \cref{ProposedModel}, we adopt the ADMM \cite{J.Eckstein1992} or the split Bregman algorithm \cite{W.Guo2014}, which can convert \cref{ProposedModel} into several subproblems with closed form solutions, together with the convergence guarantee \cite{J.F.Cai2009/10}. More precisely, let $\bc_1=\bmW_1\left(\bmD\bsv-\bq\right)$, and $\bc_2=\bmW_2\left(\bmE\bq\right)$. Then \cref{ProposedModel} can be rewritten as
\begin{align*}
&\min_{\bsv,\bq,\bc_1,\bc_2}\f{1}{2}\left\|\bmA\bsv-\bsf\right\|_2^2+\left\|\gga_1\cdot\bc_1\right\|_1+\left\|\gga_2\cdot\bc_2\right\|_1\\
&\text{subject to}~~~\bc_1=\bmW_1\left(\bmD\bsv-\bq\right),~~\text{and}~~\bc_2=\bmW_2\left(\bmE\bq\right).
\end{align*}
Under this reformulation, the overall algorithm is summarized in \cref{Alg1}.

\begin{algorithm}[t!]
\begin{algorithmic}
\STATE{\textbf{Initialization:} $\bsv^{0}$, $\bq^{0}$, $\bc_{1}^{0}$, $\bc_{2}^{0}$, $\bsd_1^0$, $\bsd_2^0$}
\FOR{$n=0$, $1$, $2$, $\cdots$}
\STATE{\textbf{(1)} Update $\bsv$ and $\bq$:
\begin{align}
\left[\begin{array}{c}\bsv^{n+1}\\
\bq^{n+1}\end{array}\right]&=\argmin_{\bsv,\bq}\f{1}{2}\left\|\bmA\bsv-\bsf\right\|_2^2+\f{\beta}{2}\left\|\bmW_1\left(\bmD\bsv-\bq\right)-\bc_1^n+\bsd_1^n\right\|_2^2+\f{\beta}{2}\left\|\bmW_2\left(\bmE\bq\right)-\bc_2^n+\bsd_2^n\right\|_2^2\label{vqsubprob}
\end{align}
\textbf{(2)} Update $\bc_1$ and $\bc_2$:
\begin{align}
\bc_1^{n+1}&=\argmin_{\bc_1}\left\|\gga_1\cdot\bc_1\right\|_1+\f{\beta}{2}\left\|\bc_1-\bmW_1\left(\bmD\bsv^{n+1}-\bq^{n+1}\right)-\bsd_1^n\right\|_2^2\label{c1sub}\\
\bc_2^{n+1}&=\argmin_{\bc_2}\left\|\gga_2\cdot\bc_2\right\|_1+\f{\beta}{2}\left\|\bc_2-\bmW_2\left(\bmE\bq^{n+1}\right)-\bsd_2^n\right\|_2^2\label{c2sub}
\end{align}
\textbf{(3)} Update $\bsd_1$ and $\bsd_2$:
\begin{align}
\bsd_1^{n+1}&=\bsd_1^n+\bmW_1\left(\bmD\bsv^{n+1}-\bq^{n+1}\right)-\bc_1^{n+1}\\
\bsd_2^{n+1}&=\bsd_2^n+\bmW_2\left(\bmE\bq^{n+1}\right)-\bc_2^{n+1}.
\end{align}}
\ENDFOR
\end{algorithmic}
\caption{Split Bregman Algorithm for \cref{ProposedModel}}\label{Alg1}
\end{algorithm}

For \cref{vqsubprob}, since $\bmW_k^*\bmW_k=\bmI$ for $k=1,2$, we solve the following system of linear equations
\begin{align}\label{vqsystem}
\left[\begin{array}{cc}
\bmK_{11}&\bmK_{21}^*\\
\bmK_{21}&\bmK_{22}
\end{array}\right]\left[\begin{array}{c}
\bsv\\
\bq
\end{array}\right]=\left[\begin{array}{c}
\bsf_1\\
\bsf_2
\end{array}\right]
\end{align}
where
\begin{align*}
\begin{array}{rl}
\bmK_{11}\hspace{-0.8em}&=\bmA^*\bmA+\beta\bmD^*\bmD\\
\bmK_{21}\hspace{-0.8em}&=-\beta\bmD\\
\bmK_{22}\hspace{-0.8em}&=\beta\bmI+\beta\bmE^*\bmE
\end{array}~~~\text{and}~~~\begin{array}{rl}
\bsf_1\hspace{-0.8em}&=\bmA^*\bsf+\beta\bmD^*\left[\bmW_1^*\left(\bc_1^n-\bsd_1^n\right)\right]\vspace{0.88em}\\
\bsf_2\hspace{-0.8em}&=-\beta\bmW_1^*\left(\bc_1^n-\bsd_1^n\right)+\beta\bmE^*\left[\bmW_2^*\left(\bc_2^n-\bsd_2^n\right)\right].
\end{array}
\end{align*}
Depending on the formulation of $\bmA$, various methods can be used to solve \cref{vqsystem} efficiently. For example, when $\bmA$ is a pointwise multiplication in the frequency domain (e.g. image denoising, image deblurring, and CS restoration), so are the constituent operators $\bmK_{11},\bmK_{21},$ and $\bmK_{22}$, and thus we can solve \cref{vqsystem} by the pointwise Cramer's rule. For a more generalized $\bmA$, we can use a distributed optimization based method \cite{S.Boyd2011}.

The closed form solutions for \cref{c1sub,c2sub} is expressed in terms of the soft thresholding:
\begin{align}
\bc_1^{n+1}&=\bmT_{\gga_1/\beta}\left(\bmW_1\left(\bmD\bsv^{n+1}-\bq^{n+1}\right)+\bsd_1^n\right),\label{c1explicit}\\
\bc_2^{n+1}&=\bmT_{\gga_2/\beta}\left(\bmW_2\left(\bmE\bq^{n+1}\right)+\bsd_2^n\right).\label{c2explicit}
\end{align}
More precisely, the soft thresholding operator $\bmT_{\gga}\left(\bc\right)$ for $\bc\in\msV_k^{M_2}$ and $\gga=\left[\begin{array}{ccc}
\gamma_1&\cdots&\gamma_{M_2}
\end{array}\right]^T$ is defined as the following componentwise manner:
\begin{align*}
\bmT_{\gga}\left(\bc\right)_{l,m}(\bk)=\max\left\{\left|\bc_{l,m}(\bk)\right|-\gamma_l,0\right\}\f{\bc_{l,m}(\bk)}{|\bc_{l,m}(\bk)|}
\end{align*}
for $\bk\in\OO$, $l=1,\ldots,M_2$, and $m=1,\ldots,2^k$, with the convention that $0/0=0$.

\section{Numerical results}\label{Experiments}

In this section, we conduct some numerical simulations in the context of restoration from the partial Fourier samples as a proof-of-concept study.  Specifically, to compare the performance of the direct rank minimization \cref{RankMinimization} and the wavelet frame relaxation \cref{ProposedModel}, we choose to compare the relaxation model
\begin{align}\label{ProposedCSMRI}
\min_{\bsv,\bq}\f{1}{2}\left\|\bmR_{\MM}\bsv-\bsf\right\|_2^2+\left\|\gga_1\cdot\bmW_1\left(\bmD\bsv-\bq\right)\right\|_1+\left\|\gga_2\cdot\bmW_2\left(\bmE\bq\right)\right\|_1
\end{align}
solved by \cref{Alg1} with the following Schatten $0$-norm relaxation \cite{M.Fornasier2011,Y.Hu2019,K.Mohan2012,G.Ongie2017} (SLRM model) of \cref{RankMinimization}:
\begin{align}\label{LRHTGVIRLSCSMRI}
\min_{\bsv,\bq}\f{1}{2}\left\|\bmR_{\MM}\bsv-\bsf\right\|_2^2&+\gamma_1\left\|\bmH\left(\bmD\bsv-\bq\right)\right\|_0+\gamma_2\left\|\bmH\left(\bmE\bq\right)\right\|_0,
\end{align}
where $\bmR_{\MM}$ denotes a restriction onto the known sample grid $\MM$. In addition, since it has been demonstrated in \cite{Y.Hu2019} that the GSLR framework outperforms the structured low rank matrix approaches which consider either the first order derivatives or the second derivatives, we also compare with the following GSLR model in \cite{Y.Hu2019}:
\begin{align}\label{LRHInfConvIRLSCSMRI}
\min_{\bsv_1,\bsv_2}\f{1}{2}\left\|\bmR_{\MM}\left(\bsv_1+\bsv_2\right)-\bsf\right\|_2^2&+\gamma_1\left\|\bmH\left(\bmD_1\bsv_1\right)\right\|_0+\gamma_2\left\|\bmH\left(\bmD_2\bsv_2\right)\right\|_0.
\end{align}
In \cref{LRHTGVIRLSCSMRI,LRHInfConvIRLSCSMRI}, $\left\|\bZ\right\|_0$ is the Schatten $0$-norm of a matrix $\bZ$ defined as
\begin{align*}
\left\|\bZ\right\|_0=\ln\det\left(\left(\bZ^*\bZ\right)^{1/2}+\eps\bsI\right)
\end{align*}
with a small constant $\eps>0$, and \cref{LRHTGVIRLSCSMRI,LRHInfConvIRLSCSMRI} are solved by the iterative reweighted least squares method described in \cite{Y.Hu2019}.

Finally, to further study the improvements over the conventional on-the-grid approaches, we compare with the piecewise linear framelet (Fra) model (e.g. \cite{J.F.Cai2009/10})
\begin{align}\label{FrameCSMRI}
\min_{\bu}\f{1}{2}\left\|\bmR_{\MM}\bmsF\bu-\bsf\right\|_2^2+\left\|\gga\cdot\bmW\bu\right\|_1,
\end{align}
the TGV model \cite{K.Bredies2010,F.Knoll2011}
\begin{align}\label{TGVCSMRI}
\min_{\bu,\bp}\f{1}{2}\left\|\bmR_{\MM}\bmsF\bu-\bsf\right\|_2^2+\gamma_1\left\|\Na\bu-\bp\right\|_1+\gamma_2\left\|\Na_s\bp\right\|_1,
\end{align}
and the inf-convolution (IF) model \cite{A.Chambolle1997}
\begin{align}\label{InfConvCSMRI}
\min_{\bu_1,\bu_2}\f{1}{2}\left\|\bmR_{\MM}\bmsF\left(\bu_1+\bu_2\right)-\bsf\right\|_2^2+\gamma_1\left\|\Na\bu_1\right\|_1+\gamma_2\left\|\Na^2\bu_2\right\|_1,
\end{align}
where $\bmsF$ denotes the two dimensional discrete Fourier transform. Throughout this paper, we use the split Bregman algorithm (e.g. \cite{J.F.Cai2009/10,T.Goldstein2009,W.Guo2014}) to solve \cref{FrameCSMRI,TGVCSMRI,InfConvCSMRI}. All experiments are implemented on MATLAB $\mathrm{R}2014\mathrm{a}$ running on a laptop with $64\mathrm{GB}$ RAM and Intel(R) Core(TM) CPU $\mathrm{i}7$-$8750\mathrm{H}$ at $2.20\mathrm{GHz}$ with $6$ cores.

\begin{figure}[t]
\centering
\subfloat[Ellipses]{\label{Logan}\includegraphics[width=3.50cm]{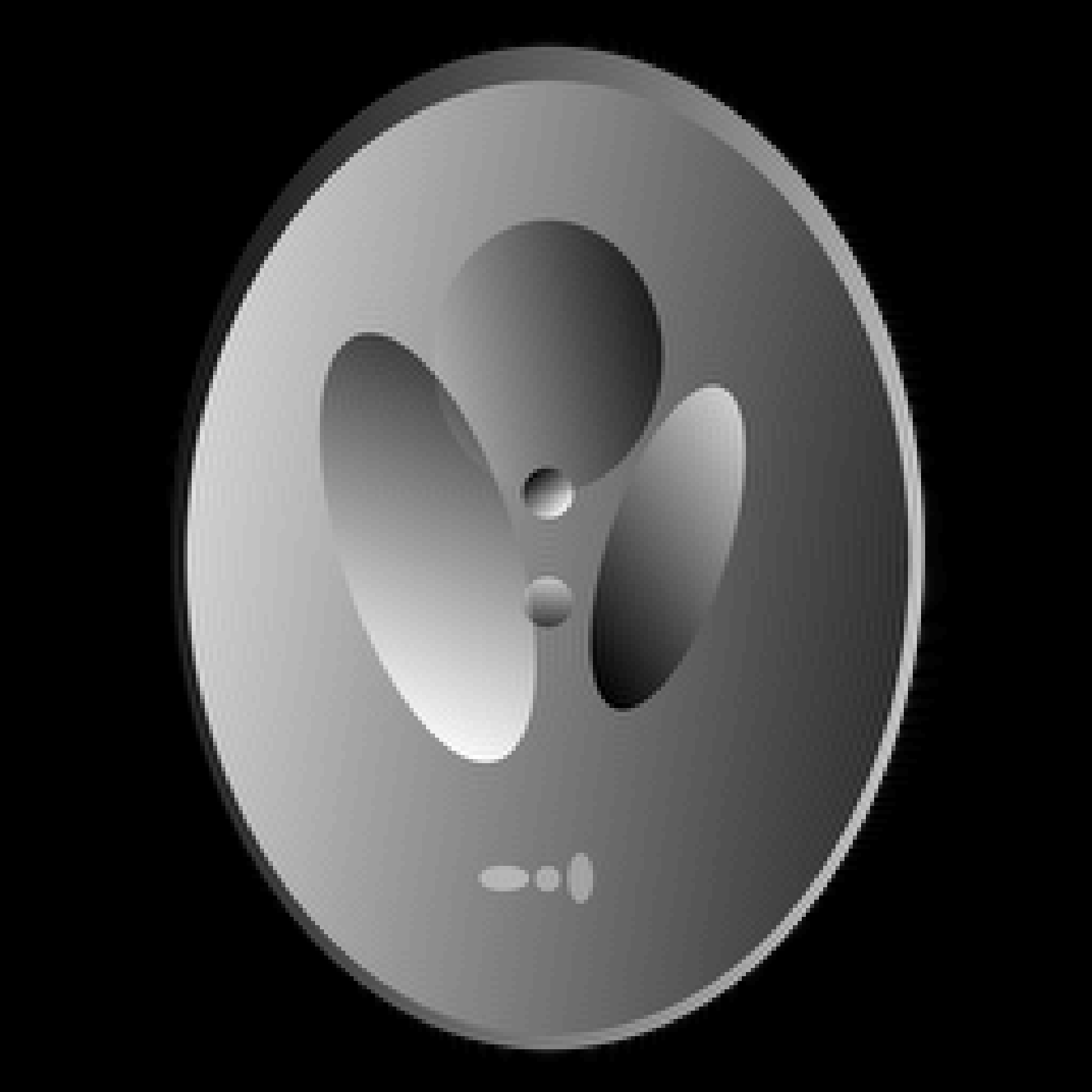}}\hspace{0.001cm}
\subfloat[Rectangles]{\label{Rectangle}\includegraphics[width=3.50cm]{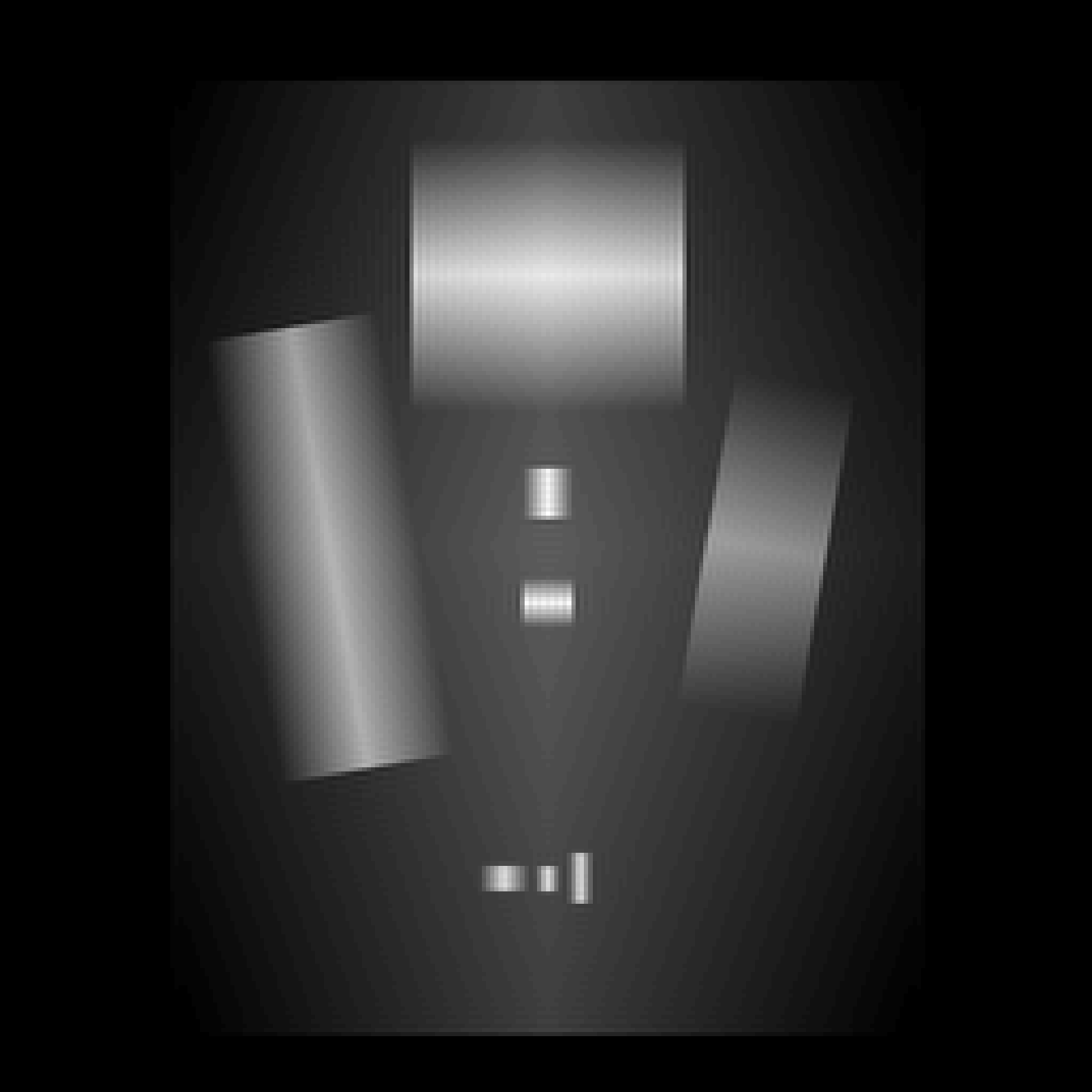}}\hspace{0.001cm}
\subfloat[Airplane]{\label{Airplane}\includegraphics[width=3.50cm]{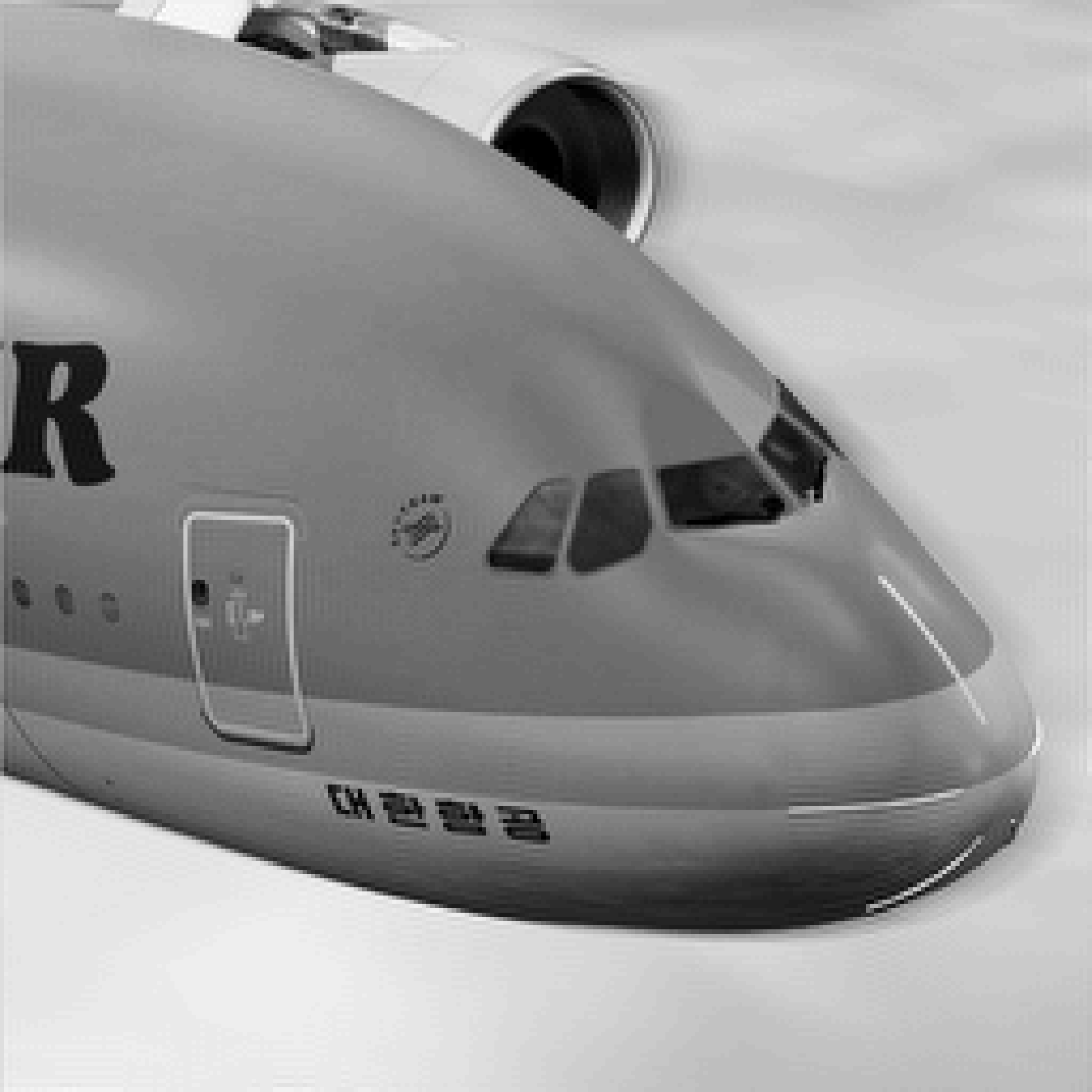}}\hspace{0.001cm}
\subfloat[Car]{\label{Car}\includegraphics[width=3.50cm]{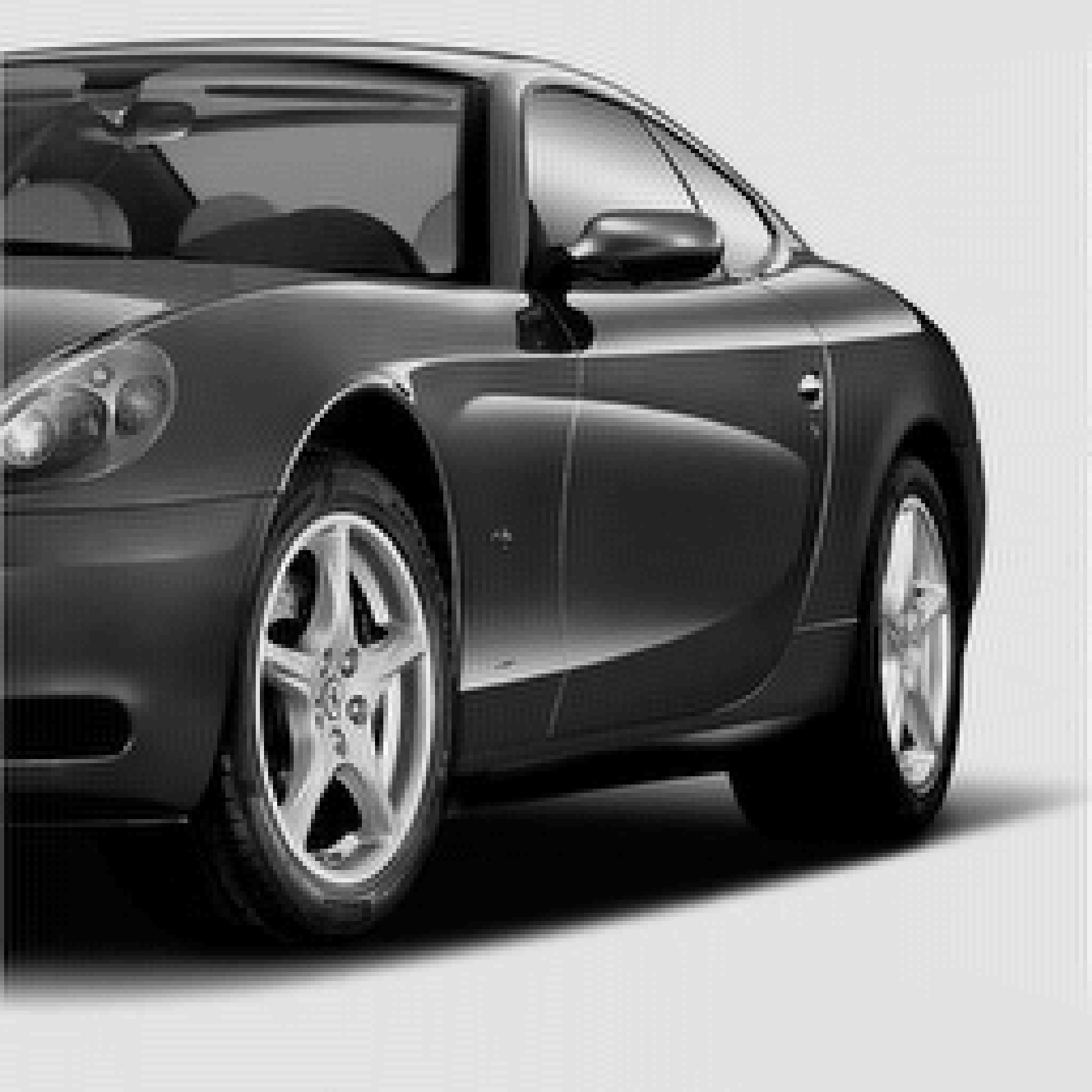}}\vspace{-0.25cm}\\
\subfloat[Ellipses-observed]{\label{LoganZeroPad}\includegraphics[width=3.50cm]{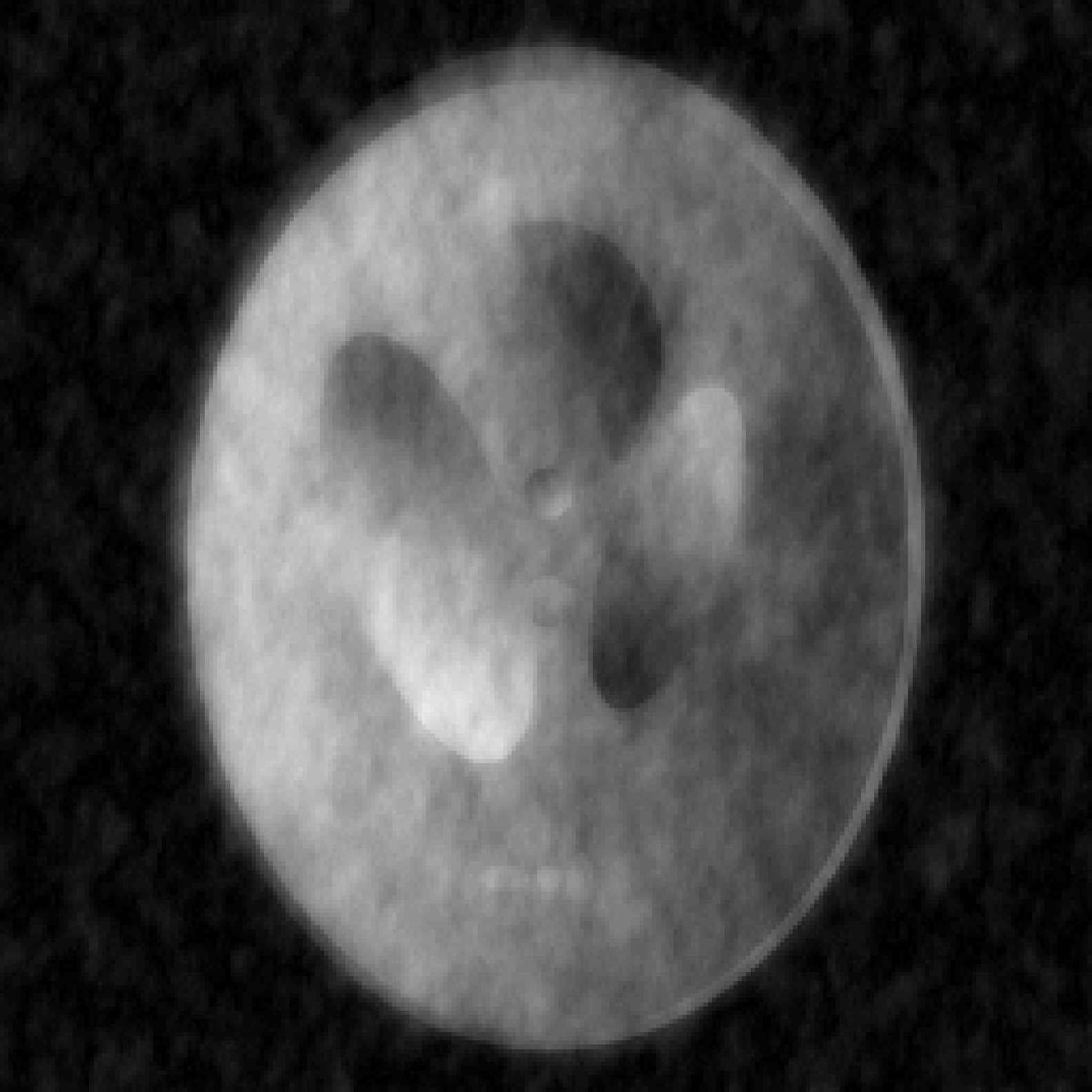}}\hspace{0.001cm}
\subfloat[Rectangles-observed]{\label{RectangleZeroPad}\includegraphics[width=3.50cm]{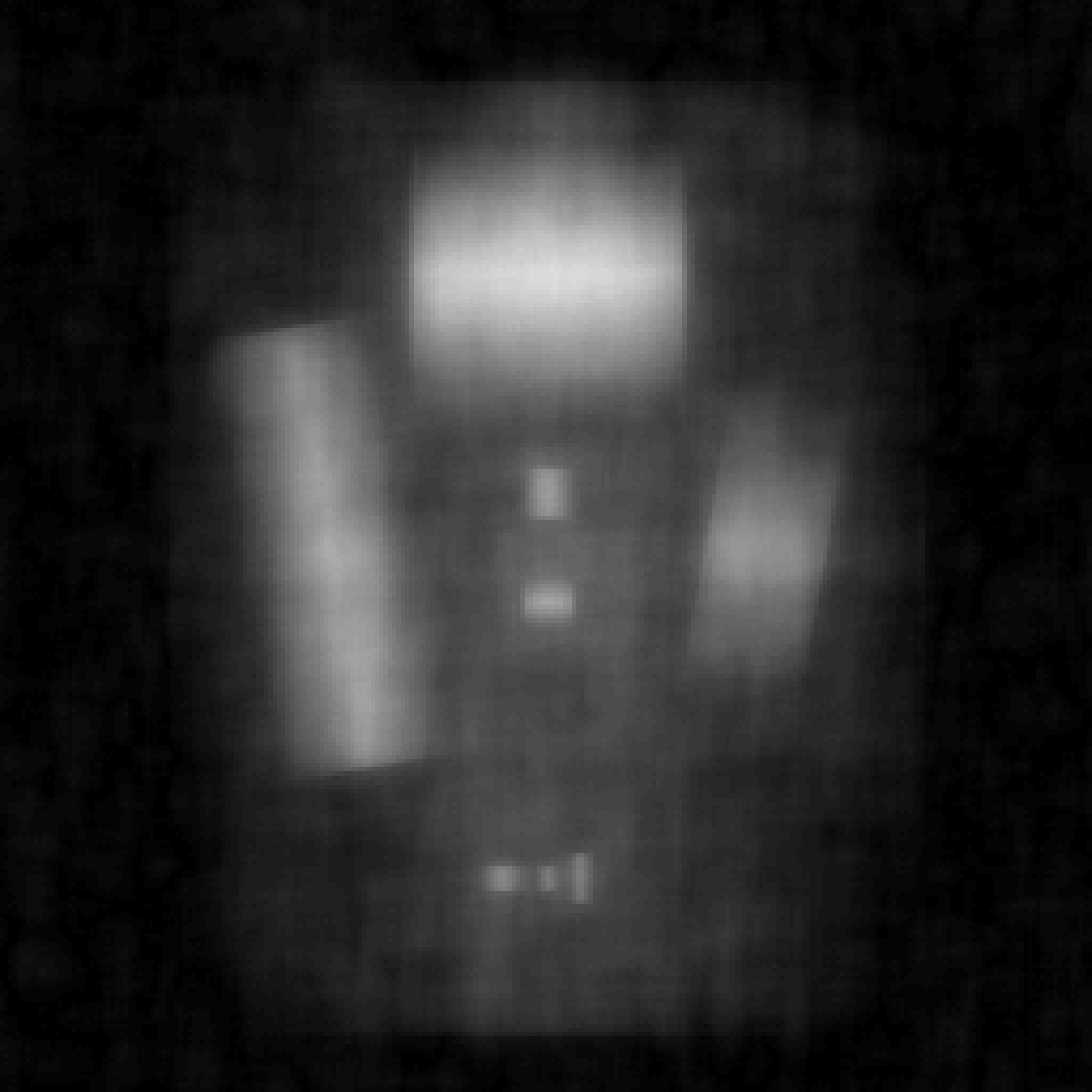}}\hspace{0.001cm}
\subfloat[Airplane-observed]{\label{AirplaneZeroPad}\includegraphics[width=3.50cm]{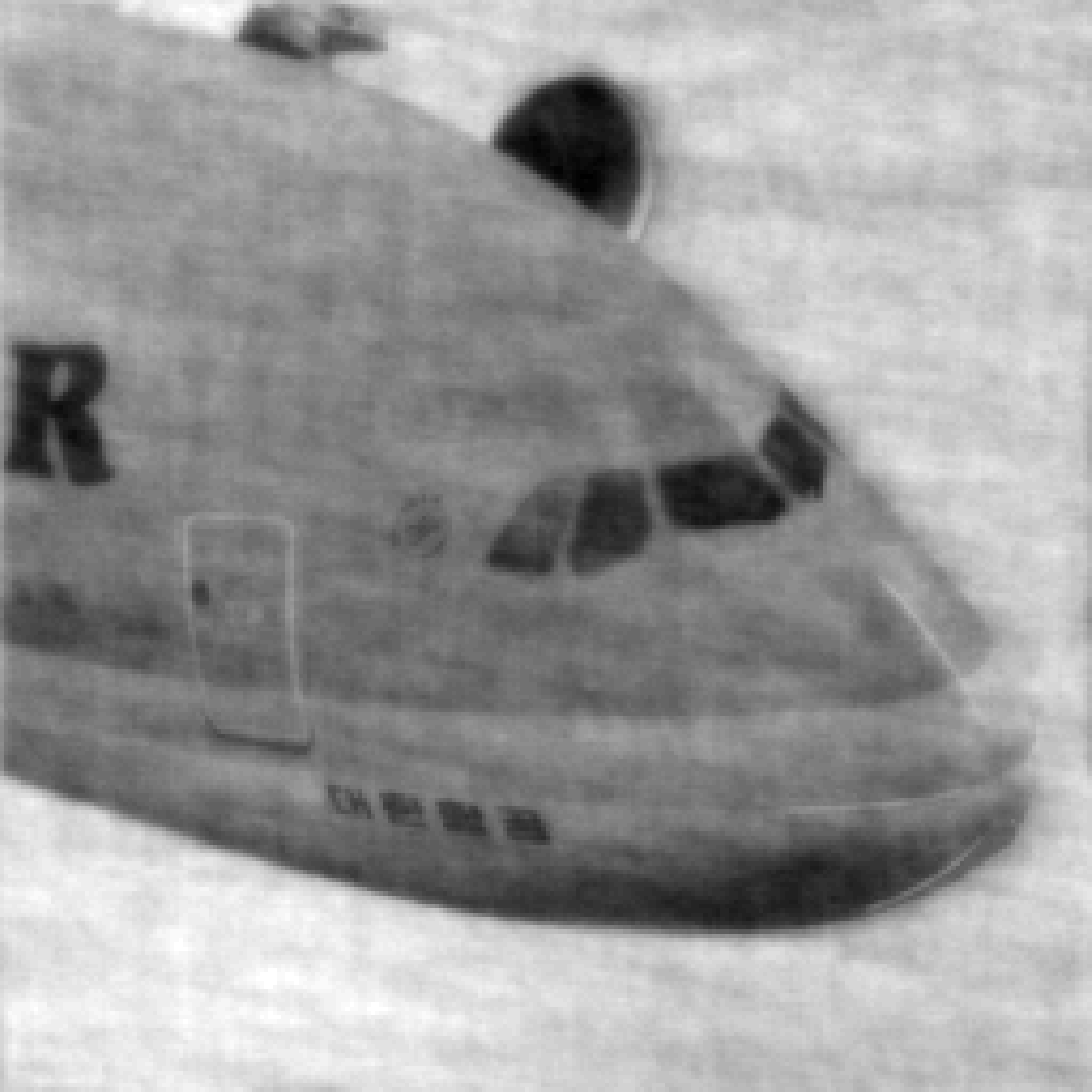}}\hspace{0.001cm}
\subfloat[Car-observed]{\label{CarZeroPad}\includegraphics[width=3.50cm]{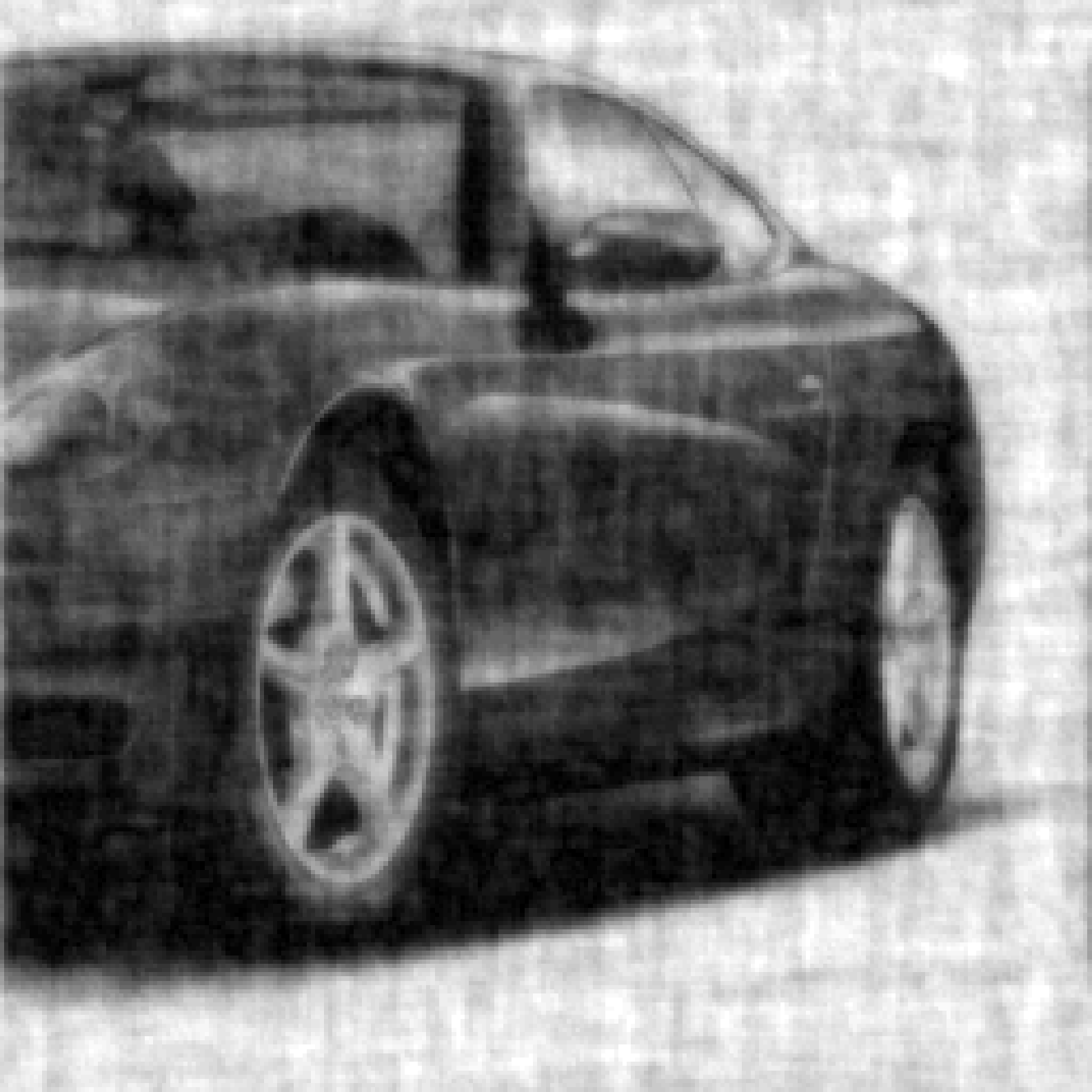}}
\caption{Visualization of test images and observed images. All images are displayed in the window level $[0,1]$ for the fair comparisons.}\label{OriginalImages}
\end{figure}

Throughout the experiments, we test two synthetic images (``Ellipses'' and ``Rectangles''), and two natural images (``Airplane'', and ``Car'') taking the values in $[0,1]$, as shown in \cref{OriginalImages}. The data $\bsf$ is synthesized by randomly sampling $20\%$ of $256\times256$ Fourier samples via the variable density sampling method described in \cite{M.Lustig2007}. A complex white Gaussian noise with the standard deviation $1$ is also added to generate a noisy partial sampling. For \cref{ProposedCSMRI,LRHTGVIRLSCSMRI,LRHInfConvIRLSCSMRI}, we use the $K\times K$ square patch for simplicity. Specifically, we choose $K=31$ for the ``Ellipses'', $K=25$ for the ``Rectangles'', and $K=51$ for ``Airplane'' and ``Car'', depending on the geometry of the target images. More precisely, in \cref{ProposedCSMRI,LRHTGVIRLSCSMRI,LRHInfConvIRLSCSMRI}, we use the discrete Fourier transform of the images restored by the TGV model \cref{TGVCSMRI} and the inf-convolution model \cref{InfConvCSMRI} respectively, to compute the SVDs of multi-fold Hankel matrices using the $K\times K$ patches, which will be used for $\bmW_k$'s in \cref{ProposedCSMRI} and initializations in \cref{LRHTGVIRLSCSMRI,LRHInfConvIRLSCSMRI}. For \cref{FrameCSMRI}, we choose $\bmW$ to be the undecimated tensor product piecewise linear B-spline framelet transform with $1$ level of decomposition \cite{B.Dong2013}. For \cref{TGVCSMRI,InfConvCSMRI}, we use the forward difference with the periodic boundary condition for the difference operators $\Na$, $\Na_s$, and $\Na^2$. In all of the experiments, we have manually tuned the regularization parameters to achieve the optimal restoration results in each scenario. For the quantitative comparison, we compute the signal-to-noise ratio (SNR), the high frequency error norm (HFEN) \cite{S.Ravishankar2011}, and the structure similarity index map (SSIM) \cite{Z.Wang2004}. Note that for \cref{ProposedCSMRI,LRHTGVIRLSCSMRI,LRHInfConvIRLSCSMRI}, the restored image is computed via the inverse DFT of the restored Fourier samples.

\begin{table}[t]
\centering
\begin{tabular}{|c|c|c|c|c|c|c|c|c|c|}
\hline
Images&Indices&Zero fill&Model \cref{ProposedCSMRI}&SLRM \cref{LRHTGVIRLSCSMRI}&GSLR \cref{LRHInfConvIRLSCSMRI}&Fra \cref{FrameCSMRI}&TGV \cref{TGVCSMRI}&IF \cref{InfConvCSMRI}\\ \hline
\multirow{3}{*}{Ellipses}&SNR&$13.75$&$\textbf{36.13}$&$32.18$&$30.93$&$27.41$&$29.67$&$28.36$\\ \cline{2-9}
&HFEN&$0.7746$&$\textbf{0.0375}$&$0.0713$&$0.0789$&$0.1533$&$0.1029$&$0.1351$\\ \cline{2-9}
&SSIM&$0.4162$&$\textbf{0.9943}$&$0.9784$&$0.9746$&$0.8966$&$0.9547$&$0.9335$\\ \hline
\multirow{3}{*}{Rectangles}&SNR&$15.68$&$\textbf{35.00}$&$30.74$&$30.42$&$24.30$&$29.28$&$27.43$\\ \cline{2-9}
&HFEN&$0.7914$&$\textbf{0.0876}$&$0.1456$&$0.1770$&$0.4444$&$0.1878$&$0.2366$\\ \cline{2-9}
&SSIM&$0.6130$&$\textbf{0.9803}$&$\textbf{0.9815}$&$0.9742$&$0.9391$&$0.9552$&$0.9211$\\ \hline
\multirow{3}{*}{Airplane}&SNR&$21.55$&$\textbf{36.23}$&$33.84$&$33.12$&$31.24$&$33.40$&$32.92$\\ \cline{2-9}
&HFEN&$0.4402$&$\textbf{0.0530}$&$0.0724$&$0.0753$&$0.1135$&$0.0803$&$0.0889$\\ \cline{2-9}
&SSIM&$0.7386$&$\textbf{0.9822}$&$0.9706$&$0.9693$&$0.9564$&$0.9669$&$0.9638$\\ \hline
\multirow{3}{*}{Car}&SNR&$15.88$&$\textbf{27.42}$&$25.88$&$24.81$&$22.72$&$25.69$&$25.16$\\ \cline{2-9}
&HFEN&$0.5528$&$\textbf{0.1057}$&$0.1447$&$0.1460$&$0.2290$&$0.1530$&$0.1669$\\ \cline{2-9}
&SSIM&$0.5627$&$\textbf{0.9569}$&$0.9291$&$0.9282$&$0.8761$&$0.9218$&$0.9112$\\ \hline
\end{tabular}
\caption{Comparison of SNRs, HFENs, and SSIMs.}\label{TableResults}
\end{table}

\begin{figure}[t]
\centering
\subfloat[]{\label{LoganOriginal}\includegraphics[width=2.25cm]{LoganOriginal.pdf}}\hspace{0.001cm}
\subfloat[]{\label{LoganLRHTGV}\includegraphics[width=2.25cm]{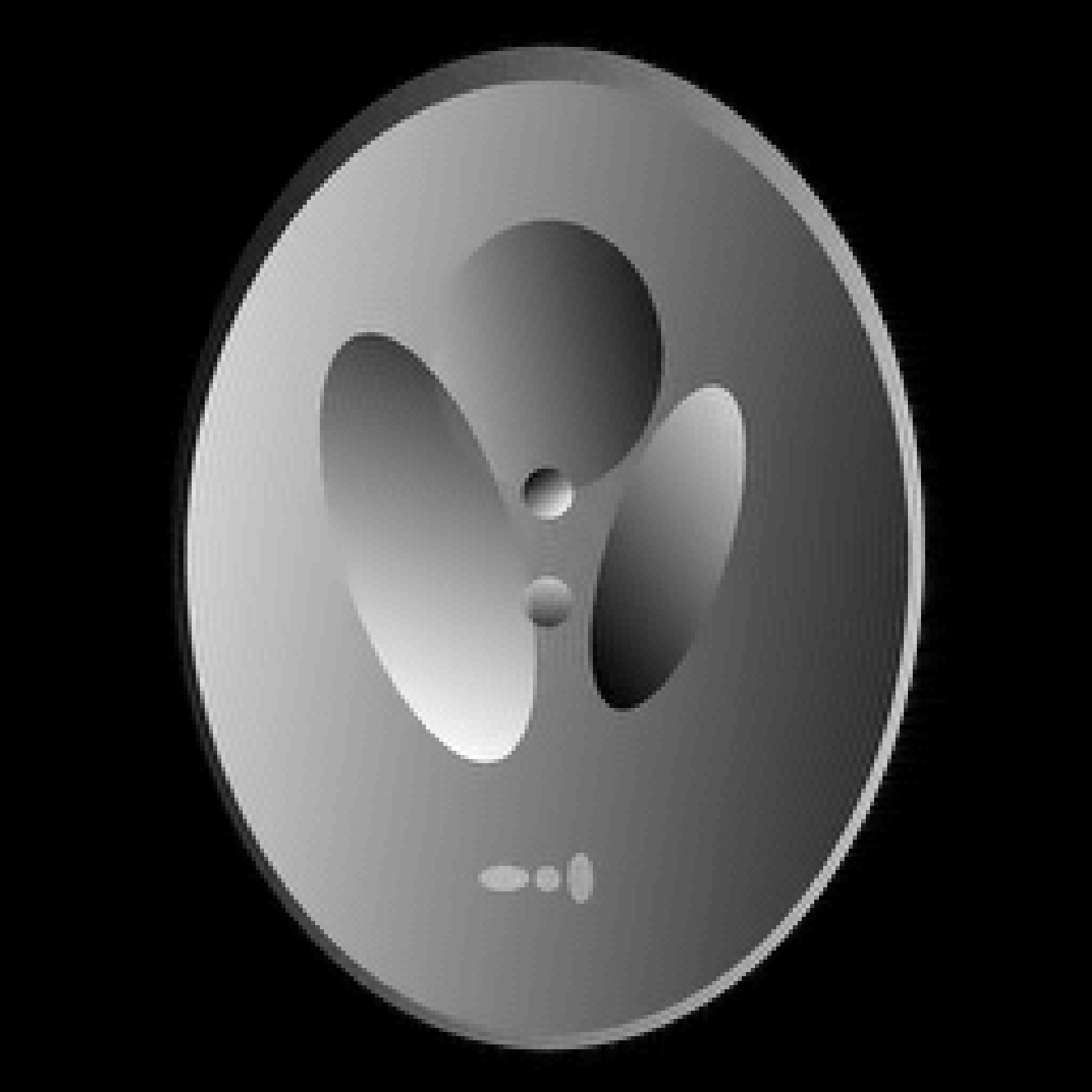}}\hspace{0.001cm}
\subfloat[]{\label{LoganLRHTGVIRLS}\includegraphics[width=2.25cm]{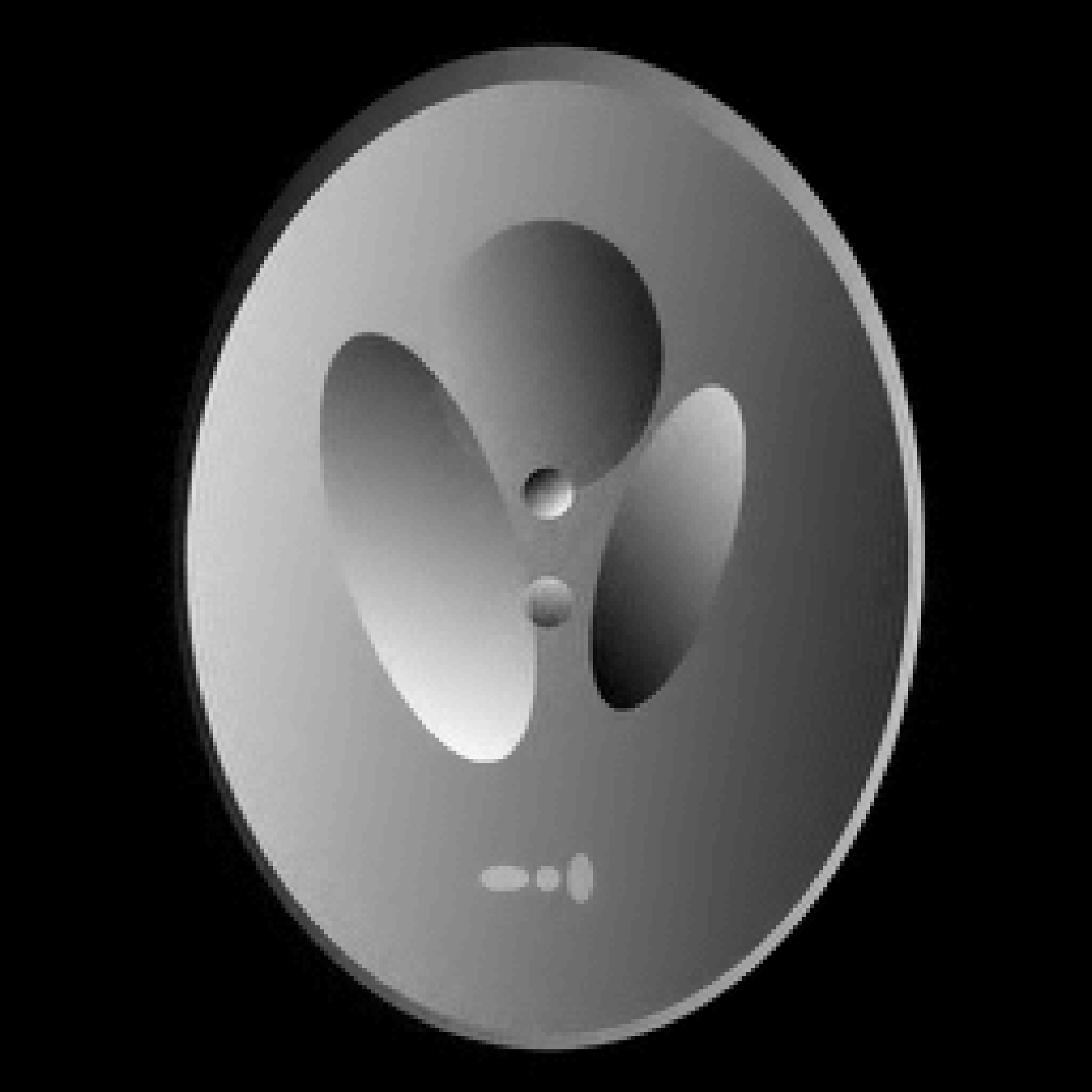}}\hspace{0.001cm}
\subfloat[]{\label{LoganLRHInfConvIRLS}\includegraphics[width=2.25cm]{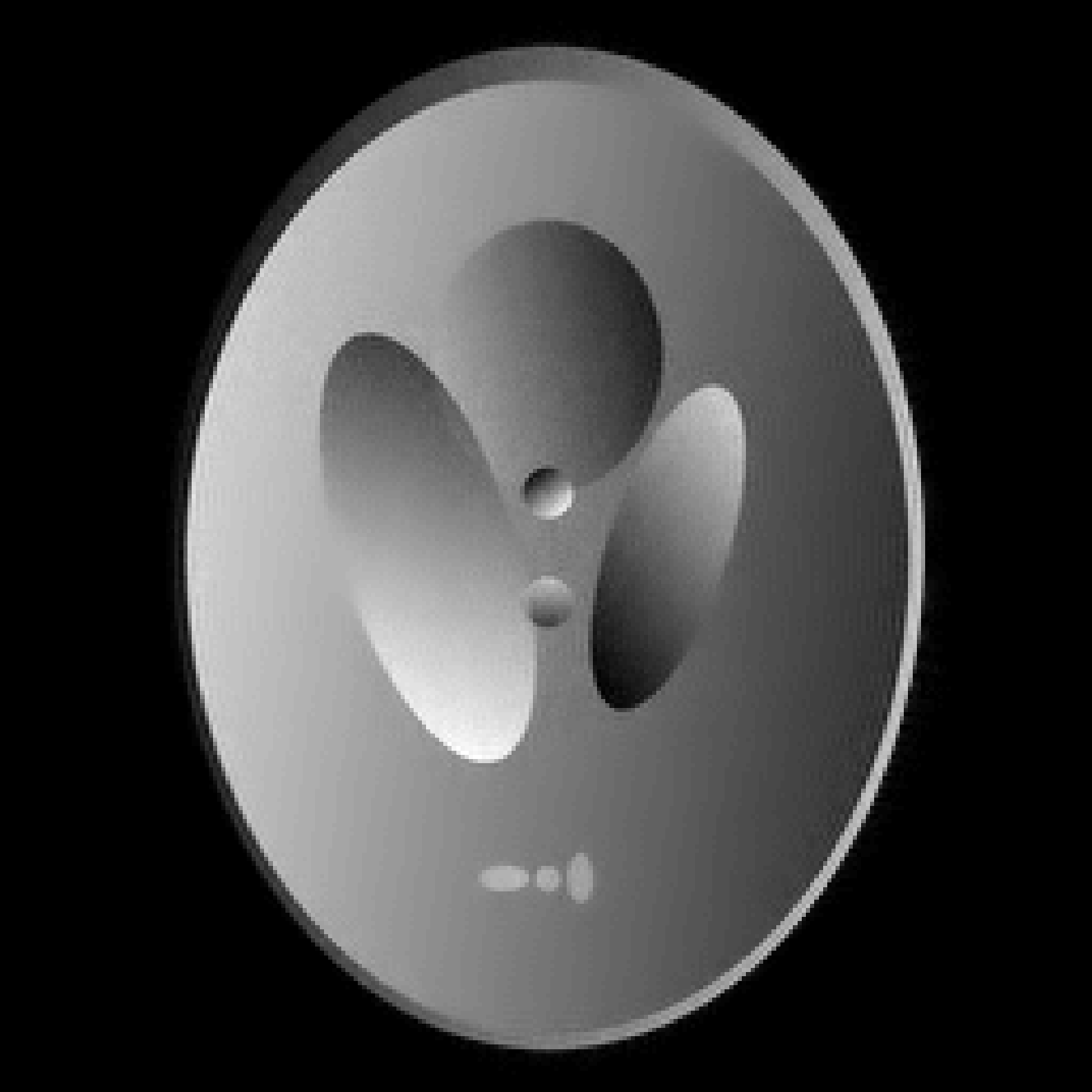}}\hspace{0.001cm}
\subfloat[]{\label{LoganFra}\includegraphics[width=2.25cm]{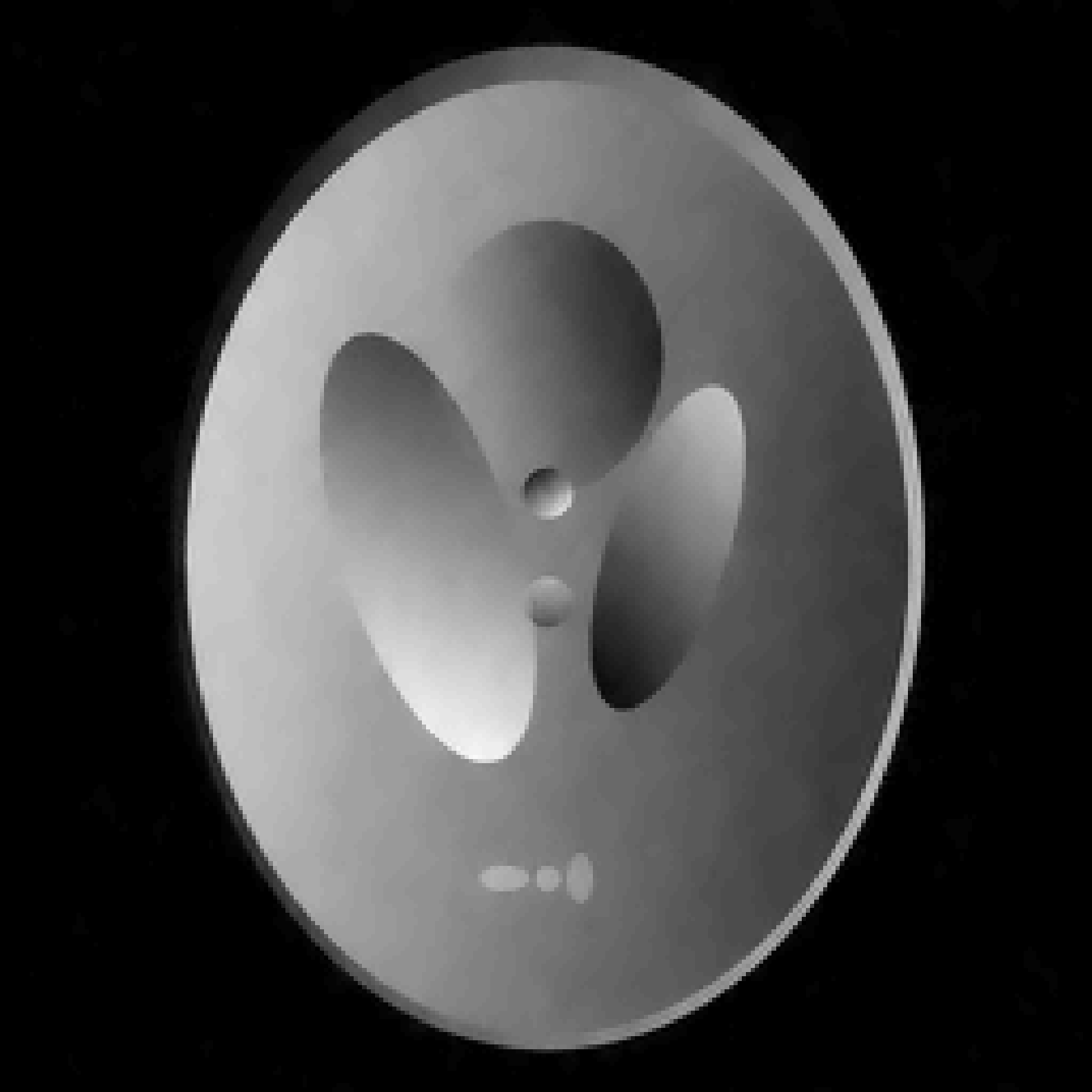}}\hspace{0.001cm}
\subfloat[]{\label{LoganTGV}\includegraphics[width=2.25cm]{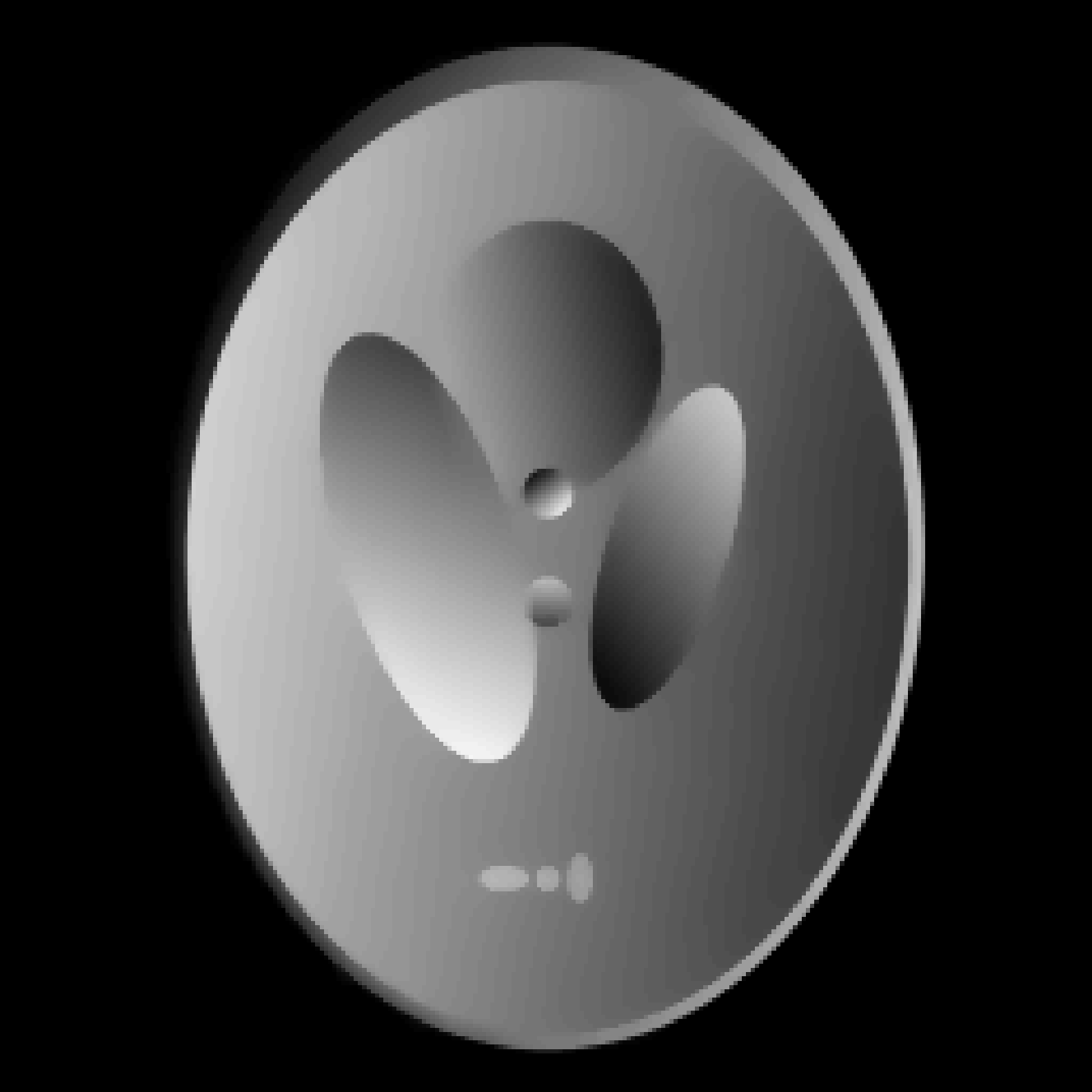}}\hspace{0.001cm}
\subfloat[]{\label{LoganInfConv}\includegraphics[width=2.25cm]{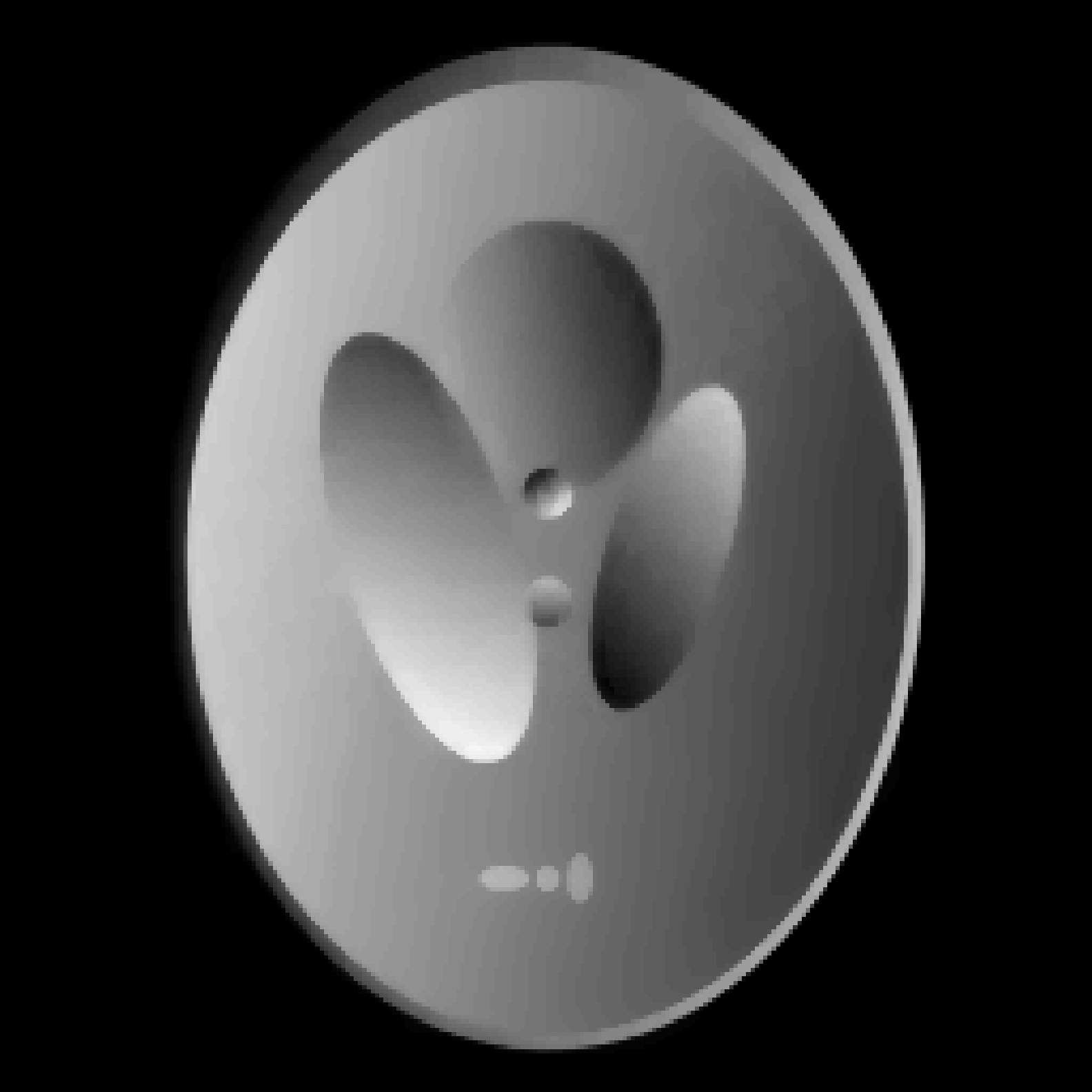}}\vspace{-0.75em}\\
\subfloat[]{\label{LoganOriginalZoom}\includegraphics[width=2.25cm]{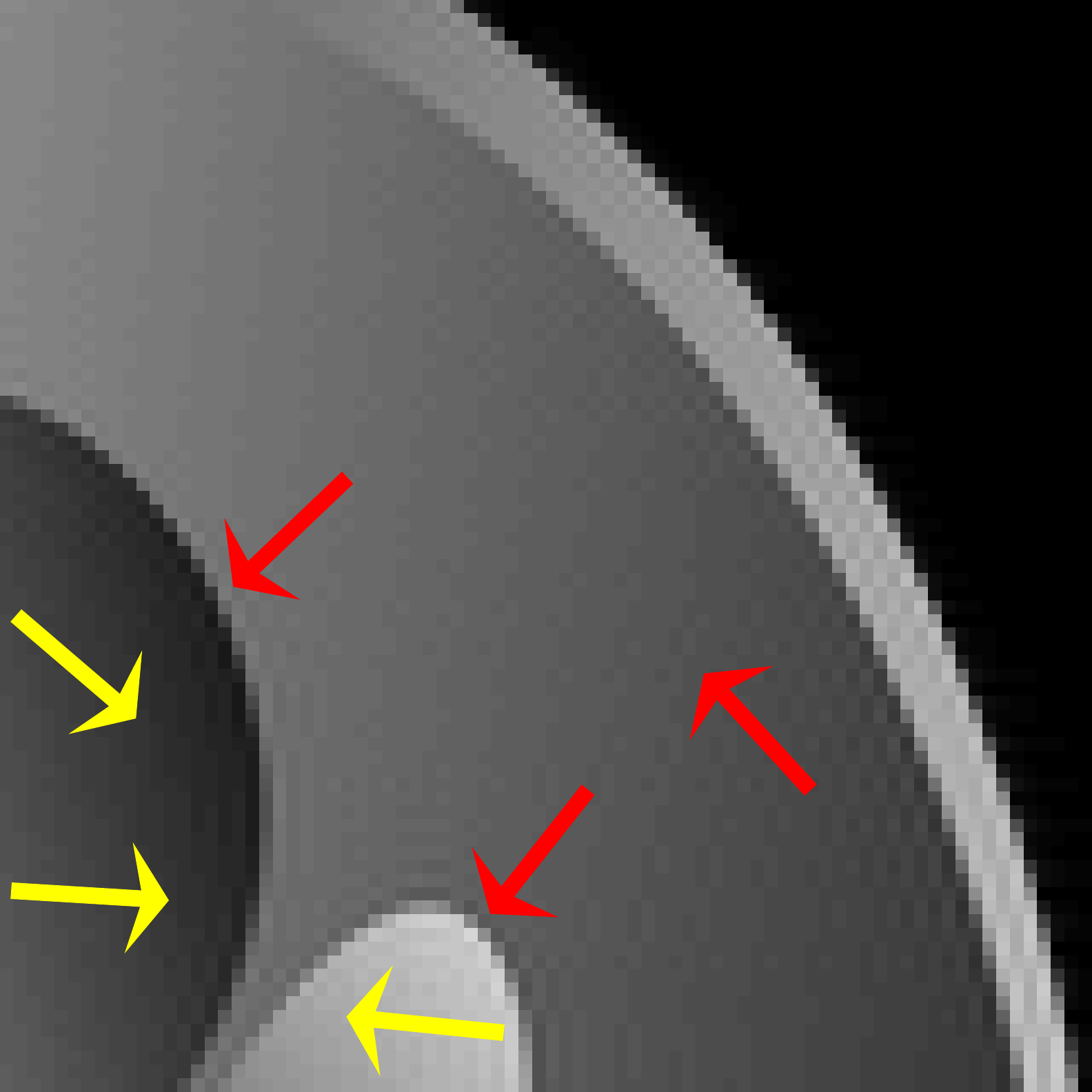}}\hspace{0.001cm}
\subfloat[]{\label{LoganLRHTGVZoom}\includegraphics[width=2.25cm]{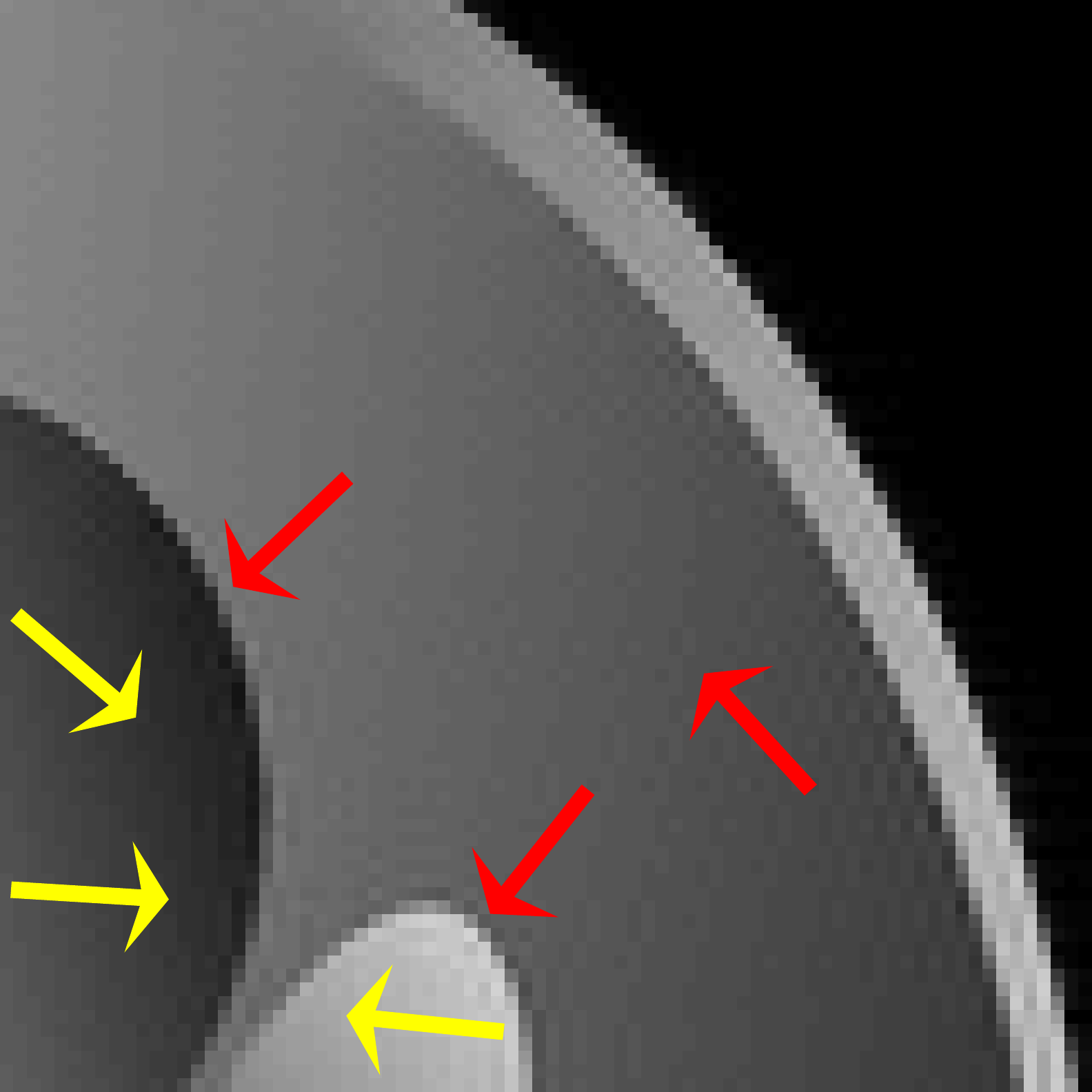}}\hspace{0.001cm}
\subfloat[]{\label{LoganLRHTGVIRLSZoom}\includegraphics[width=2.25cm]{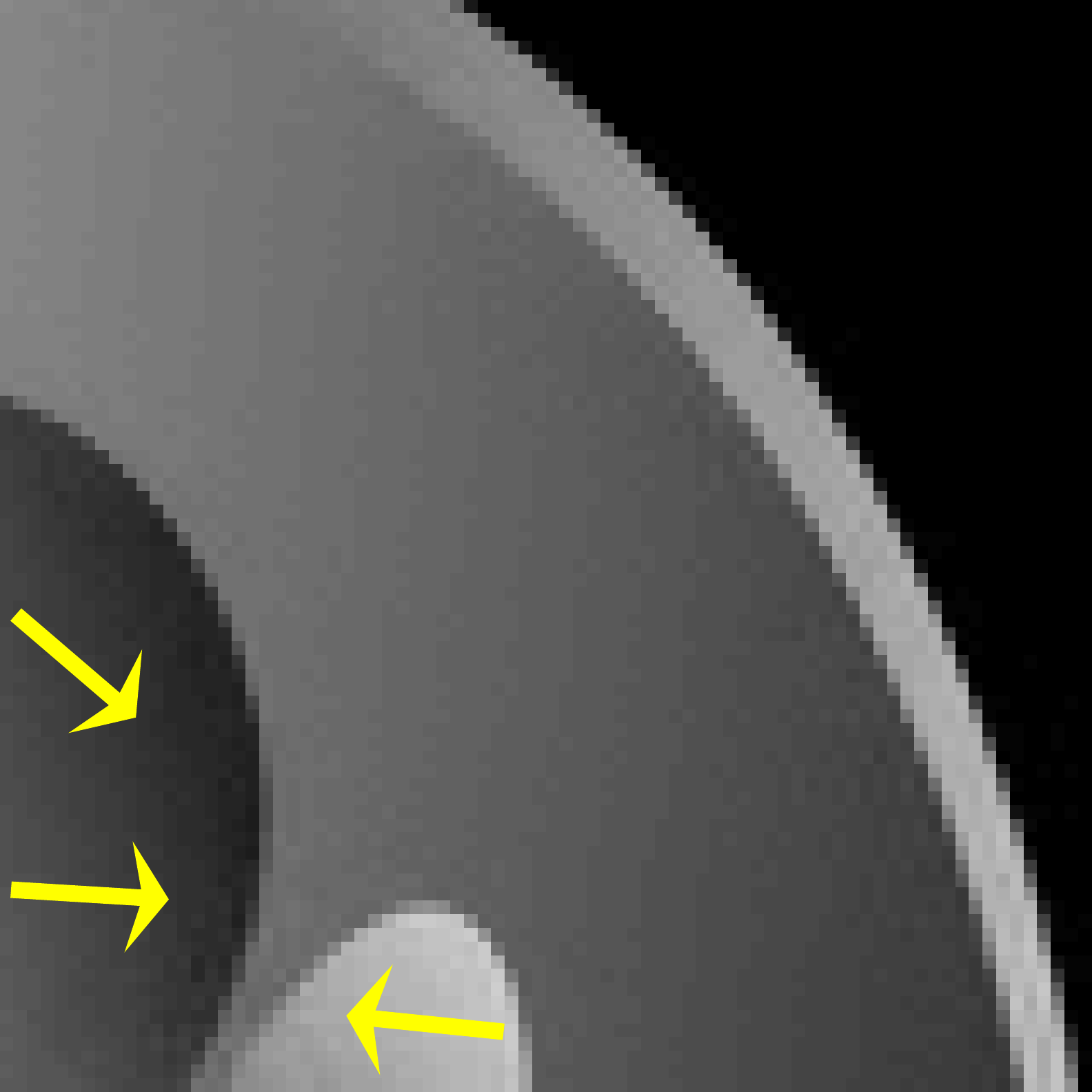}}\hspace{0.001cm}
\subfloat[]{\label{LoganLRHInfConvIRLSZoom}\includegraphics[width=2.25cm]{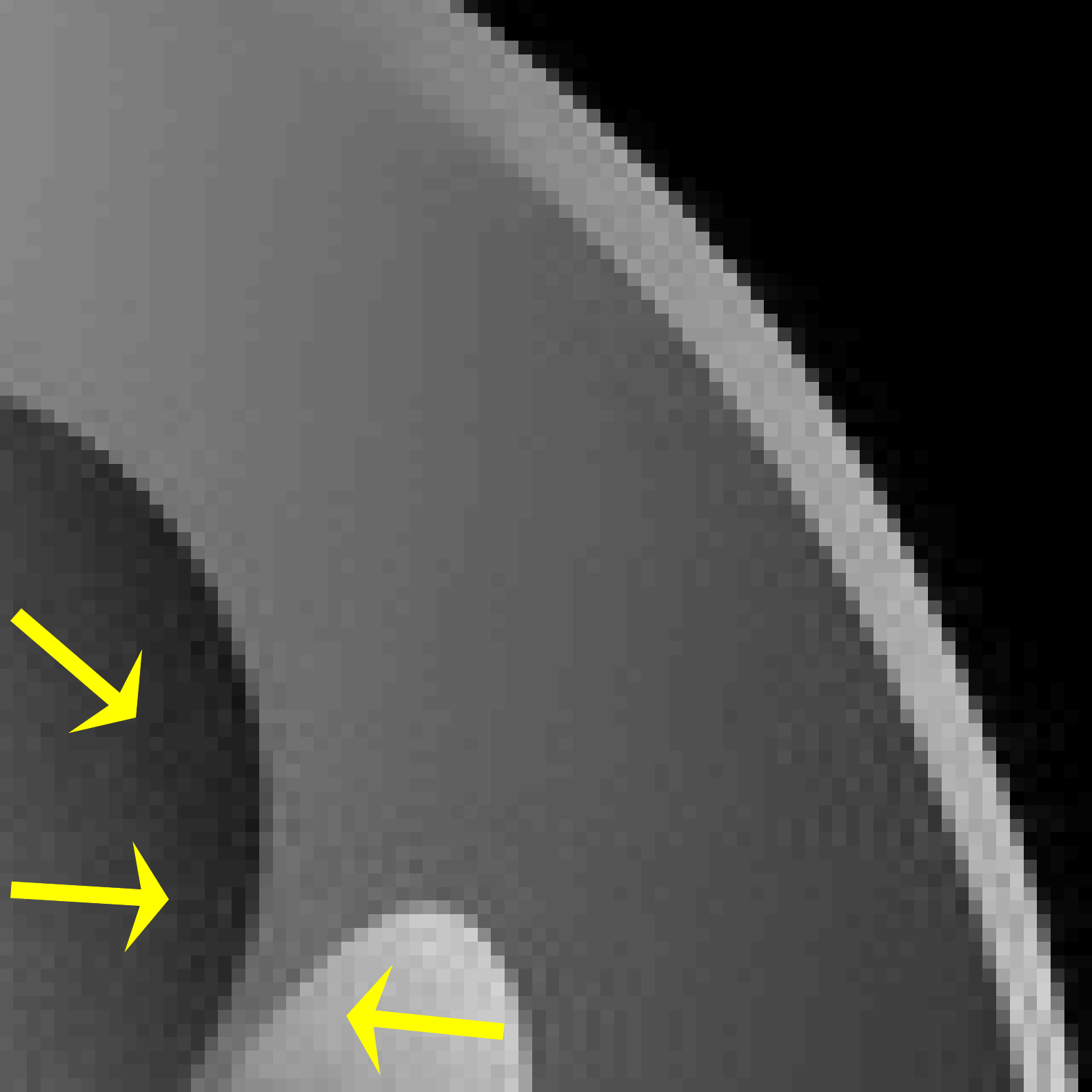}}\hspace{0.001cm}
\subfloat[]{\label{LoganFraZoom}\includegraphics[width=2.25cm]{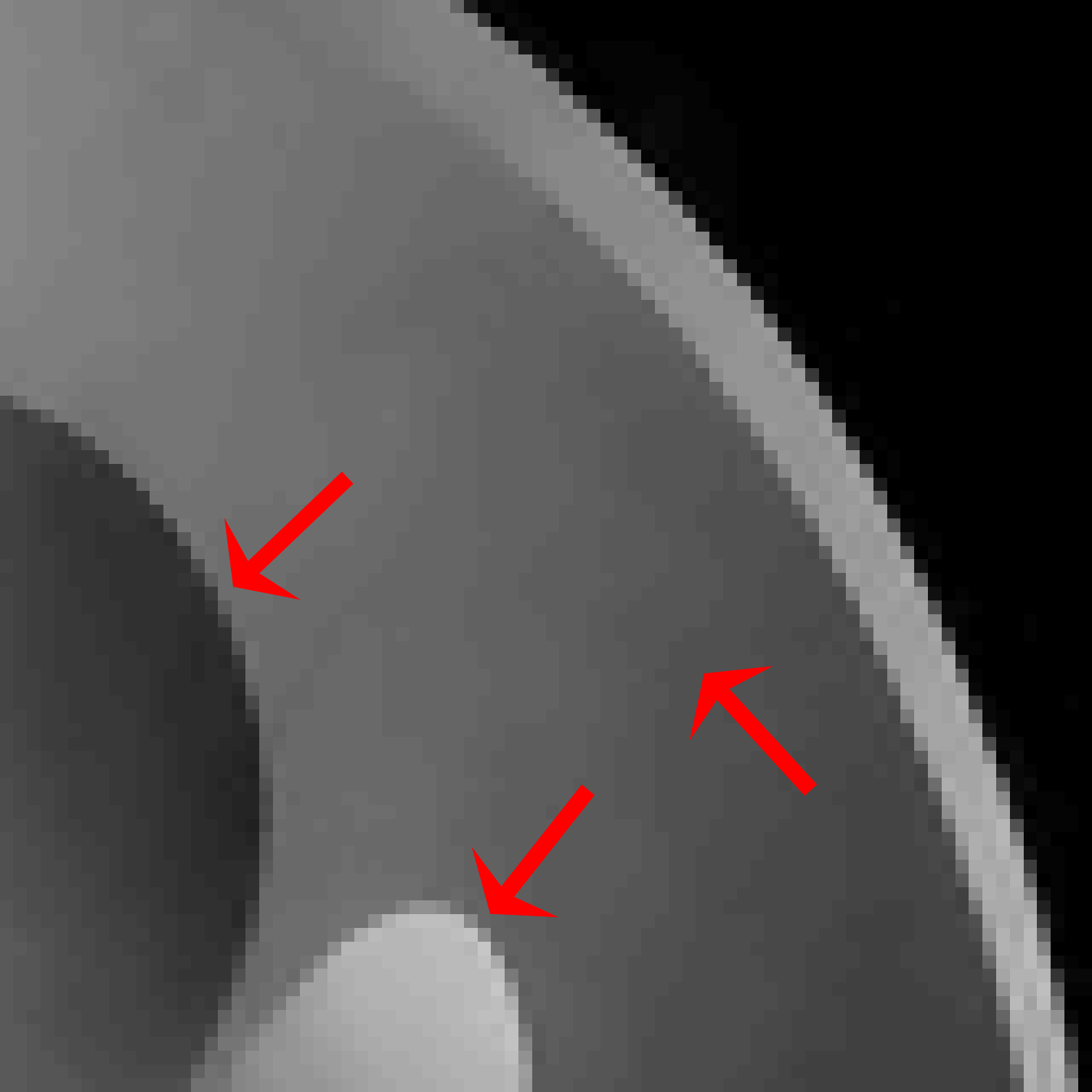}}\hspace{0.001cm}
\subfloat[]{\label{LoganTGVZoom}\includegraphics[width=2.25cm]{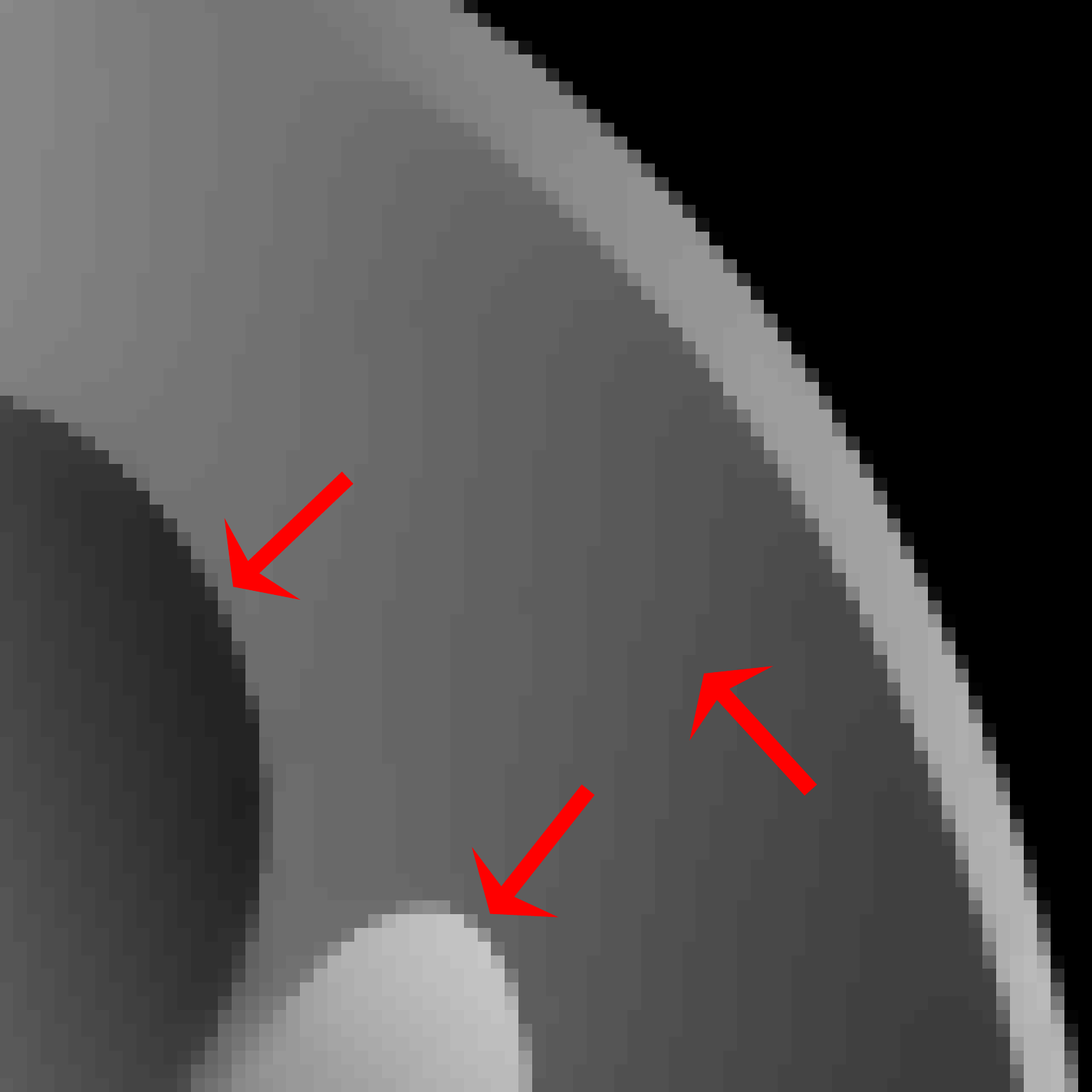}}\hspace{0.001cm}
\subfloat[]{\label{LoganInfConvZoom}\includegraphics[width=2.25cm]{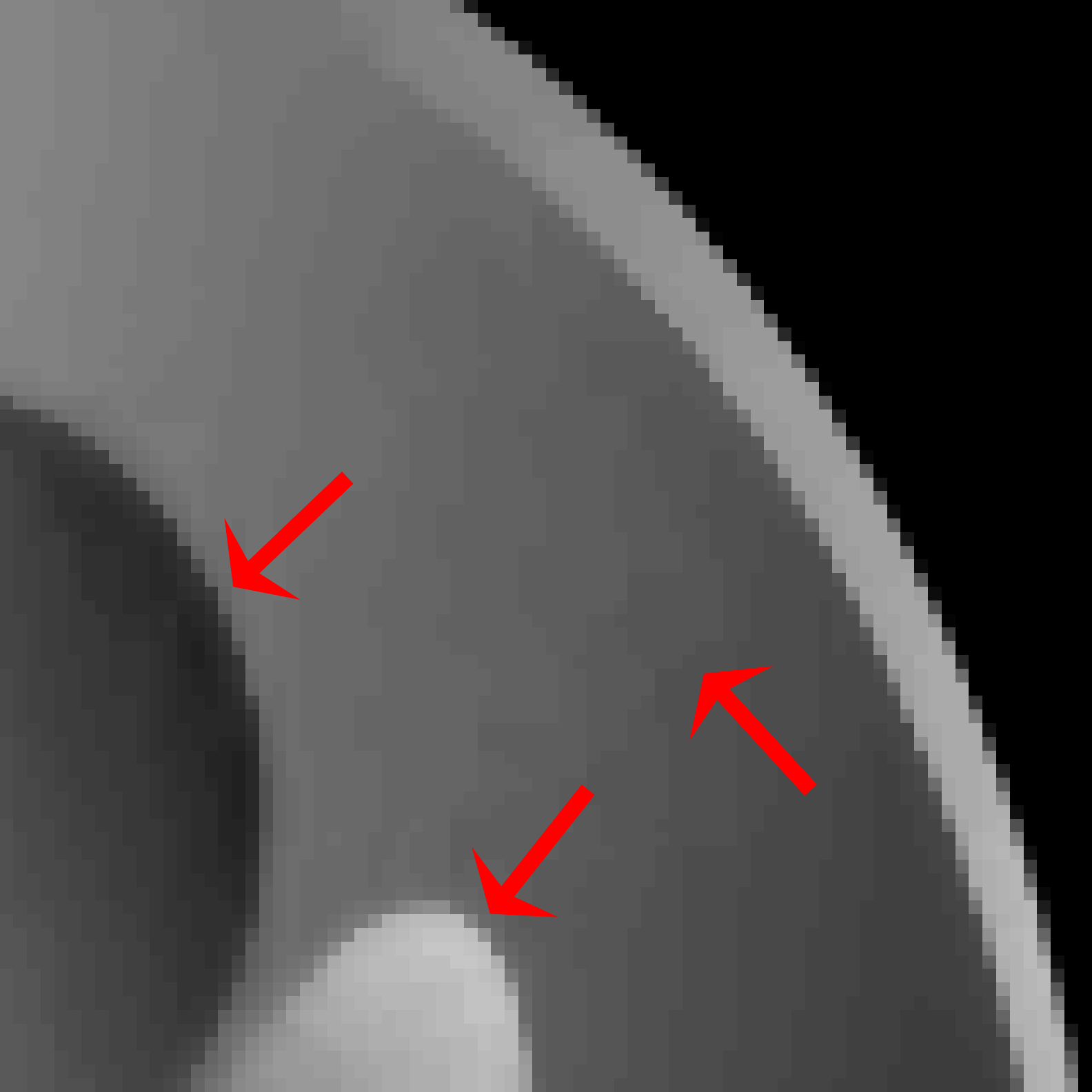}}\vspace{-0.75em}\\
\subfloat[]{\label{LoganMask}\includegraphics[width=2.25cm]{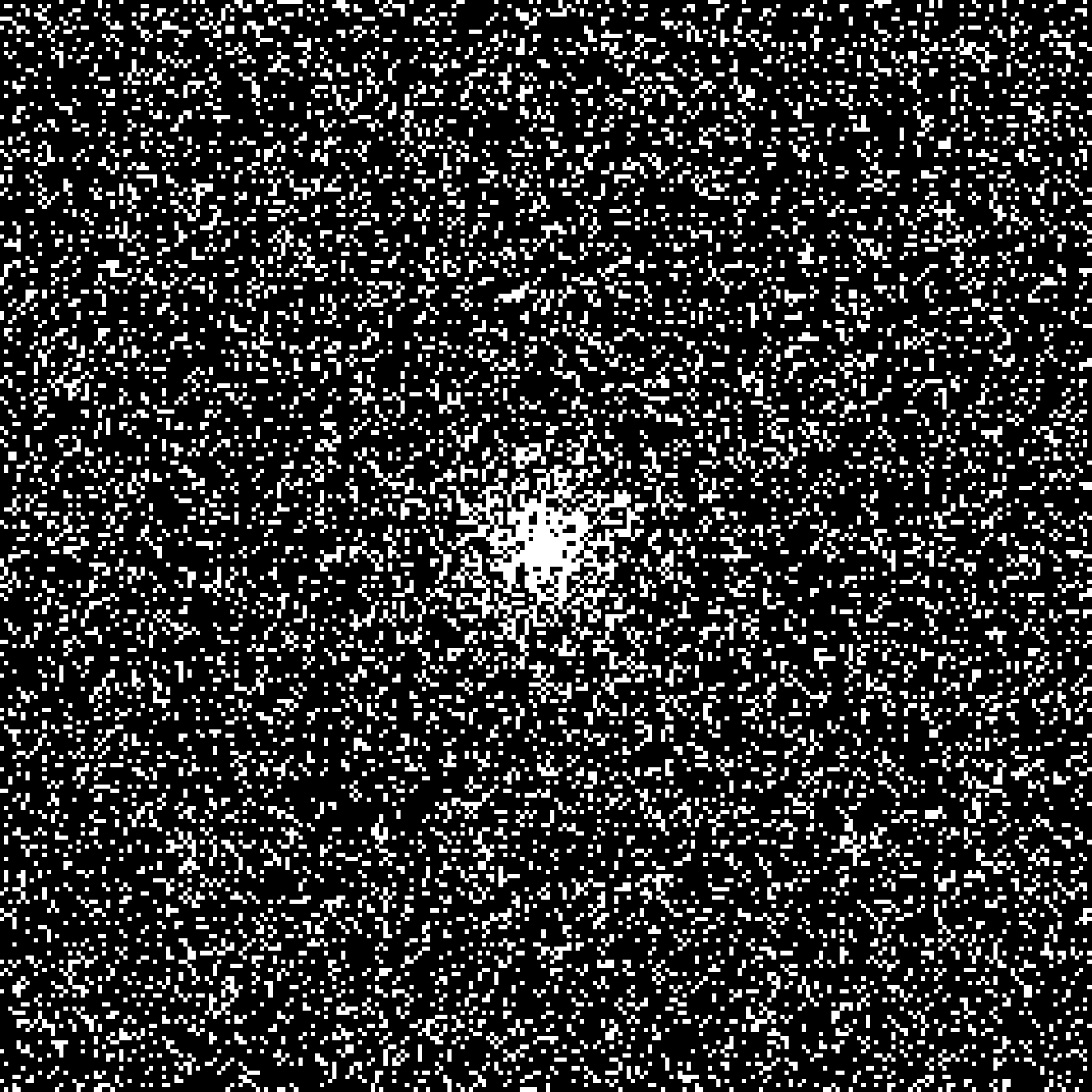}}\hspace{0.001cm}
\subfloat[]{\label{LoganLRHTGVError}\includegraphics[width=2.25cm]{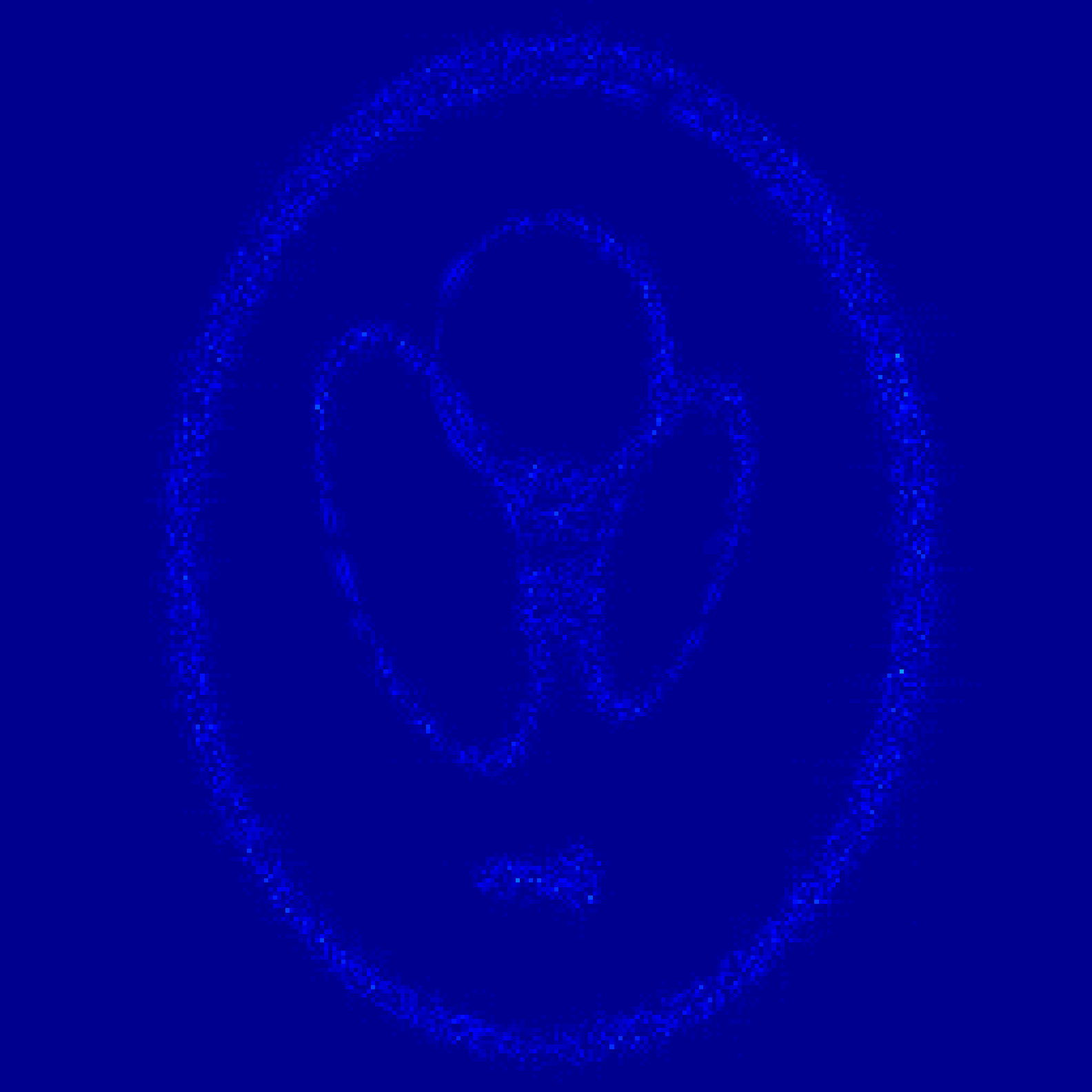}}\hspace{0.001cm}
\subfloat[]{\label{LoganLRHTGVIRLSError}\includegraphics[width=2.25cm]{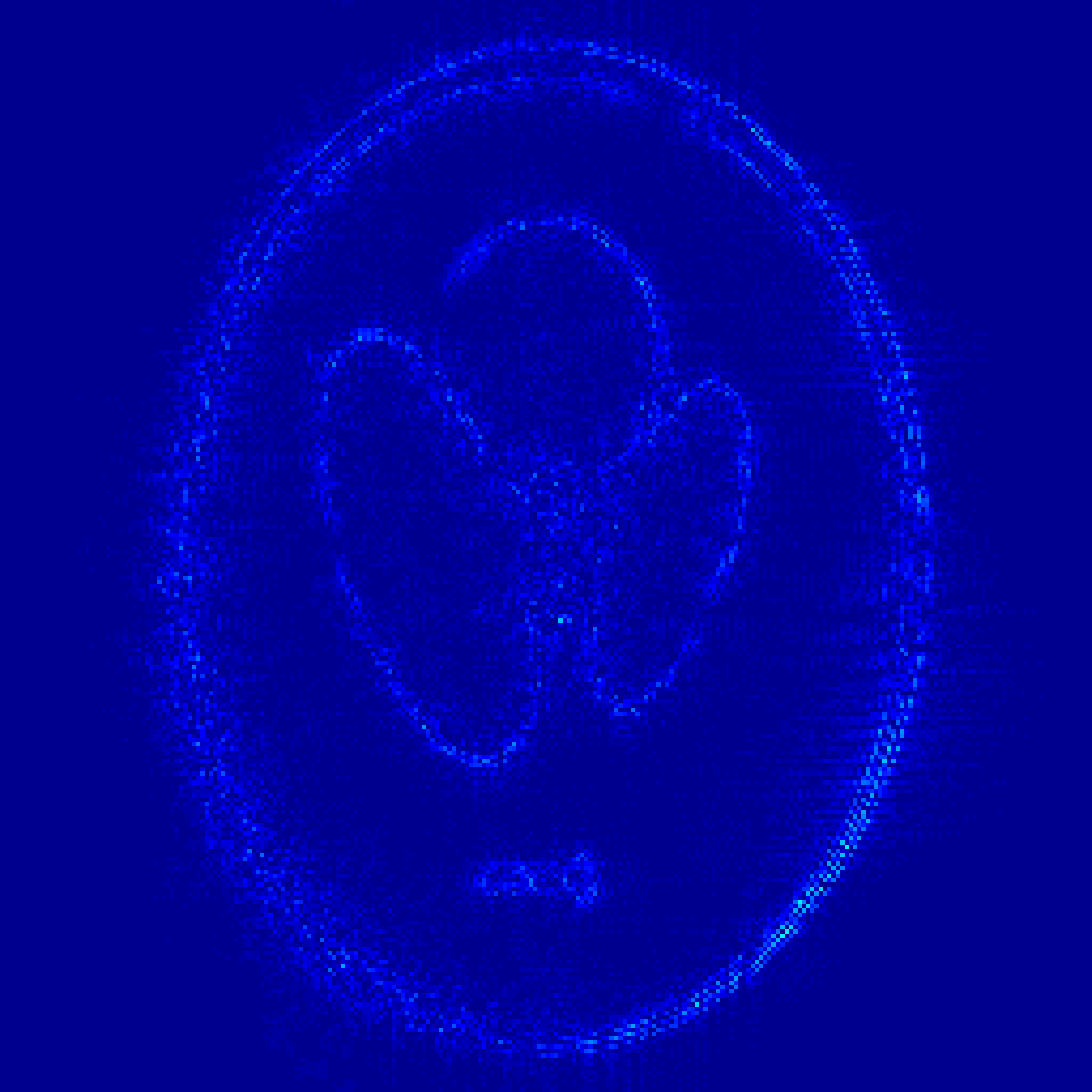}}\hspace{0.001cm}
\subfloat[]{\label{LoganLRHInfConvIRLSError}\includegraphics[width=2.25cm]{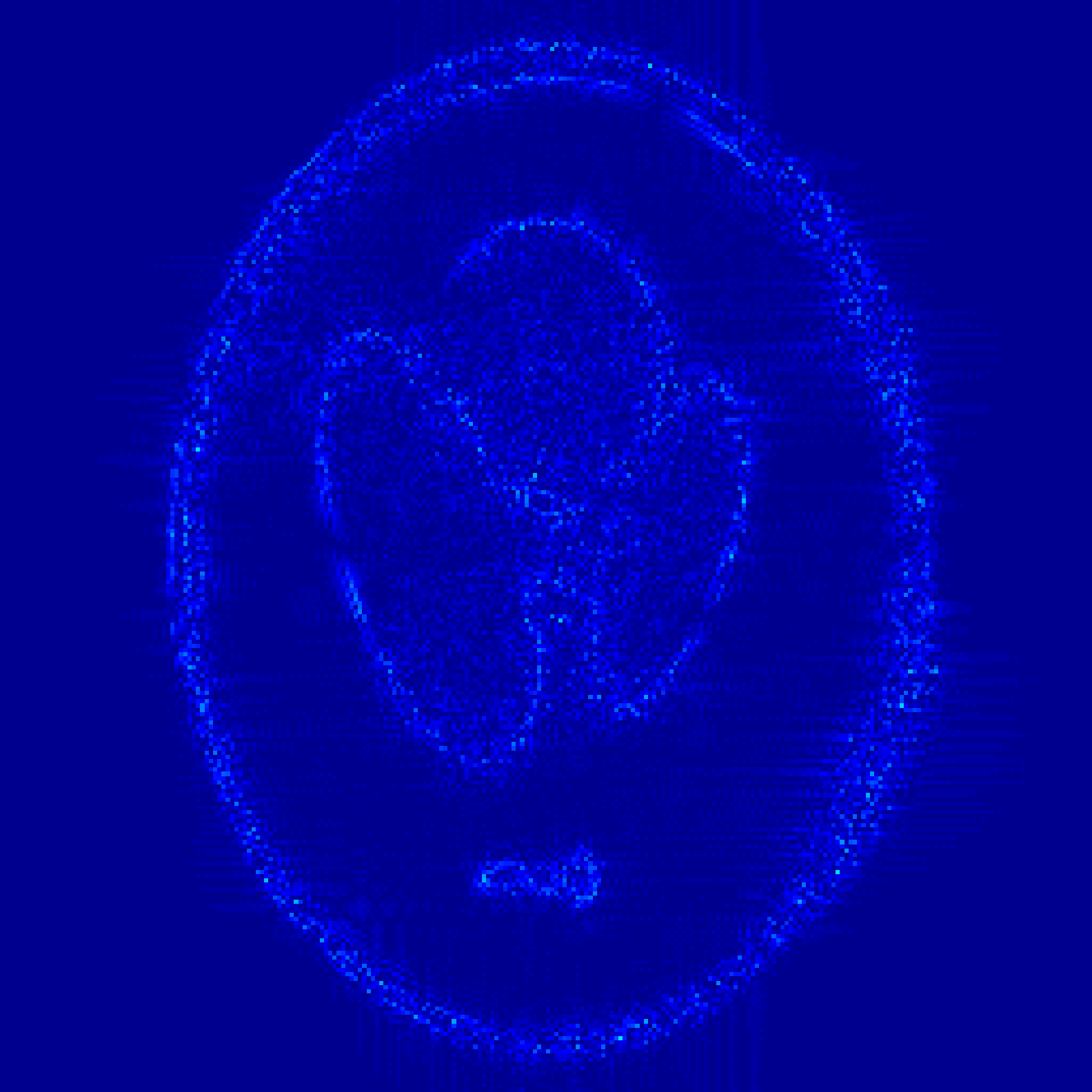}}\hspace{0.001cm}
\subfloat[]{\label{LoganFraError}\includegraphics[width=2.25cm]{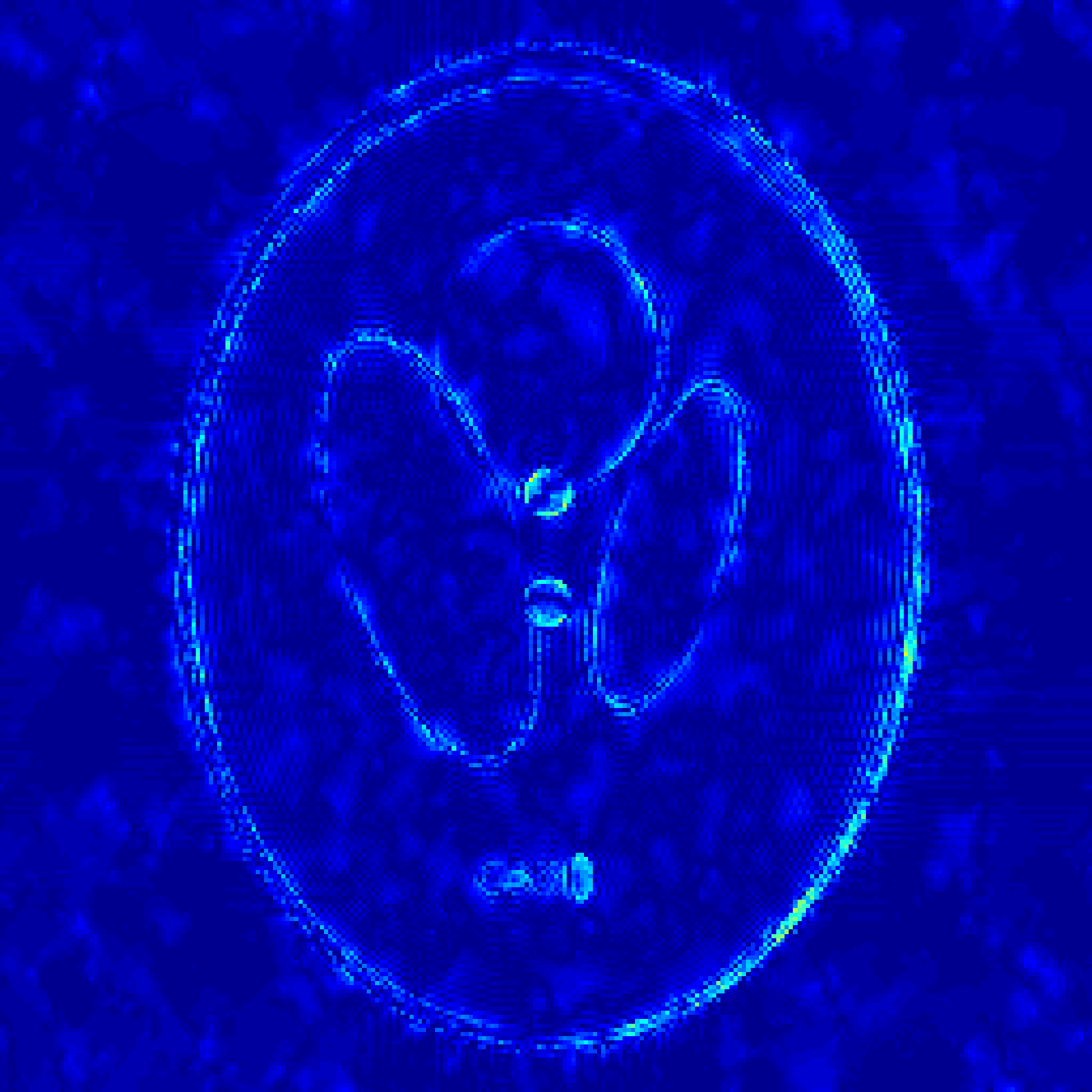}}\hspace{0.001cm}
\subfloat[]{\label{LoganTGVError}\includegraphics[width=2.25cm]{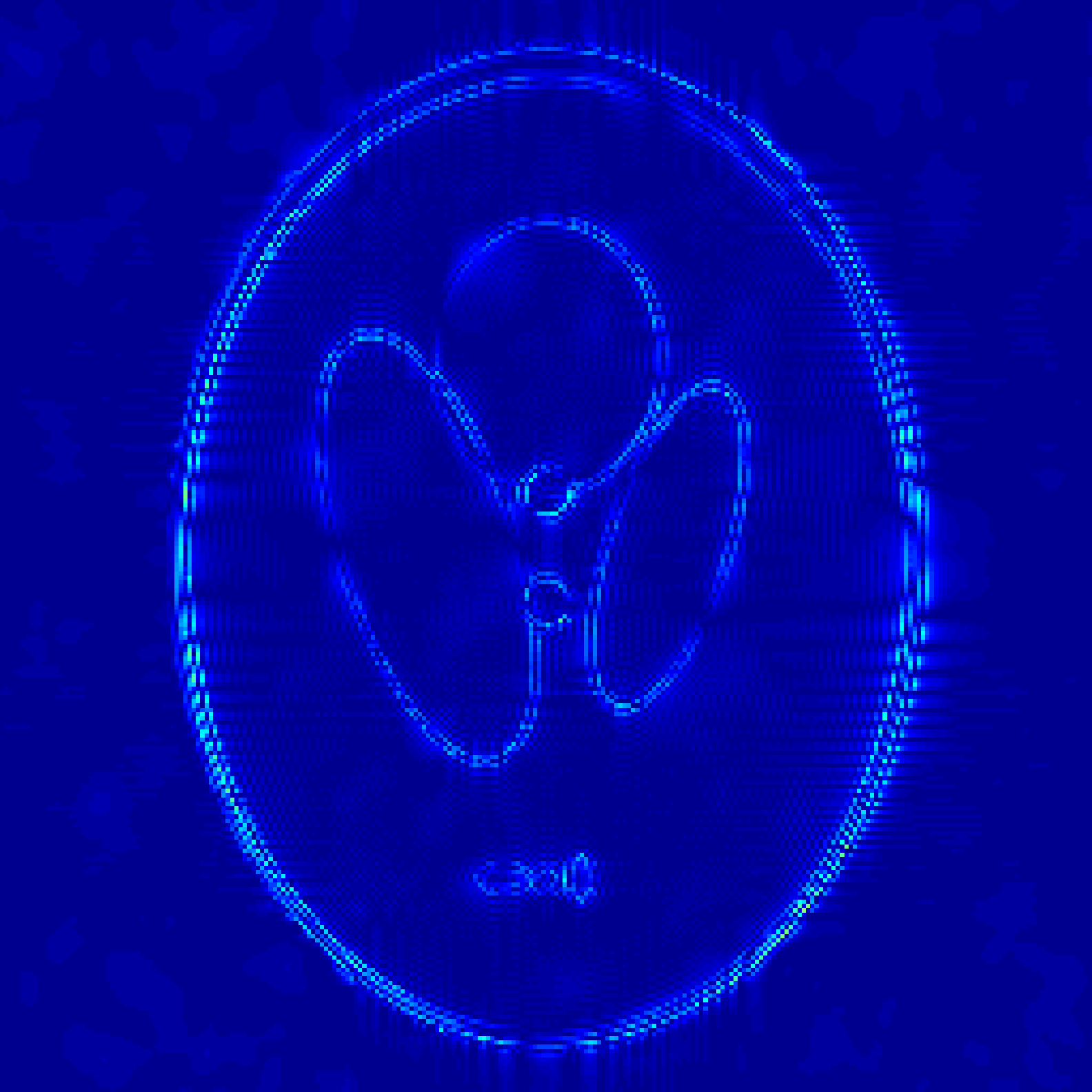}}\hspace{0.001cm}
\subfloat[]{\label{LoganInfConvError}\includegraphics[width=2.25cm]{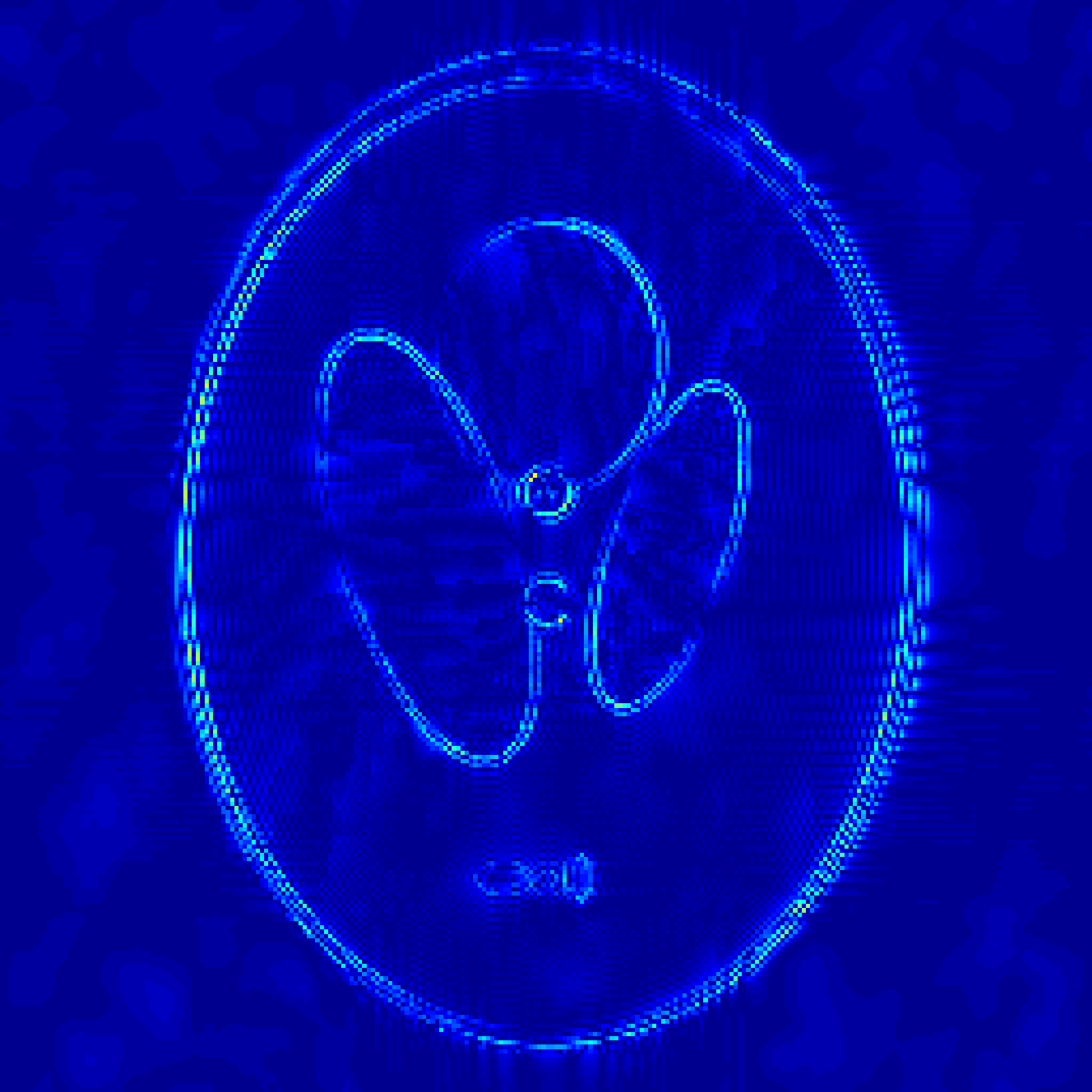}}
\caption{Visual comparisons for ``Ellipses''. \cref{LoganOriginal}: original image, \cref{LoganLRHTGV}: model \cref{ProposedCSMRI}, \cref{LoganLRHTGVIRLS}: SLRM \cref{LRHTGVIRLSCSMRI}, \cref{LoganLRHInfConvIRLS}: GSLR \cref{LRHInfConvIRLSCSMRI}, \cref{LoganFra}: framelet \cref{FrameCSMRI}, \cref{LoganTGV}: TGV \cref{TGVCSMRI}, \cref{LoganInfConv}: infimal convolution \cref{InfConvCSMRI}. \cref{LoganOriginalZoom,LoganLRHTGVZoom,LoganLRHTGVIRLSZoom,LoganLRHInfConvIRLSZoom,LoganFraZoom,LoganTGVZoom,LoganInfConvZoom}: zoom-in views of \cref{LoganOriginal,LoganLRHTGV,LoganLRHTGVIRLS,LoganLRHInfConvIRLS,LoganFra,LoganTGV,LoganInfConv}. Yellow arrows indicate the region worth noting for comparisons among \cref{ProposedCSMRI,LRHTGVIRLSCSMRI,LRHInfConvIRLSCSMRI}, and red arrows indicate the region worth noting for comparisons among \cref{ProposedCSMRI,FrameCSMRI,TGVCSMRI,InfConvCSMRI}. \cref{LoganMask}: sample region, \cref{LoganLRHTGVError,LoganLRHTGVIRLSError,LoganLRHInfConvIRLSError,LoganFraError,LoganTGVError,LoganInfConvError}: error maps of \cref{LoganLRHTGV,LoganLRHTGVIRLS,LoganLRHInfConvIRLS,LoganFra,LoganTGV,LoganInfConv}.}\label{LoganResults}
\end{figure}

\begin{figure}[t]
\centering
\subfloat[]{\label{RectangleOriginal}\includegraphics[width=2.25cm]{RectangleOriginal.pdf}}\hspace{0.001cm}
\subfloat[]{\label{RectangleLRHTGV}\includegraphics[width=2.25cm]{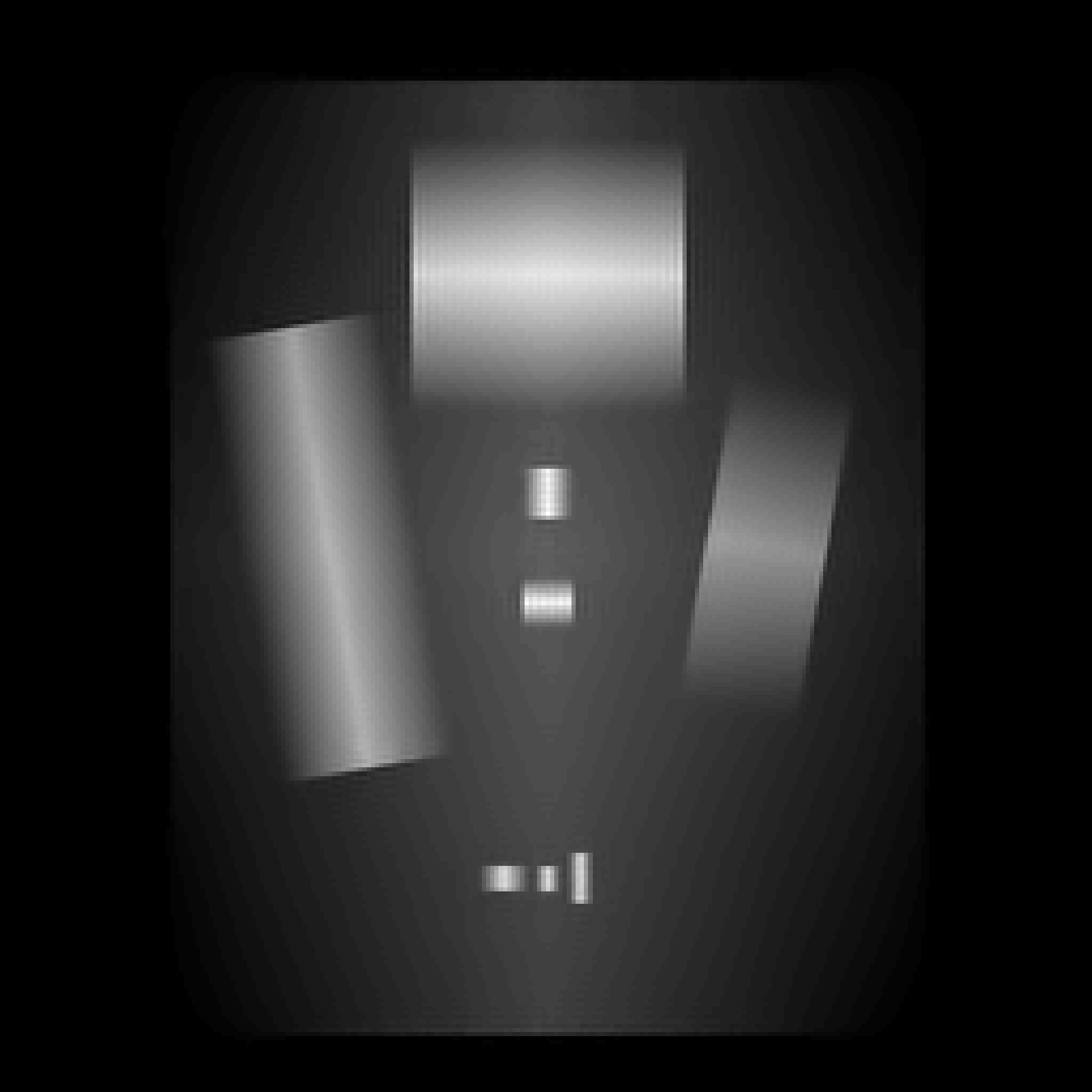}}\hspace{0.001cm}
\subfloat[]{\label{RectangleLRHTGVIRLS}\includegraphics[width=2.25cm]{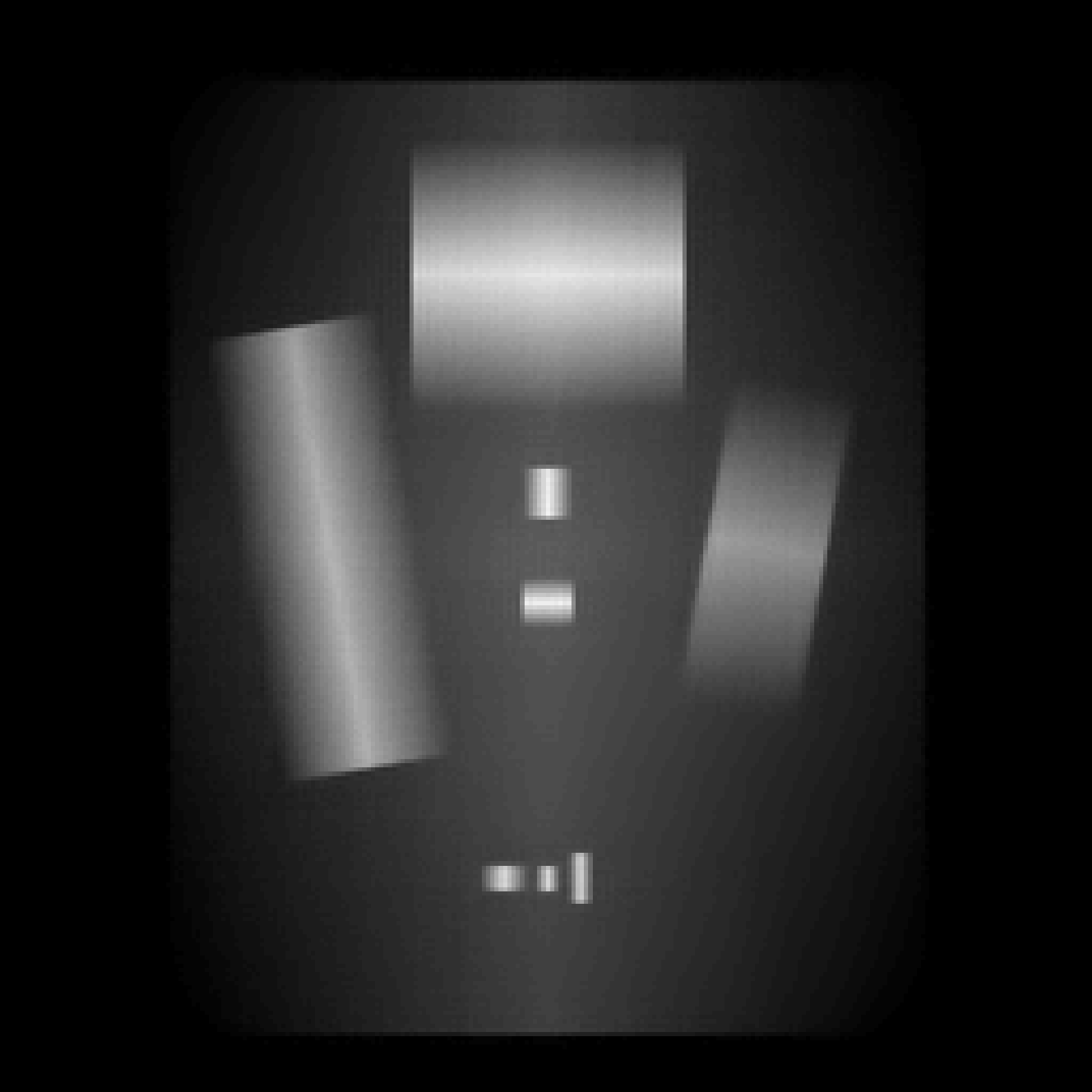}}\hspace{0.001cm}
\subfloat[]{\label{RectangleLRHInfConvIRLS}\includegraphics[width=2.25cm]{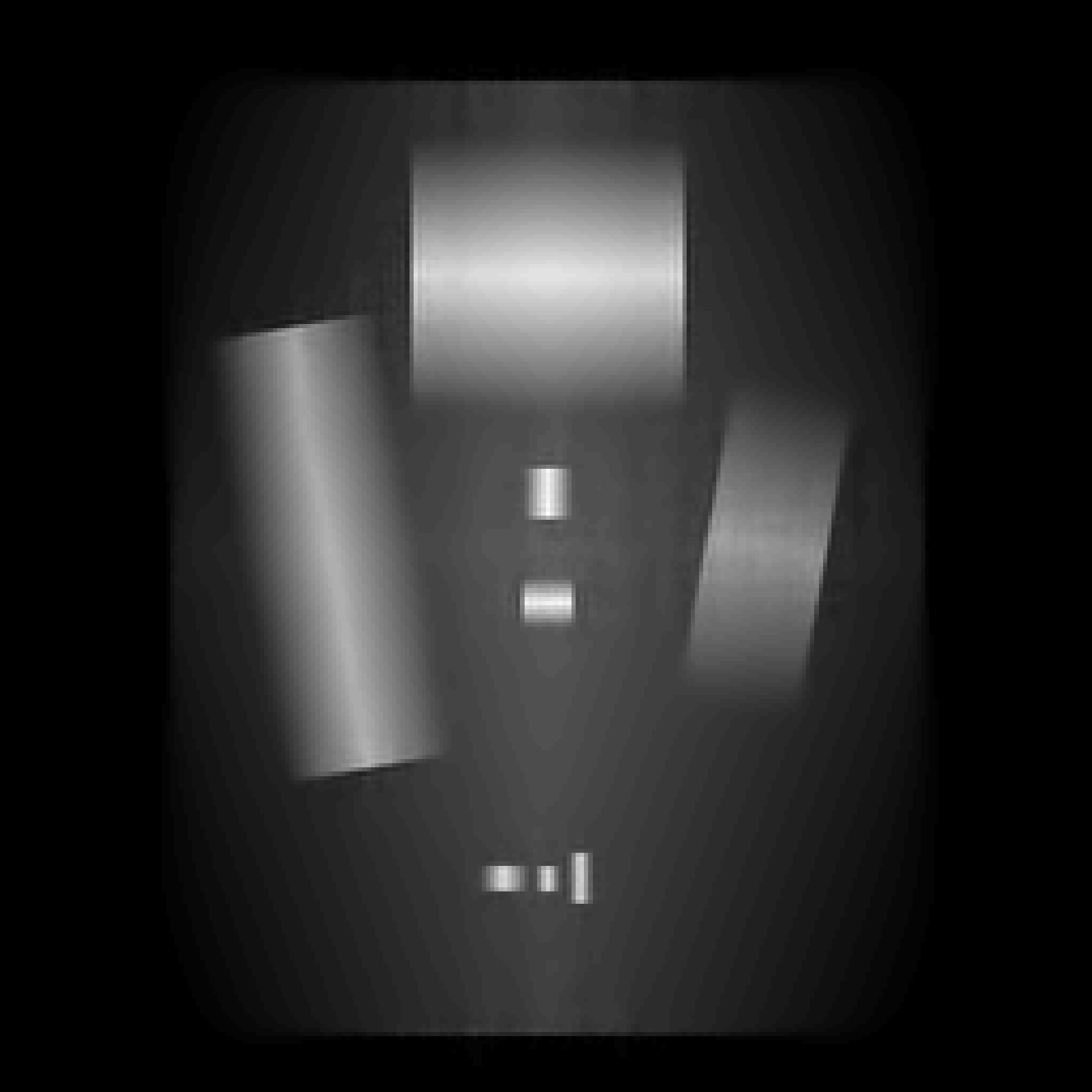}}\hspace{0.001cm}
\subfloat[]{\label{RectangleFra}\includegraphics[width=2.25cm]{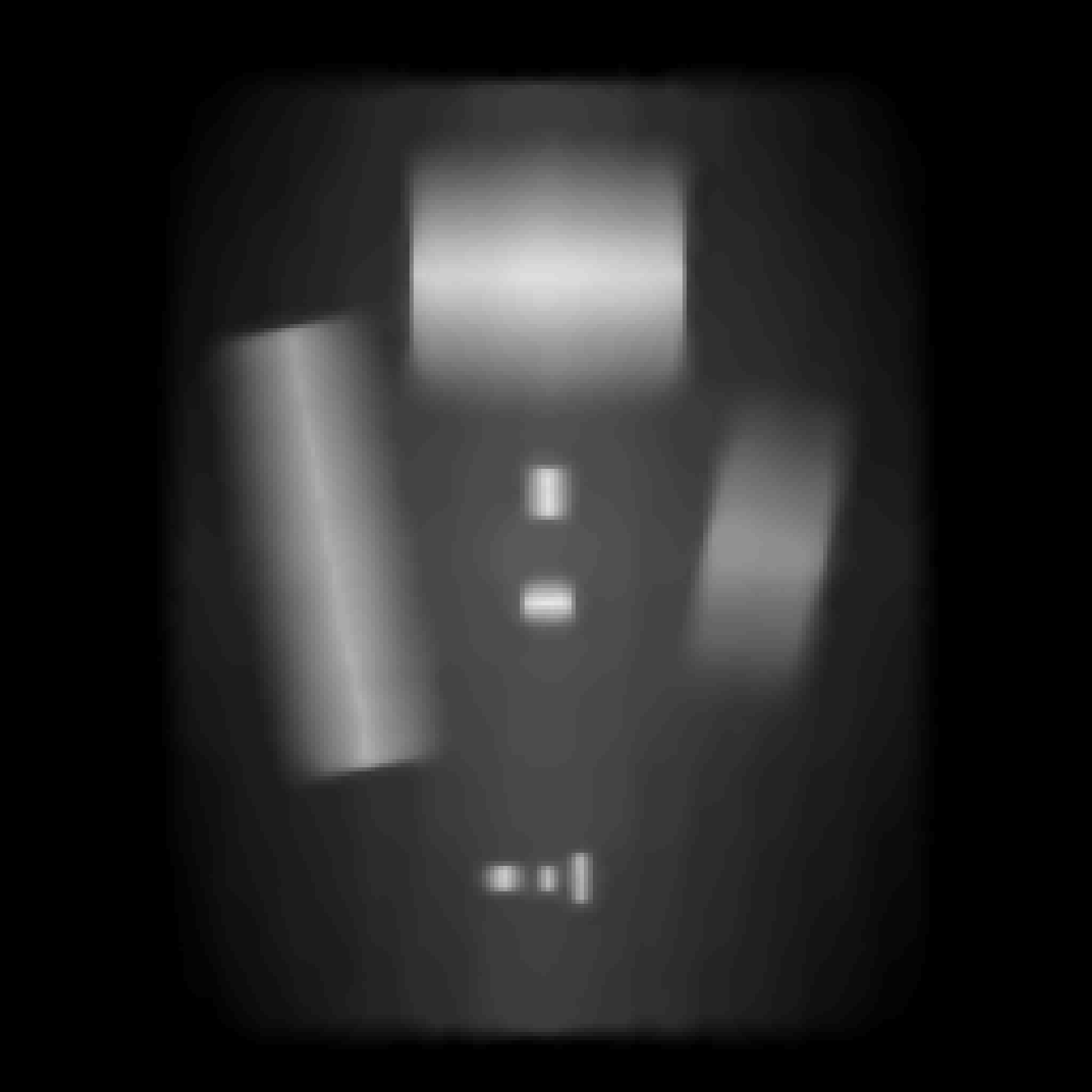}}\hspace{0.001cm}
\subfloat[]{\label{RectangleTGV}\includegraphics[width=2.25cm]{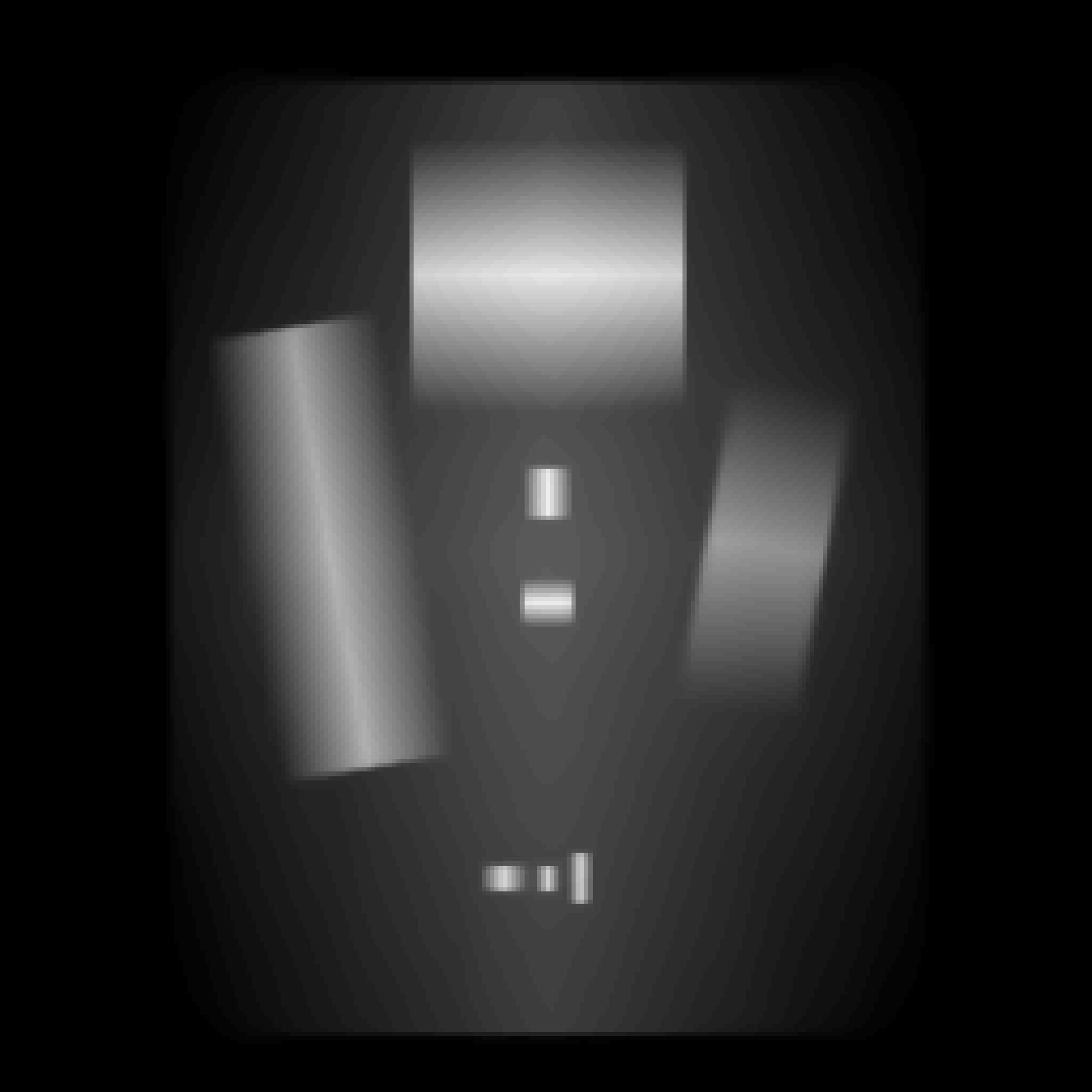}}\hspace{0.001cm}
\subfloat[]{\label{RectangleInfConv}\includegraphics[width=2.25cm]{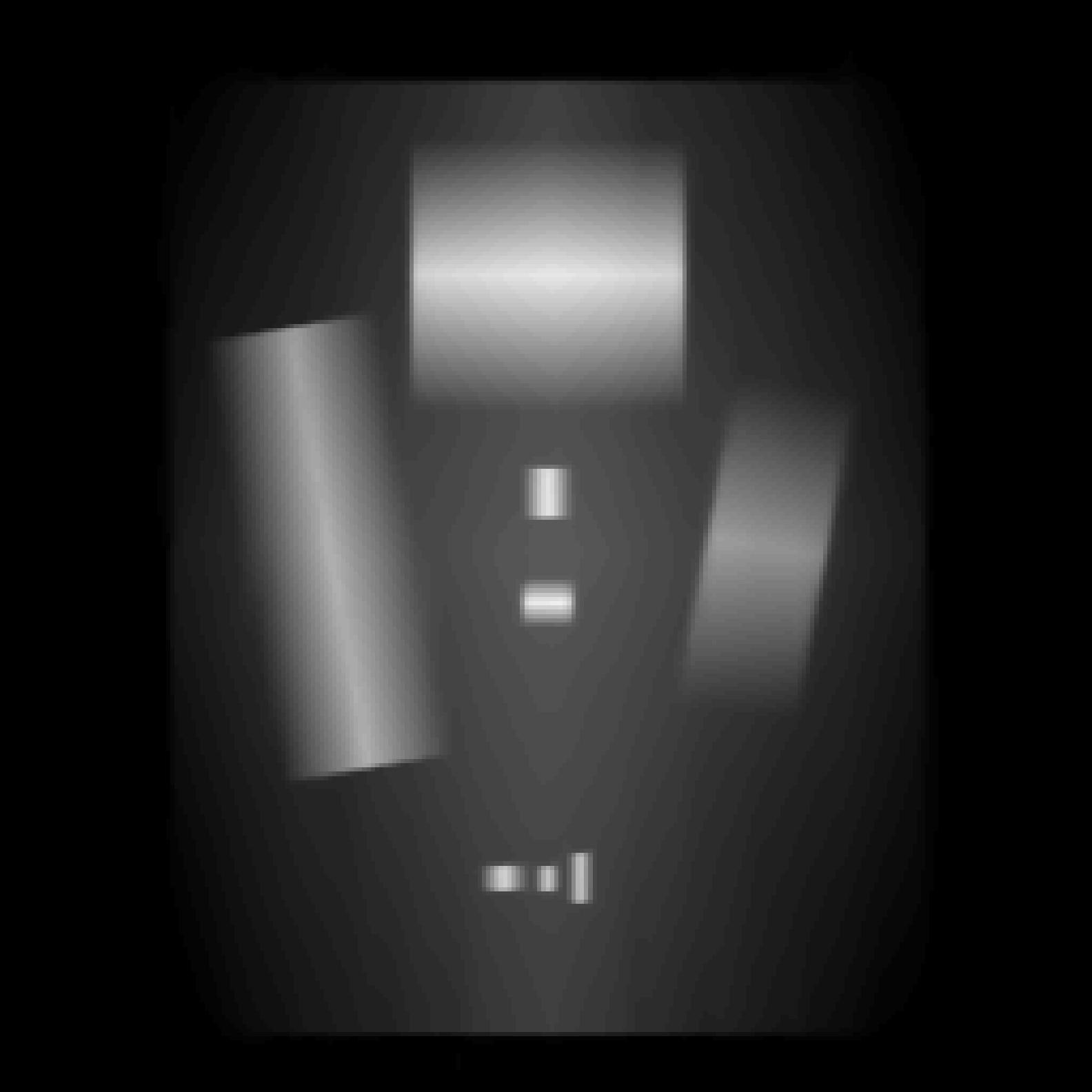}}\vspace{-0.75em}\\
\subfloat[]{\label{RectangleOriginalZoom}\includegraphics[width=2.25cm]{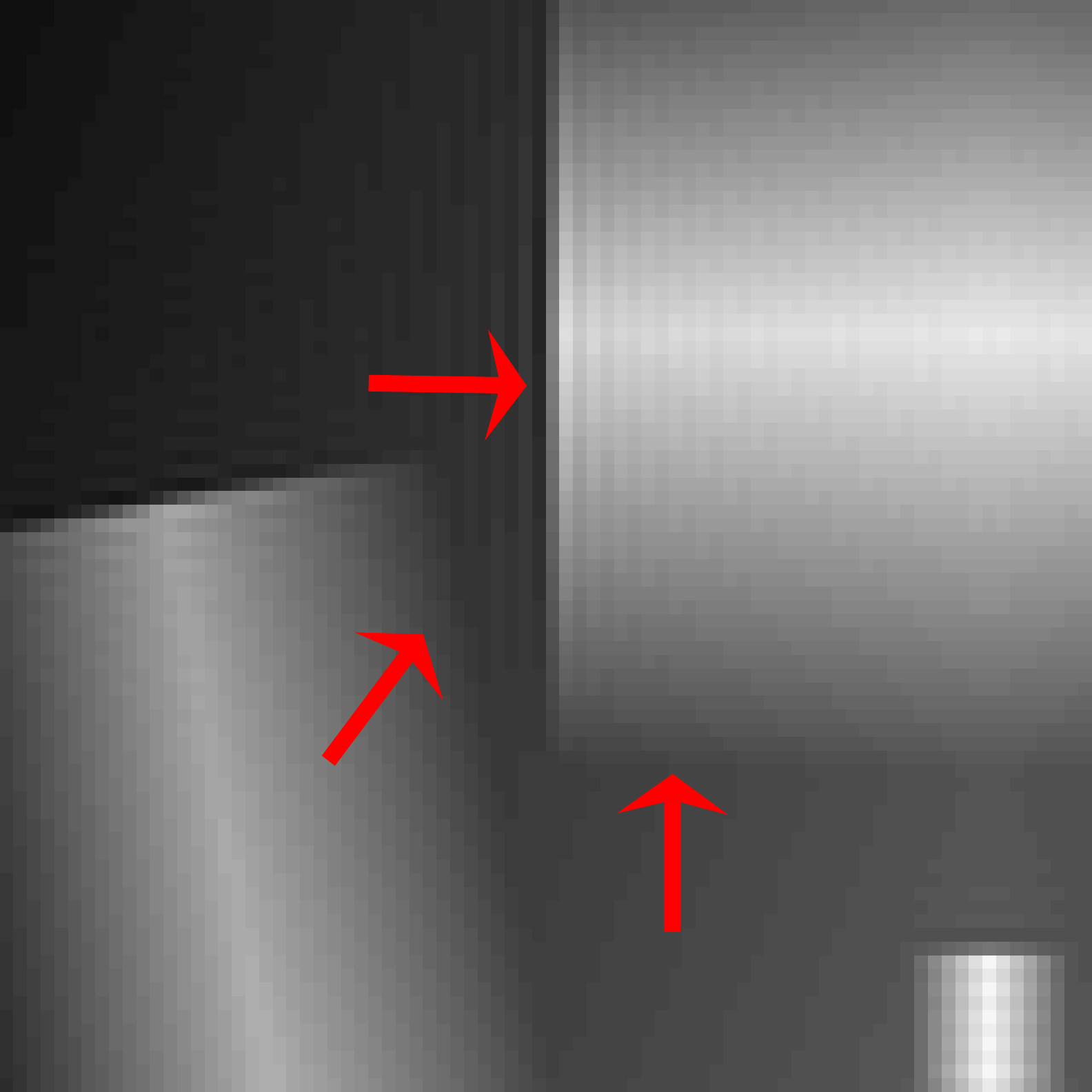}}\hspace{0.001cm}
\subfloat[]{\label{RectangleLRHTGVZoom}\includegraphics[width=2.25cm]{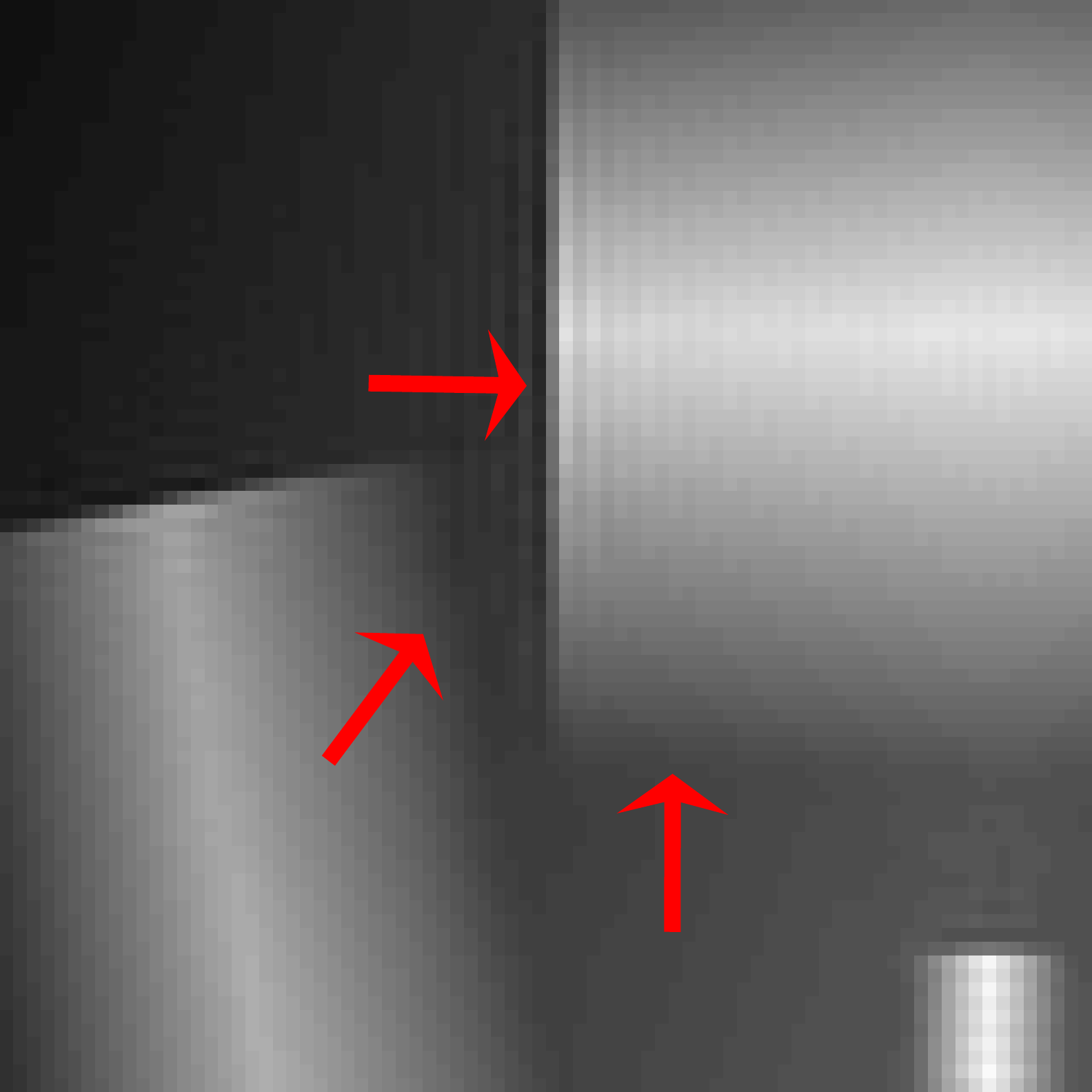}}\hspace{0.001cm}
\subfloat[]{\label{RectangleLRHTGVIRLSZoom}\includegraphics[width=2.25cm]{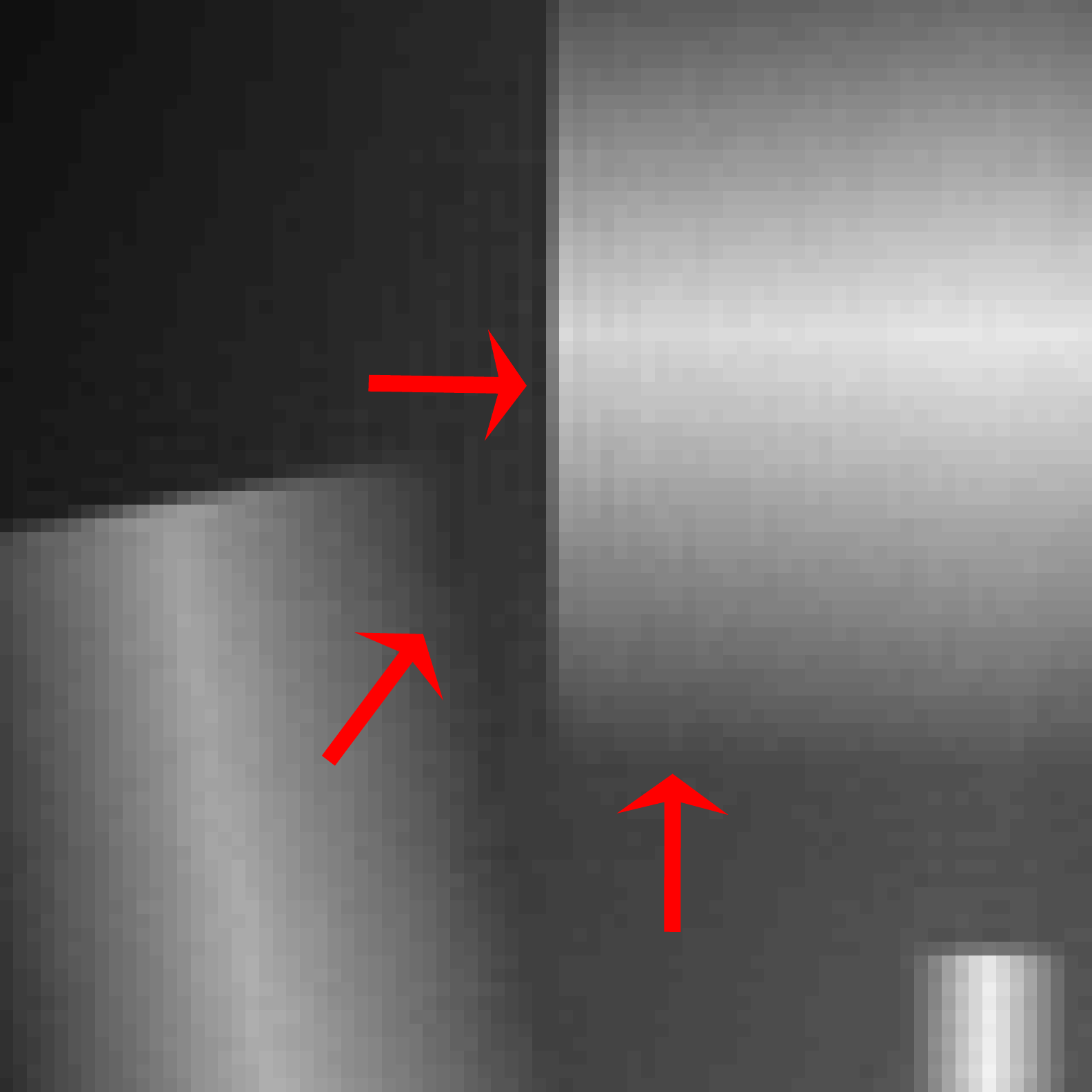}}\hspace{0.001cm}
\subfloat[]{\label{RectangleLRHInfConvIRLSZoom}\includegraphics[width=2.25cm]{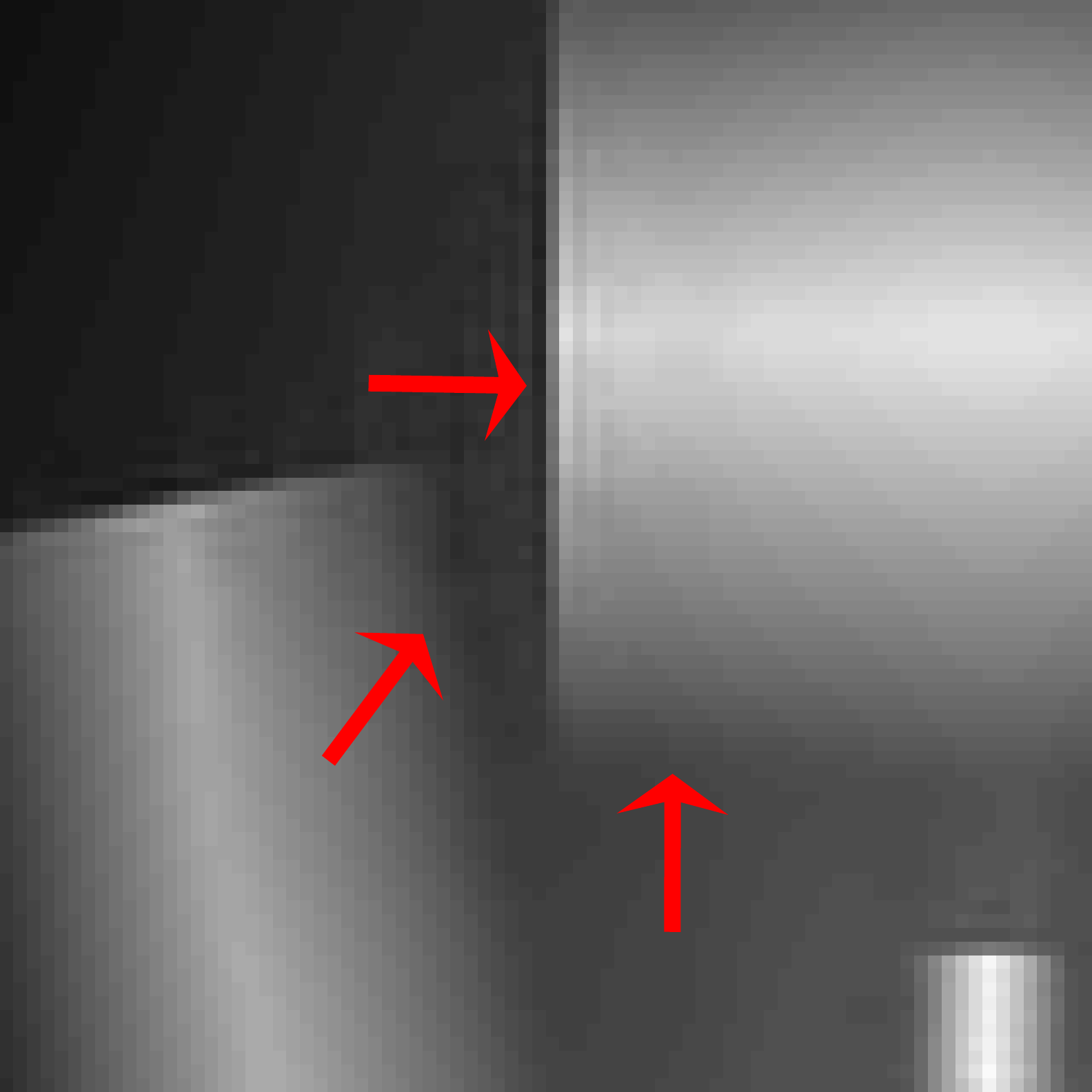}}\hspace{0.001cm}
\subfloat[]{\label{RectangleFraZoom}\includegraphics[width=2.25cm]{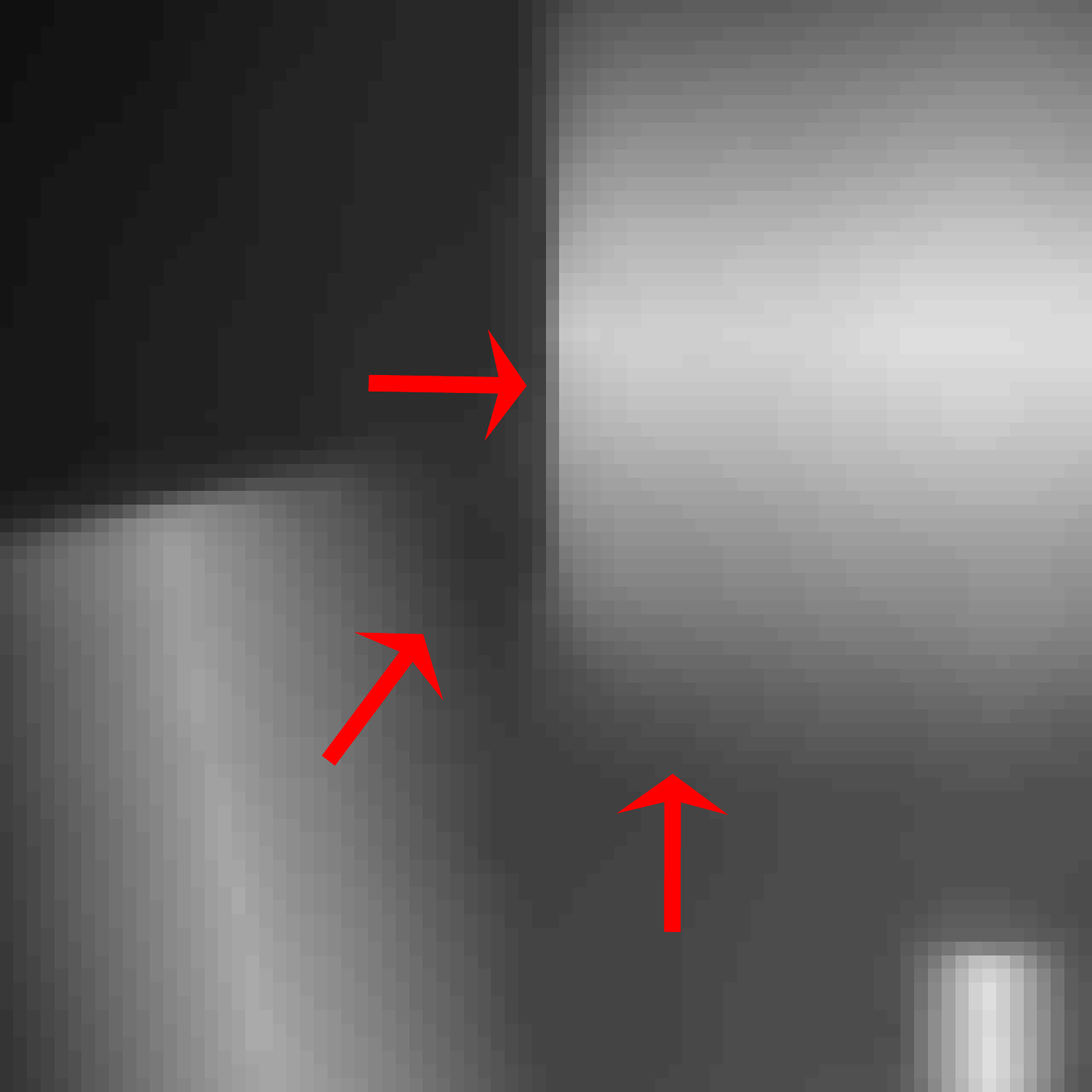}}\hspace{0.001cm}
\subfloat[]{\label{RectangleTGVZoom}\includegraphics[width=2.25cm]{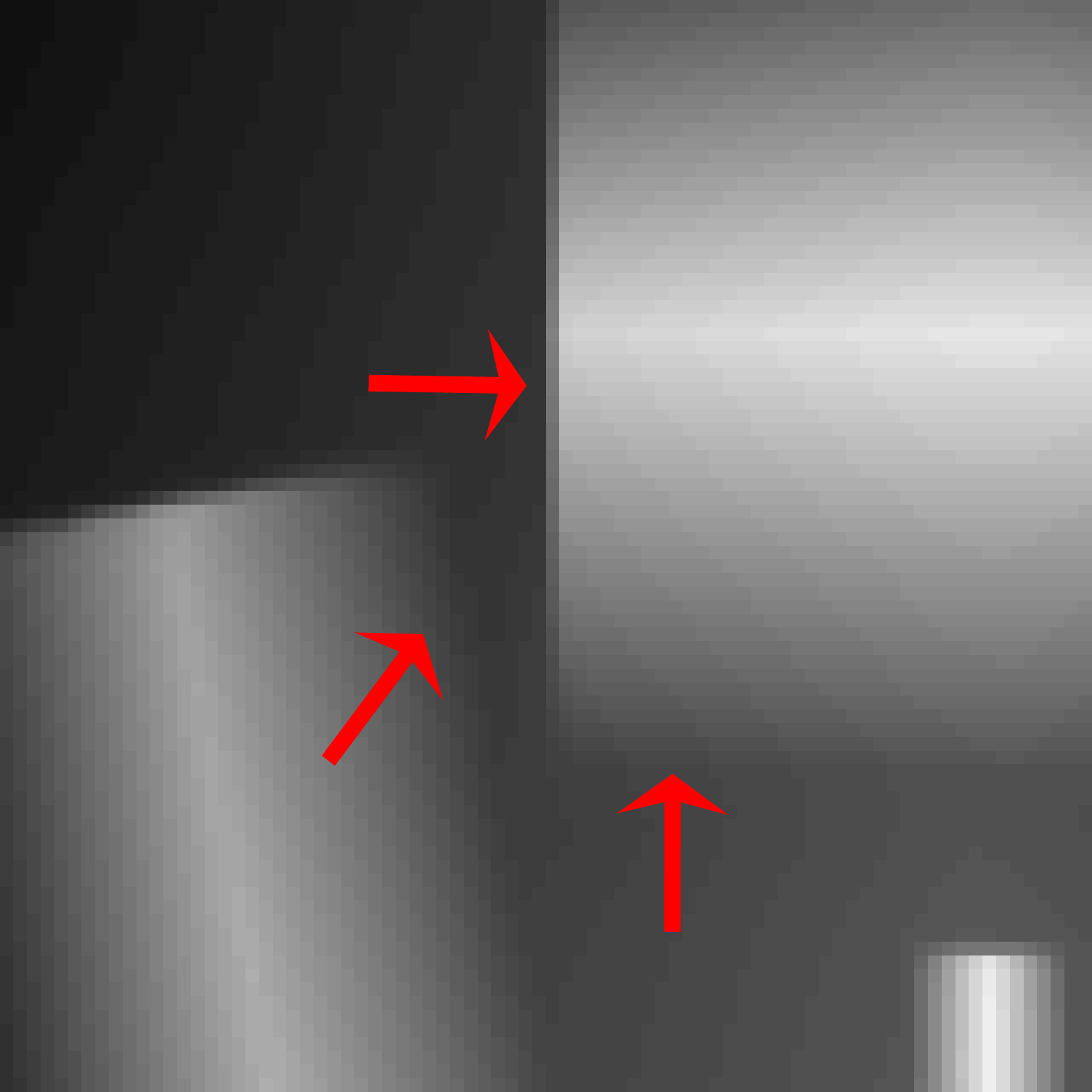}}\hspace{0.001cm}
\subfloat[]{\label{RectangleInfConvZoom}\includegraphics[width=2.25cm]{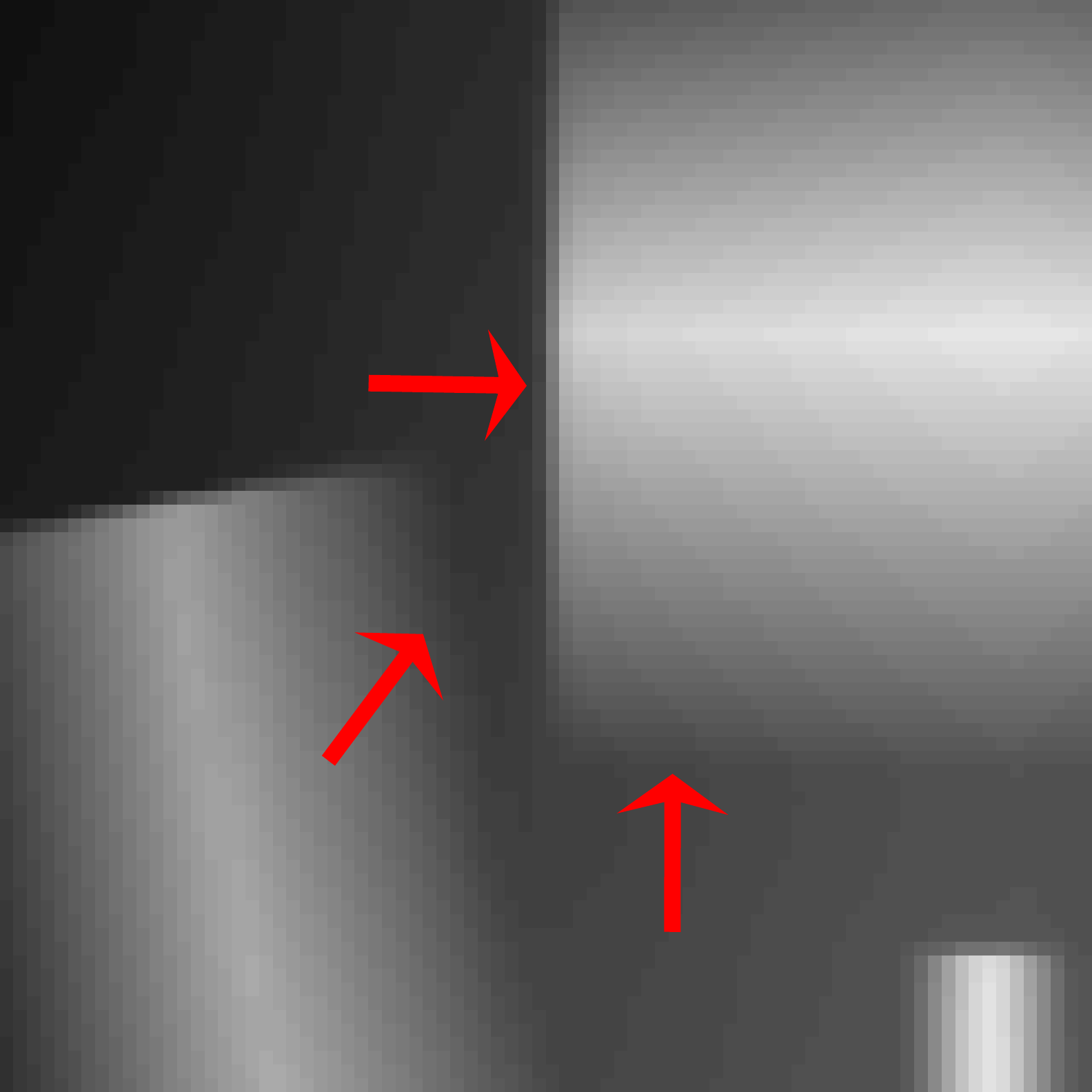}}\vspace{-0.75em}\\
\subfloat[]{\label{RectangleMask}\includegraphics[width=2.25cm]{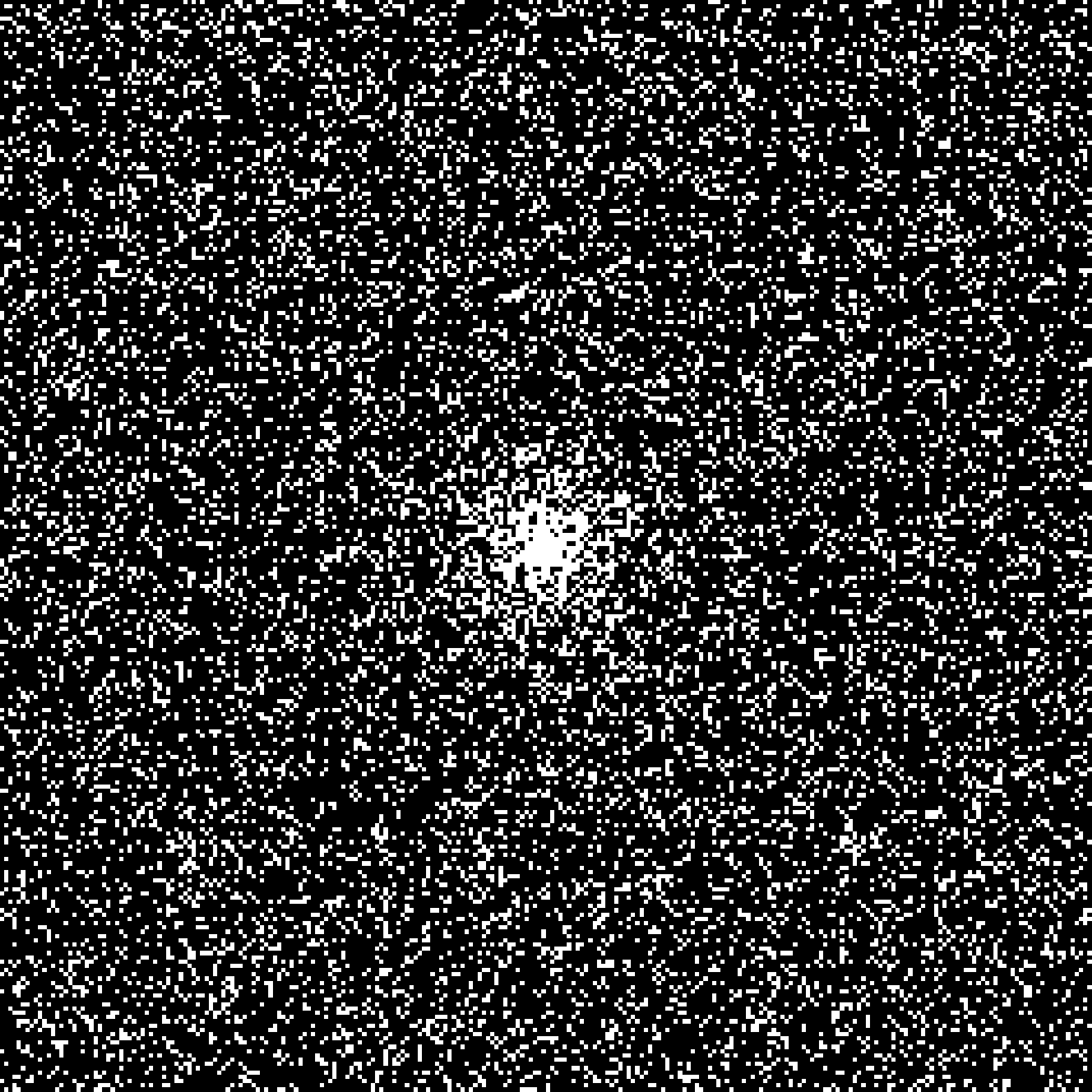}}\hspace{0.001cm}
\subfloat[]{\label{RectangleLRHTGVError}\includegraphics[width=2.25cm]{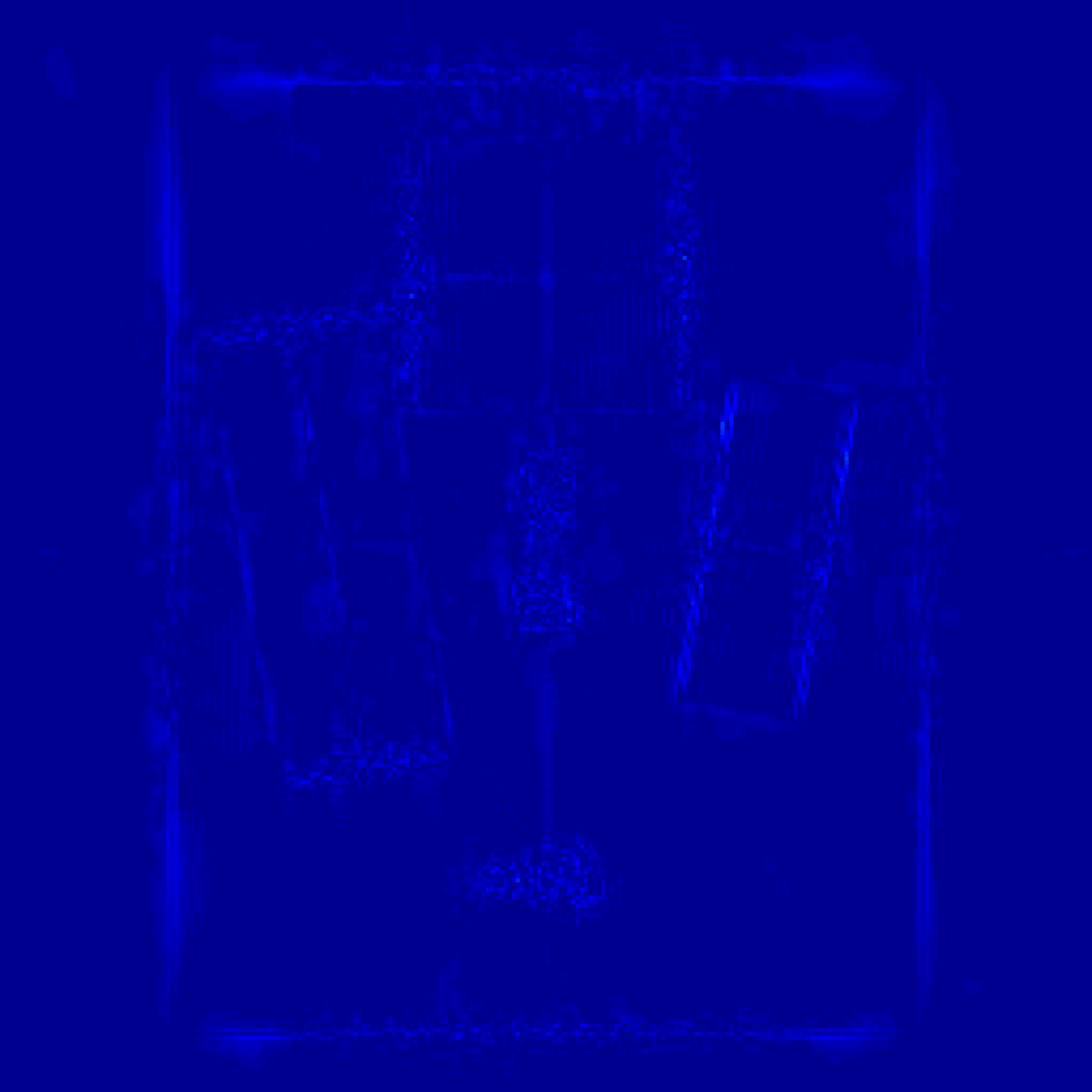}}\hspace{0.001cm}
\subfloat[]{\label{RectangleLRHTGVIRLSError}\includegraphics[width=2.25cm]{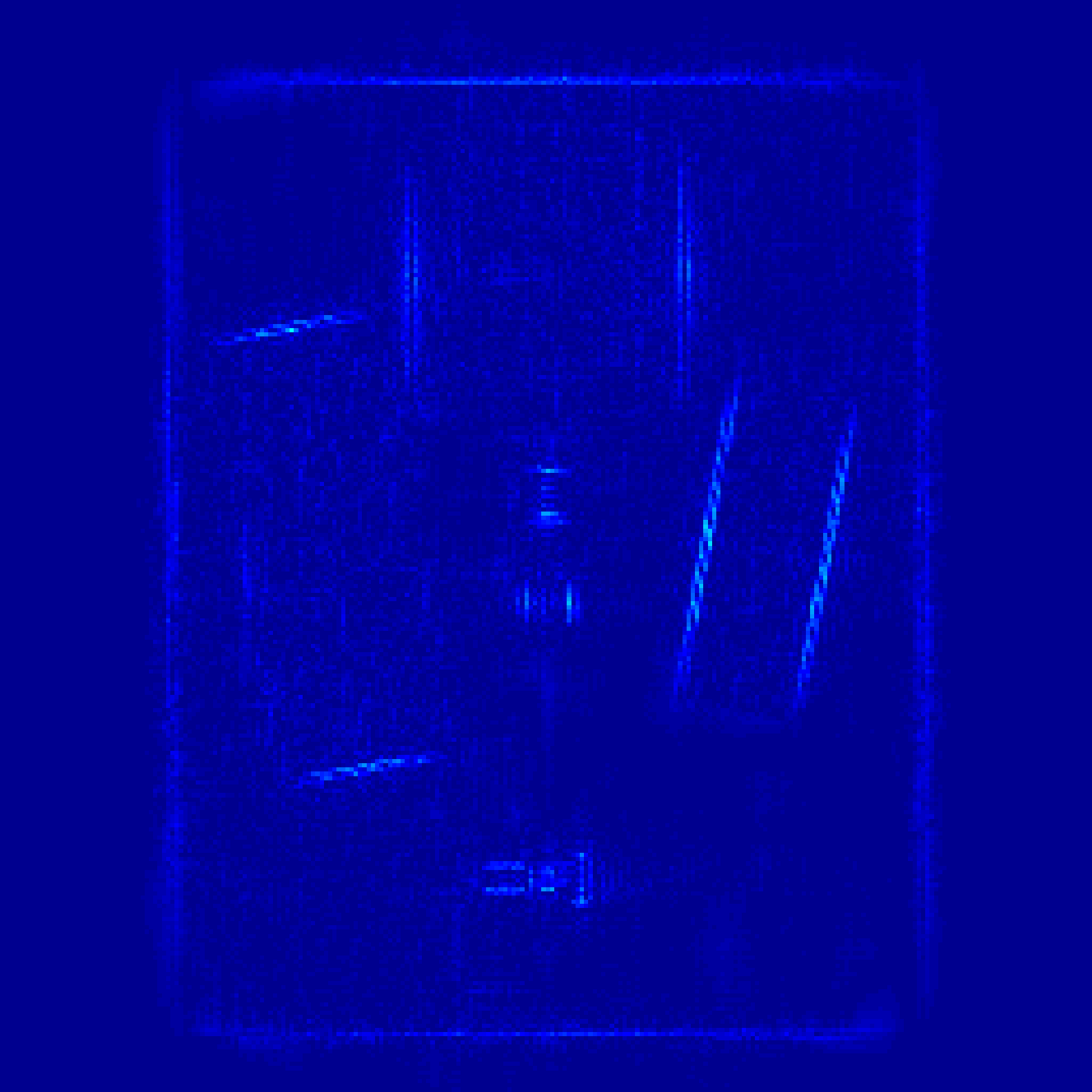}}\hspace{0.001cm}
\subfloat[]{\label{RectangleLRHInfConvIRLSError}\includegraphics[width=2.25cm]{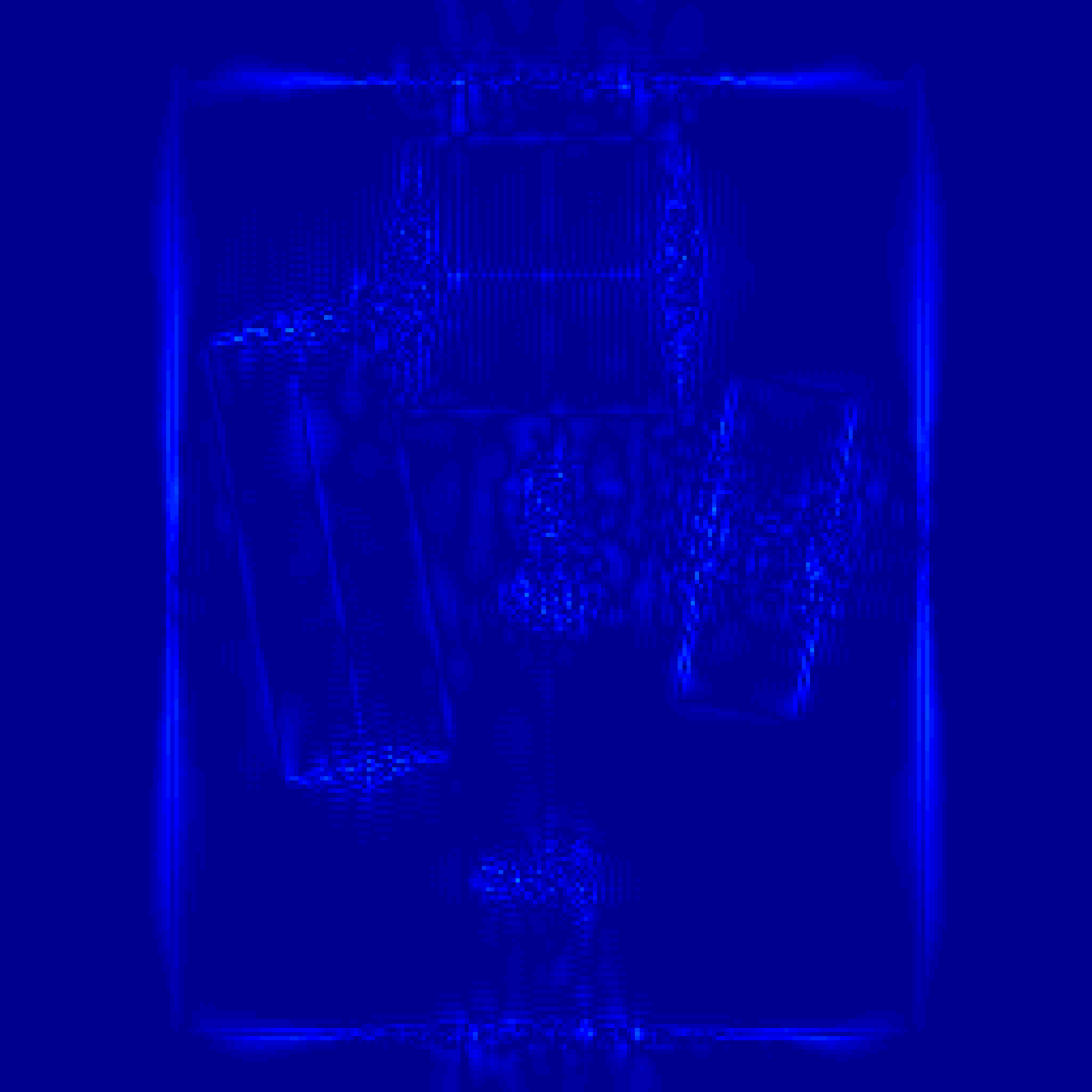}}\hspace{0.001cm}
\subfloat[]{\label{RectangleFraError}\includegraphics[width=2.25cm]{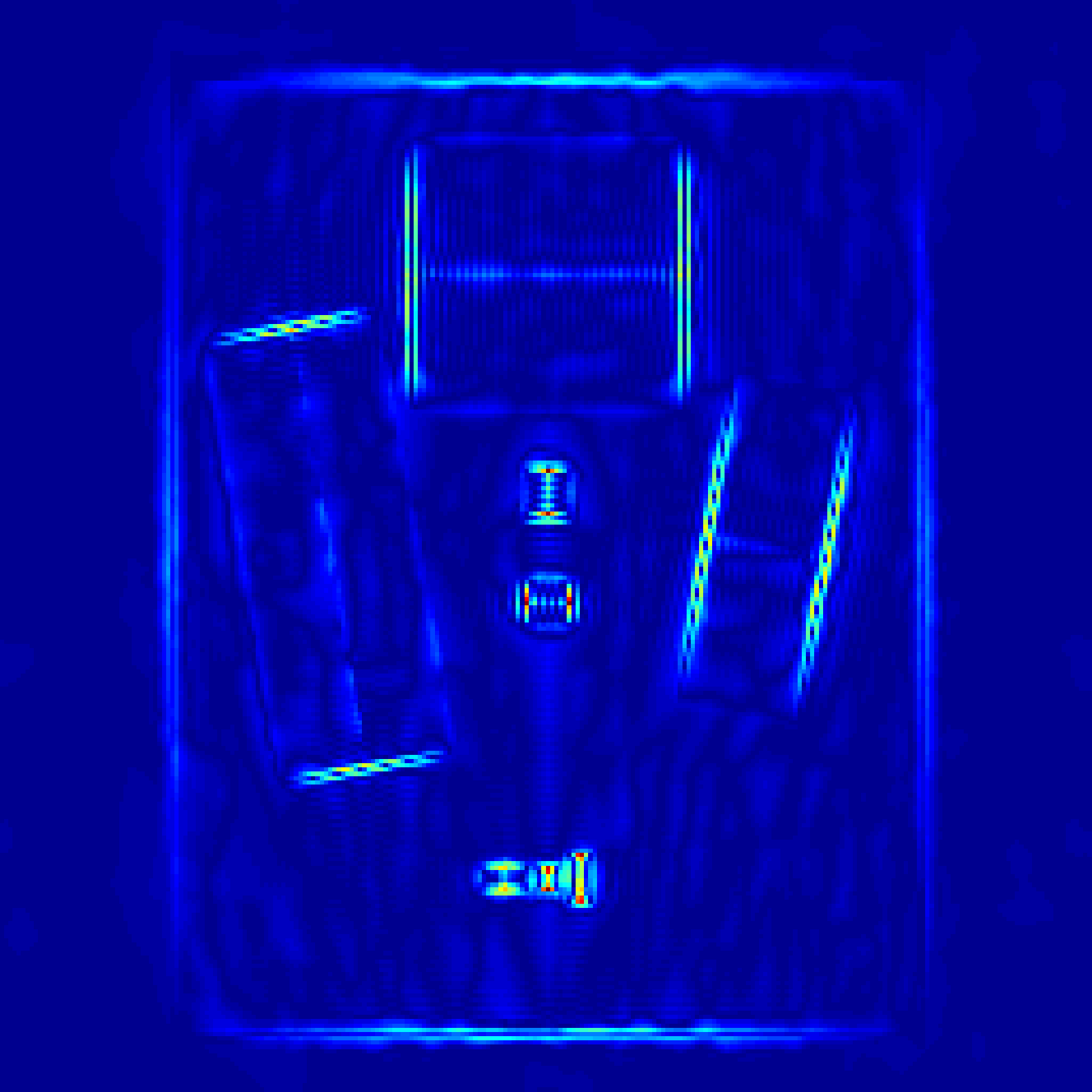}}\hspace{0.001cm}
\subfloat[]{\label{RectangleTGVError}\includegraphics[width=2.25cm]{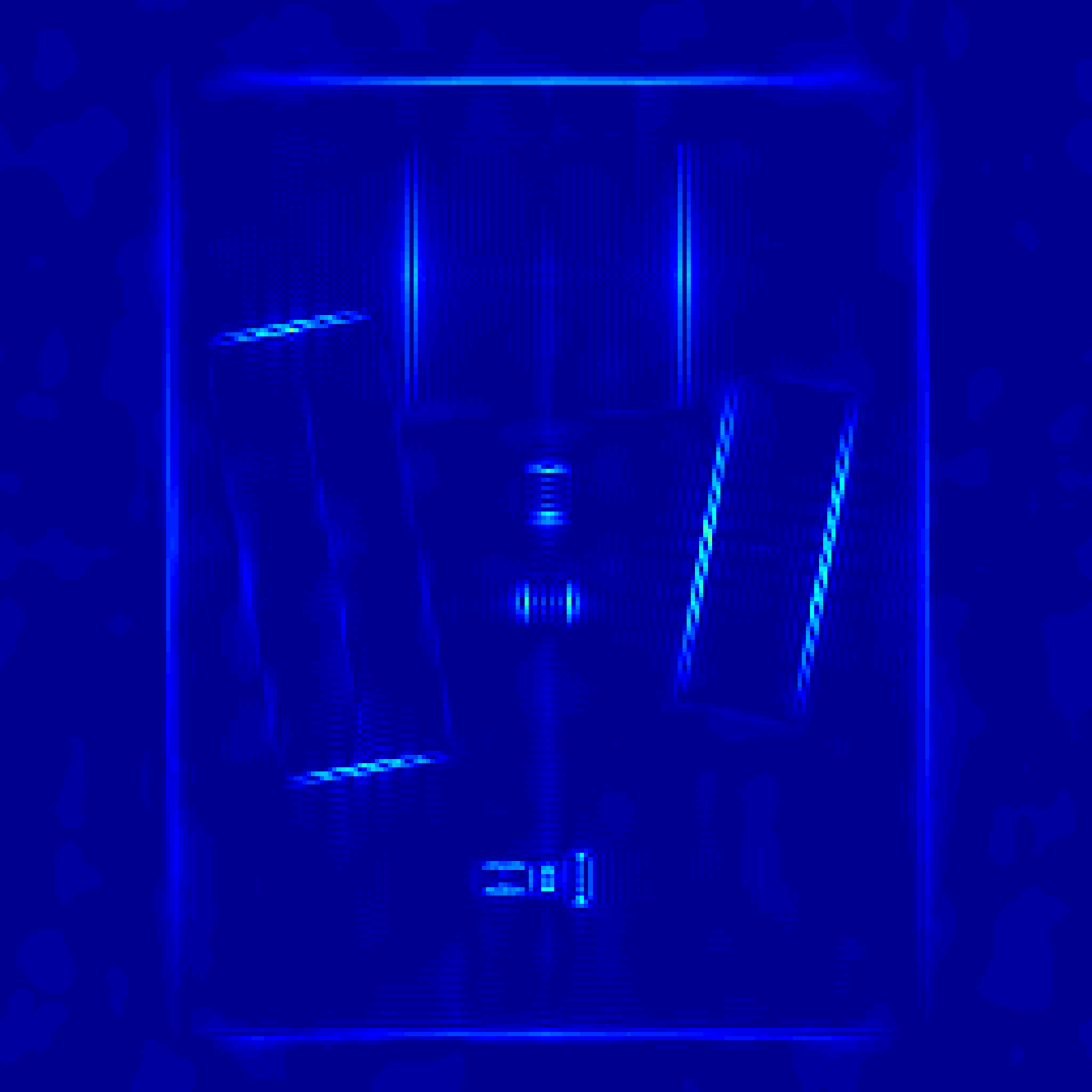}}\hspace{0.001cm}
\subfloat[]{\label{RectangleInfConvError}\includegraphics[width=2.25cm]{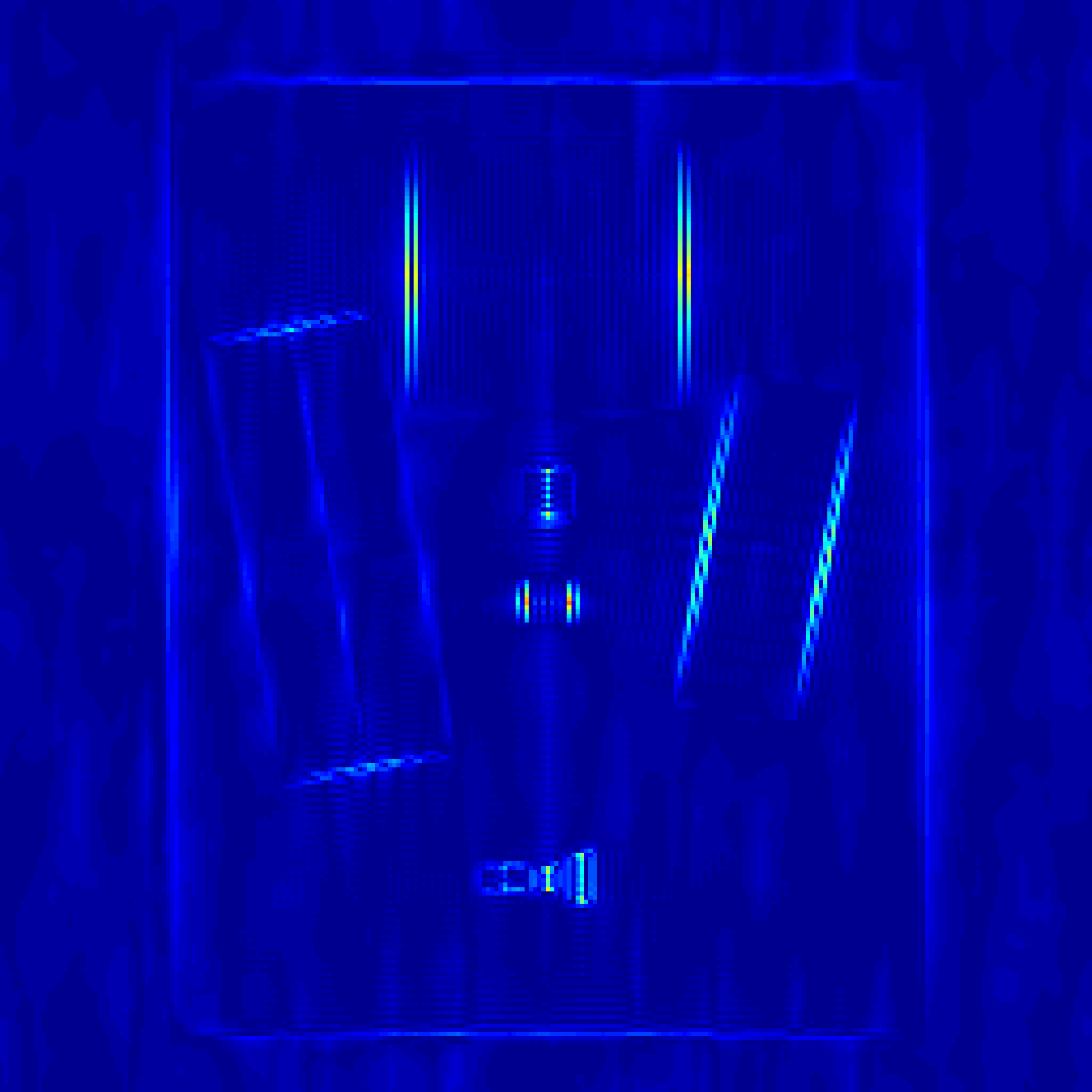}}
\caption{Visual comparisons for ``Rectangles''. \cref{RectangleOriginal}: original image, \cref{RectangleLRHTGV}: model \cref{ProposedCSMRI}, \cref{RectangleLRHTGVIRLS}: SLRM \cref{LRHTGVIRLSCSMRI}, \cref{RectangleLRHInfConvIRLS}: GSLR \cref{LRHInfConvIRLSCSMRI}, \cref{RectangleFra}: framelet \cref{FrameCSMRI}, \cref{RectangleTGV}: TGV \cref{TGVCSMRI}, \cref{RectangleInfConv}: infimal convolution \cref{InfConvCSMRI}. \cref{RectangleOriginalZoom,RectangleLRHTGVZoom,RectangleLRHTGVIRLSZoom,RectangleLRHInfConvIRLSZoom,RectangleFraZoom,RectangleTGVZoom,RectangleInfConvZoom}: zoom-in views of \cref{RectangleOriginal,RectangleLRHTGV,RectangleLRHTGVIRLS,RectangleLRHInfConvIRLS,RectangleFra,RectangleTGV,RectangleInfConv}. Red arrows indicate the region worth noting. \cref{RectangleMask}: sample region, \cref{RectangleLRHTGVError,RectangleLRHTGVIRLSError,RectangleLRHInfConvIRLSError,RectangleFraError,RectangleTGVError,RectangleInfConvError}: error maps of \cref{RectangleLRHTGV,RectangleLRHTGVIRLS,RectangleLRHInfConvIRLS,RectangleFra,RectangleTGV,RectangleInfConv}.}\label{RectangleResults}
\end{figure}

\begin{figure}[t]
\centering
\subfloat[]{\label{AirplaneOriginal}\includegraphics[width=2.25cm]{AirplaneOriginal.pdf}}\hspace{0.001cm}
\subfloat[]{\label{AirplaneLRHTGV}\includegraphics[width=2.25cm]{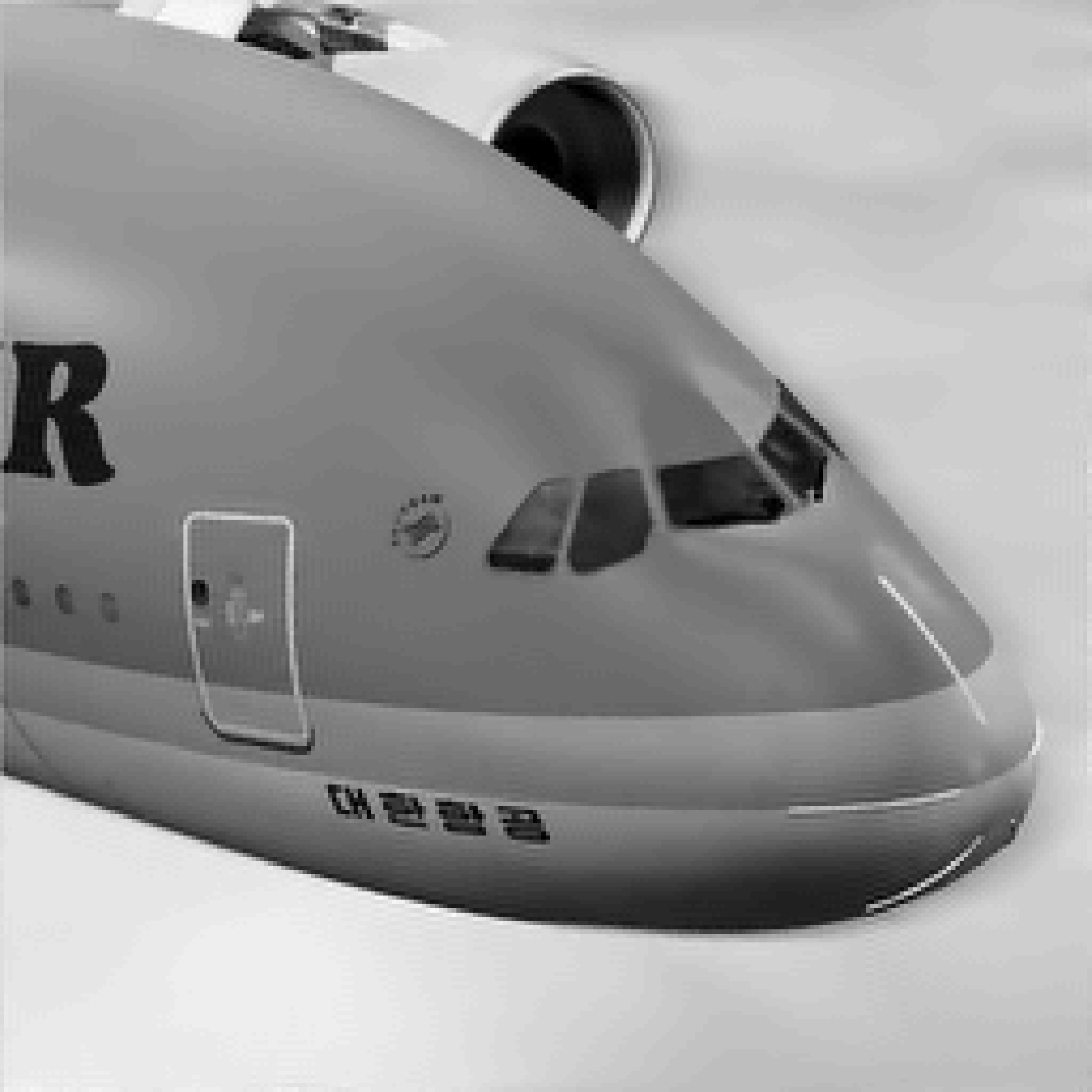}}\hspace{0.001cm}
\subfloat[]{\label{AirplaneLRHTGVIRLS}\includegraphics[width=2.25cm]{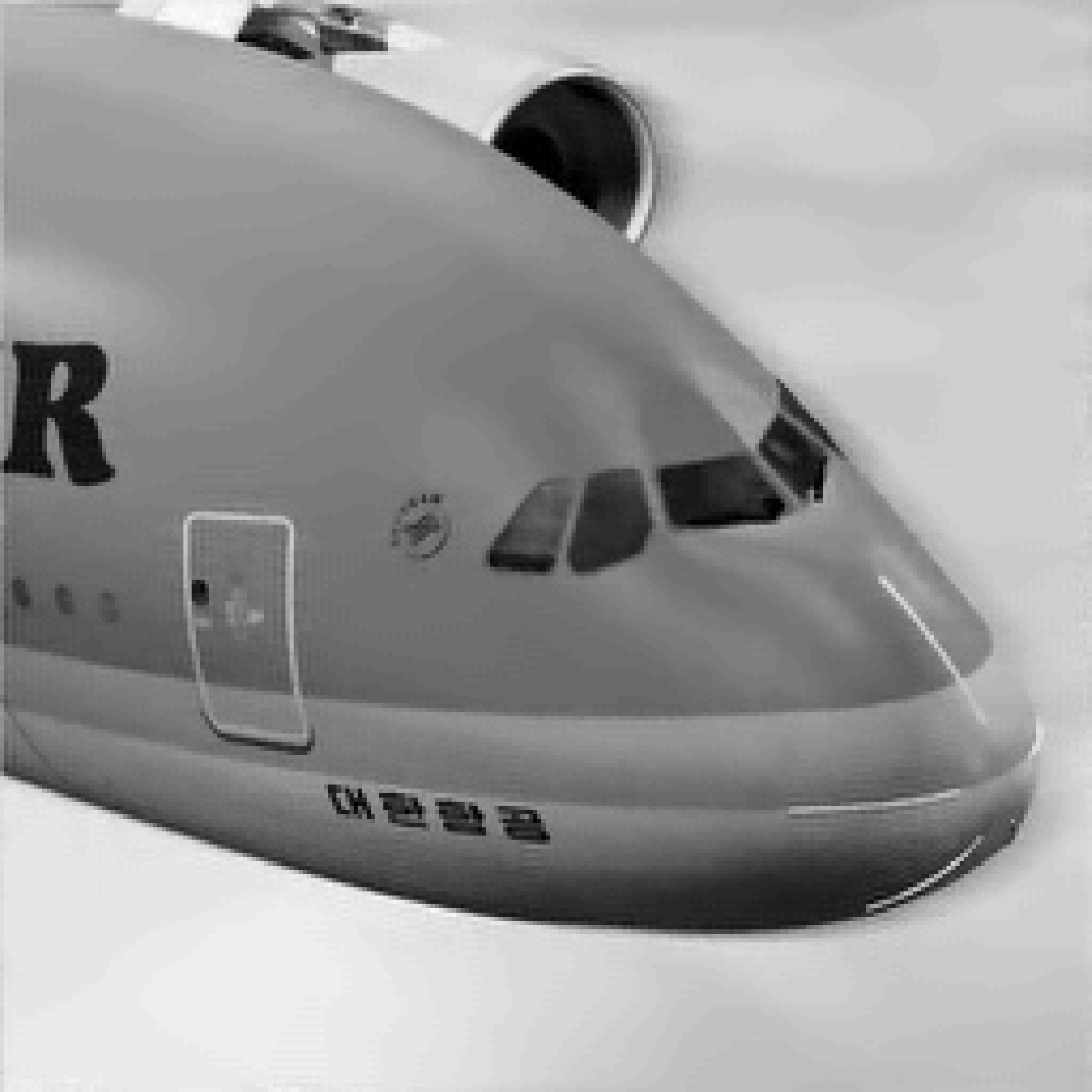}}\hspace{0.001cm}
\subfloat[]{\label{AirplaneLRHInfConvIRLS}\includegraphics[width=2.25cm]{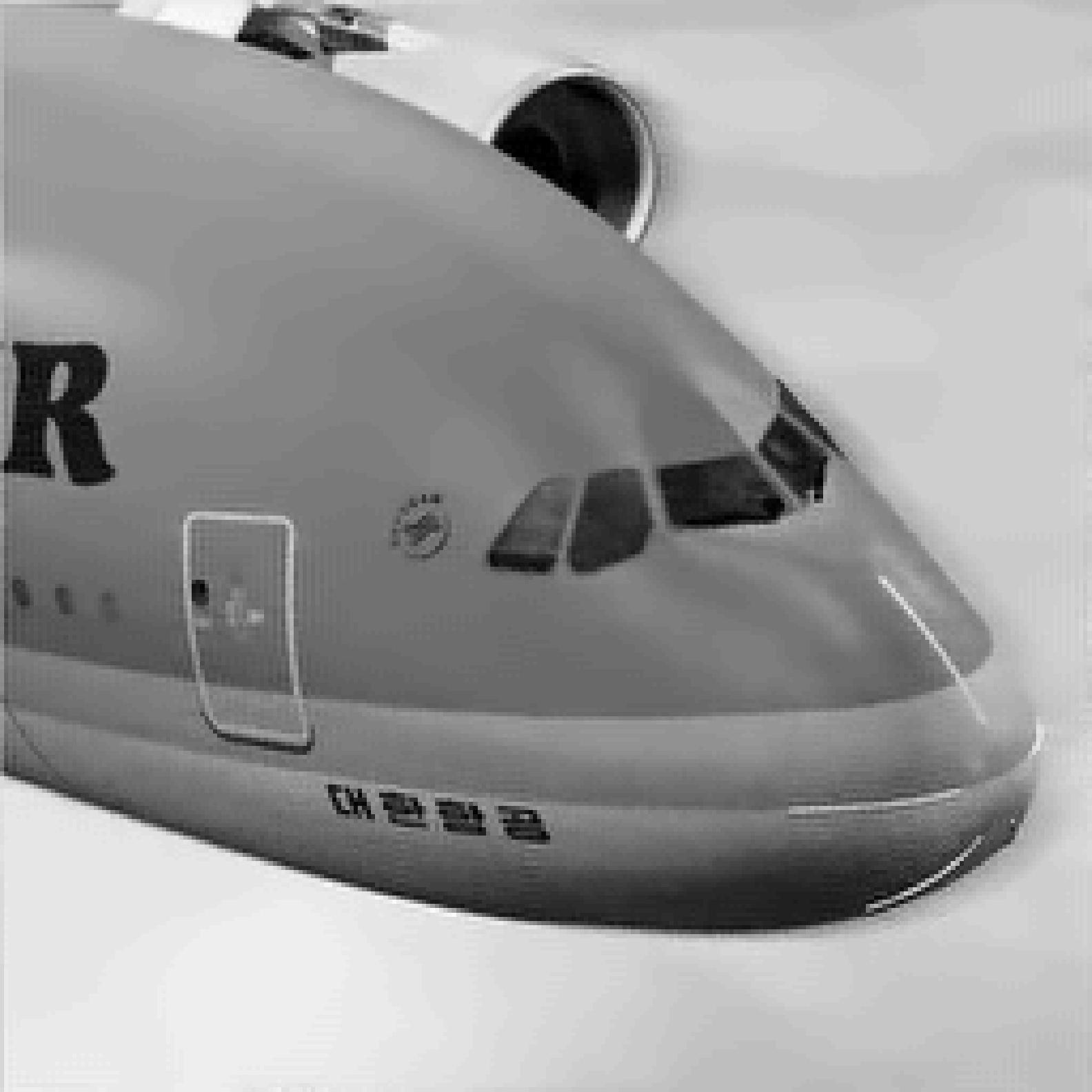}}\hspace{0.001cm}
\subfloat[]{\label{AirplaneFra}\includegraphics[width=2.25cm]{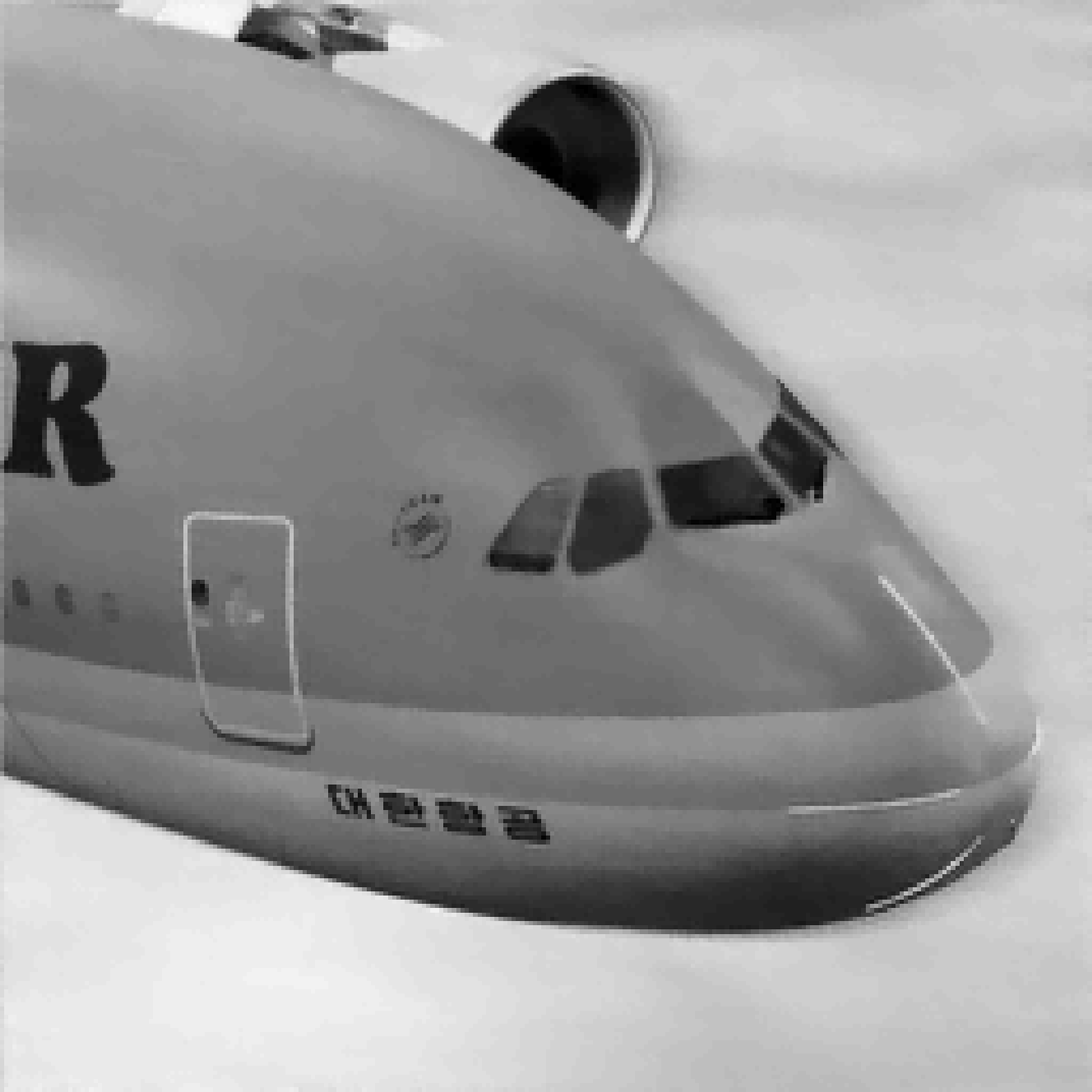}}\hspace{0.001cm}
\subfloat[]{\label{AirplaneTGV}\includegraphics[width=2.25cm]{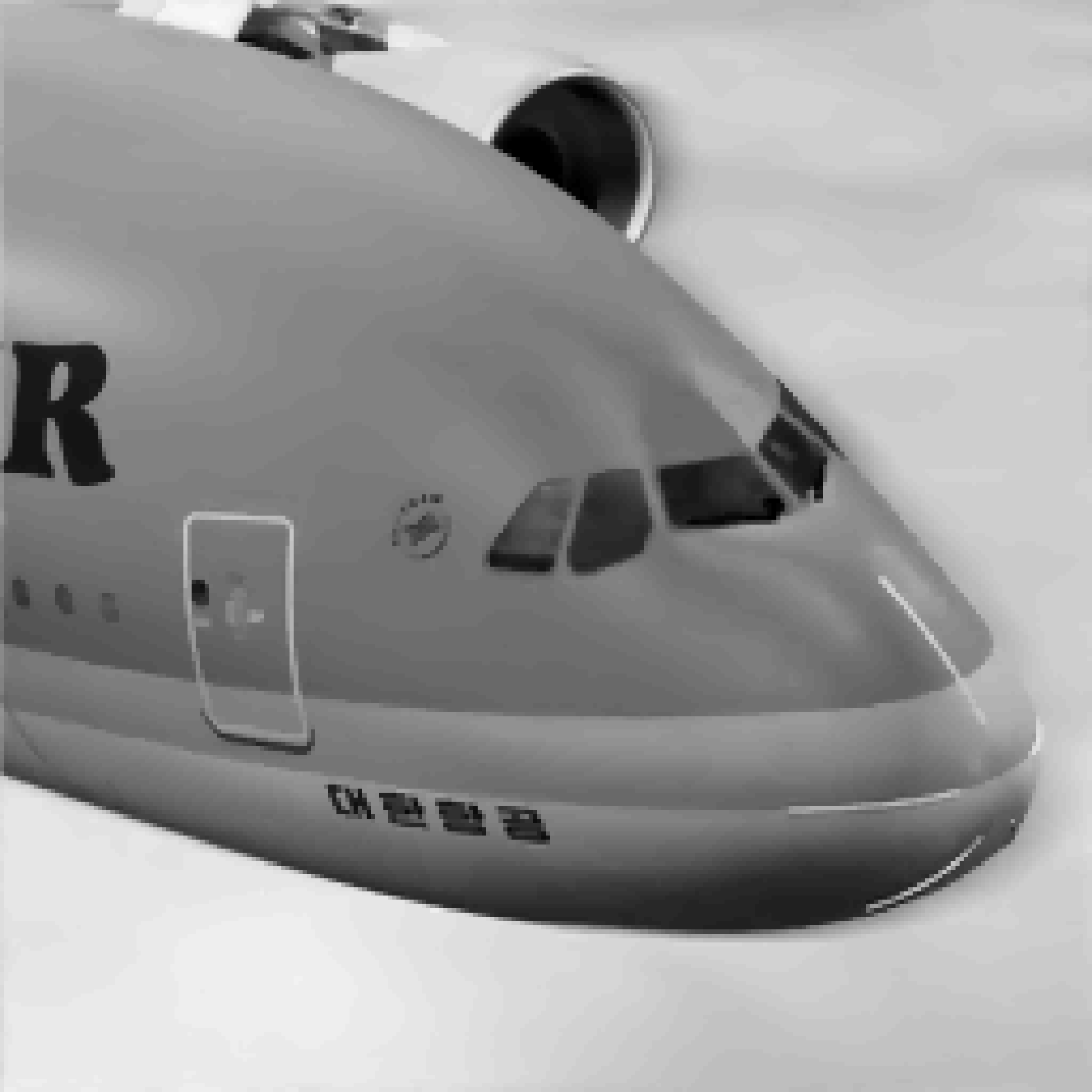}}\hspace{0.001cm}
\subfloat[]{\label{AirplaneInfConv}\includegraphics[width=2.25cm]{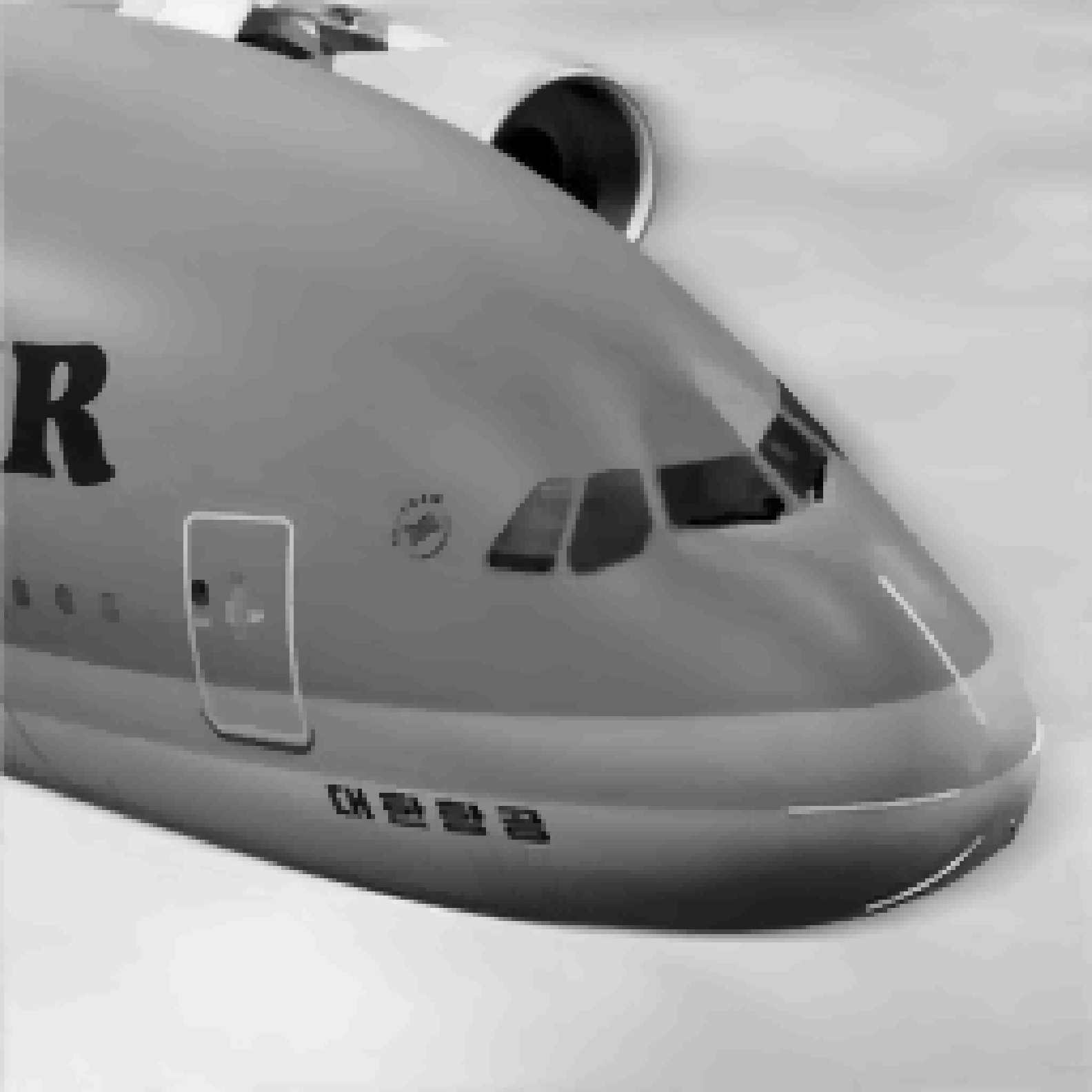}}\vspace{-0.75em}\\
\subfloat[]{\label{AirplaneOriginalZoom}\includegraphics[width=2.25cm]{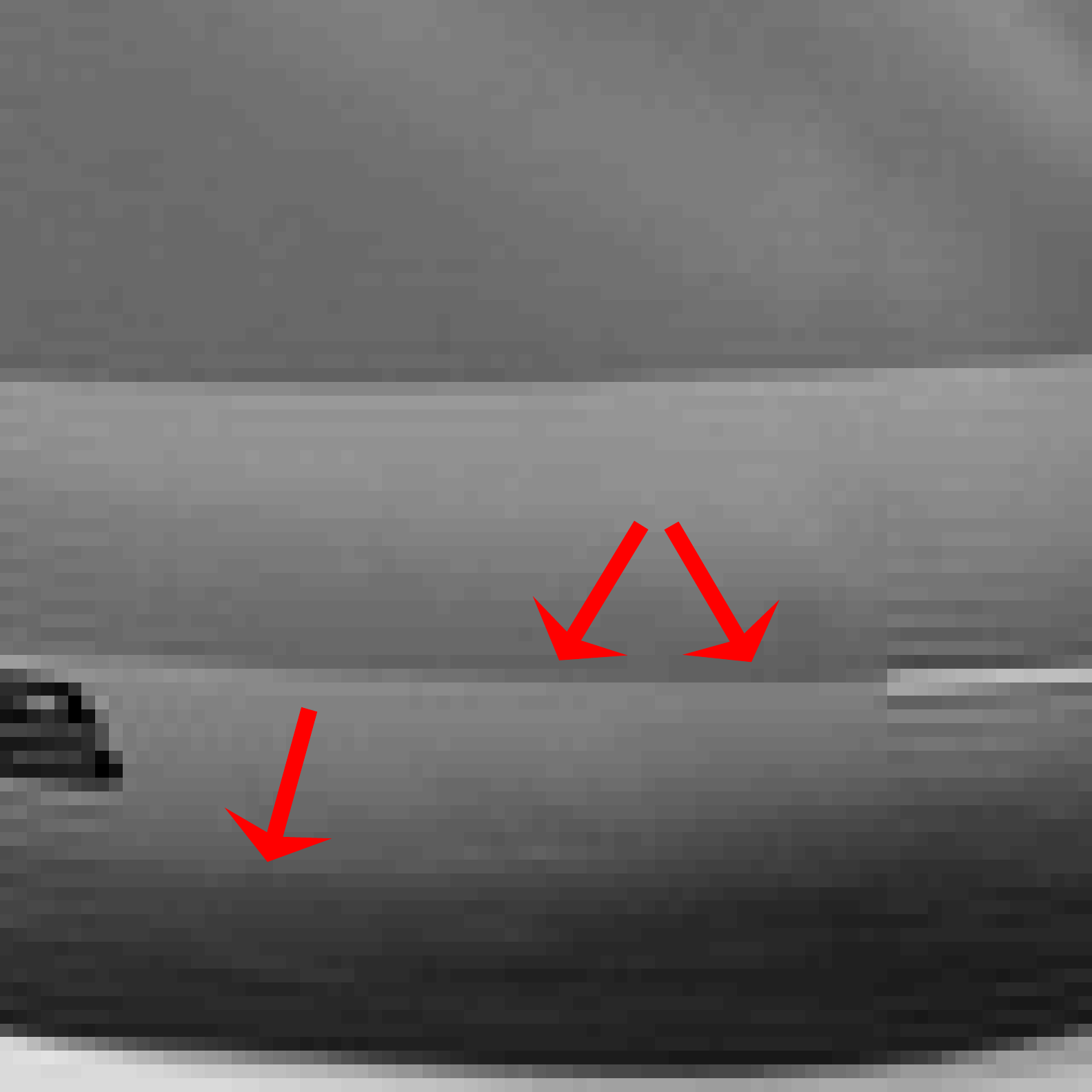}}\hspace{0.001cm}
\subfloat[]{\label{AirplaneLRHTGVZoom}\includegraphics[width=2.25cm]{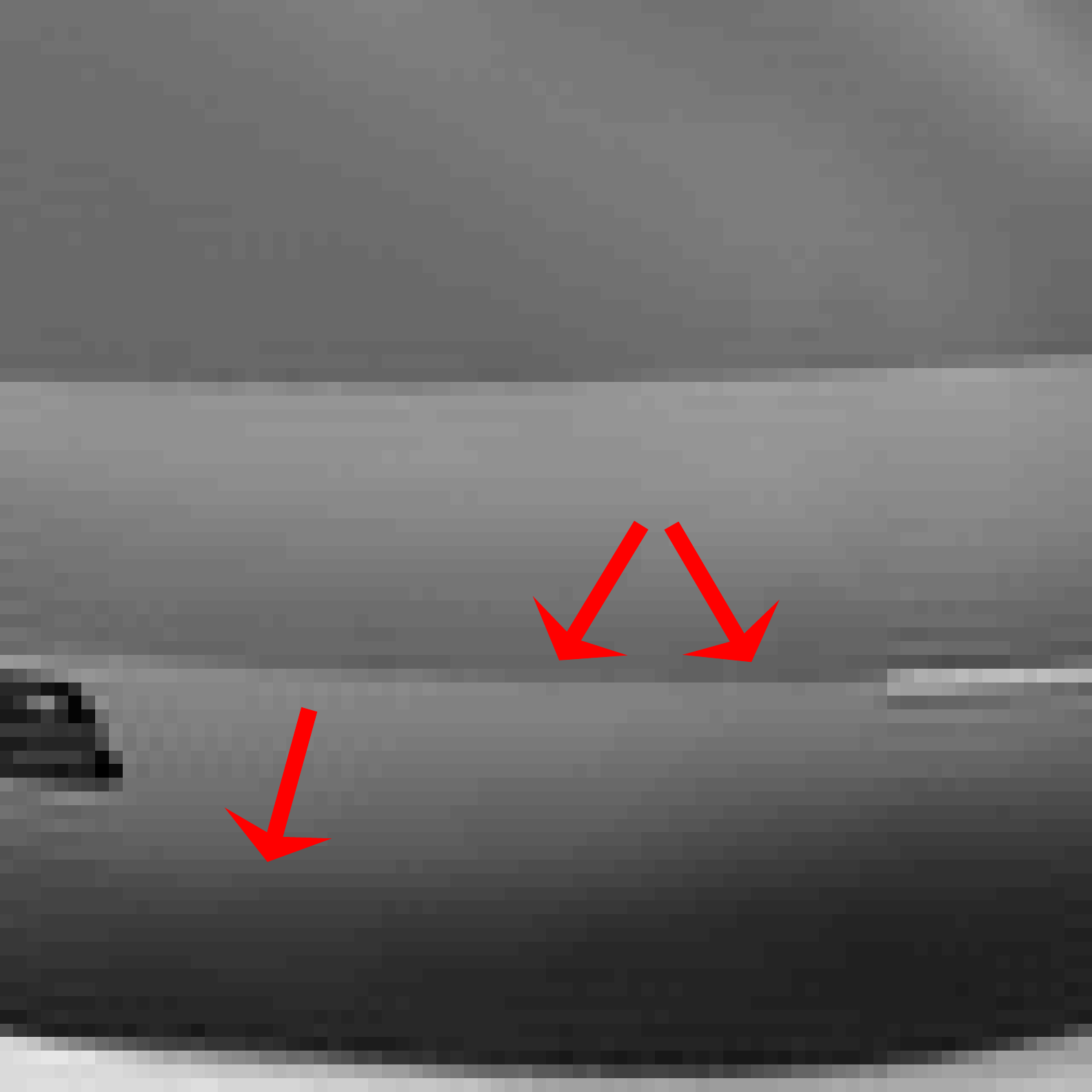}}\hspace{0.001cm}
\subfloat[]{\label{AirplaneLRHTGVIRLSZoom}\includegraphics[width=2.25cm]{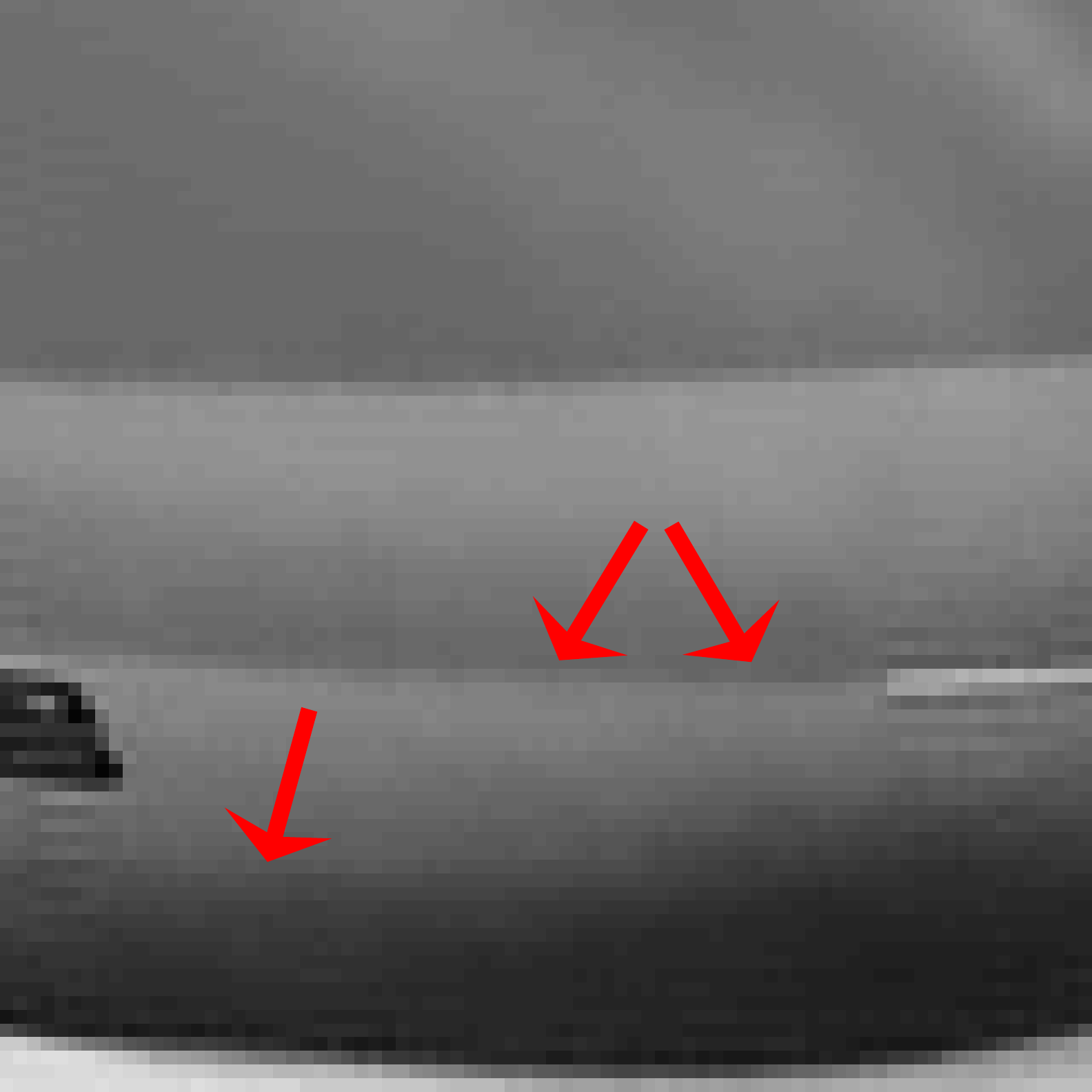}}\hspace{0.001cm}
\subfloat[]{\label{AirplaneLRHInfConvIRLSZoom}\includegraphics[width=2.25cm]{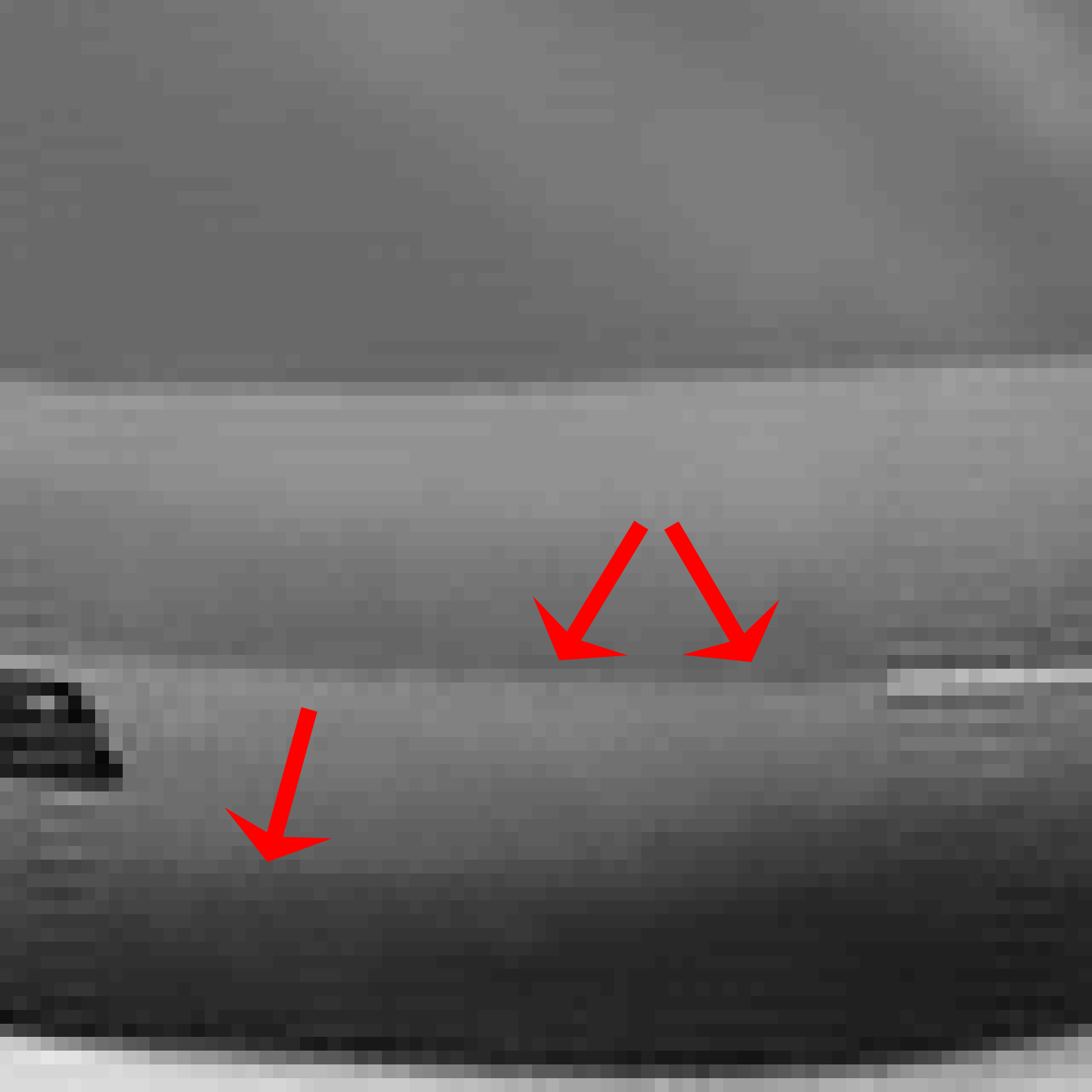}}\hspace{0.001cm}
\subfloat[]{\label{AirplaneFraZoom}\includegraphics[width=2.25cm]{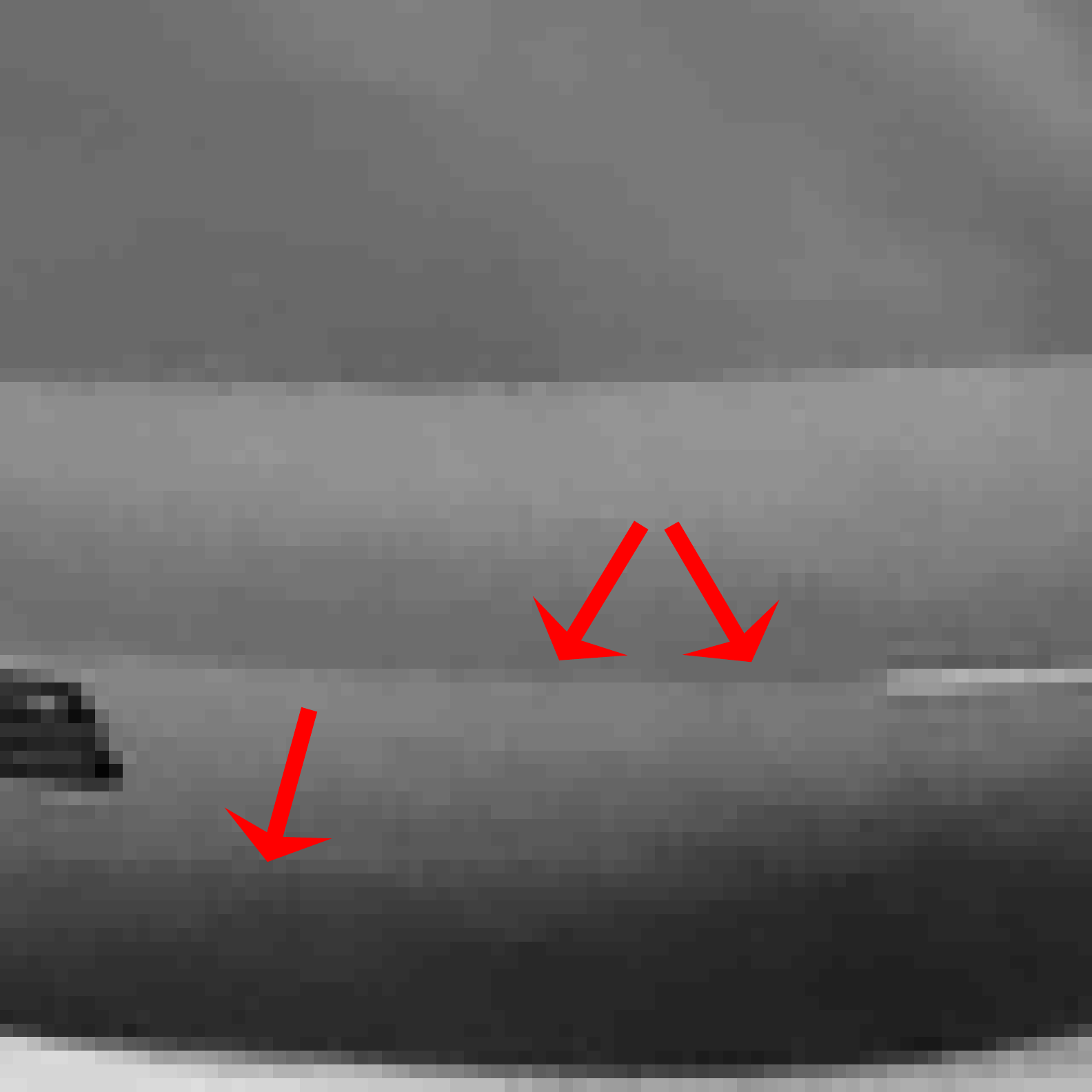}}\hspace{0.001cm}
\subfloat[]{\label{AirplaneTGVZoom}\includegraphics[width=2.25cm]{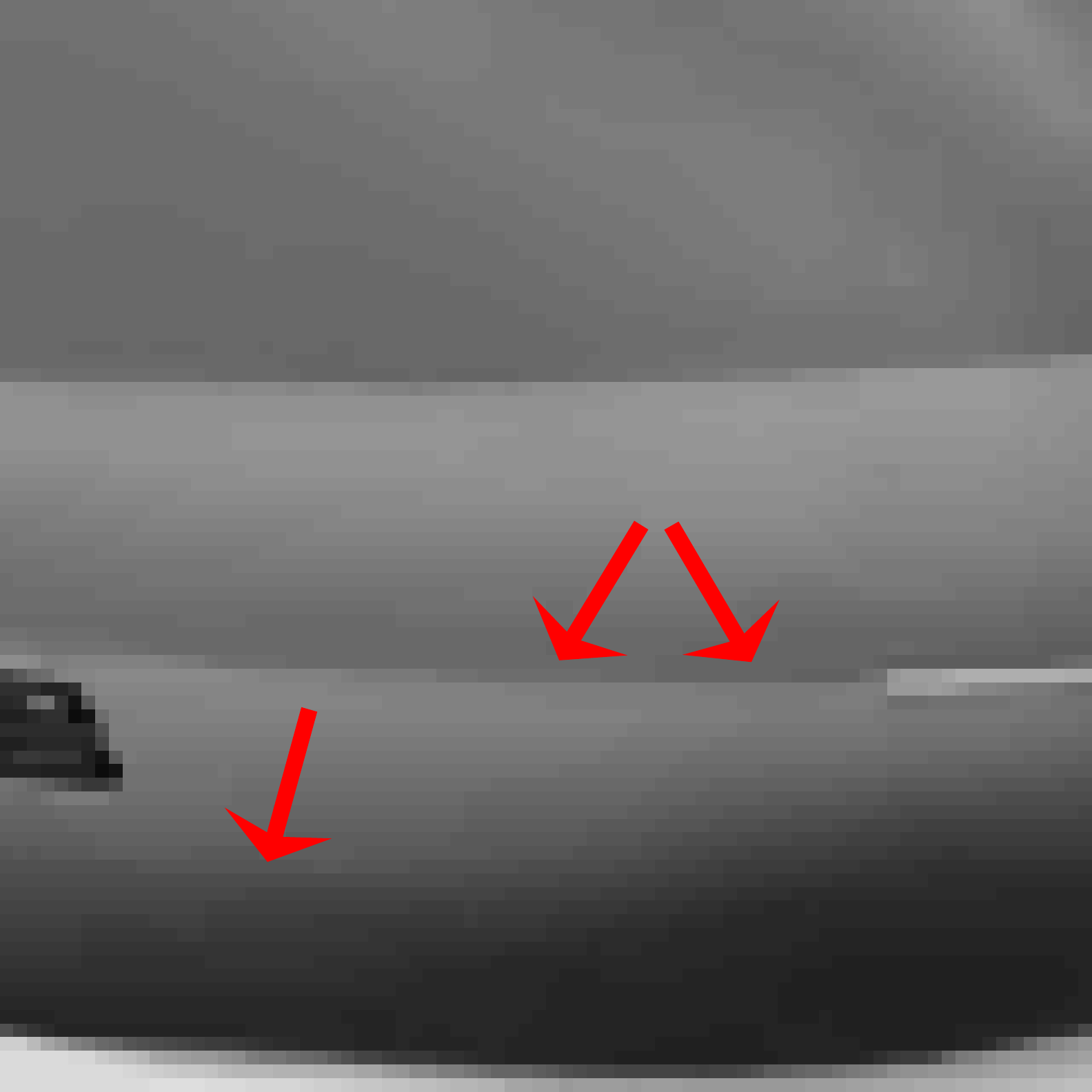}}\hspace{0.001cm}
\subfloat[]{\label{AirplaneInfConvZoom}\includegraphics[width=2.25cm]{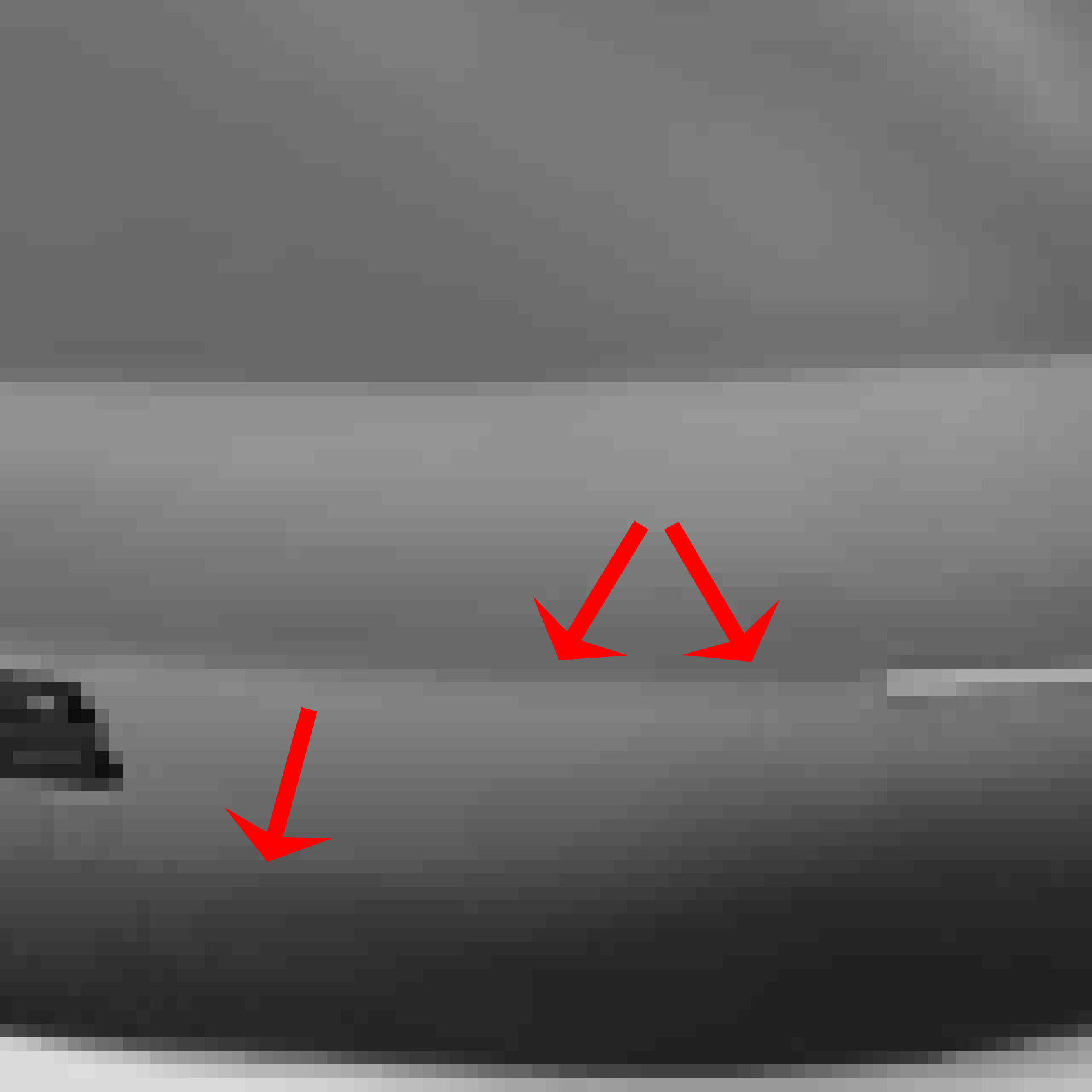}}\vspace{-0.75em}\\
\subfloat[]{\label{AirplaneMask}\includegraphics[width=2.25cm]{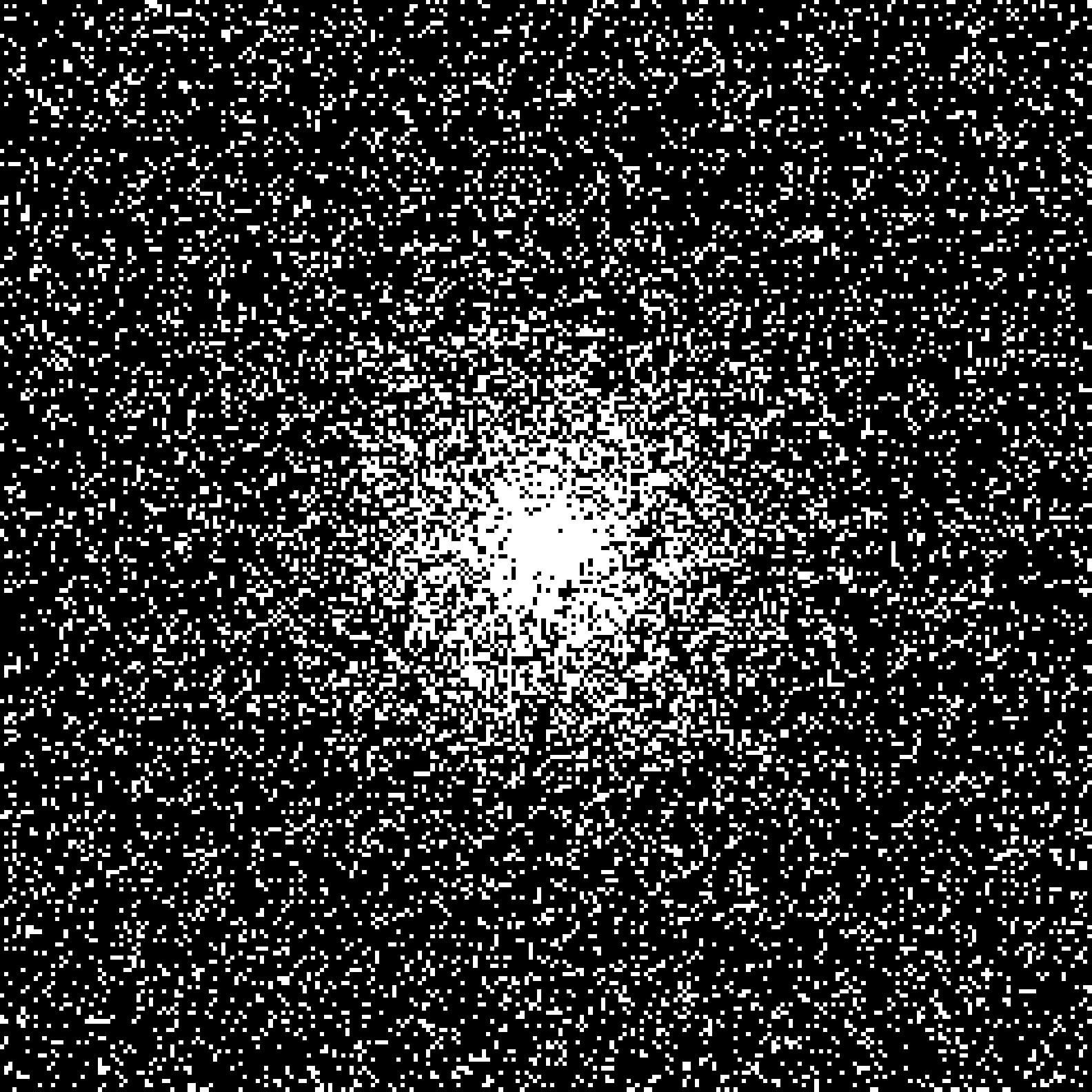}}\hspace{0.001cm}
\subfloat[]{\label{AirplaneLRHTGVError}\includegraphics[width=2.25cm]{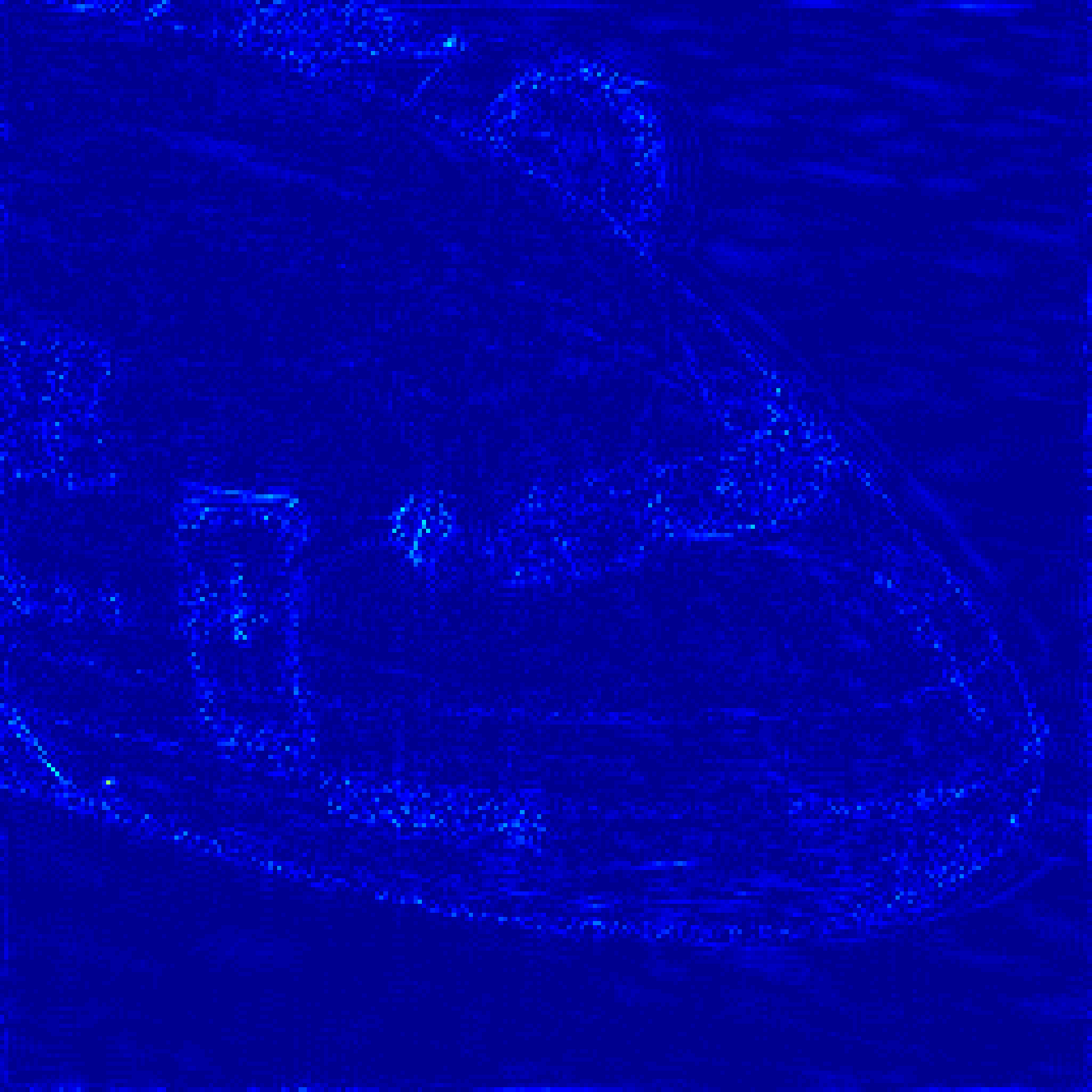}}\hspace{0.001cm}
\subfloat[]{\label{AirplaneLRHTGVIRLSError}\includegraphics[width=2.25cm]{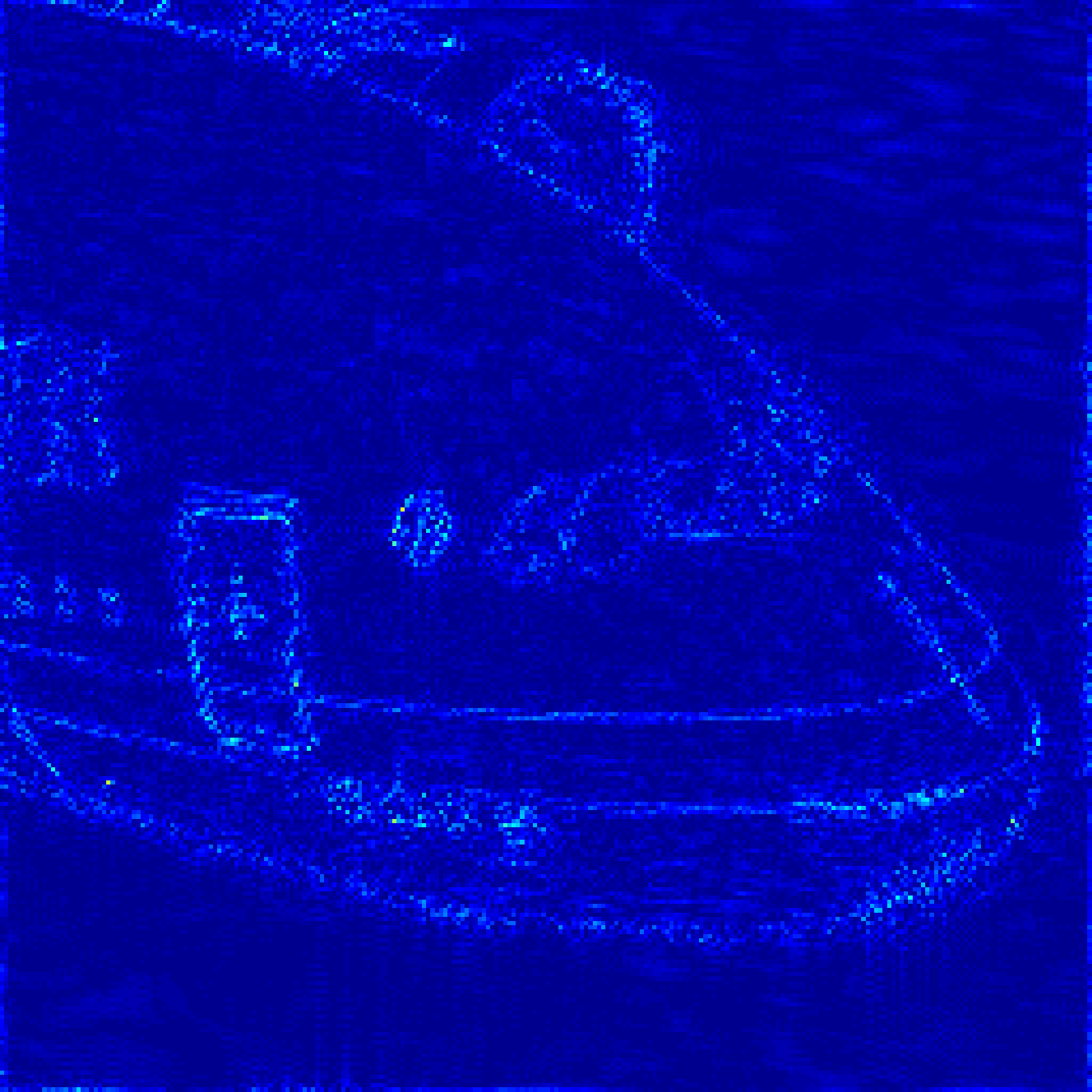}}\hspace{0.001cm}
\subfloat[]{\label{AirplaneLRHInfConvIRLSError}\includegraphics[width=2.25cm]{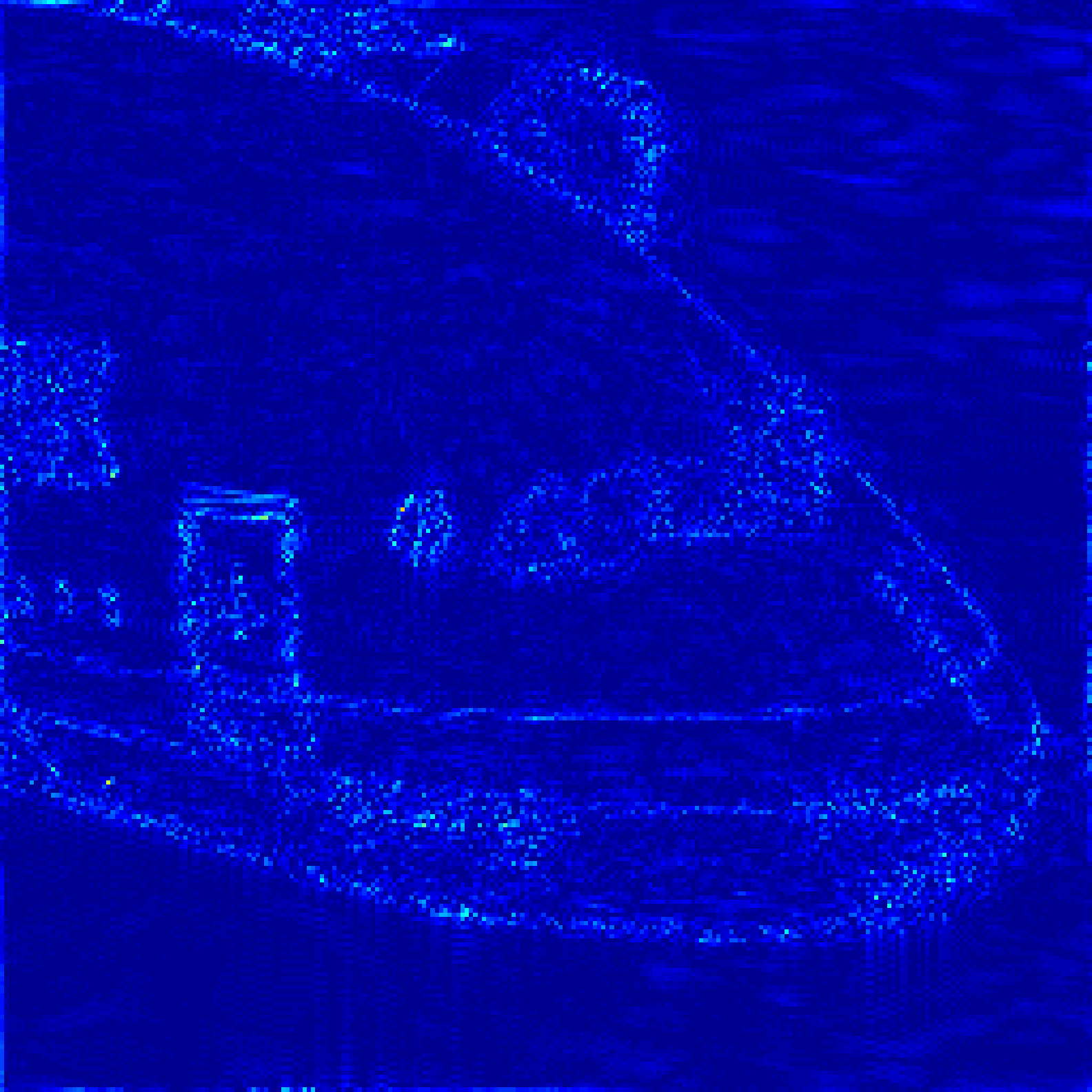}}\hspace{0.001cm}
\subfloat[]{\label{AirplaneFraError}\includegraphics[width=2.25cm]{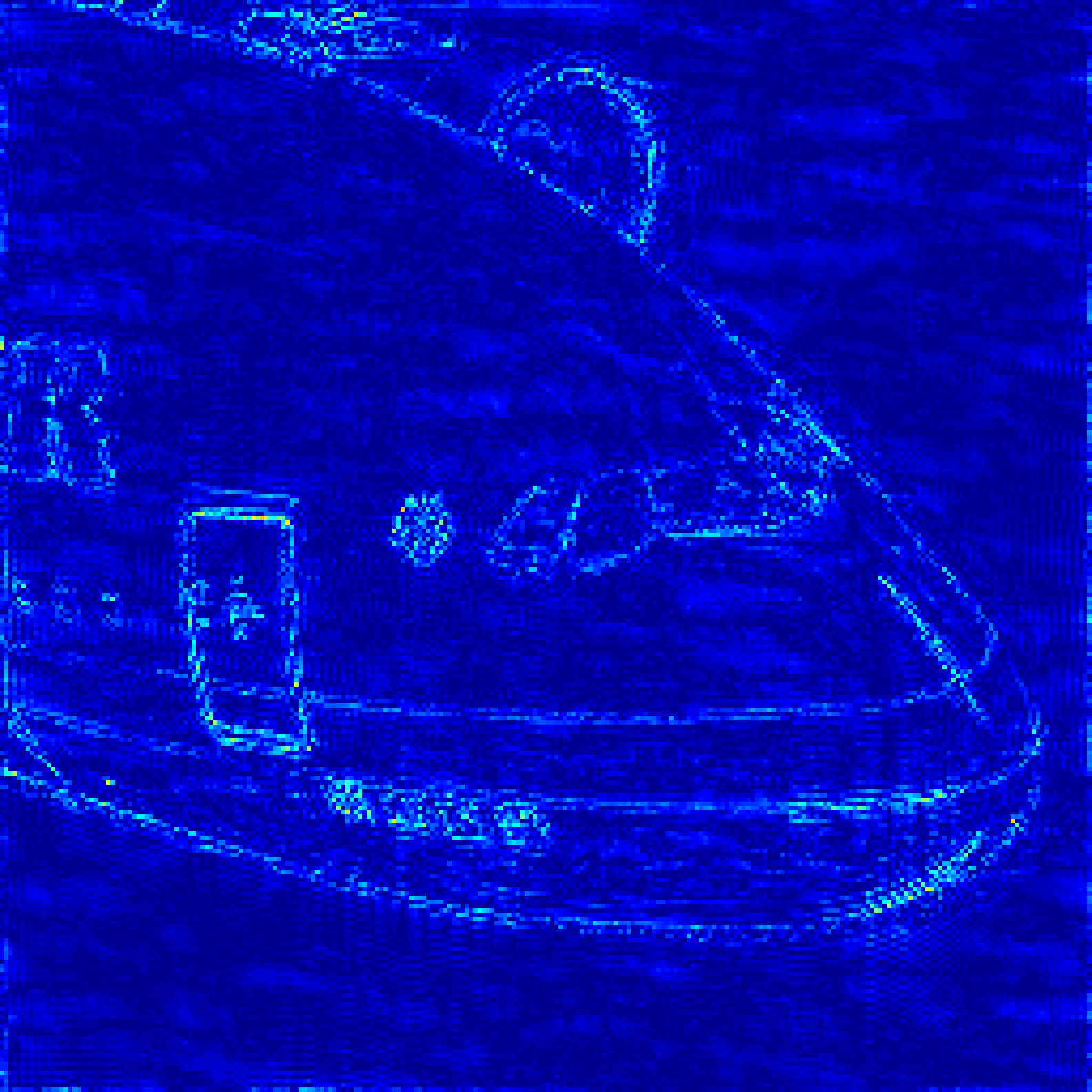}}\hspace{0.001cm}
\subfloat[]{\label{AirplaneTGVError}\includegraphics[width=2.25cm]{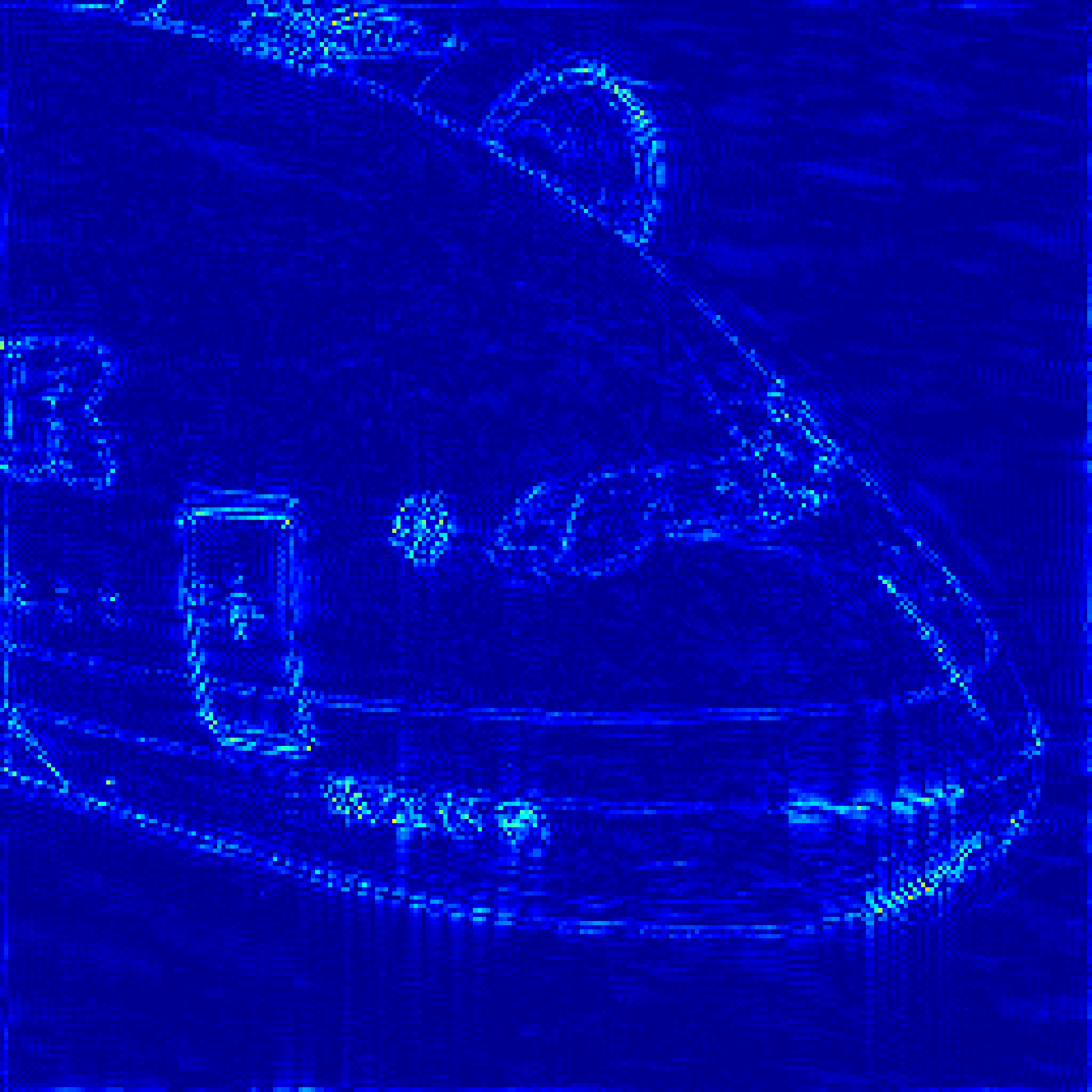}}\hspace{0.001cm}
\subfloat[]{\label{AirplaneInfConvError}\includegraphics[width=2.25cm]{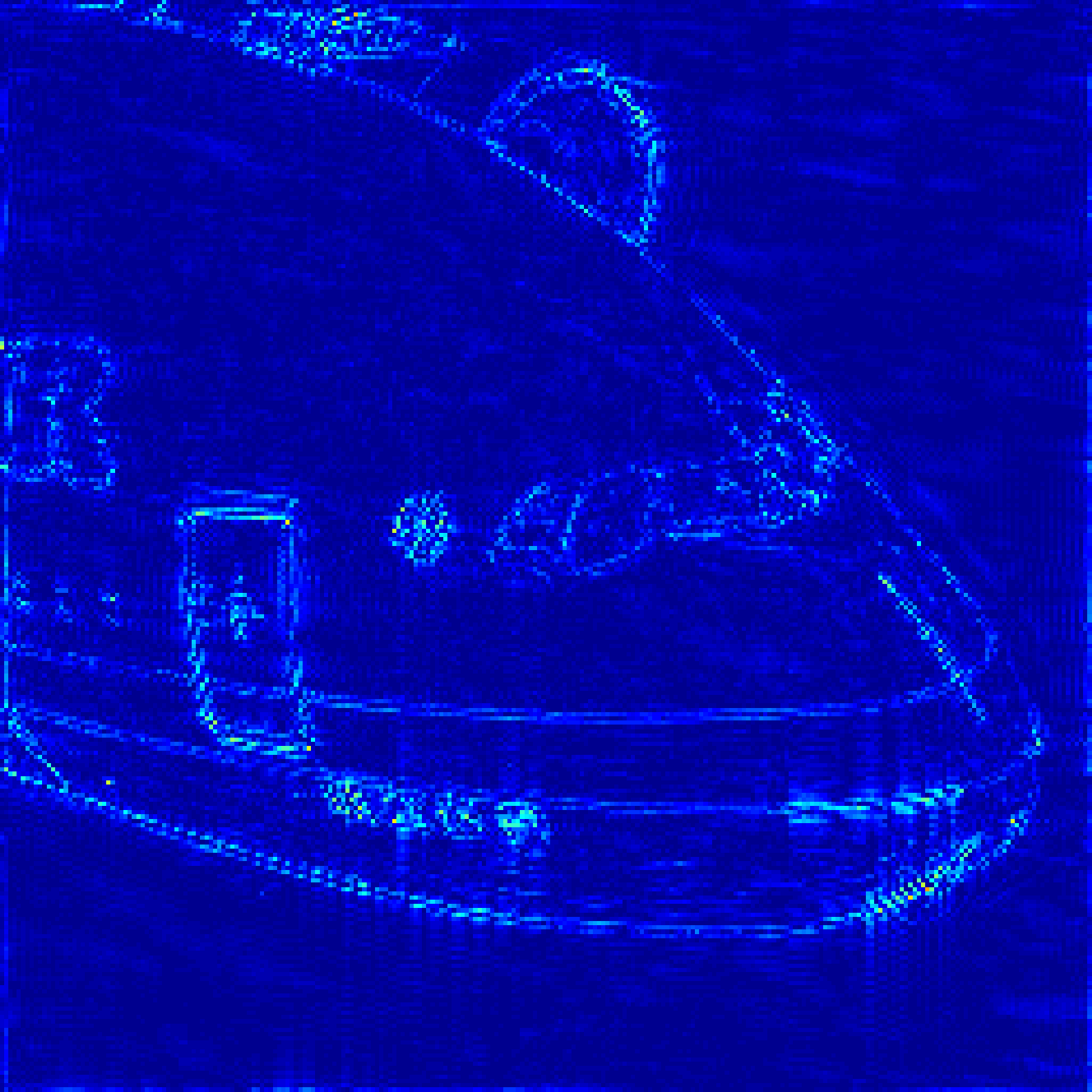}}
\caption{Visual comparisons for ``Airplane''. \cref{AirplaneOriginal}: original image, \cref{AirplaneLRHTGV}: model \cref{ProposedCSMRI}, \cref{AirplaneLRHTGVIRLS}: SLRM \cref{LRHTGVIRLSCSMRI}, \cref{AirplaneLRHInfConvIRLS}: GSLR \cref{LRHInfConvIRLSCSMRI}, \cref{AirplaneFra}: framelet \cref{FrameCSMRI}, \cref{AirplaneTGV}: TGV \cref{TGVCSMRI}, \cref{AirplaneInfConv}: infimal convolution \cref{InfConvCSMRI}. \cref{AirplaneOriginalZoom,AirplaneLRHTGVZoom,AirplaneLRHTGVIRLSZoom,AirplaneLRHInfConvIRLSZoom,AirplaneFraZoom,AirplaneTGVZoom,AirplaneInfConvZoom}: zoom-in views of \cref{AirplaneOriginal,AirplaneLRHTGV,AirplaneLRHTGVIRLS,AirplaneLRHInfConvIRLS,AirplaneFra,AirplaneTGV,AirplaneInfConv}. Red arrows indicate the region worth noting. \cref{AirplaneMask}: sample region, \cref{AirplaneLRHTGVError,AirplaneLRHTGVIRLS,AirplaneLRHInfConvIRLS,AirplaneFraError,AirplaneTGVError,AirplaneInfConvError,AirplaneLRHTGVIRLSError,AirplaneLRHInfConvIRLSError}: error maps of \cref{AirplaneLRHTGV,AirplaneFra,AirplaneTGV,AirplaneInfConv}.}\label{AirplaneResults}
\end{figure}

\begin{figure}[t]
\centering
\subfloat[]{\label{CarOriginal}\includegraphics[width=2.25cm]{CarOriginal.pdf}}\hspace{0.001cm}
\subfloat[]{\label{CarLRHTGV}\includegraphics[width=2.25cm]{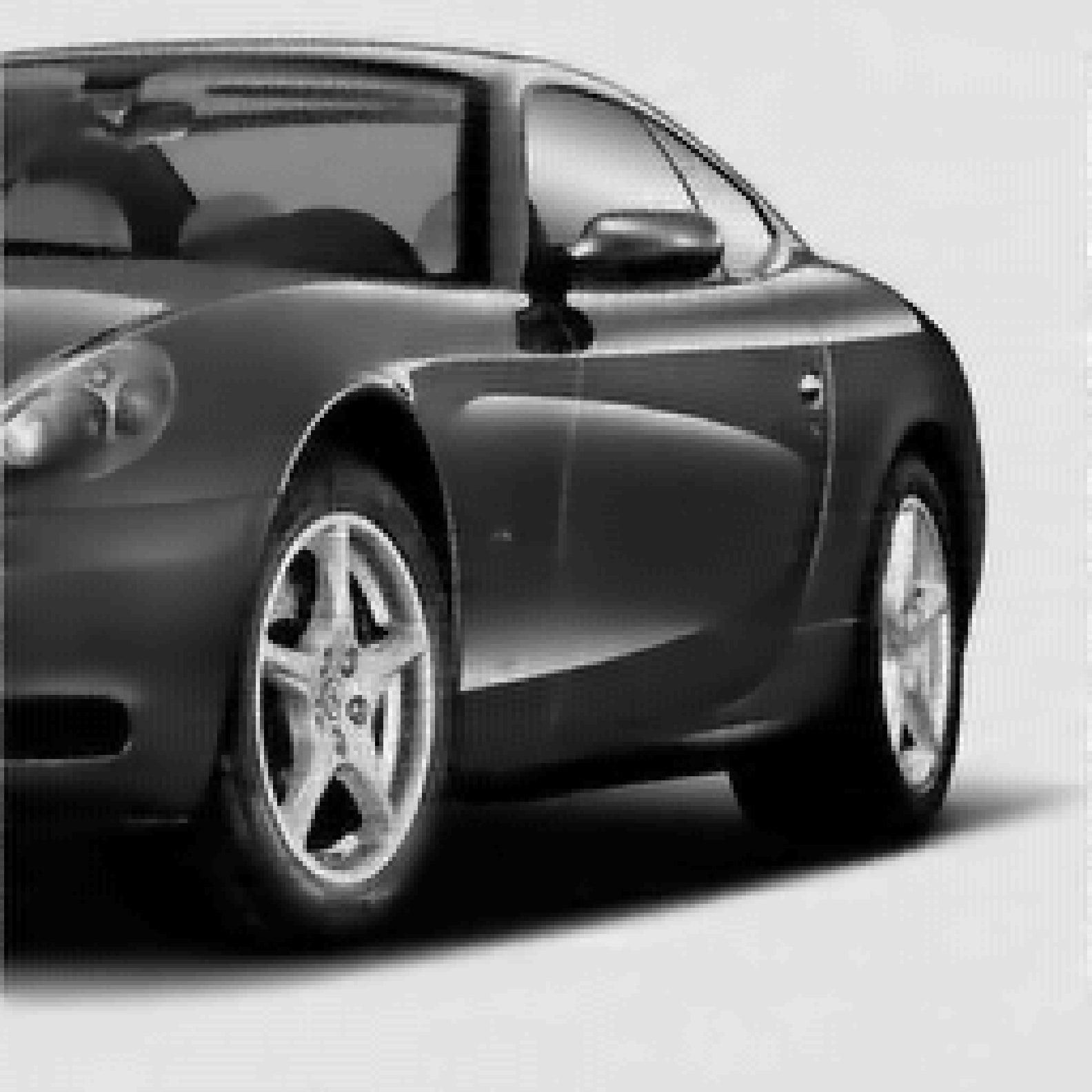}}\hspace{0.001cm}
\subfloat[]{\label{CarLRHTGVIRLS}\includegraphics[width=2.25cm]{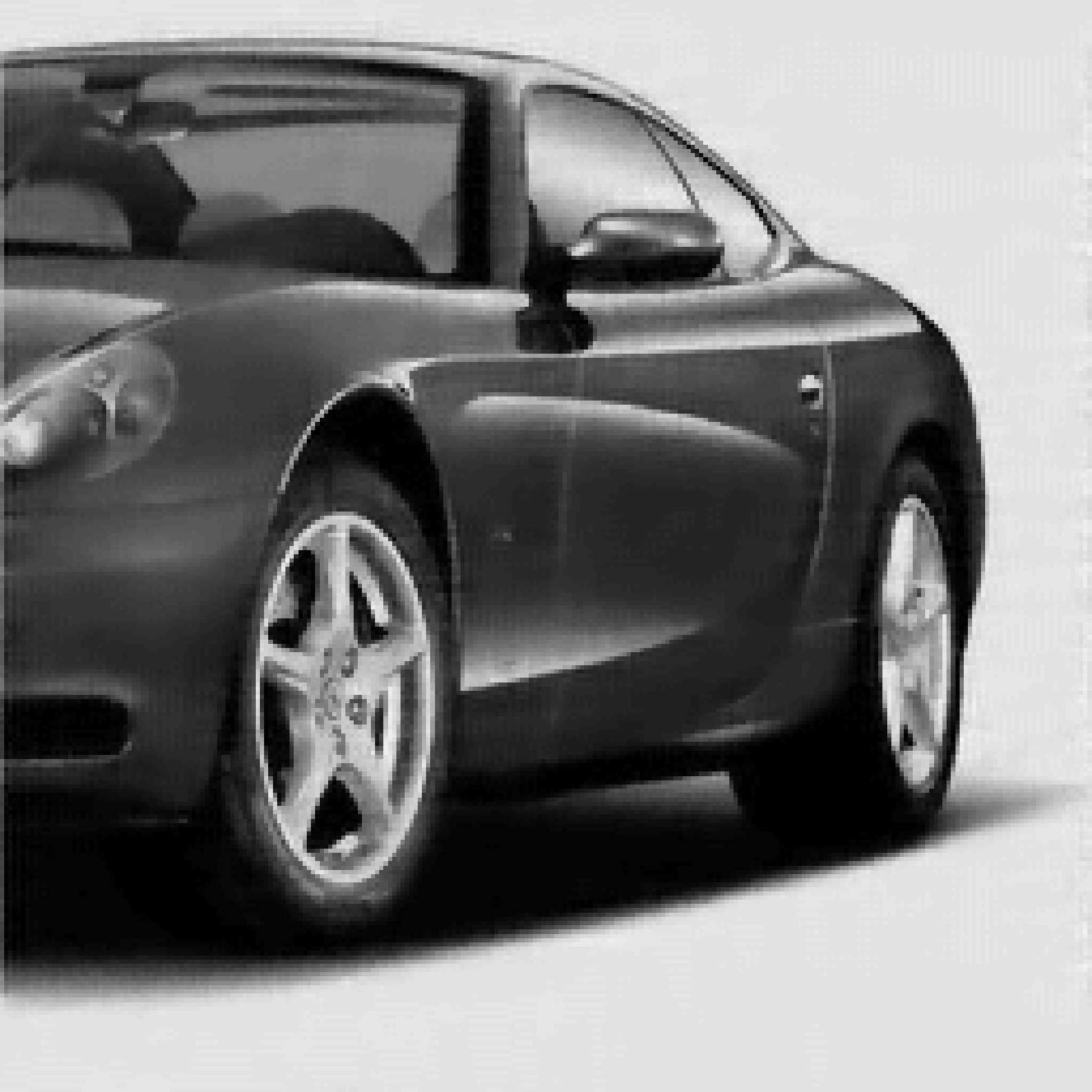}}\hspace{0.001cm}
\subfloat[]{\label{CarLRHInfConvIRLS}\includegraphics[width=2.25cm]{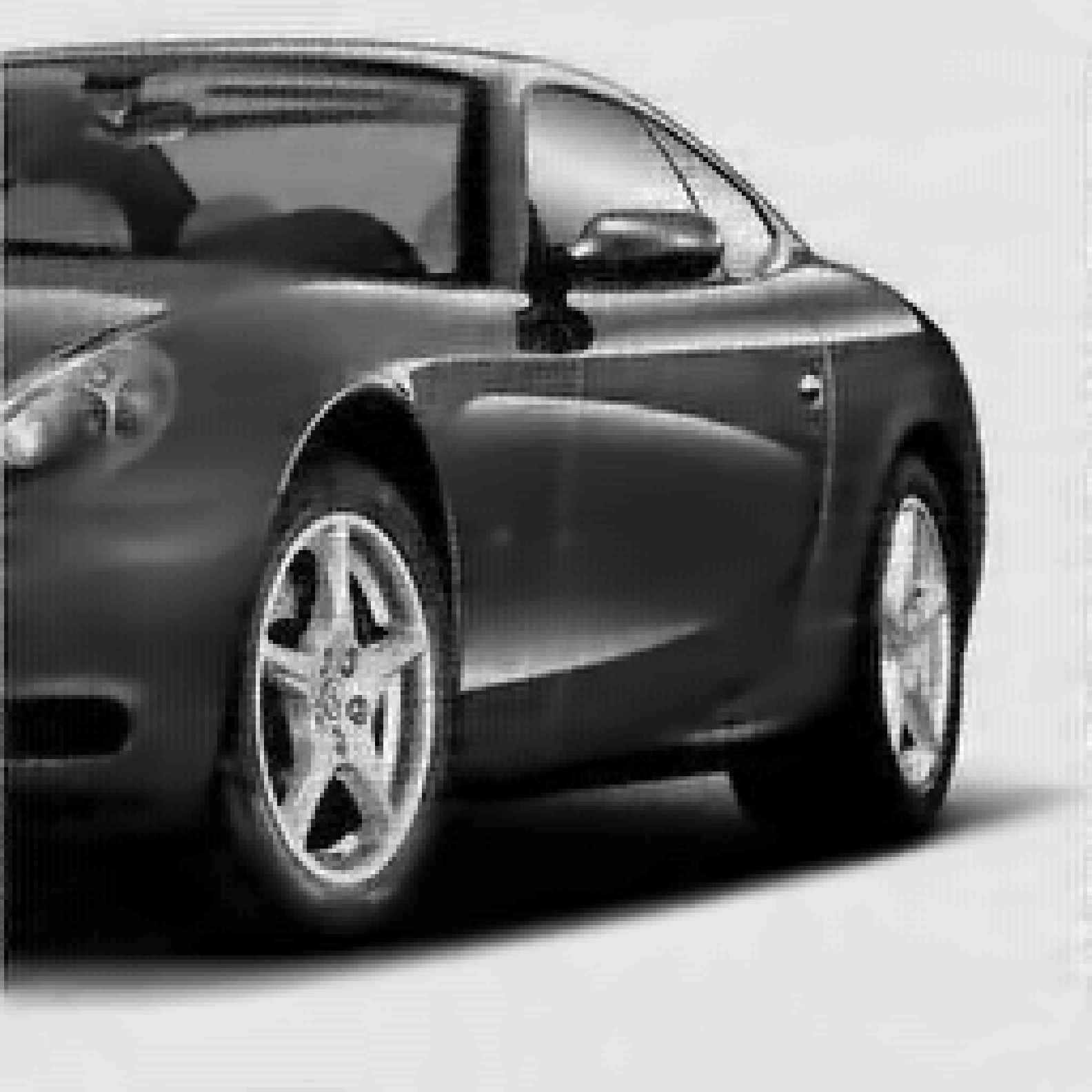}}\hspace{0.001cm}
\subfloat[]{\label{CarFra}\includegraphics[width=2.25cm]{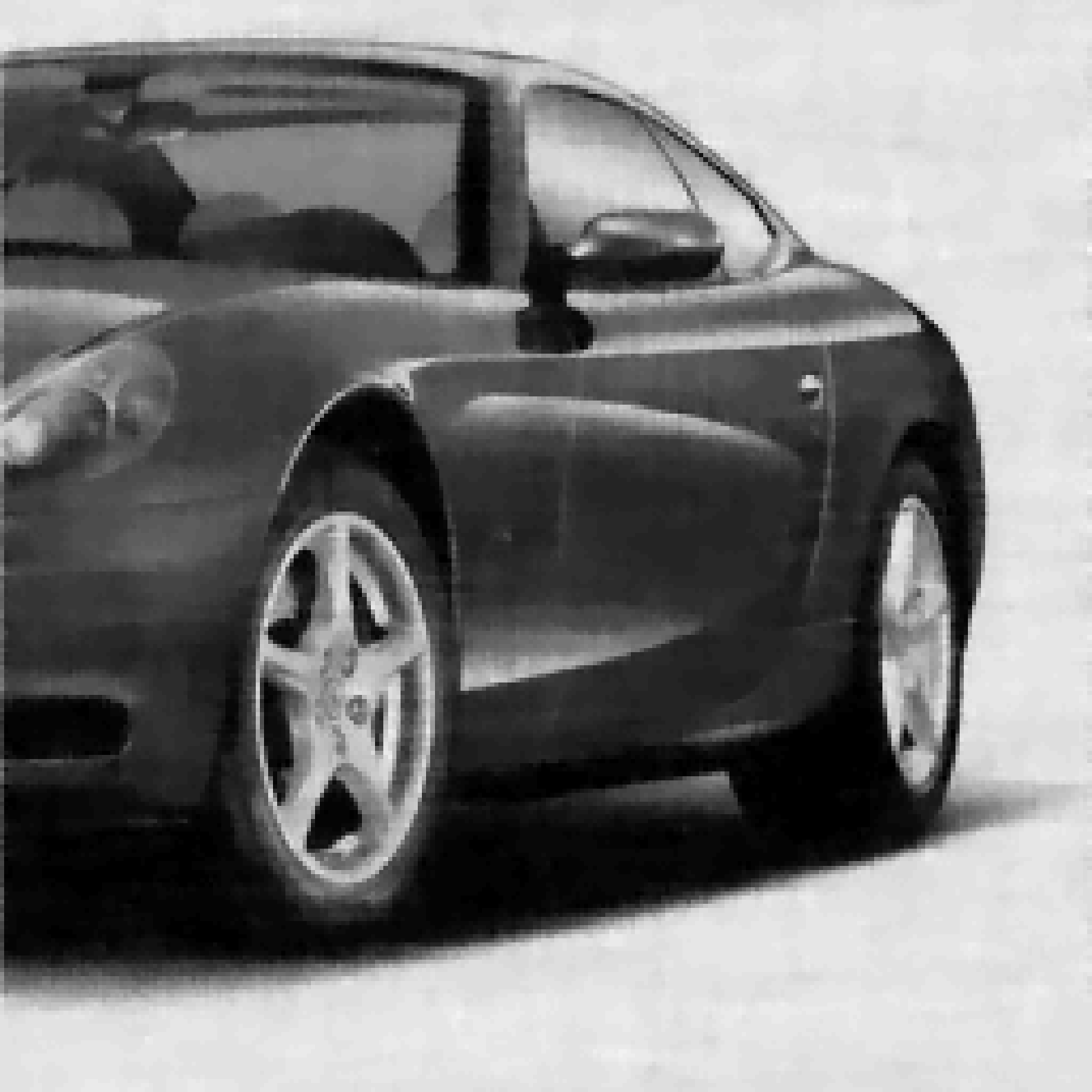}}\hspace{0.001cm}
\subfloat[]{\label{CarTGV}\includegraphics[width=2.25cm]{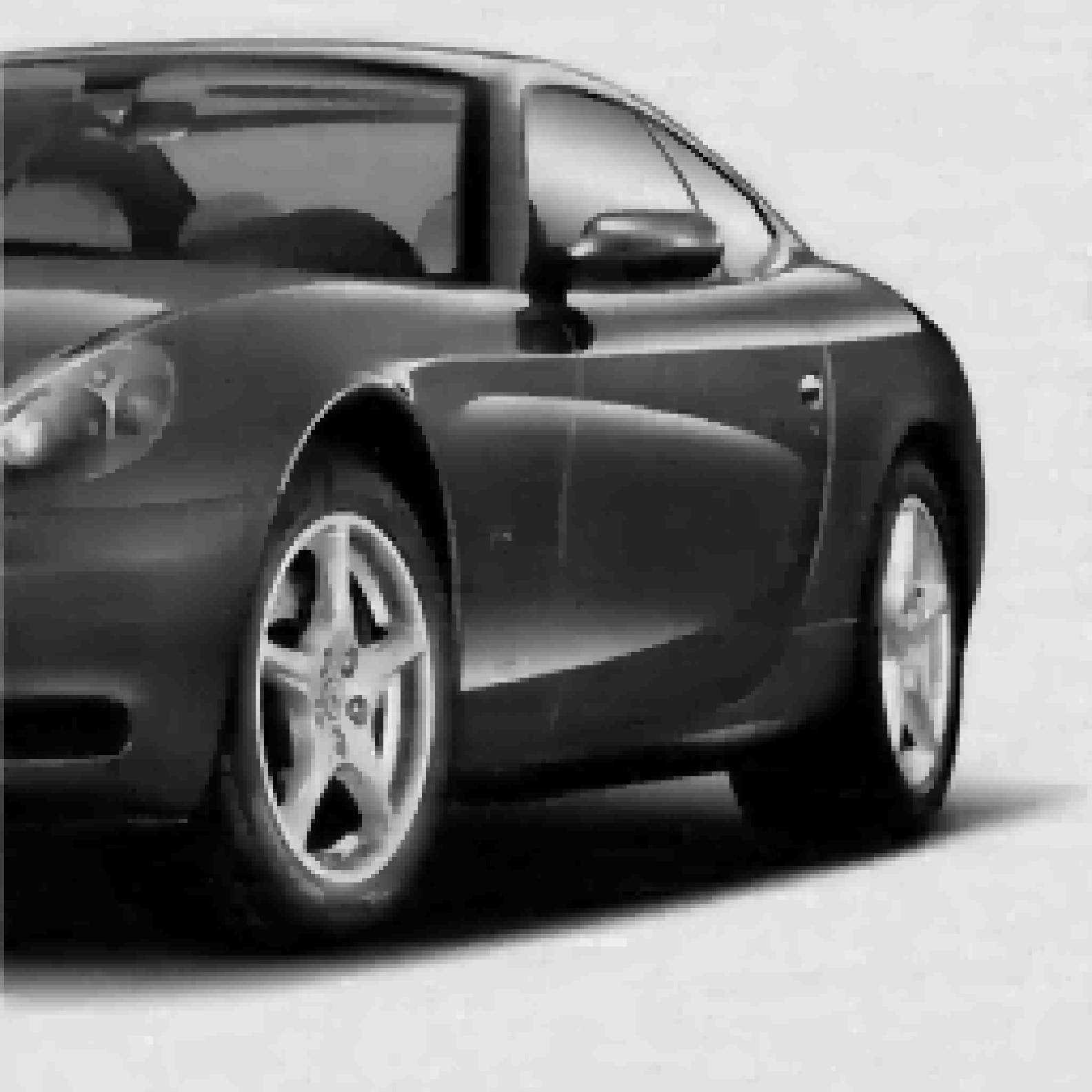}}\hspace{0.001cm}
\subfloat[]{\label{CarInfConv}\includegraphics[width=2.25cm]{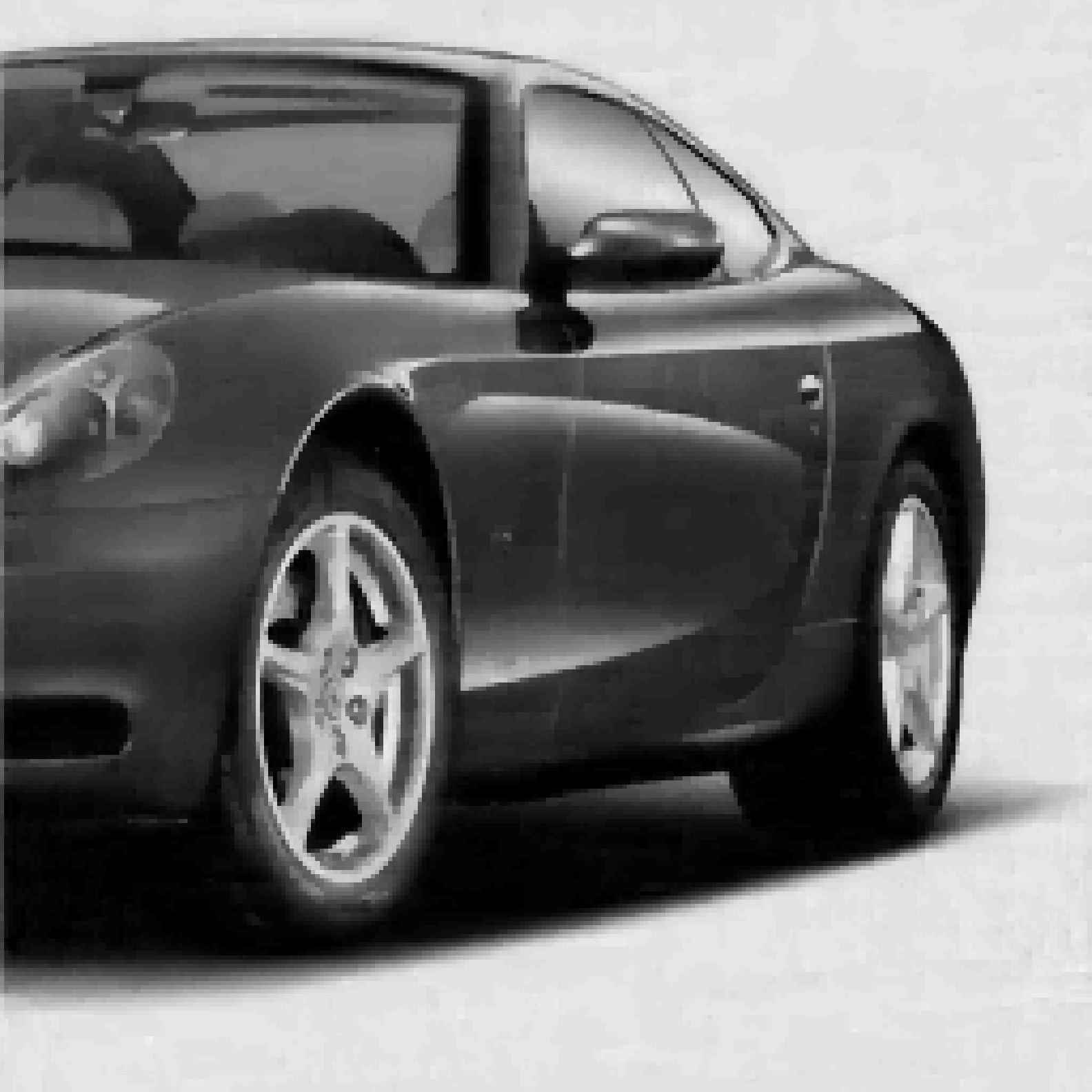}}\vspace{-0.75em}\\
\subfloat[]{\label{CarOriginalZoom}\includegraphics[width=2.25cm]{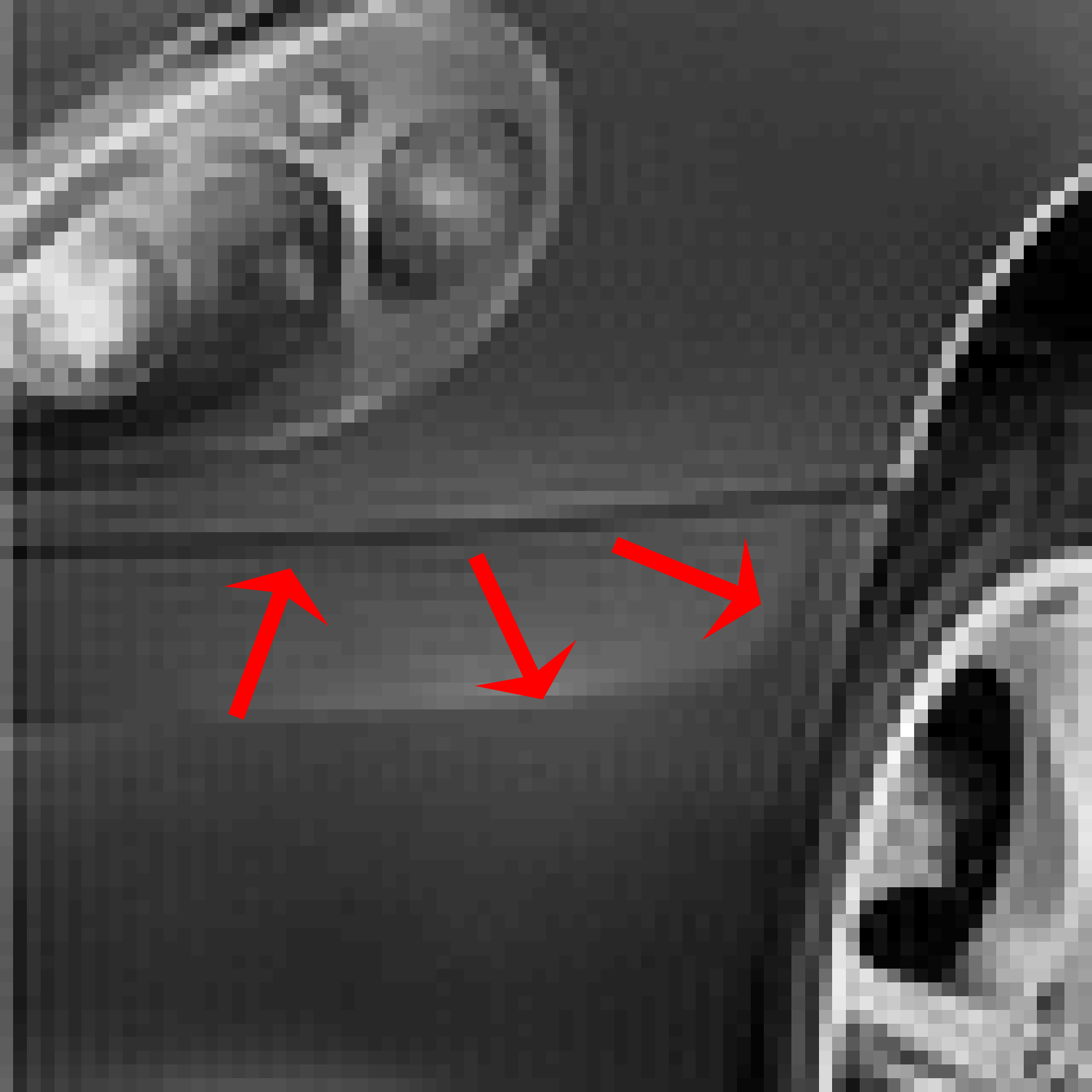}}\hspace{0.001cm}
\subfloat[]{\label{CarLRHTGVZoom}\includegraphics[width=2.25cm]{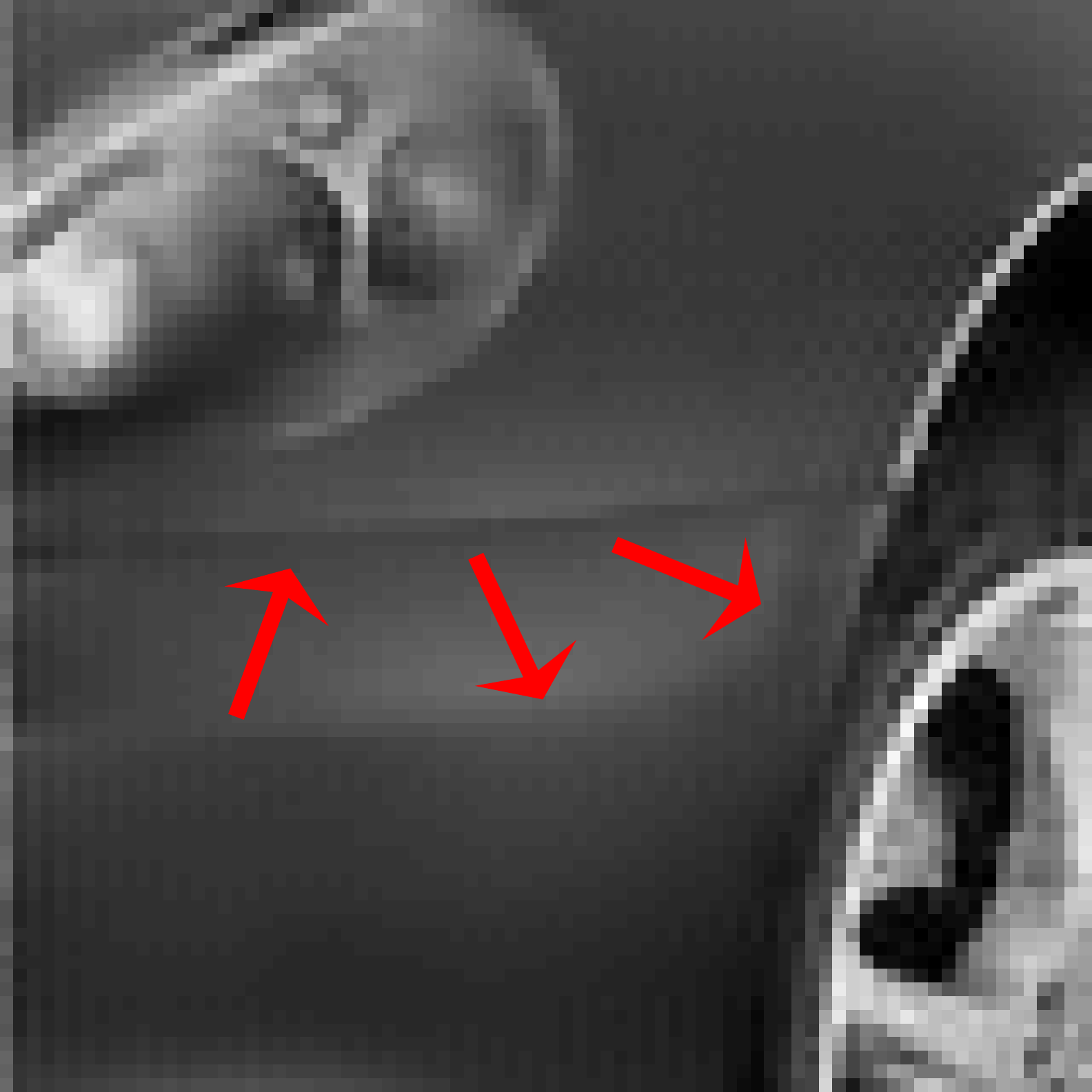}}\hspace{0.001cm}
\subfloat[]{\label{CarLRHTGVIRLSZoom}\includegraphics[width=2.25cm]{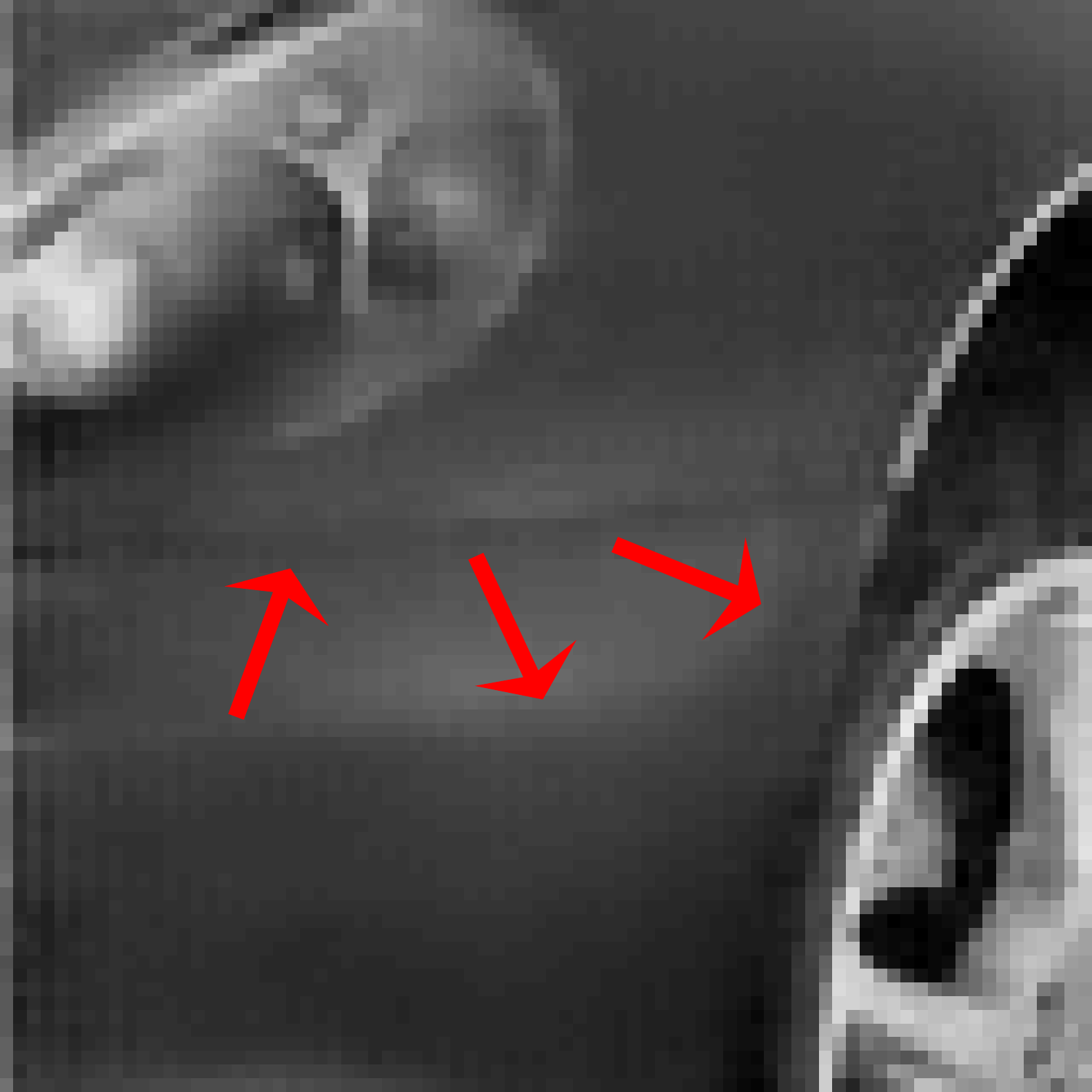}}\hspace{0.001cm}
\subfloat[]{\label{CarLRHInfConvIRLSZoom}\includegraphics[width=2.25cm]{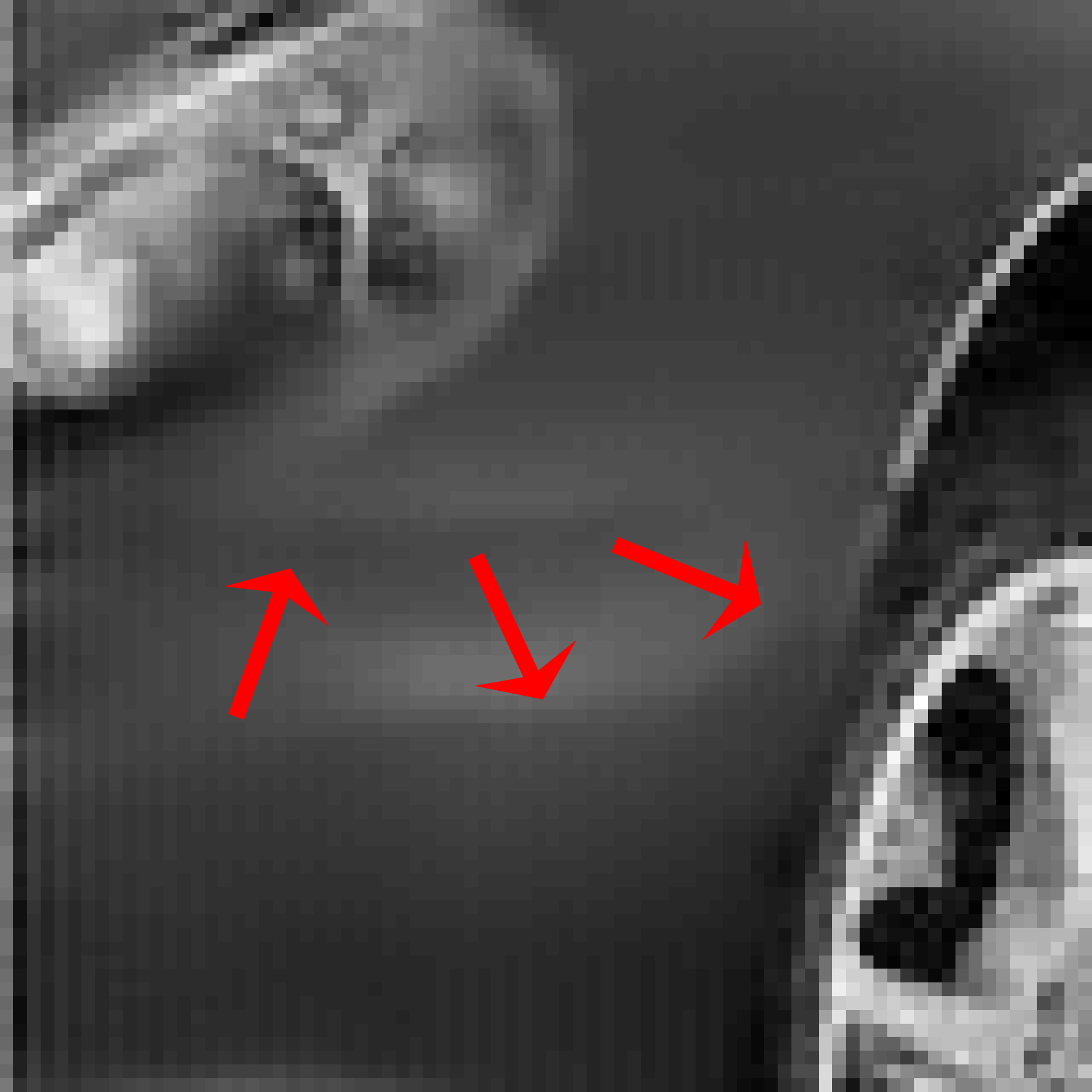}}\hspace{0.001cm}
\subfloat[]{\label{CarFraZoom}\includegraphics[width=2.25cm]{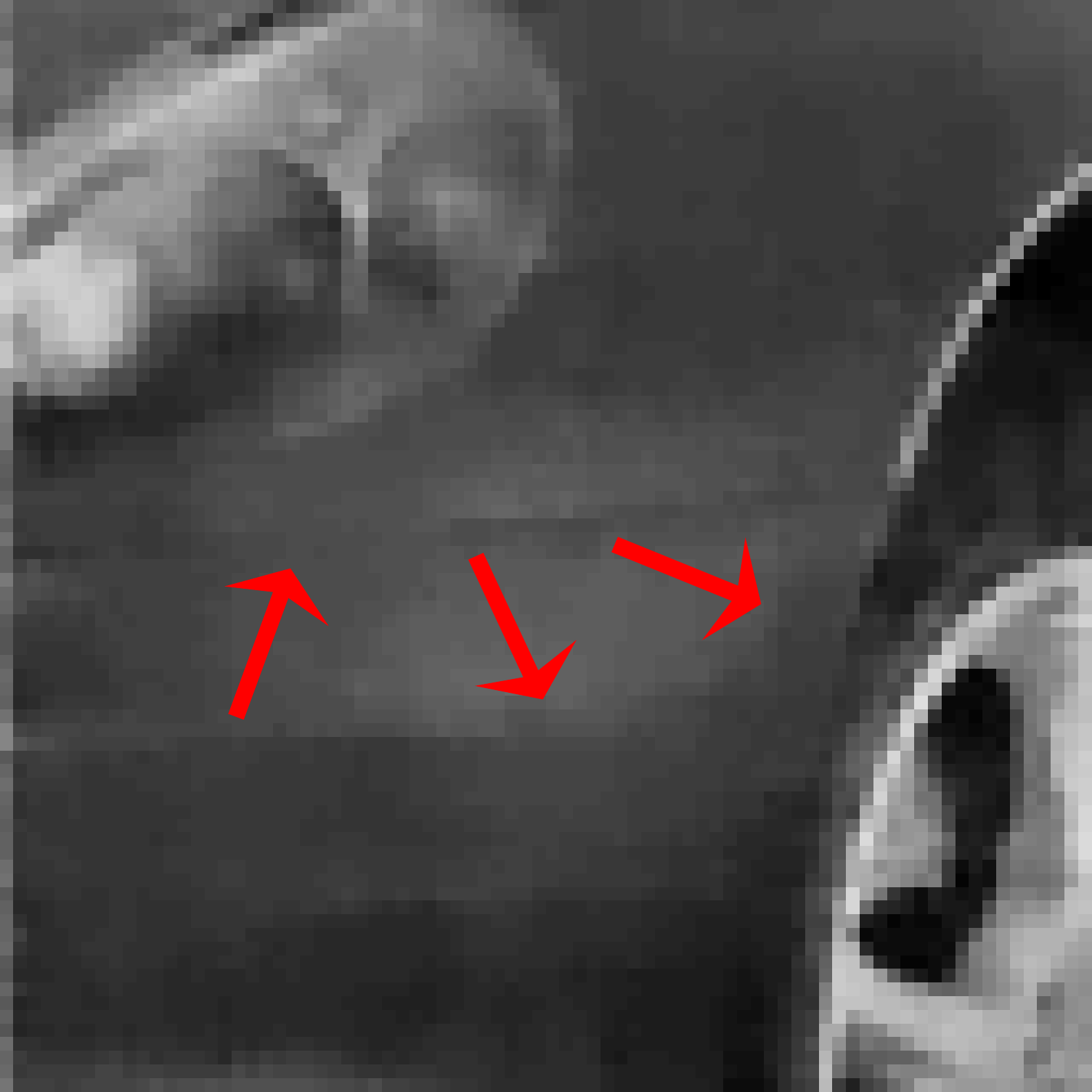}}\hspace{0.001cm}
\subfloat[]{\label{CarTGVZoom}\includegraphics[width=2.25cm]{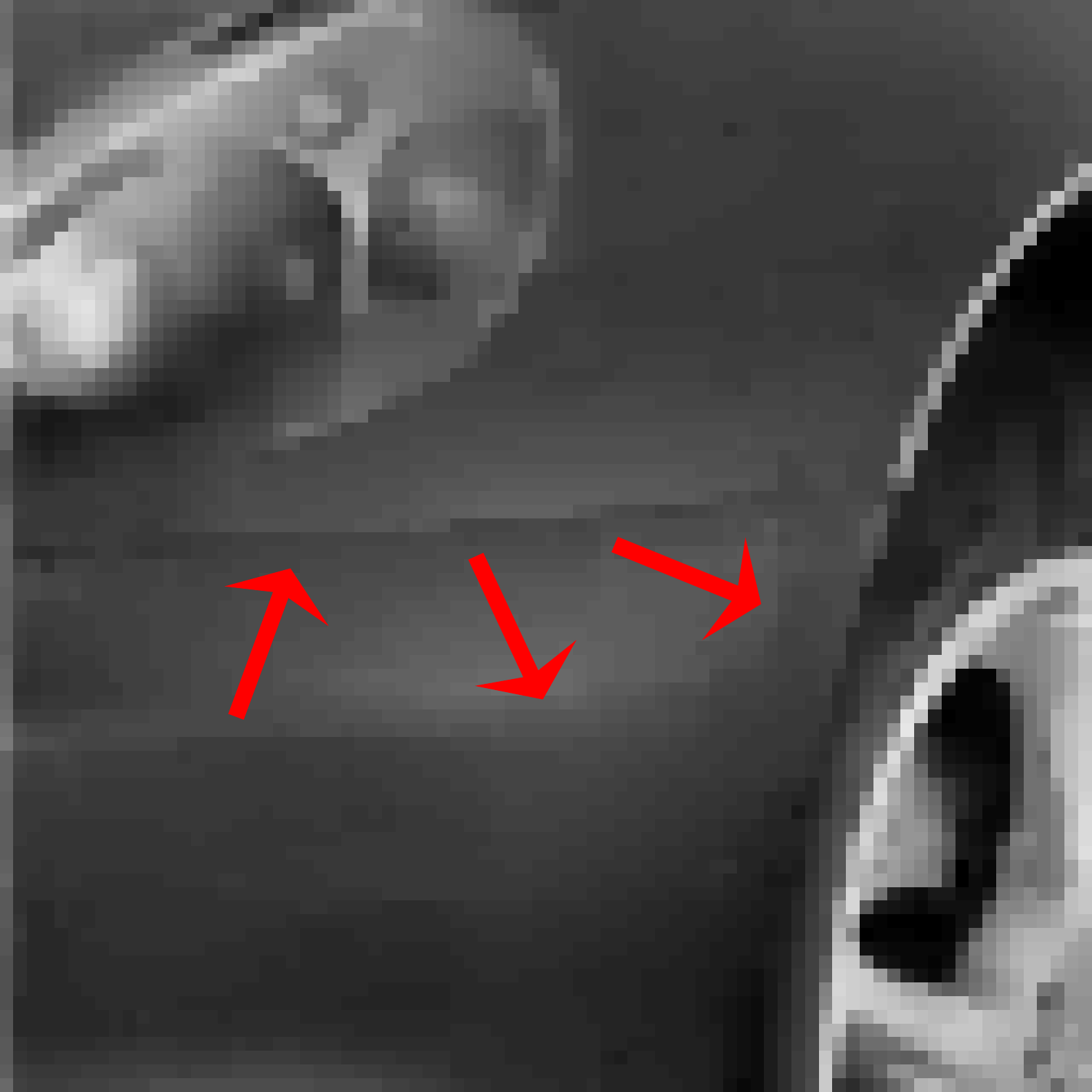}}\hspace{0.001cm}
\subfloat[]{\label{CarInfConvZoom}\includegraphics[width=2.25cm]{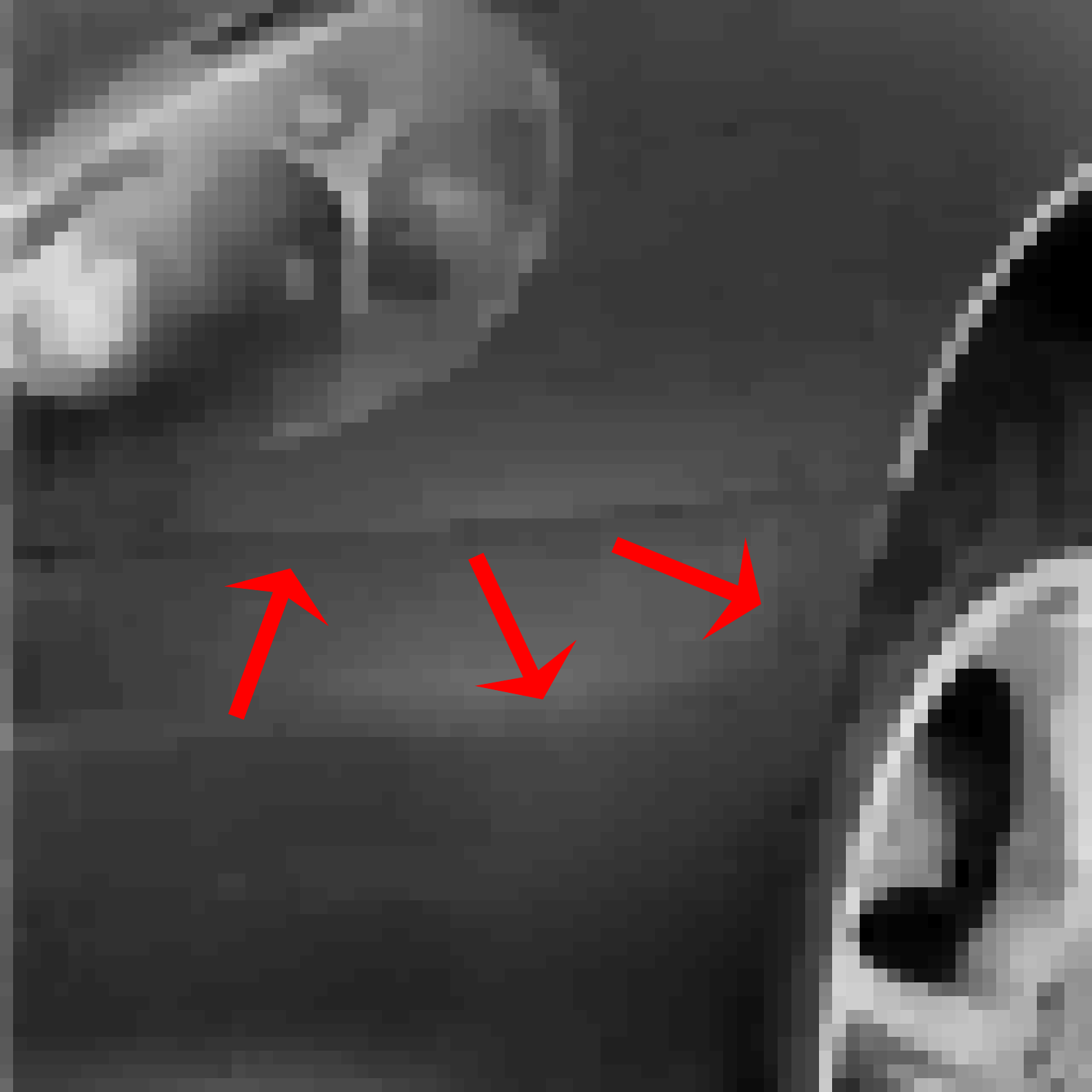}}\vspace{-0.75em}\\
\subfloat[]{\label{CarMask}\includegraphics[width=2.25cm]{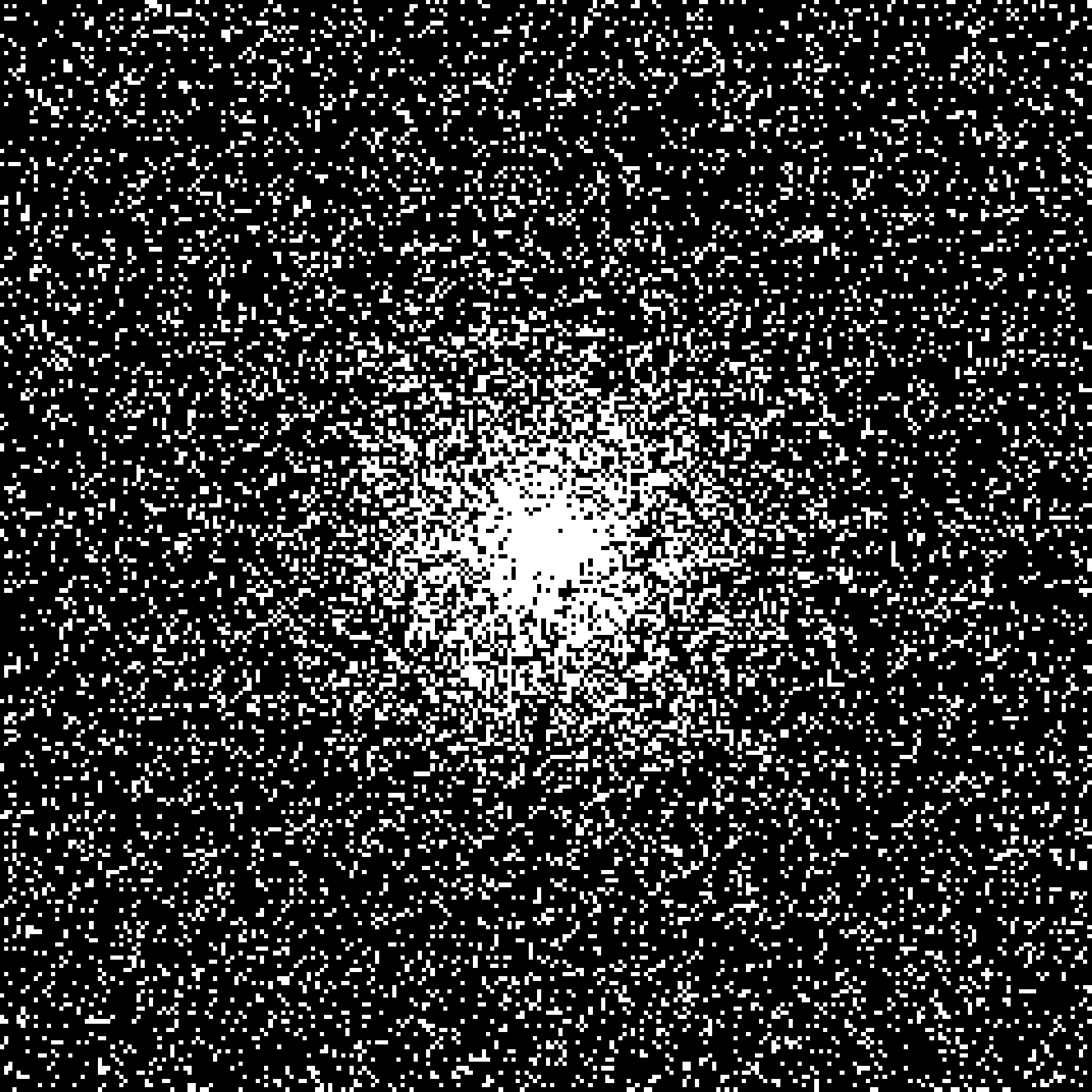}}\hspace{0.001cm}
\subfloat[]{\label{CarLRHTGVError}\includegraphics[width=2.25cm]{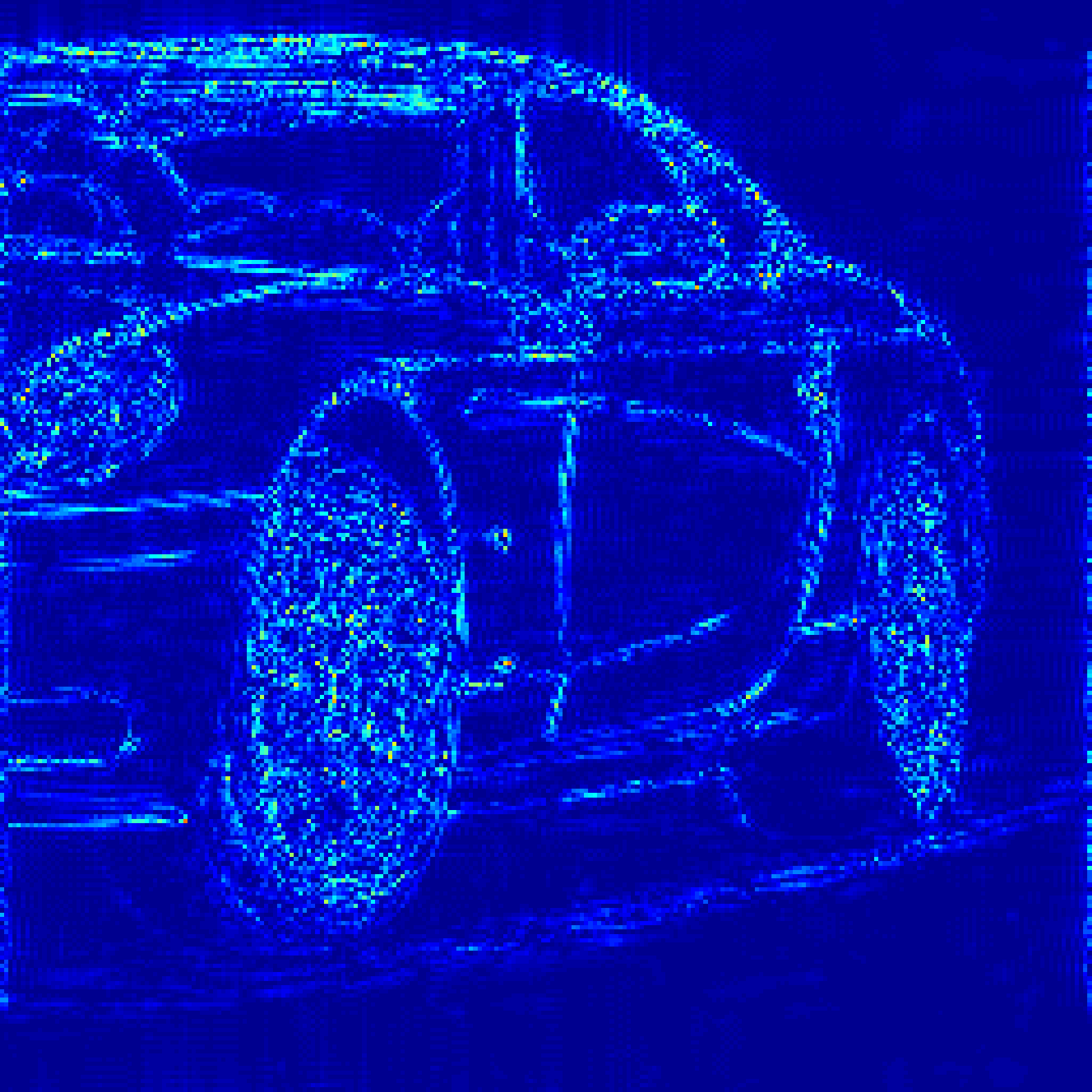}}\hspace{0.001cm}
\subfloat[]{\label{CarLRHTGVIRLSError}\includegraphics[width=2.25cm]{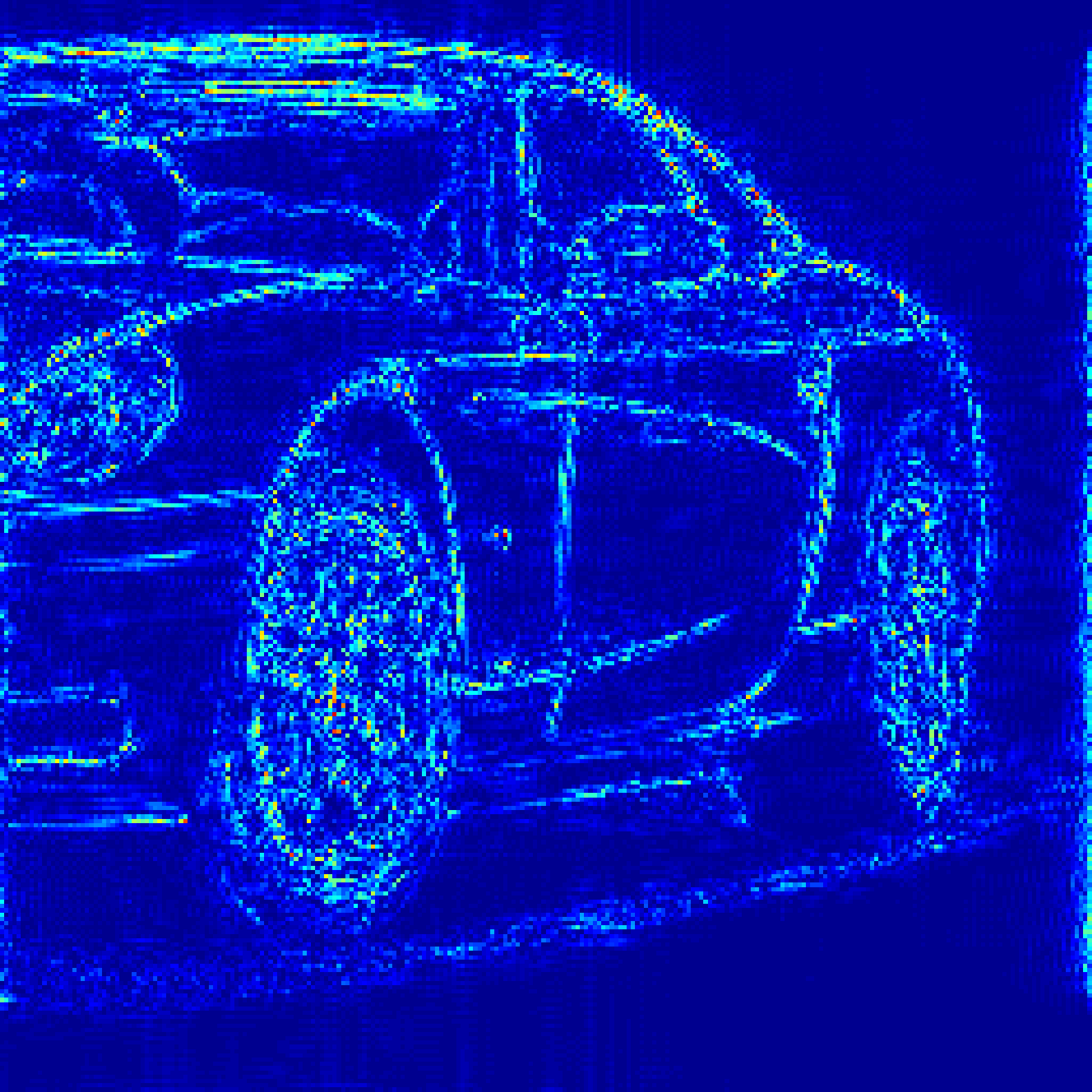}}\hspace{0.001cm}
\subfloat[]{\label{CarLRHInfConvIRLSError}\includegraphics[width=2.25cm]{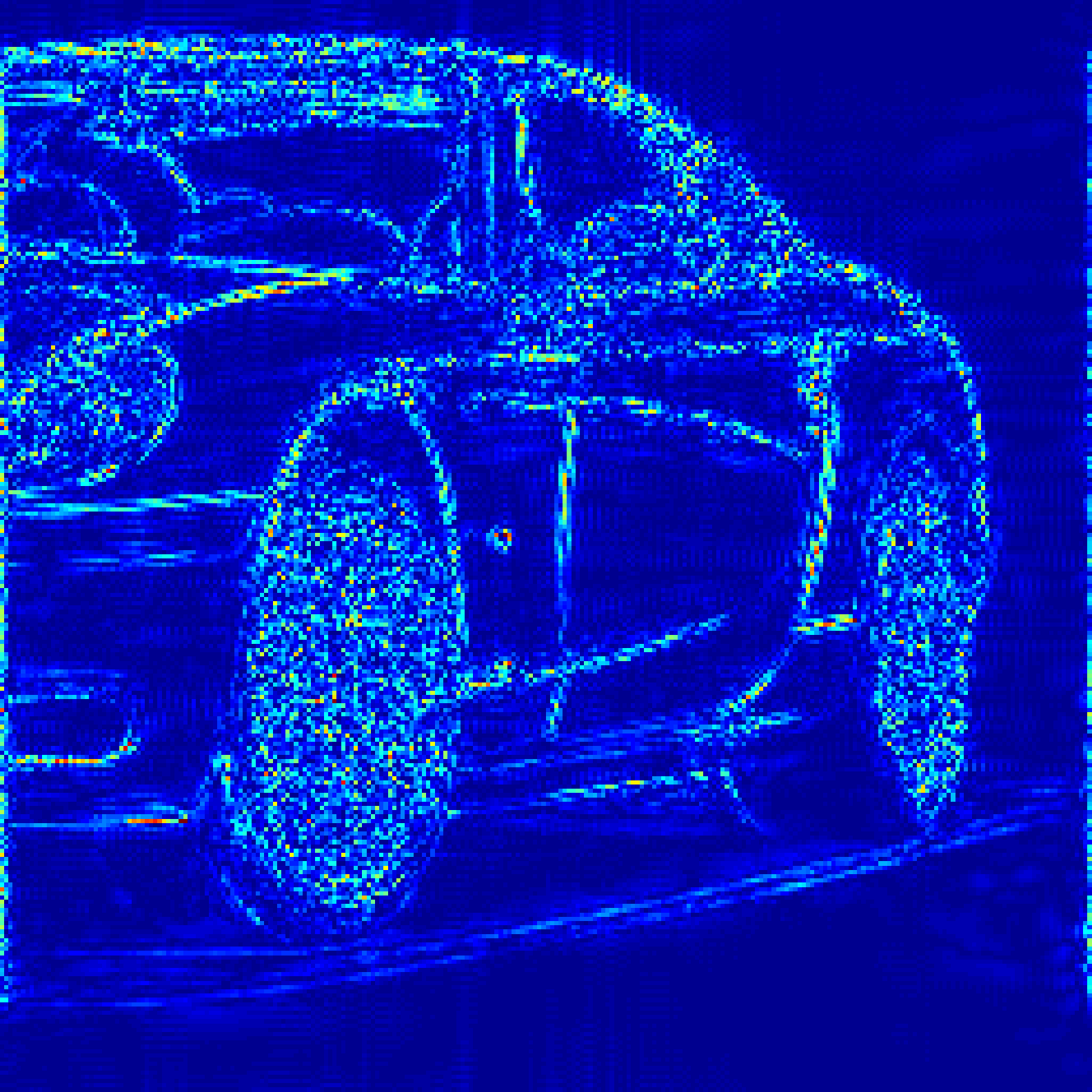}}\hspace{0.001cm}
\subfloat[]{\label{CarFraError}\includegraphics[width=2.25cm]{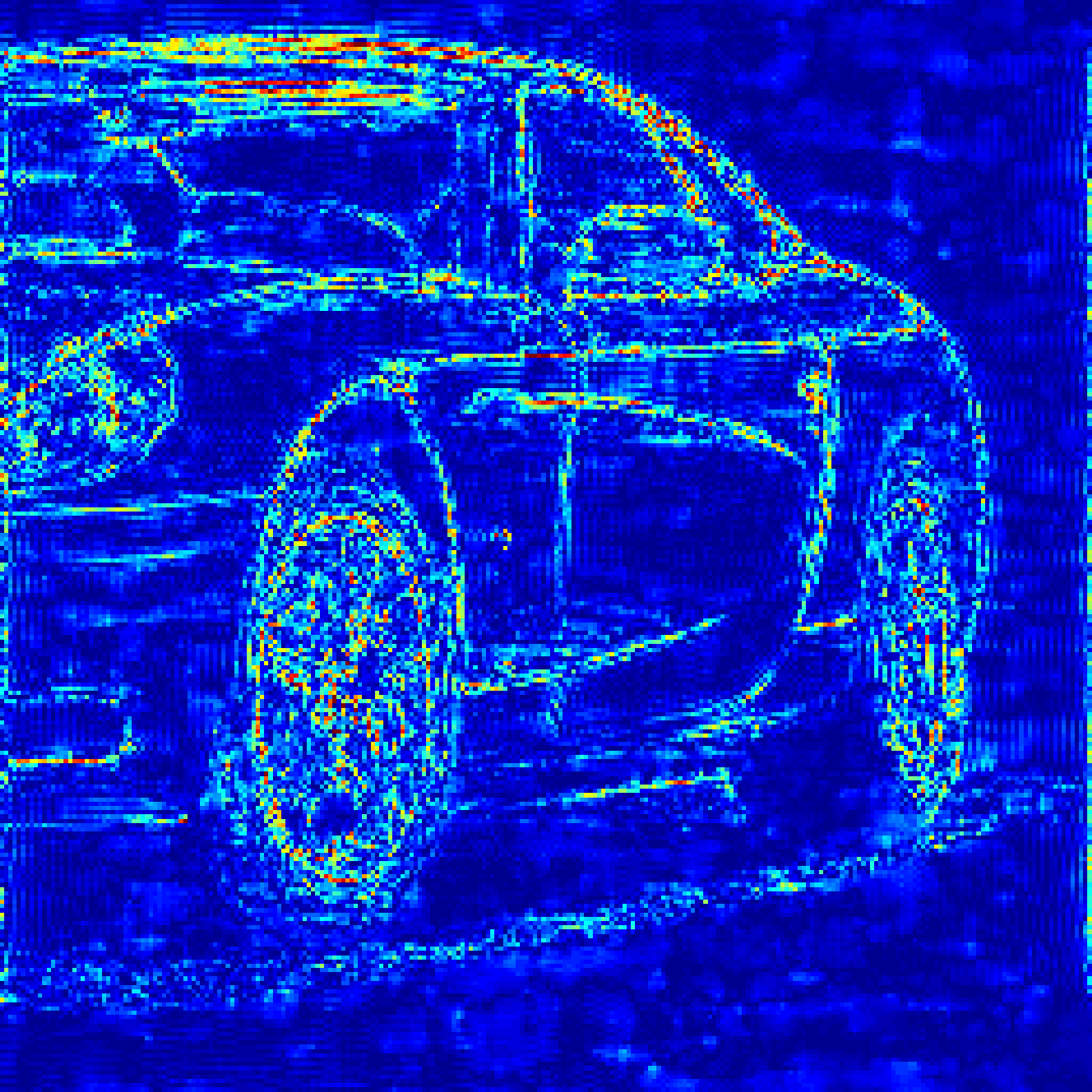}}\hspace{0.001cm}
\subfloat[]{\label{CarTGVError}\includegraphics[width=2.25cm]{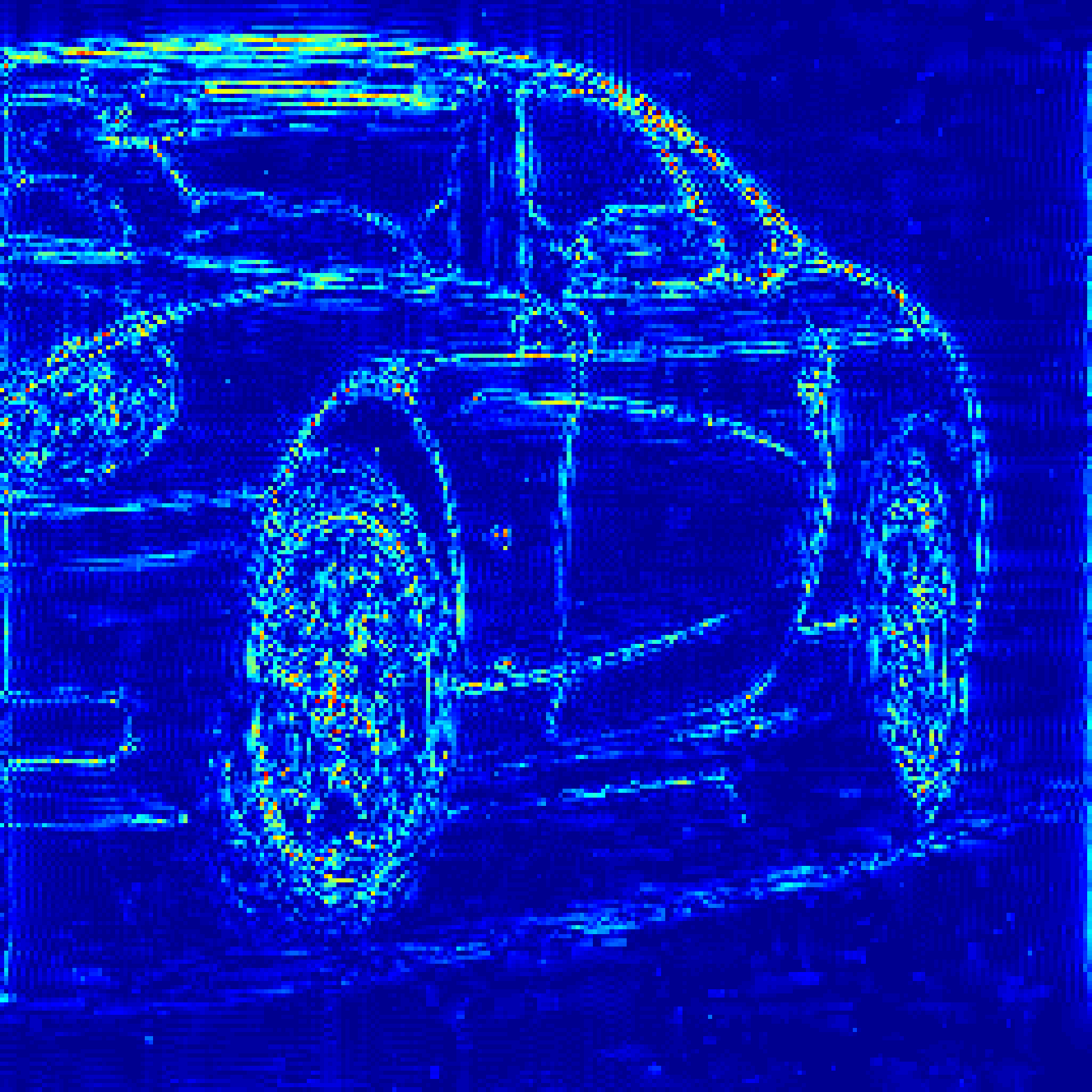}}\hspace{0.001cm}
\subfloat[]{\label{CarInfConvError}\includegraphics[width=2.25cm]{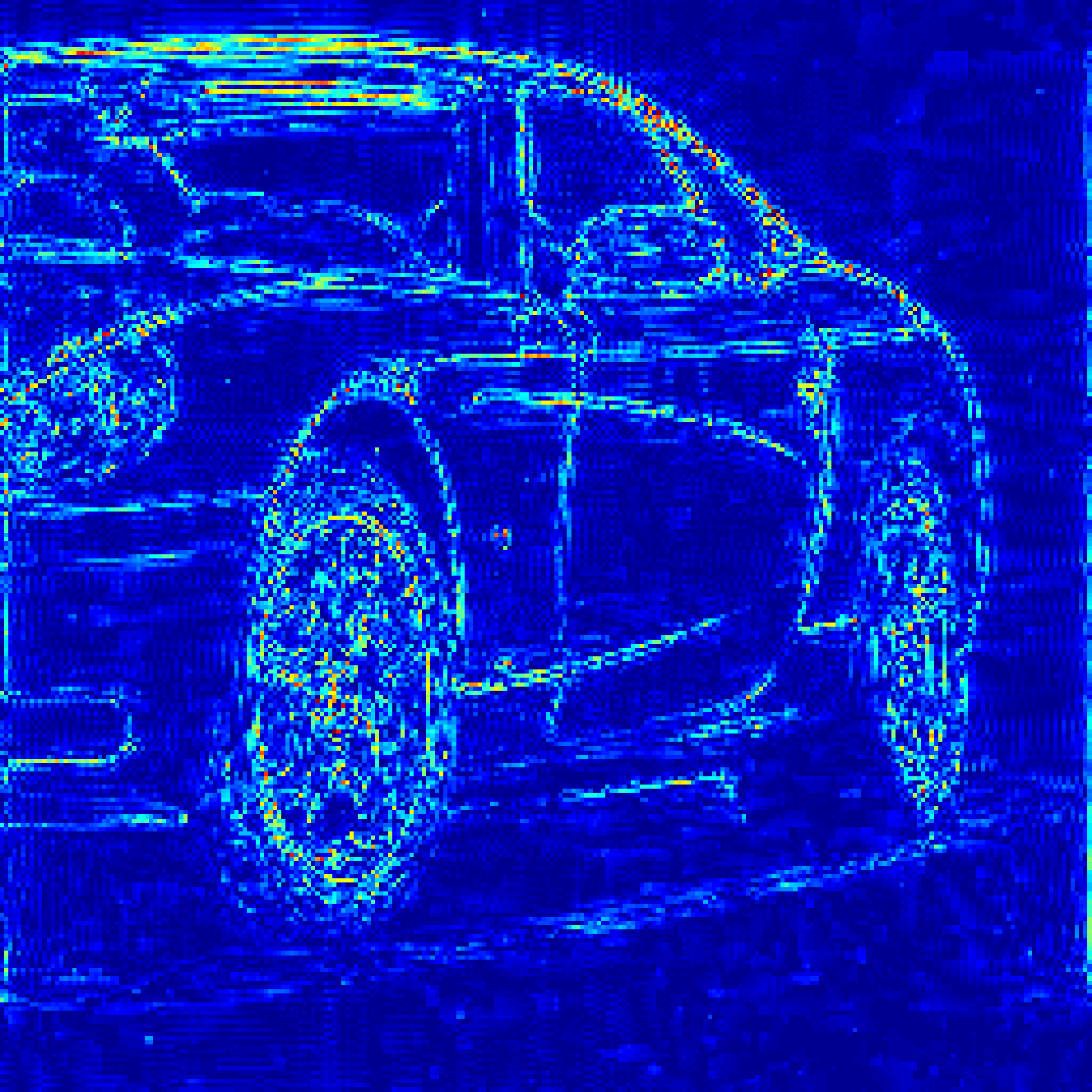}}
\caption{Visual comparisons for ``Car''. \cref{CarOriginal}: original image, \cref{CarLRHTGV}: model \cref{ProposedCSMRI}, \cref{CarLRHTGVIRLS}: SLRM \cref{LRHTGVIRLSCSMRI}, \cref{CarLRHInfConvIRLS}: GSLR \cref{LRHInfConvIRLSCSMRI}, \cref{CarFra}: framelet \cref{FrameCSMRI}, \cref{CarTGV}: TGV \cref{TGVCSMRI}, \cref{CarInfConv}: infimal convolution \cref{InfConvCSMRI}. \cref{CarOriginalZoom,CarLRHTGVZoom,CarLRHTGVIRLSZoom,CarLRHInfConvIRLSZoom,CarFraZoom,CarTGVZoom,CarInfConvZoom}: zoom-in views of \cref{CarOriginal,CarLRHTGV,CarLRHTGVIRLS,CarLRHInfConvIRLS,CarFra,CarTGV,CarInfConv}. Red arrows indicate the region worth noting. \cref{CarMask}: sample region, \cref{CarLRHTGVError,CarLRHTGVIRLSError,CarLRHInfConvIRLSError,CarFraError,CarTGVError,CarInfConvError}: error maps of \cref{CarLRHTGV,CarLRHTGVIRLS,CarLRHInfConvIRLS,CarFra,CarTGV,CarInfConv}.}\label{CarResults}
\end{figure}

\cref{TableResults} summarizes the SNR, the HFEN, and the SSIM of the aforementioned approaches. \cref{LoganResults,RectangleResults,AirplaneResults,CarResults} display the visual comparisons (the first row) together with the zoom-in views (the second row) and the error maps (the third row) of \cref{ProposedCSMRI} against \cref{FrameCSMRI,TGVCSMRI,InfConvCSMRI,LRHTGVIRLSCSMRI,LRHInfConvIRLSCSMRI}. Throughout this paper, all restored images are displayed in the window level $[0,1]$, and all error maps are displayed in the window level $[0,0.2]$ for fair comparisons. We can see that the relaxed model \cref{ProposedCSMRI} consistently performs best, and the SLRM model \cref{LRHTGVIRLSCSMRI} performs the second in almost every index in all scenarios, and the improvements are visually observable as well, both of which demonstrate that our proposed SLRM framework performs well in the piecewise smooth image restoration.

At first glance, since the models \cref{ProposedCSMRI,LRHTGVIRLSCSMRI} can be regarded as different relaxations of the structured low rank matrix framework for the continuous domain regularization, the results also suggest that the off-the-grid regularization performs better than the on-the-grid regularization as it can reduce the basis mismatch between the true support (or the true singularity) in continuum and the discrete grid, leading to the improvements in the indices. In fact, due to such a basis mismatch, the conventional on-the-grid approaches (\cref{FrameCSMRI,TGVCSMRI,InfConvCSMRI}) suffers from the errors concentrated near the image singularities as well as the distorted shapes compared to \cref{ProposedCSMRI}, which can be seen in the error maps and the zoom-in views of \cref{LoganResults,RectangleResults,AirplaneResults,CarResults}.

Most importantly, the numerical results illustrate the followings. First of all, from the comparisons between \cref{LRHTGVIRLSCSMRI,LRHInfConvIRLSCSMRI}, we can see that the proposed SLRM framework performs better than the GSLR framework in \cite{Y.Hu2019} for the piecewise smooth image restoration. Second, from the comparisons between \cref{ProposedCSMRI,LRHTGVIRLSCSMRI}, we can further see that the relaxation into the wavelet frame analysis approach can achieve further improvements over the direct rank minimization. As previously mentioned, since the proposed SLRM framework takes the GSLR framework in \cite{Y.Hu2019} as a special case, we can restore broader range of piecewise smooth images. In addition, it should be noted that the annihilation of the Fourier samples of the second order derivatives in general require larger annihilating filter than the first order derivatives \cite{G.Ongie2015a}, which may degrade the low rank structures of the corresponding multi-fold Hankel matrix. As a consequence, the GSLR framework can introduce staircase artifacts in the smooth region compared to the proposed SLRM framework, as can be seen in the zoom-in views of \cref{LoganResults,RectangleResults,AirplaneResults,CarResults}. Finally, since the weights corresponding to the derivatives in the image domain amplifies the noise in the high frequencies, it is likely that the direct rank minimization approaches (\cref{LRHTGVIRLSCSMRI,LRHInfConvIRLSCSMRI}) lead to the artifacts near the image singularities corresponding to the high frequency components, whereas the relaxation into the wavelet frame analysis approach (\cref{ProposedCSMRI}) can achieve the denoising effect in spite of the amplified noise, leading to the better restoration results with less artifacts near the edges, as can be seen in the error maps in \cref{LoganResults,RectangleResults,AirplaneResults,CarResults}.

\section{Conclusion and future directions}\label{Conclusion}

In this paper, we have introduced a new structured low rank matrix framework for the piecewise smooth image restoration. Following the previous structured low rank matrix framework for the piecewise constant image \cite{G.Ongie2016}, we assume that the image singularities lie in the zero level set of a band-limited periodic function. Inspired by the total generalized variation \cite{K.Bredies2010}, we derive the annihilation relation of the Fourier samples of the gradient and the (symmetric) gradient tensor, which in turn leads to the low rank multi-fold Hankel matrices to balance between the first order derivatives and the second order derivatives. In addition, as a by-product of the proposed structured low rank matrix framework, we further introduce a wavelet frame based sparse regularization model for the piecewise image restoration via the continuous domain regularization, from the SVDs of low rank multi-fold Hankel matrices. Finally, the numerical experiments show that the proposed wavelet frame based model based on the proposed SLRM framework outperforms both the conventional on-the-grid approaches and the existing the SLRM framework as well as the existing rank minimization approaches for the piecewise smooth image restoration.

For the future work, we plan to generalize the frameworks in \cite{G.Ongie2016} to the piecewise smooth image restoration. More precisely, we need to develop sampling guarantees for the unique restoration of a continuous domain piecewise smooth image from its uniform low pass Fourier samples, based on the proposed annihilation relations. Unfortunately, it is not clear for us at this moment under what conditions we can restore the singularity set from uniform Fourier samples as the annihilation relations involve a blind deconvolution problem due to the additional unknown vector field $p$ in \cref{Th1}. Nevertheless, this is definitely a future direction we would like to work on. We are also interested in proving the restoration guarantee of \cref{RankMinimization} (or its convex nuclear norm relaxation) given that $\bmA=\bmR_{\MM}$, i.e. the sample grid $\MM$ is drawn from the image grid $\OO$ uniformly at random, which generalizes the theoretical frameworks in \cite{G.Ongie2018} to the piecewise smooth image restoration.

Finally, to broaden the scope of applications, it is also likely to apply the idea in this paper to the various image restoration tasks such as the blind deconvolution (especially that in the frequency domain which is related to the retinex \cite{E.H.Land1971}) and the sparse angle CT. In fact, it is not clear at this moment how to extend this approach to the sparse angle CT restoration, as the measurement can be viewed as the Fourier samples on the discrete radial grid. Nevertheless, this is also an important future work for the broader scope of applications. Finally, observing the successive structure in our SLRM framework, we are also interested in developing a deep learning framework, motivated by the recent works on the deep learning techniques combined with the structured low rank matrix approaches \cite{M.Jacob2019,T.H.Kim2019,A.Pramanik2019a,A.Pramanik2019}.

{\appendix

\section{Preliminaries on tight wavelet frames}\label{PreliminariesTightFrame}

Provided here is a brief introduction on tight wavelet frames and the wavelet frame based image restoration. Interested readers may consult \cite{J.F.Cai2008,B.Dong2013,B.Dong2015,A.Ron1997,Shen2010} for detailed surveys on tight wavelet frames, and \cite{J.F.Cai2012} for detailed reviews on the wavelet frame based image restoration. For the sake of simplicity, we only discuss the real valued wavelet tight frame systems, but note that it is not difficult to extend the idea to the complex case.

Let $\msH$ be a Hilbert space equipped with an inner product $\la\cdot,\cdot\ra$. A (Bessel) sequence $\{\bvphi_n:n\in\Z\}\subseteq\msH$ is called a tight frame on $\msH$ if
\begin{align}\label{TightFrame}
\|\bu\|^2=\sum_{n\in\Z}|\la\bu,\bvphi_n\ra|^2~~~~~\text{for all}~~~\bu\in\msH.
\end{align}
Given $\{\bvphi_n:n\in\Z\}\subseteq\msH$, we define the analysis operator $\bmW:\msH\to\ell_2(\Z)$ as
\begin{align*}
\bu\in\msH\mapsto \bmW\bu=\{\la\bu,\bvphi_n\ra:n\in\Z\}\in\ell_2(\Z).
\end{align*}
The synthesis operator $\bmW^T:\ell_2(\Z)\to\msH$ is defined as the adjoint of $\bmW$:
\begin{align*}
\bc\in\ell_2(\Z)\mapsto\bmW^T\bc=\sum_{n\in\Z}\bc(n)\bvphi_n\in\msH.
\end{align*}
Then $\{\bvphi_n:n\in\Z\}$ is a tight frame on $\msH$ if and only if $\bmW^T\bmW=\bmI$. It follows that, for a given tight frame $\{\bvphi_n:n\in\Z\}$, we have the following canonical expression:
\begin{align*}
\bu=\sum_{n\in\Z}\la\bu,\bvphi_n\ra\bvphi_n,
\end{align*}
with $\bmW\bu=\{\la\bu,\bvphi_n\ra:n\in\Z\}$ being called the canonical tight frame coefficients. Hence, the tight frames are extensions of orthonormal bases to the redundant systems. In fact, a tight frame is an orthonormal basis if and only if $\|\bvphi_n\|=1$ for all $n\in\Z$.

One of the most widely used class of tight frames is the discrete wavelet frame generated by a set of finitely supported filters $\{\a_1,\cdots,\a_m\}$. Throughout this paper, we only discuss the two dimensional undecimated wavelet frames on $\ell_2(\Z^2)$, but note that it is not difficult to extend to $\ell_2(\Z^d)$ with $d\geq3$. For $\a\in\ell_1(\Z^2)$, define a convolution operator $\bmS_{\a}:\ell_2(\Z^2)\to\ell_2(\Z^2)$ by
\begin{align}\label{DiscreteConv}
(\bmS_{\a}\bu)(\bk)=(\a\ast\bu)(\bk)=\sum_{\bsl\in\Z^2}\a(\bk-\bsl)\bu(\bsl)~~~~~\text{for}~~~\bu\in\ell_2(\Z^2).
\end{align}
Given a set of finitely supported filters $\{\a_1,\cdots,\a_m\}$, define the analysis operator $\bmW$ and the synthesis operator $\bmW^T$ by
\begin{align}
\bmW&=\left[\bmS_{\a_1(-\cdot)}^T,\bmS_{\a_2(-\cdot)}^T,\cdots,\bmS_{\a_m(-\cdot)}^T\right]^T,\label{AnalConv2}\\
\bmW^T&=\big[\bmS_{\a_1},\bmS_{\a_2},\cdots,\bmS_{\a_m}\big],\label{SyntConv2}
\end{align}
respectively. Then, the direct computation can show that the rows of $\bmW$ form a tight frame on $\ell_2(\Z^2)$ (i.e. $\bmW^T\bmW=\bmI$) if and only if the filters $\{\a_1,\cdots,\a_m\}$ satisfy
\begin{align}\label{SomeUEP}
\sum_{j=1}^m\sum_{\bsl\in\Z^2}\a_j(\bk+\bsl)\a_j(\bsl)=\dde(\bk)=\left\{\begin{array}{cl}
1~&\text{if}~\bk=\0,\\
0~&\text{if}~\bk\neq\0,
\end{array}\right.
\end{align}
called the {\emph{unitary extension principle}} (UEP) condition \cite{B.Han2011}. Finally, for a two dimensional discrete image on the finite grid, $\bmS_{\a}$ in \cref{DiscreteConv} with a finitely supported $\a$ denotes the discrete convolution under the periodic boundary condition throughout this paper.

For the image restoration, we assume that there exists a tight wavelet frame $\bmW$ defined as \cref{AnalConv2} with filters $\left\{\a_1,\cdots,\a_m\right\}$ such that $\bu$ is sparse under $\bmW$. This leads us to solve
\begin{align}\label{FrameModelGeneric}
\min_{\bc}\f{1}{2}\|\bmA\bmW^T\bc-\bsf\|_2^2+\f{1}{2}\|(\bmI-\bmW\bmW^T)\bc\|_2^2+\|\gga\cdot\bc\|_1
\end{align}
and restore $\bu=\bmW^T\bc$. In \cref{FrameModelGeneric}, the first term is a data fidelity term, the third $\ell_1$ norm is defined as
\begin{align}\label{Frameell1norm}
\|\gga\cdot\bc\|_1=\sum_{j=1}^m\gamma_j\left\|\bc_j\right\|_1
\end{align}
with the parameters $\gga=\left[\gamma_1,\ldots,\gamma_m\right]^T$ and $\gamma_j>0$, to promote the sparsity of $\bc$. The second term penalizes the distance between $\bc$ and its projection onto $\msR(\bmW)$. Since we obtain the image via $\bu=\bmW^T\bc$, the second term in fact forces the sparse coefficient $\bc$ close to the canonical coefficient of the restored image. Since the magnitude/decay of the canonical coefficient reflects the regularity of the image under some mild conditions on the tight frame system $\bmW$ \cite{L.Borup2004}, we can see that the second and third terms promote a balance between the sparsity of coefficient and the regularity of image.

We can further impose the flexibility on \cref{FrameModelGeneric} by introducing the weight on the second term:
\begin{align}\label{Balanced}
\min_{\bc}\f{1}{2}\|\bmA\bmW^T\bc-\bsf\|_2^2+\f{\mu}{2}\|(\bmI-\bmW\bmW^T)\bc\|_2^2+\|\gga\cdot\bc\|_1,
\end{align}
called the balanced approach \cite{J.F.Cai2008,R.H.Chan2003}. If $\mu=0$, then \cref{Balanced} becomes the synthesis approach \cite{I.Daubechies2007,M.J.Fadili2007,M.A.T.Figueiredo2003}
\begin{align}\label{Synthesis}
\min_{\bc}\f{1}{2}\|\bmA\bmW^T\bc-\bsf\|_2^2+\|\gga\cdot\bc\|_1,
\end{align}
as we aim to find the sparsest coefficient that synthesizes the image. On the other hand, if $\mu=\infty$, then it must be $(\bmI-\bmW\bmW^T)\bc=\0$. Since we then have $\bc\in\msR(\bmW)$, it follows that $\bc=\bmW\bu$ for some $\bu$. Hence, \cref{FrameModelGeneric} becomes the analysis approach \cite{J.F.Cai2009/10}
\begin{align}\label{Analysis}
\min_{\bu}\f{1}{2}\|\bmA\bu-\bsf\|_2^2+\|\gga\cdot\bmW\bu\|_1.
\end{align}
Since the analysis approach promotes the sparse canonical coefficient, it emphasizes the regularity of the restored image.

These three approaches are equivalent if and only if $\bmW^T\bmW=\bmW\bmW^T=\bmI$, i.e. $\bmW$ forms an orthonormal basis. However, since we consider the redundant system, each approach will lead to the different restoration results. Notice that the coefficient vector $\bc$ obtained from the synthesis approach \cref{Synthesis} is much sparser than the canonical coefficient $\bmW\bu$ obtained from the analysis approach \cref{Analysis}. It is empirically observed that the synthesis approach \cref{Synthesis} tends to yield some artifacts in the restored image. In contrast, the analysis approach \cref{Analysis} usually have less artifacts as it promotes the regularity of the image along the singularities. Theoretically, it has been proved in \cite{J.F.Cai2012} that, under a suitable choice of parameters, the analysis approach \cref{Analysis} can be seen as sophisticated discretization of a variational model including the Rudin-Osher-Fatemi model \cite{L.I.Rudin1992}, which enables a geometric interpretation on \cref{Analysis}.

\section{Proof of \cref{Th1}}\label{ProofTh1} Let $u(\x)$ be defined as in \cref{uModel} where the singularity curves $\Gamma=\bigcup_{j=1}^J\p\Om_j$ satisfy \cref{MajorAssumption}. In the sense of distribution, the gradient of $u$ satisfies
\begin{align}\label{Gradu}
\na u(\x)=\left(\p_1 u(\x),\p_2 u(\x)\right)=\sum_{j=1}^J\left(\aal_j1_{\Om_j}(\x)+\left(\aal_j^T\x+\beta_j\right)\bn_j(\x)\rd\ssi(\x)\big|_{\p\Om_j}\right)
\end{align}
where $\bn_j=\left(n_{j1},n_{j2}\right)$ is the outward normal vector on $\p\Om_j$, and $\ssi$ is the surface measure. We define
\begin{align*}
p(\x)=\left(p_1(\x),p_2(\x)\right)=\sum_{j=1}^J\aal_j1_{\Om_j}(\x)
\end{align*}
and
\begin{align*}
\rd\nnu(\x)=\sum_{j=1}^J\rd\nnu_j(\x)=\sum_{j=1}^J\left(\aal_j^T\x+\beta_j\bn_j(\x)\rd\ssi(\x)\big|_{\p\Om_j}\right).
\end{align*}
Then it is obvious that we have
\begin{align}\label{Gradump}
\na u-p=\rd\nnu=\sum_{j=1}^J\rd\nnu_j.
\end{align}
Since $|\bn_j|=1$ almost everywhere on $\p\Om_j$, $\aal_j^T\x+\beta_j$ is continuous on a compact set $\p\Om_j$, $\na u-p$ in \cref{Gradump} defines a finite Radon vector measure on $\R^2$ supported on $\Gamma$. In addition, since \cref{Gradu} as well as \cref{Gradump} holds in the sense of tempered distribution, we can compute $\msF(\na u-p)(\xxi)$ as a Fourier transform of a measure (e.g. \cite{Folland1999}). Namely, we have
\begin{align}\label{FouGradump}
\msF(\na u-p)(\xxi)=\sum_{j=1}^J\int_{\p\Om_j}e^{-2\pi i\xxi\cdot\x}\rd\nnu_j(\x).
\end{align}
Let $\varphi\in\msS$, where $\msS$ is the space of Schwartz functions (e.g. \cite{Folland1999}). Then the direct computation shows that $\left(\msF(\na u-p)\ast\wh{\varphi}\right)(\xxi)$ satisfies
\begin{align}\label{FouGradumpConv}
\begin{split}
\left(\msF(\na u-p)\ast\wh{\varphi}\right)(\xxi)=\sum_{j=1}^J\int_{\p\Om_j}e^{-2\pi i\xxi\cdot\x}\varphi(\x)\rd\nnu_j(\x).
\end{split}
\end{align}
From \cref{FouGradumpConv}, it is easy to see that if $\Gamma$ lies in the zero level set of $\varphi$, then $\msF(\na u-p)\ast\wh{\varphi}=0$. In particular, if $\varphi$ is defined as \cref{MajorAssumption}, then
\begin{align*}
\wh{\varphi}(\xxi)=\sum_{\bk\in\KK}\a(\bk)\delta(\xxi+\bk),
\end{align*}
which leads to
\begin{align*}
\left(\msF(\na u-p)\ast\wh{\varphi}\right)(\xxi)=\sum_{\bk\in\KK}\msF(\na u-p)(\xxi+\bk)\a(\bk)=0.
\end{align*}
This proves \cref{FirstAnnihil}.

To complete the proof, we further note that in the sense of distribution, $\na_sp$ satisfies
\begin{align}\label{SymGradp}
\na_sp=-\sum_{j=1}^J\left[\begin{array}{cc}
\alpha_{j1}n_{j1}&\displaystyle{\f{1}{2}\left(\alpha_{j1}n_{j2}+\alpha_{j2}n_{j1}\right)}\\
\displaystyle{\f{1}{2}\left(\alpha_{j1}n_{j2}+\alpha_{j2}n_{j1}\right)}&\alpha_{j2}n_{j2}
\end{array}\right]\rd\ssi\big|_{\p\Om_j}:=\sum_{j=1}^J\rd\Ph_j.
\end{align}
Following the similar arguments, it is not difficult to verify that $\rd\Ph_j$ in \cref{SymGradp} defines a finite Radon tensor measure on $\R^2$ supported on $\Gamma$. which means that $\msF(\na_sp)$ is defined as
\begin{align}\label{FouSymGradp}
\msF(\na_sp)(\xxi)=\sum_{j=1}^J\int_{\p\Om_j}e^{-2\pi i\xxi\cdot\x}\rd\Ph_j(\x).
\end{align}
Therefore, for $\varphi\in\msS$, we again have
\begin{align*}
\left(\msF(\na_sp)\ast\wh{\varphi}\right)(\xxi)=\sum_{j=1}^J\int_{\p\Om_j}e^{-2\pi i\xxi\cdot\x}\varphi(\x)\rd\Ph_j(\x),
\end{align*}
and in particular, for $\varphi$ satisfying \cref{MajorAssumption}, we have
\begin{align*}
\left(\msF(\na_sp)\ast\wh{\varphi}\right)(\xxi)=\sum_{\bk\in\KK}\msF(\na_sp)(\xxi+\bk)\a(\bk)=0,
\end{align*}
and this completes the proof.

\section{Proof of \cref{Th2}}\label{ProofTh2}

Let $u(\x)$ be a piecewise linear function defined as \cref{uModel}, and $\bsv=\msF(u)\big|_{\OO}$. Assume \cref{Assumption}, and we consider the full SVDs
\begin{align}\label{HankelSVD1}
\bmH\left(\bmD\bsv-\bq\right)&=\sum_{l=1}^{M_2}\Sig_1^{(l,l)}\bX_1^{(:,l)}\left(\bY_1^{(:,l)}\right)^*,
\end{align}
and
\begin{align}\label{HankelSVD2}
\bmH\left(\bmE\bq\right)&=\sum_{l=1}^{M_2}\Sig_2^{(l,l)}\bX_2^{(:,l)}\left(\bY_2^{(:,l)}\right)^*,
\end{align}
with
\begin{align*}
\Sig_1^{(1,1)}&\geq\cdots\geq\Sig_1^{(r_1,r_1)}>0,~~~\text{and}~~~\Sig_1^{(l,l)}=0~~\text{for}~~l>r_1,
\end{align*}
and
\begin{align*}
\Sig_2^{(1,1)}&\geq\cdots\geq\Sig_2^{(r_1,r_1)}>0,~~~\text{and}~~~\Sig_2^{(l,l)}=0~~\text{for}~~l>r_2.
\end{align*}
Using the right singular vectors in \cref{HankelSVD1,HankelSVD2}, we define
\begin{align*}
\a_{1l}=M_2^{-1/2}\bY_1^{(:,l)}~~~~~\text{and}~~~~~\a_{2l}=M_2^{-1/2}\bY_2^{(:,l)}
\end{align*}
by reformulating $M_2\times 1$ vectors into $K_1\times K_2$ filters supported on $\KK$.

Firstly, we consider $\a_{11},\ldots,\a_{1M_2}$. Note that we have
\begin{align}\label{Stiefel}
\sum_{l=1}^{M_2}\left(M_2^{-1/2}\bY_1^{(m,l)}\right)\left(M_2^{-1/2}\overline{\bY}_1^{(n,l)}\right)=M_2^{-1}\dde(m-n)=\left\{\begin{array}{cl}
1~&\text{if}~m=n,\\
0~&\text{if}~m\neq n,
\end{array}\right.
\end{align}
for $m,n=1,\ldots,M_2.$ Then taking summation along the diagonals and rearranging into the two dimensional multi-indices, we can write \cref{Stiefel} as
\begin{align*}
\sum_{l=1}^{M_2}\sum_{\bm\in\KK}\a_{1l}(\bk+\bm)\overline{\a}_{1l}(\bm)=\dde(\bk)=\left\{\begin{array}{cl}
1~&\text{if}~\bk=\0,\\
0~&\text{if}~\bk\neq\0.
\end{array}\right.
\end{align*}
This means that $\bmW_1$ and $\bmW_1^*$ defined as \cref{OurAnalysis,OurSynthesis} using $\a_{11},\ldots,\a_{1M_2}$ satisfies
\begin{align*}
\bmW_1^*\bmW_1\left(\bmD\bsv-\bq\right)=\sum_{l=1}^{M_2}\bmS_{\overline{\a}_{1l}}\left(\bmS_{\a_{1l}(-\cdot)}\left(\bmD\bsv-\bq\right)\right)=\bmD\bsv-\bq,
\end{align*}
which shows that the filters $\a_{11},\ldots,\a_{1M_2}$ form a tight frame system.

In addition, we note that
\begin{align*}
\Sig_1^{(l,l)}\bX_1^{(:,l)}=\left(\bmH\left(\bmD\bsv-\bq\right)\right)\bY_1^{(:,l)}=\left(\bY_1^{(:,l)}(-\cdot)\ast\left(\bmD\bsv-\bq\right)\right)\big|_{\OO:\KK}
\end{align*}
where the discrete convolution $\ast$ is performed by reformulating $\bY_1^{(:,l)}\in\C^{M_2}$ into a $K_1\times K_2$ filter supported on $\KK$. Then \cref{HankelSVD1} implies that for $l=r_1+1,\ldots,M_2,$ we have
\begin{align*}
\left(\bmS_{\a_{1l}(-\cdot)}\left(\bmD\bsv-\bq\right)\right)(\bk)=0,~~~~~\bk\in\OO:\KK,
\end{align*}
which proves \cref{LowRankHankelTightFrame1,TightFrameSparse1}.

Similarly, we can prove that $\a_{21},\ldots,\a_{2M_2}$ satisfies
\begin{align*}
\sum_{l=1}^{M_2}\sum_{\bm\in\KK}\a_{2l}(\bk+\bm)\overline{\a}_{2l}(\bm)=\dde(\bk),
\end{align*}
so that
\begin{align*}
\bmW_2^*\bmW_2\left(\bmE\bq\right)=\sum_{l=1}^{M_2}\bmS_{\overline{\a}_{2l}}\left(\bmS_{\a_{2l}(-\cdot)}\left(\bmE\bq\right)\right)=\bmE\bq.
\end{align*}
Finally, by \cref{HankelSVD2}, we have, for $l=r_2+1,\ldots,M_2$,
\begin{align*}
\left(\bmS_{\a_{2l}(-\cdot)}\left(\bmE\bq\right)\right)(\bk)=0,~~~~~\bk\in\OO:\KK.
\end{align*}
This completes the proof.

\section*{Acknowledgments} The authors would like to thank Dr. Greg Ongie in the Department of Statistics at the University of Chicago, an author of \cite{G.Ongie2018,G.Ongie2015,G.Ongie2016,G.Ongie2017}, and Prof. Yue Hu in the School of Electronics and Information Engineering, Harbin Institute of Technology, an author of \cite{Y.Hu2019} for making the data sets as well as the MATLAB toolbox available so that the experiments can be implemented.}

\section*{References}

\bibliographystyle{siamplain}

\end{document}